\documentclass[twoside,twocolumn,9pt]{article}
\usepackage{extsizes}
\usepackage[super,sort&compress,comma]{natbib} 
\usepackage[version=3]{mhchem}
\usepackage[left=1.5cm, right=1.5cm, top=1.785cm, bottom=2.0cm]{geometry}
\usepackage{balance}
\usepackage{mathptmx}
\usepackage{sectsty}
\usepackage{graphicx} 
\usepackage{lastpage}
\usepackage[format=plain,justification=justified,singlelinecheck=false,font={stretch=1.125,small,sf},labelfont=bf,labelsep=space]{caption}
\usepackage{float}
\usepackage{fancyhdr}
\usepackage{fnpos}
\usepackage[english]{babel}
\addto{\captionsenglish}{%
  
}
\usepackage{array}
\usepackage{droidsans}
\usepackage{charter}
\usepackage[T1]{fontenc}
\usepackage[usenames,dvipsnames]{xcolor}
\usepackage{setspace}
\usepackage[compact]{titlesec}
\usepackage{hyperref}

\usepackage{amsmath,amssymb,amsfonts}
\usepackage{amsthm}
\usepackage{mathrsfs}
\usepackage{caption} 
\usepackage{subcaption}
\usepackage{comment}

\definecolor{cream}{RGB}{222,217,201}

\begin{document}

\pagestyle{fancy}
\thispagestyle{plain}
\fancypagestyle{plain}{
\renewcommand{\headrulewidth}{0pt}
}

\makeFNbottom
\makeatletter
\renewcommand\LARGE{\@setfontsize\LARGE{15pt}{17}}
\renewcommand\Large{\@setfontsize\Large{12pt}{14}}
\renewcommand\large{\@setfontsize\large{10pt}{12}}
\renewcommand\footnotesize{\@setfontsize\footnotesize{7pt}{10}}
\makeatother

\renewcommand{\thefootnote}{\fnsymbol{footnote}}
\renewcommand\footnoterule{\vspace*{1pt}%
\color{cream}\hrule width 3.5in height 0.4pt \color{black}\vspace*{5pt}} 
\setcounter{secnumdepth}{5}

\makeatletter 
\renewcommand\@biblabel[1]{#1}            
\renewcommand\@makefntext[1]%
{\noindent\makebox[0pt][r]{\@thefnmark\,}#1}
\makeatother 
\renewcommand{\figurename}{\small{Fig.}~}
\sectionfont{\sffamily\Large}
\subsectionfont{\normalsize}
\subsubsectionfont{\bf}
\setstretch{1.125} 
\setlength{\skip\footins}{0.8cm}
\setlength{\footnotesep}{0.25cm}
\setlength{\jot}{10pt}
\titlespacing*{\section}{0pt}{4pt}{4pt}
\titlespacing*{\subsection}{0pt}{15pt}{1pt}

\fancyfoot{}
\fancyfoot[LO,RE]{\vspace{-7.1pt}\includegraphics[height=9pt]{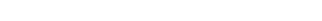}}
\fancyfoot[CO]{\vspace{-7.1pt}\hspace{13.2cm}\includegraphics{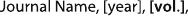}}
\fancyfoot[CE]{\vspace{-7.2pt}\hspace{-14.2cm}\includegraphics{RF}}
\fancyfoot[RO]{\footnotesize{\sffamily{1--\pageref{LastPage} ~\textbar  \hspace{2pt}\thepage}}}
\fancyfoot[LE]{\footnotesize{\sffamily{\thepage~\textbar\hspace{3.45cm} 1--\pageref{LastPage}}}}
\fancyhead{}
\renewcommand{\headrulewidth}{0pt} 
\renewcommand{\footrulewidth}{0pt}
\setlength{\arrayrulewidth}{1pt}
\setlength{\columnsep}{6.5mm}
\setlength\bibsep{1pt}

\makeatletter 
\newlength{\figrulesep} 
\setlength{\figrulesep}{0.5\textfloatsep} 

\newcommand{\topfigrule}{\vspace*{-1pt}%
\noindent{\color{cream}\rule[-\figrulesep]{\columnwidth}{1.5pt}} }

\newcommand{\botfigrule}{\vspace*{-2pt}%
\noindent{\color{cream}\rule[\figrulesep]{\columnwidth}{1.5pt}} }

\newcommand{\dblfigrule}{\vspace*{-1pt}%
\noindent{\color{cream}\rule[-\figrulesep]{\textwidth}{1.5pt}} }

\makeatother

\twocolumn[
  \begin{@twocolumnfalse}
{\includegraphics[height=30pt]{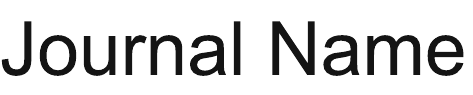}\hfill\raisebox{0pt}[0pt][0pt]{\includegraphics[height=55pt]{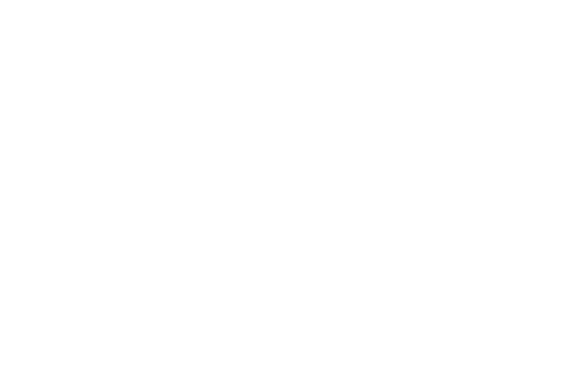}}\\[1ex]
\includegraphics[width=18.5cm]{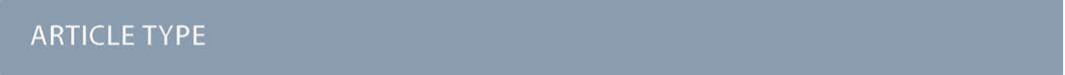}}\par
\vspace{1em}
\sffamily
\begin{tabular}{m{4.5cm} p{13.5cm} }

\includegraphics{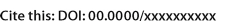} & \noindent\LARGE{\textbf{The untangling number as a measurement of entanglement complexity in 3-periodic networks}} \\

\vspace{0.3cm} & \vspace{0.3cm} \\

 & \noindent\large{Toky Andriamanalina,\textit{$^{a}$} Sonia Mahmoudi,\textit{$^{b,c}$} and Myfanwy E. Evans$^{\ast}$\textit{$^{a}$}} \\

\includegraphics{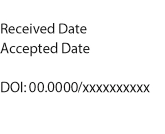} & \noindent\normalsize{Periodic networks serve as models for the structural organisation of biological and chemical crystalline systems. Single or multiple networks can have different configurations in space, where entanglement may arise due to the way the (possibly curvilinear) edges weave around each other. This entanglement influences the functional, physical, and chemical properties of the materials modelled by the networks, which highlights the need to quantify its complexity.
In this paper, we define the least tangled embeddings of 3-periodic networks that we call \textit{ground states}, through the use of knot-theoretic crossing diagrams. The concept of a ground state permits the definition of a measure of entanglement complexity called the \textit{untangling number} that quantifies the distance between a given 3-periodic structure and its least tangled version.} \\

\end{tabular}

 \end{@twocolumnfalse} \vspace{0.6cm}

  ]

\renewcommand*\rmdefault{bch}\normalfont\upshape
\rmfamily
\section*{}
\vspace{-1cm}


\footnotetext{\textit{$^{a}$~Institute for Mathematics, University of Potsdam, Karl-Liebknecht-Str. 24-25, Potsdam Golm, 14476, Germany.}}
\footnotetext{\textit{$^{b}$~Advanced Institute for Materials Research, Tohoku University, 2-1-1 Katahira, Aoba-ku, Sendai, 980-8577, Japan.}}
\footnotetext{\textit{$^{c}$~RIKEN iTHEMS, 2-1 Hirosawa, Wako, Saitama, 351-0198, Japan.}}
\footnotetext{\textit{$^{\ast}$~Corresponding author. E-mail: evans@uni-potsdam.de}}




\section{Introduction}\label{sec:1}

Tangling is a fundamental characteristic of three-dimensional structures formed by a single or multiple networks. Measuring tangling is therefore essential for understanding various entangled biological and chemical structures \cite{batten_1998,Hyde2008ASH,srs_net_nature}, including liquid crystals \cite{C4SM01226G}, DNA origami crystals \cite{dna_tensegrity_triangle, site_directed_placement_DNA_origami}, and, most particularly, coordination polymers \cite{alexandrov2011}, metal-organic frameworks \cite{mof_2001,nets_MOFs}, and covalent organic frameworks \cite{liu2016}. In fact, some physical and chemical properties of materials are influenced by their entanglement \cite{zhang2022}.

\begin{figure*}[hbtp]
    \centering
    \begin{subfigure}[b]{0.3\textwidth}
    \centering
        \includegraphics[width=\textwidth]{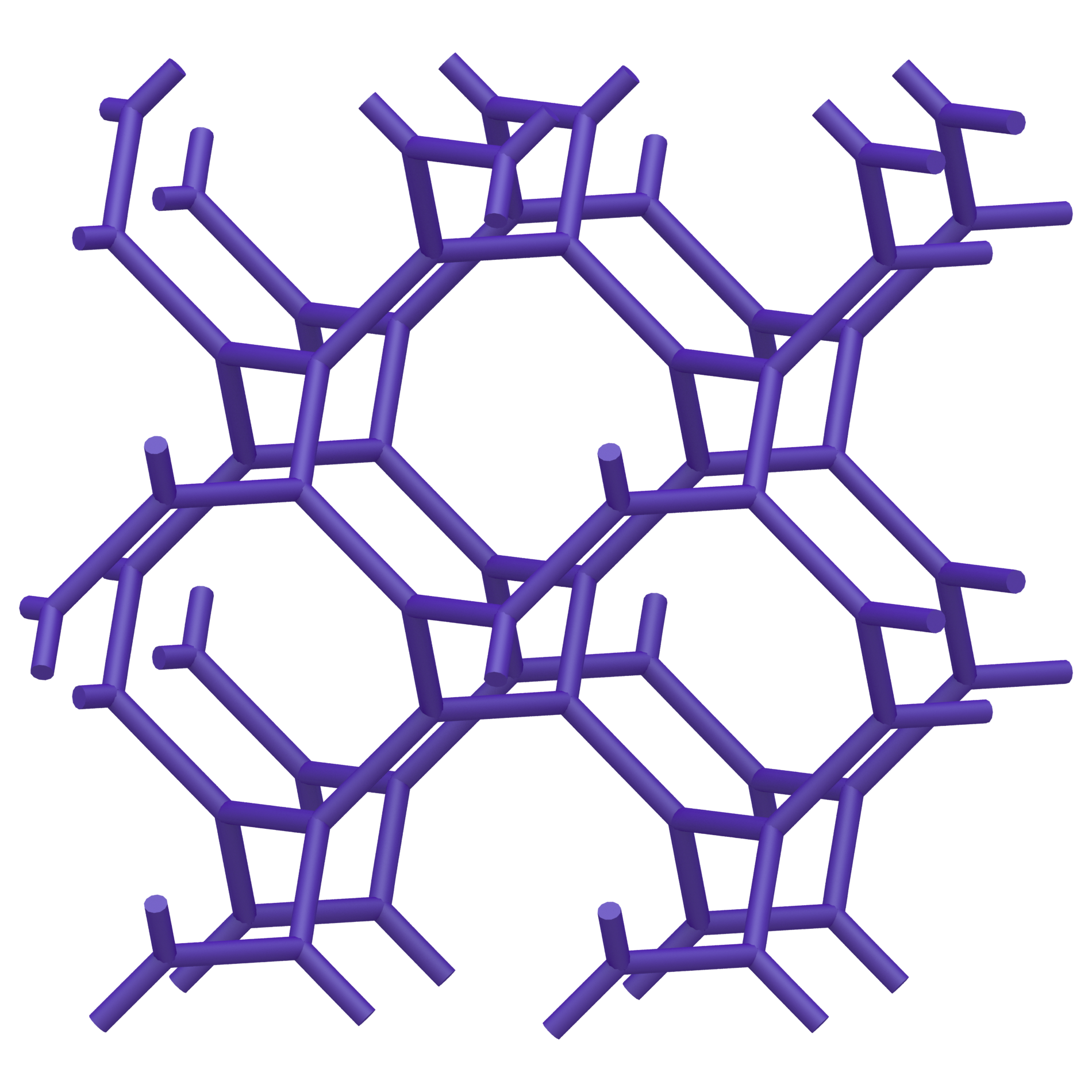}
        \caption{}
        \label{fig:srs_extended}
    \end{subfigure}
    \begin{subfigure}[b]{0.3\textwidth}
    \centering
        \includegraphics[width=\textwidth]{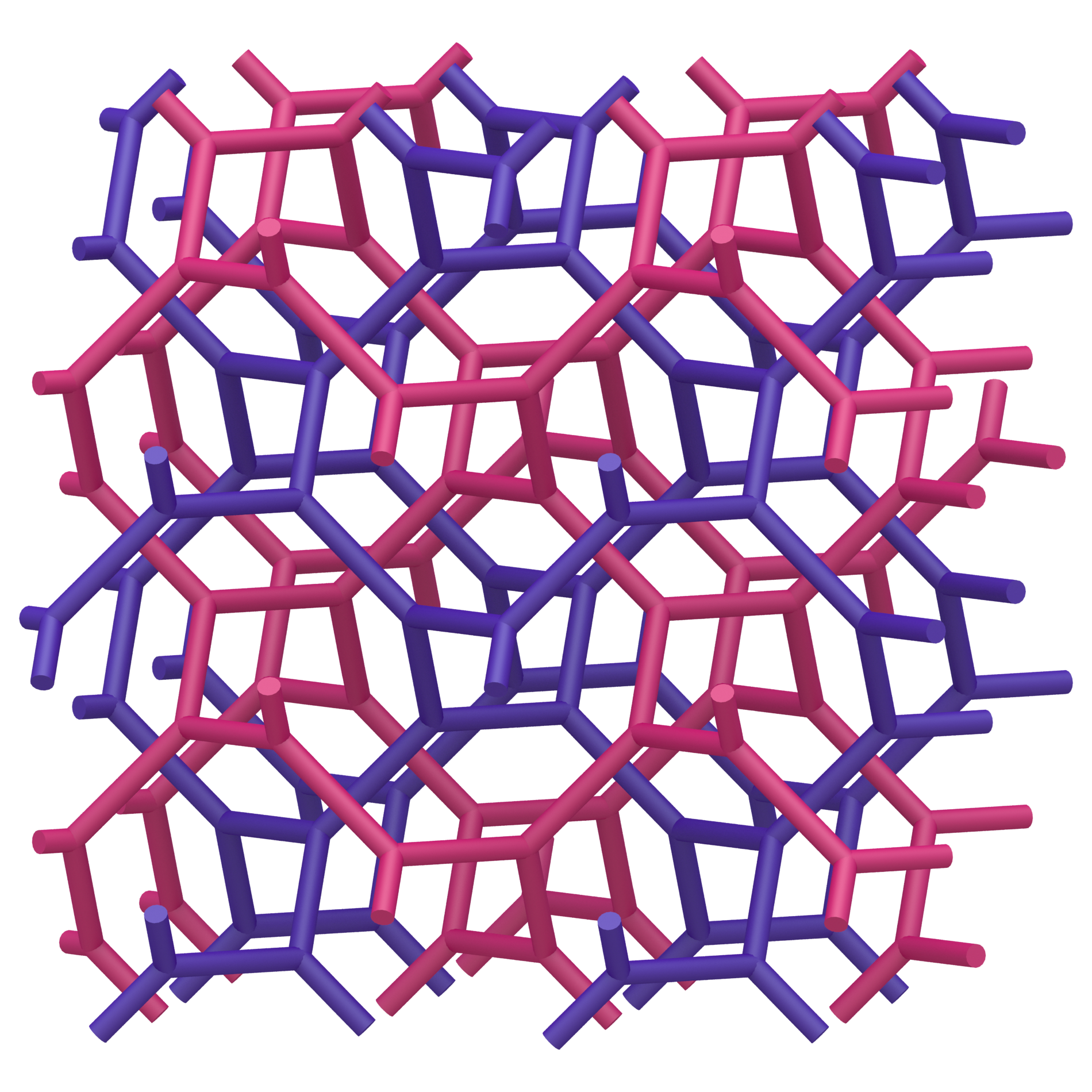}
        \caption{}
        \label{fig:interpenetrating_enantiomorphic_srs_nets_extended}
    \end{subfigure}
    \caption{Various graph embeddings comprising a single or multiple networks: (a) The standard embedding of the regular degree-3 and highly symmetric chiral 3-periodic network called \textbf{srs}. (b) An embedding of a graph comprising two \textbf{srs} networks of opposite handedness called \textbf{srs-c}.}
    \label{fig:srs_and_double_srs}
\end{figure*}

\begin{figure*}[hbtp]
    \centering
    \begin{subfigure}[b]{0.3\textwidth}
    \centering
        \includegraphics[width=\textwidth]{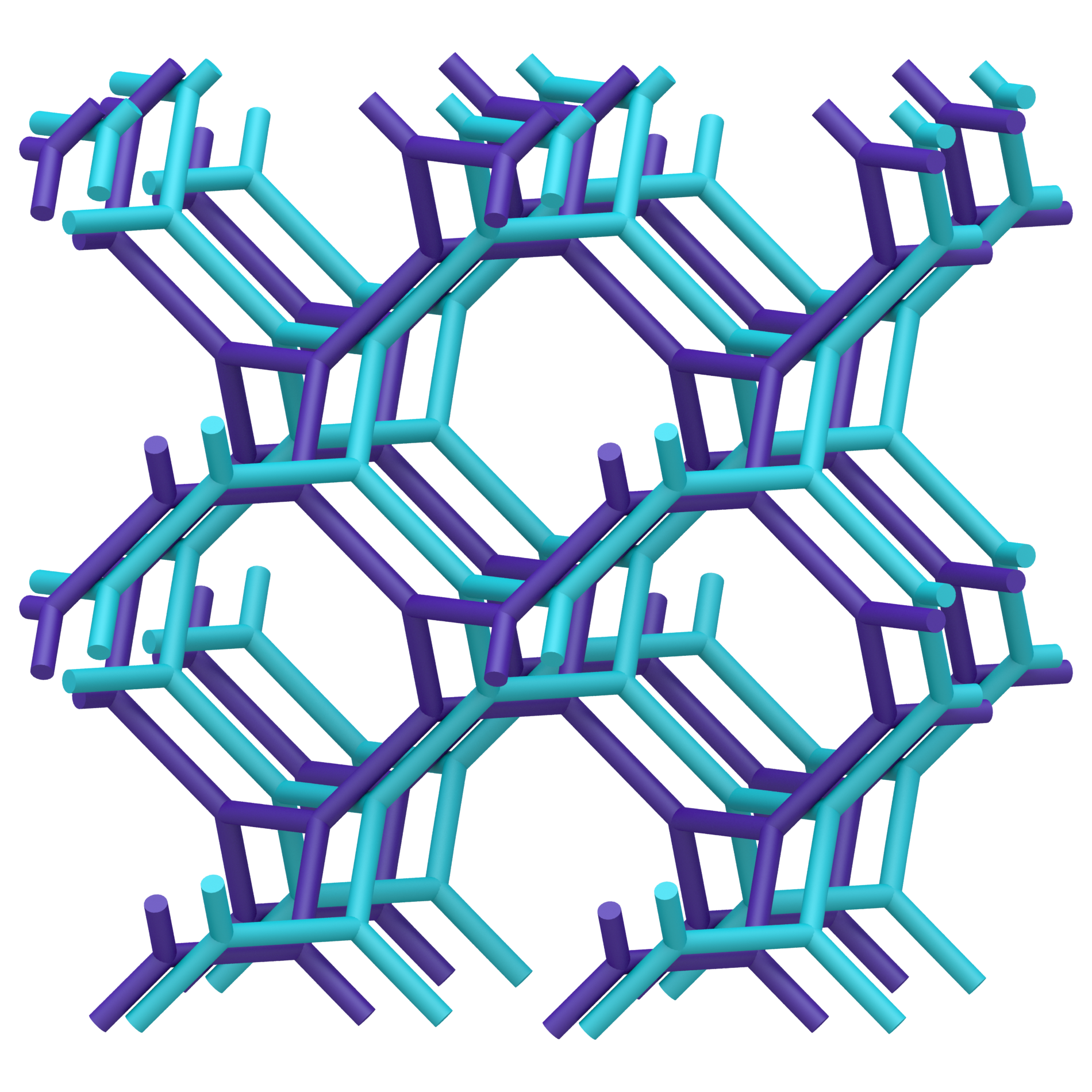}
        \caption{}
        \label{fig:interpenetrating_left_handed_srs_nets_extended}
    \end{subfigure}
    \begin{subfigure}[b]{0.3\textwidth}
    \centering
        \includegraphics[width=\textwidth]{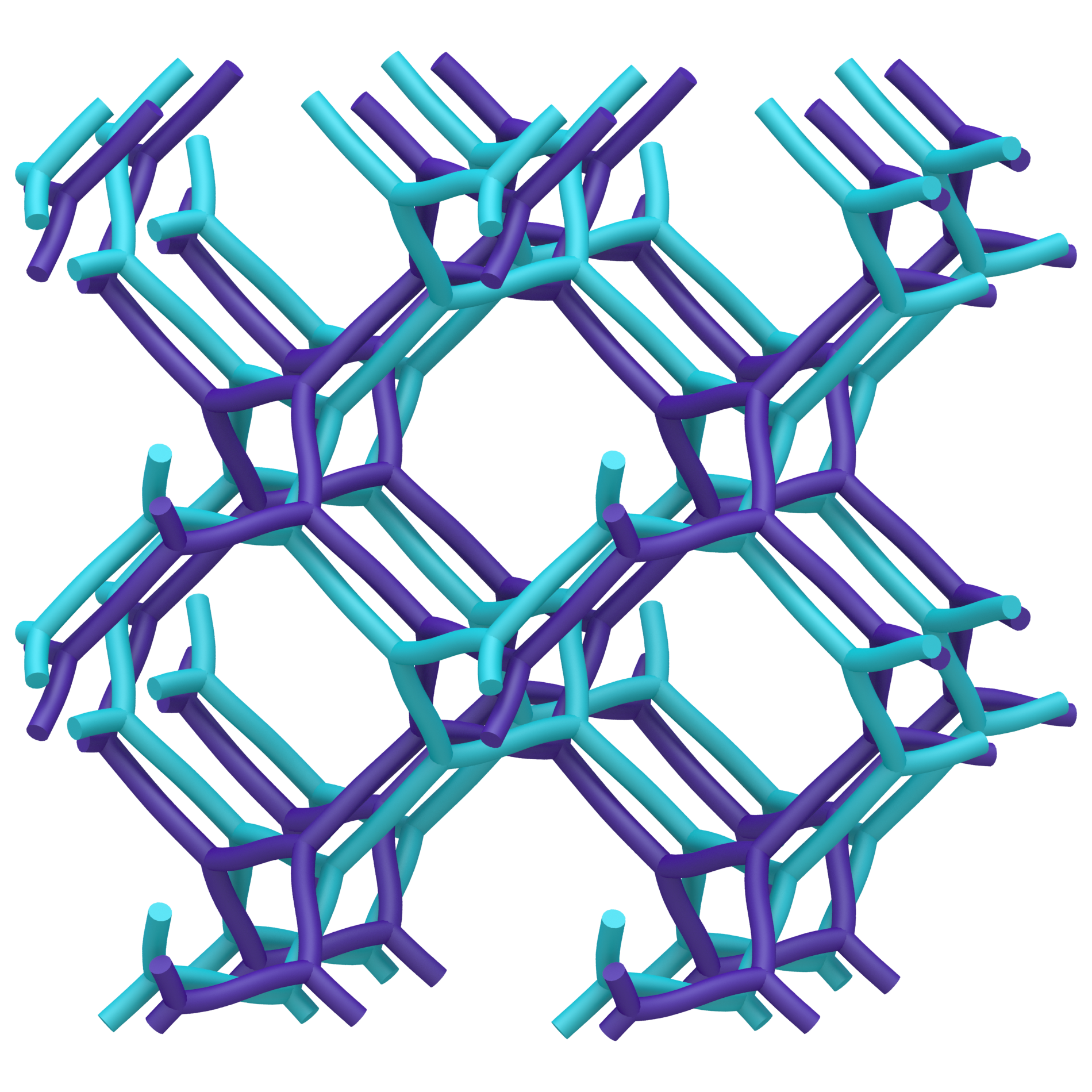}
        \caption{}
        \label{fig:n0p4on2_to_the_6_extended}
    \end{subfigure}
    \begin{subfigure}[b]{0.3\textwidth}
    \centering
        \includegraphics[width=\textwidth]{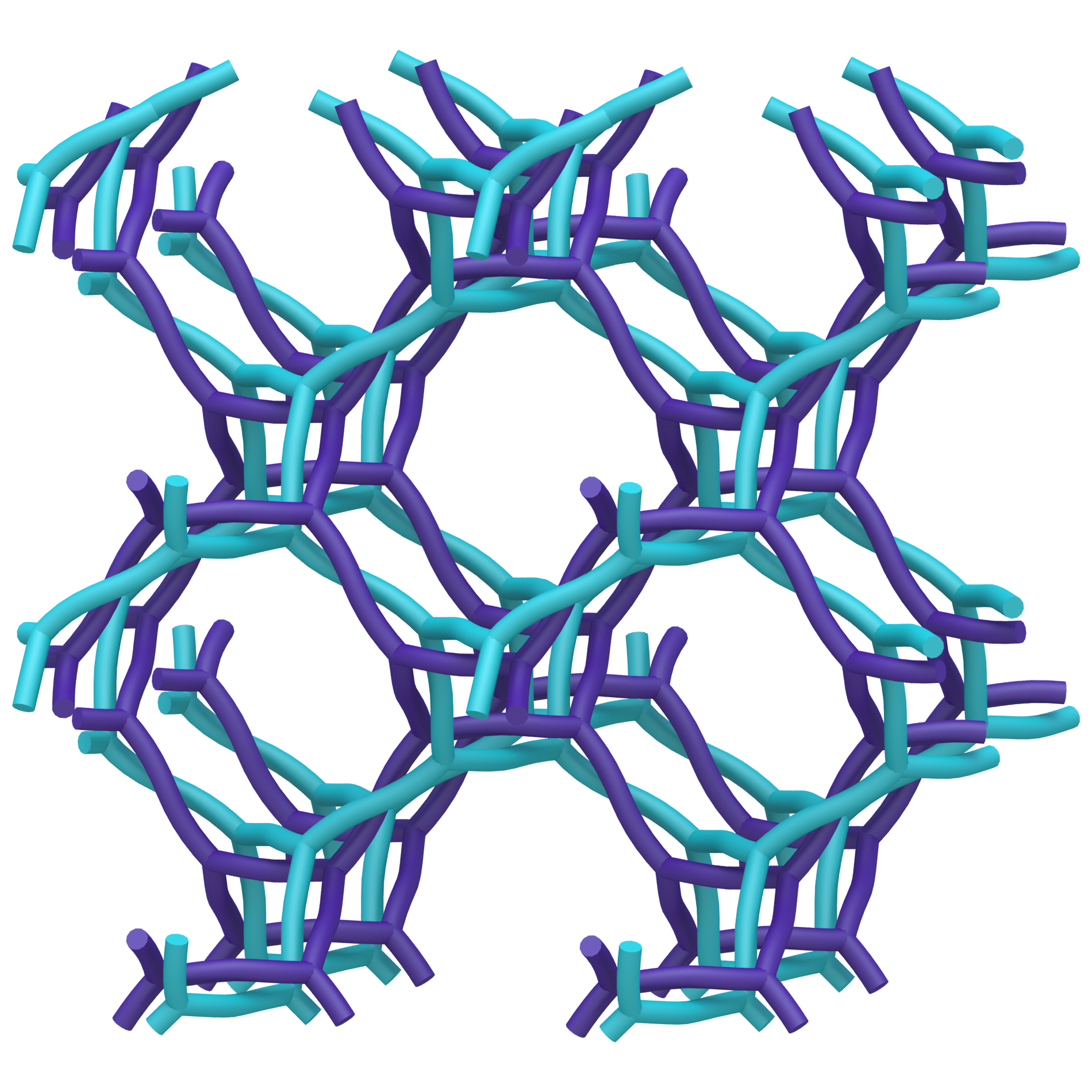}
        \caption{}
        \label{fig:0p6on2_to_the_6_extended}
    \end{subfigure}
    \caption{Distinct embeddings of the same graph: (a) An embedding of two copies of the \textbf{srs} network of the same handedness called \textbf{srs-c*}. (b,c) Two highly symmetric embeddings that are intuitively more tangled than the one in (a). The embedding in (b) is called \textbf{srs-c**}.}
    \label{fig:srs_c_star_srs_c_star_star_0p6on2_to_the_6}
\end{figure*}

Geometric and topological characterisation of 3-periodic network entanglements has been explored in various contexts \cite{CARLUCCI2003247,  DELGADOFRIEDRICHS20052533, Alexandrov:eo5016, periodic_ent_I, topospro, sym14040822}. Formally, a network is a particular type of the more general mathematical concept of a \textit{graph}. A \textit{graph} is an abstract collection of vertices connected by edges. A \textit{network} is an infinite connected graph, and multiple networks are an infinite graph with finitely many connected components \cite{Eon:ib5036,terminology_MOF_COPolymers}. An \textit{embedding} of a given graph is a geometric realisation of the graph in a given space, in which the edges need not be straight, but must not intersect. Geometric properties such as periodicity are thus properties of embeddings rather than of abstract graphs. Nevertheless, in this paper we sometimes abuse terminology by using terms such as `3-periodic graph' or `3-periodic network' for simplicity. In Fig. \ref{fig:srs_extended} we show an embedding of the \textbf{srs} network in the Euclidean space $\mathbb{R}^3$ (For the names of 3-periodic graphs, we use the terminology of RCSR \cite{rcsr} whenever possible). The \textbf{srs} is a network of degree-3, that is, each vertex of the network connects three edges; regular, that is, the convex hull of the neighbours of each vertex is a regular polygon or polyhedron, which in this case is an equilateral triangle; highly symmetric and chiral. Another example of a graph embedding is shown in Fig. \ref{fig:interpenetrating_enantiomorphic_srs_nets_extended}, which is that of two \textbf{srs} networks of opposite handedness called \textbf{srs-c}.
Two embeddings are said to be ambient isotopic if one can be continuously transformed into the other by deforming the structure through space, allowing edges to bend and stretch but not to pass through one another or through vertices.
A fundamental challenge that arises is to determine how tangled an embedding is, up to ambient isotopy, and which embedding is considered the least entangled state \cite{HYDE2011676}. For example, Fig. \ref{fig:srs_c_star_srs_c_star_star_0p6on2_to_the_6} shows three isotopically distinct embeddings of the same graph comprising two left-handed \textbf{srs} networks. Intuitively, the embedding in Fig. \ref{fig:interpenetrating_left_handed_srs_nets_extended} called \textbf{srs-c*} is the least tangled, and the embeddings in Fig. \ref{fig:n0p4on2_to_the_6_extended}, called \textbf{srs-c**}, and Fig. \ref{fig:0p6on2_to_the_6_extended} are more tangled variants. It is this intuitive idea that we formalise rigorously in this paper. We note in particular that the embedding in Fig. \ref{fig:0p6on2_to_the_6_extended} constitutes one of the chiral subspaces of a 3-periodic structure recently observed in a liquid crystal polymer \cite{4NG}, which suggests that embeddings exhibiting progressively more intricate entanglement are increasingly being observed in nature or realised through synthesis.

\begin{figure}[hbtp]
    \centering
    \begin{subfigure}[b]{0.5\textwidth}
        \centering
        \includegraphics[width=0.3125\textwidth]{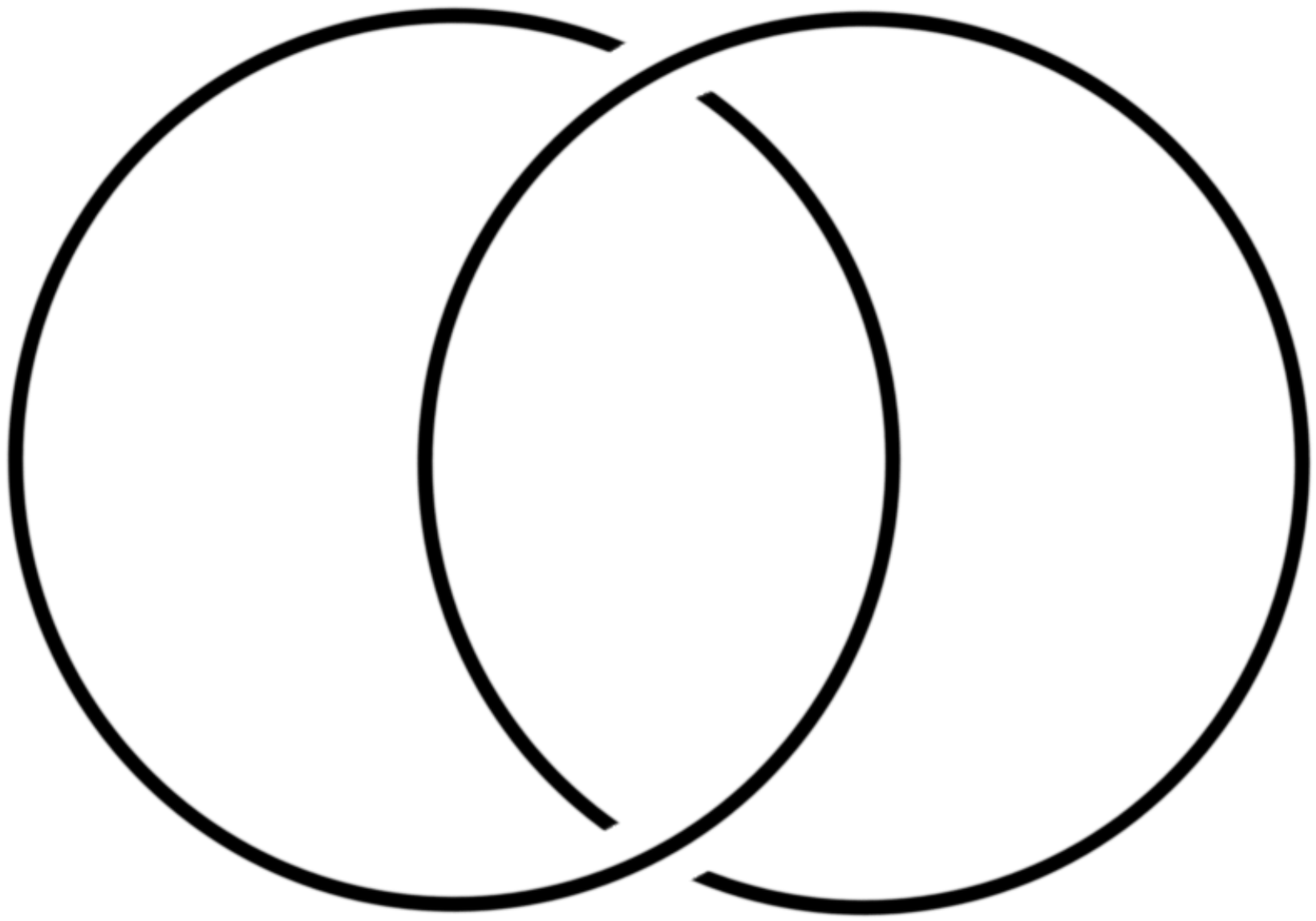}
        \caption{}
        \label{fig:hopf_link_a}
    \end{subfigure}
    \vskip\baselineskip    
    \begin{subfigure}[b]{0.5\textwidth}
        \centering
        \includegraphics[width=0.3125\textwidth]{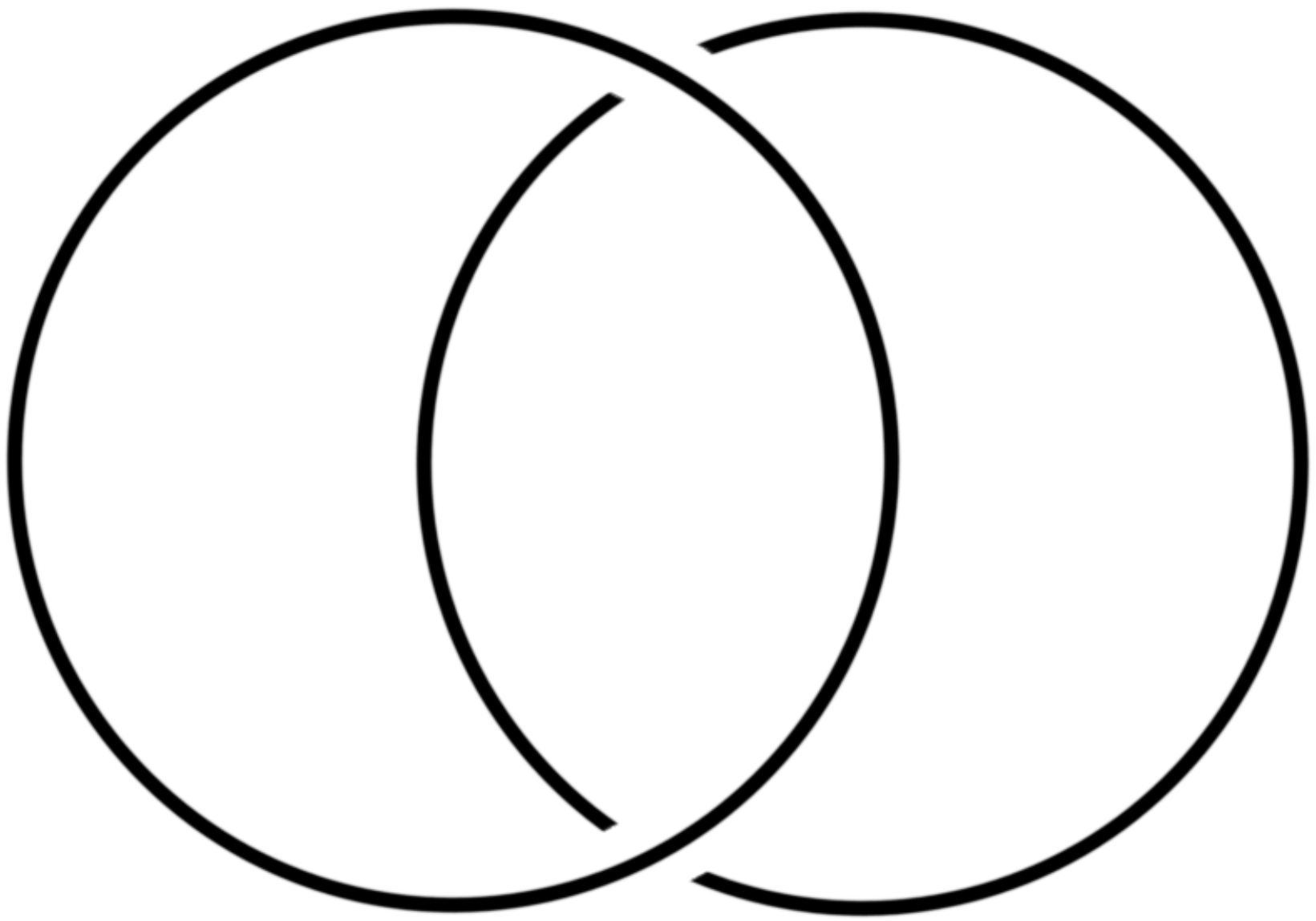}
        \hspace{0.5cm}
        \includegraphics[width=0.4625\textwidth]{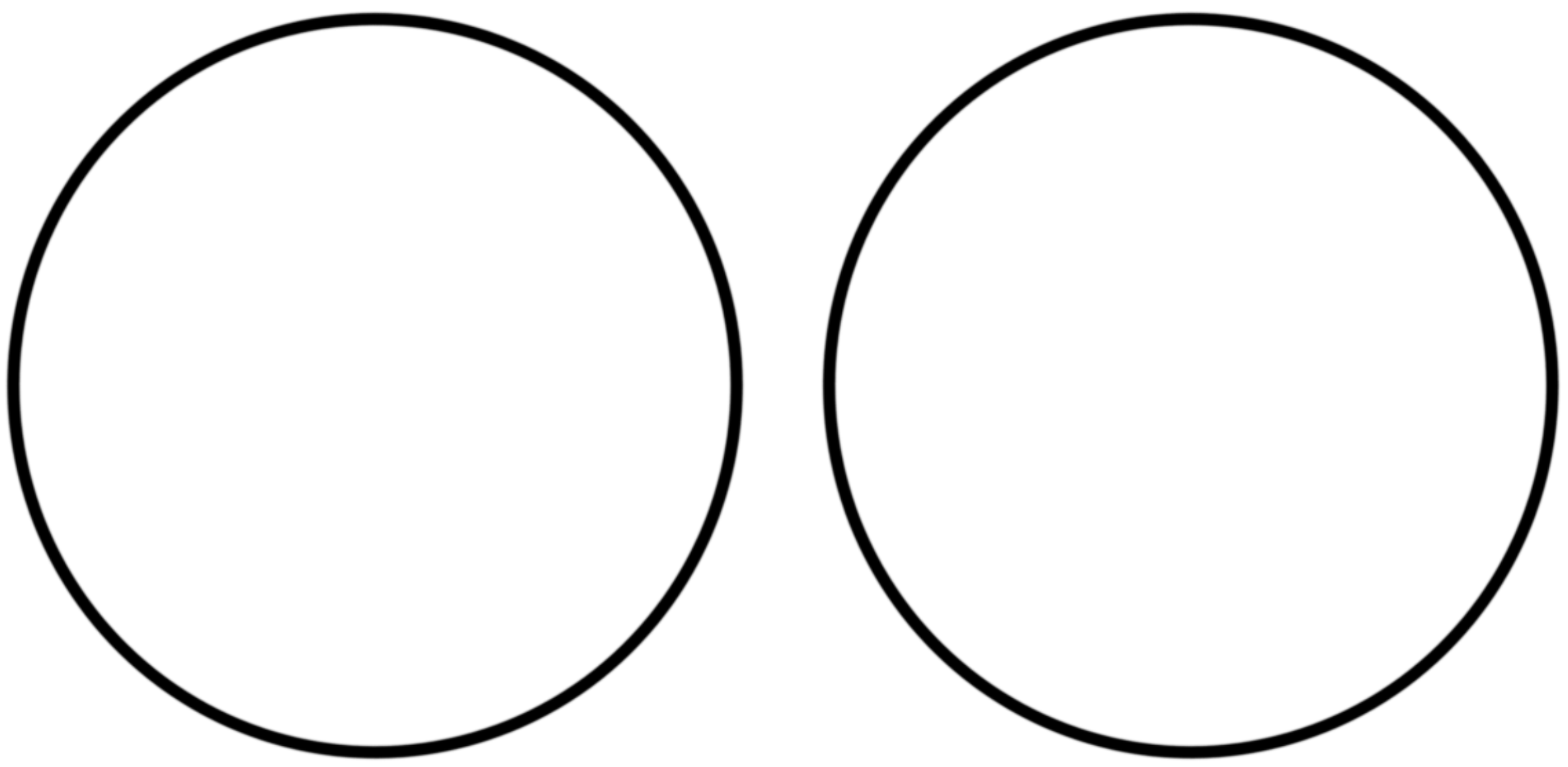}
        \caption{}
        \label{fig:hopf_link_b}
    \end{subfigure}
    \caption{(a) A crossing diagram of the Hopf link, a link with two components. (b) Changing one crossing of the diagram of the Hopf link transforms it into the unlink with two components, a disjoint union of two circles.}
    \label{fig:hopf_link}
\end{figure}

For a much simpler system, knotting and linking in finitely many closed loops is characterised by the mathematical theory of knots and links \cite{Adams.book,Murasugi1996}. A link is a collection of finitely many loops in space, and a knot is a link with a single loop. Links are also classified up to ambient isotopy, meaning that two links are equivalent if one can be transformed into the other without cutting and gluing. The simplest link with $m$ components is the unlink, a disjoint collection of $m$ circles.
A link diagram is a two-dimensional projection of a link together with the crossing information indicating which strand goes over which. Fig. \ref{fig:hopf_link_a}, for instance, displays a diagram of the Hopf link. The complexity of a link is quantified by the use of invariants, which are quantities that remain unchanged under isotopy transformations of a diagram. For example, the \textit{crossing number} is defined as the minimum number of crossings in a diagram among all possible diagrams of a link. Another example is the \textit{unknotting number}, which is the minimum number of crossing changes, similar to allowing curves to pass through one another, that are required to transform a given knot into the unknot \cite{Murasugi1996}. The unknotting number also applies to links, as shown by the example in Fig. \ref{fig:hopf_link_b}, where the Hopf link is transformed into the unlink with two components by changing one crossing of its diagram.

\begin{figure}[hbtp]
    \centering
    \includegraphics[width=0.25\textwidth]{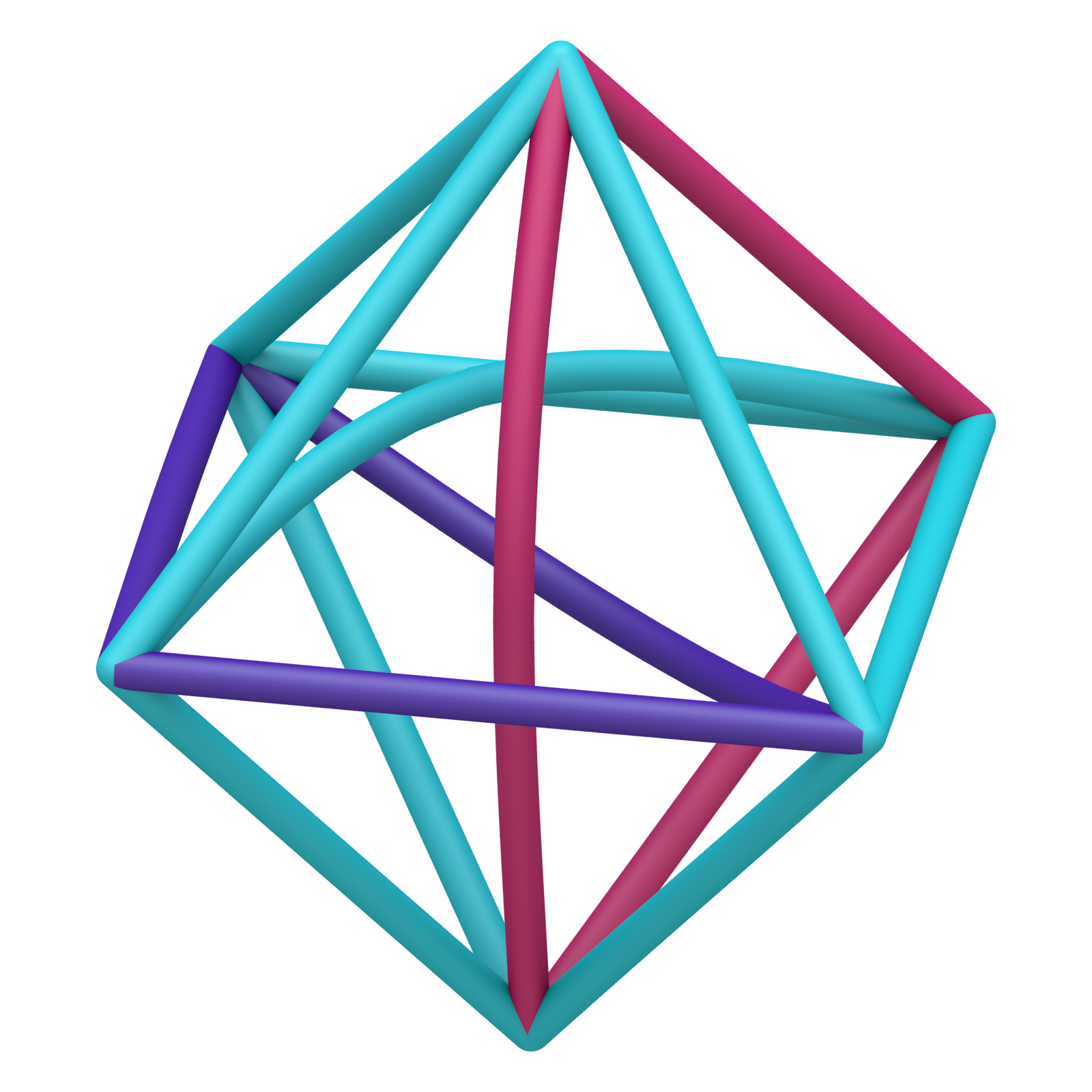}
    \caption{An embedding of the complete graph $K_6$: Two cycles of edges of the embedding form a Hopf link (compare with Fig. \ref{fig:hopf_link_a}). Any embedding of this graph must contain linked cycles.}
    \label{fig:K6}
\end{figure}

\begin{figure}[hbtp]
    \centering
    \begin{subfigure}[b]{0.22\textwidth}
    \centering
        \includegraphics[width=0.55\textwidth]{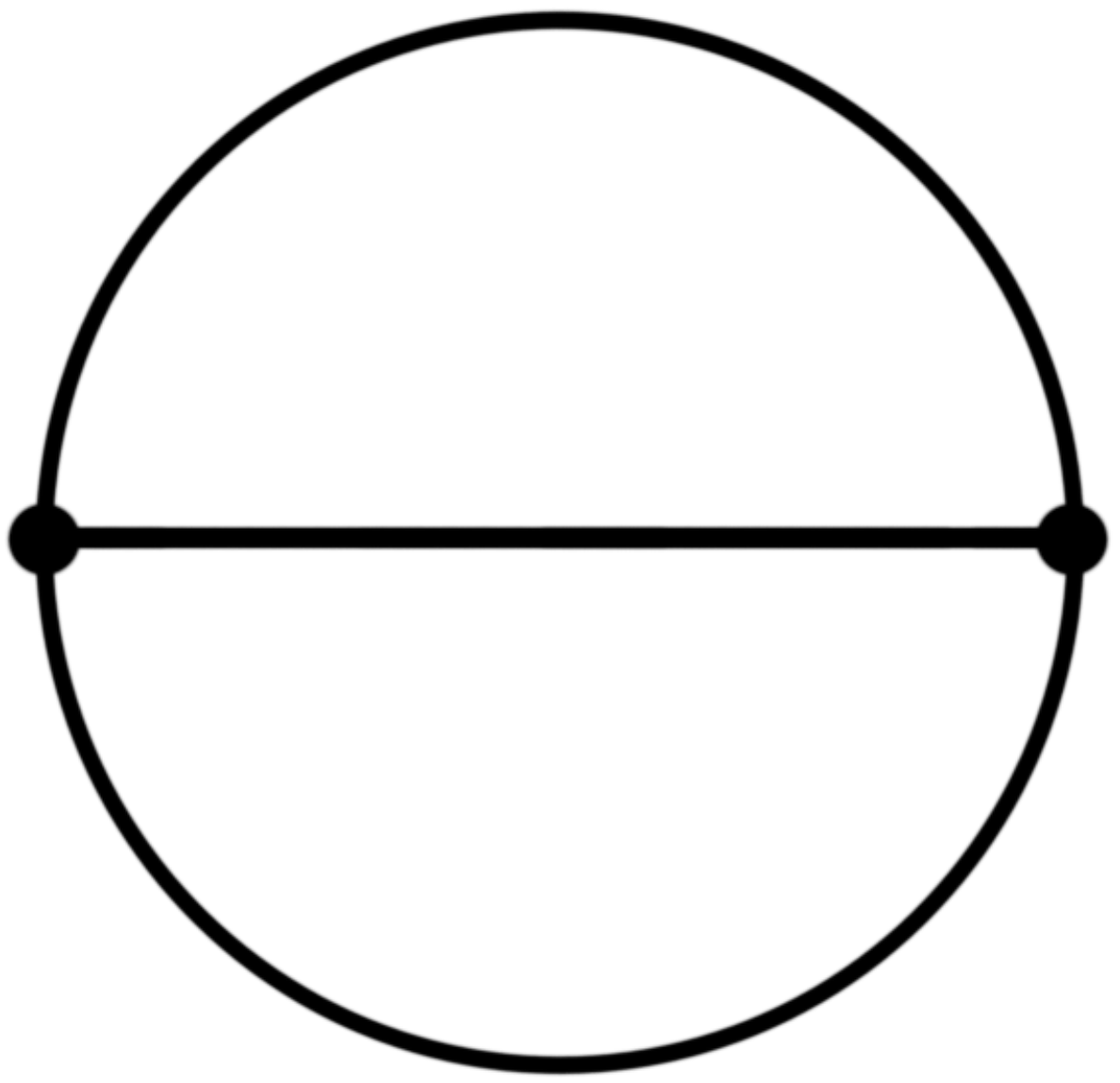}
        \caption{$0_1$}
        \label{fig:unknotted_theta_graph}
    \end{subfigure}
    \begin{subfigure}[b]{0.22\textwidth}
    \centering
        \includegraphics[width=0.55\textwidth]{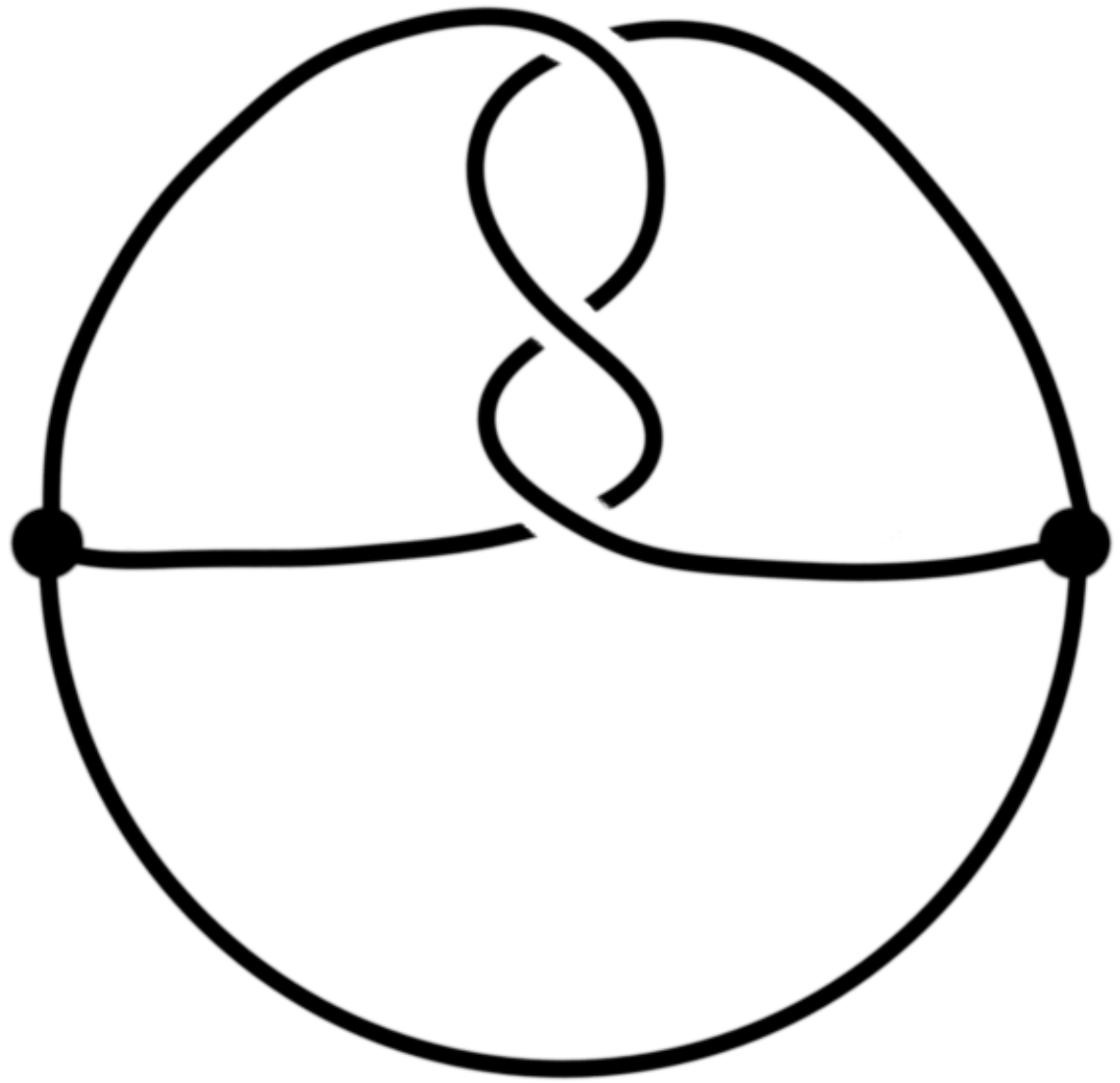}
        \caption{$3_1$}
        \label{fig:3_1}
    \end{subfigure}
    \vskip\baselineskip
    \begin{subfigure}[b]{0.22\textwidth}
    \centering
        \includegraphics[width=0.55\textwidth]{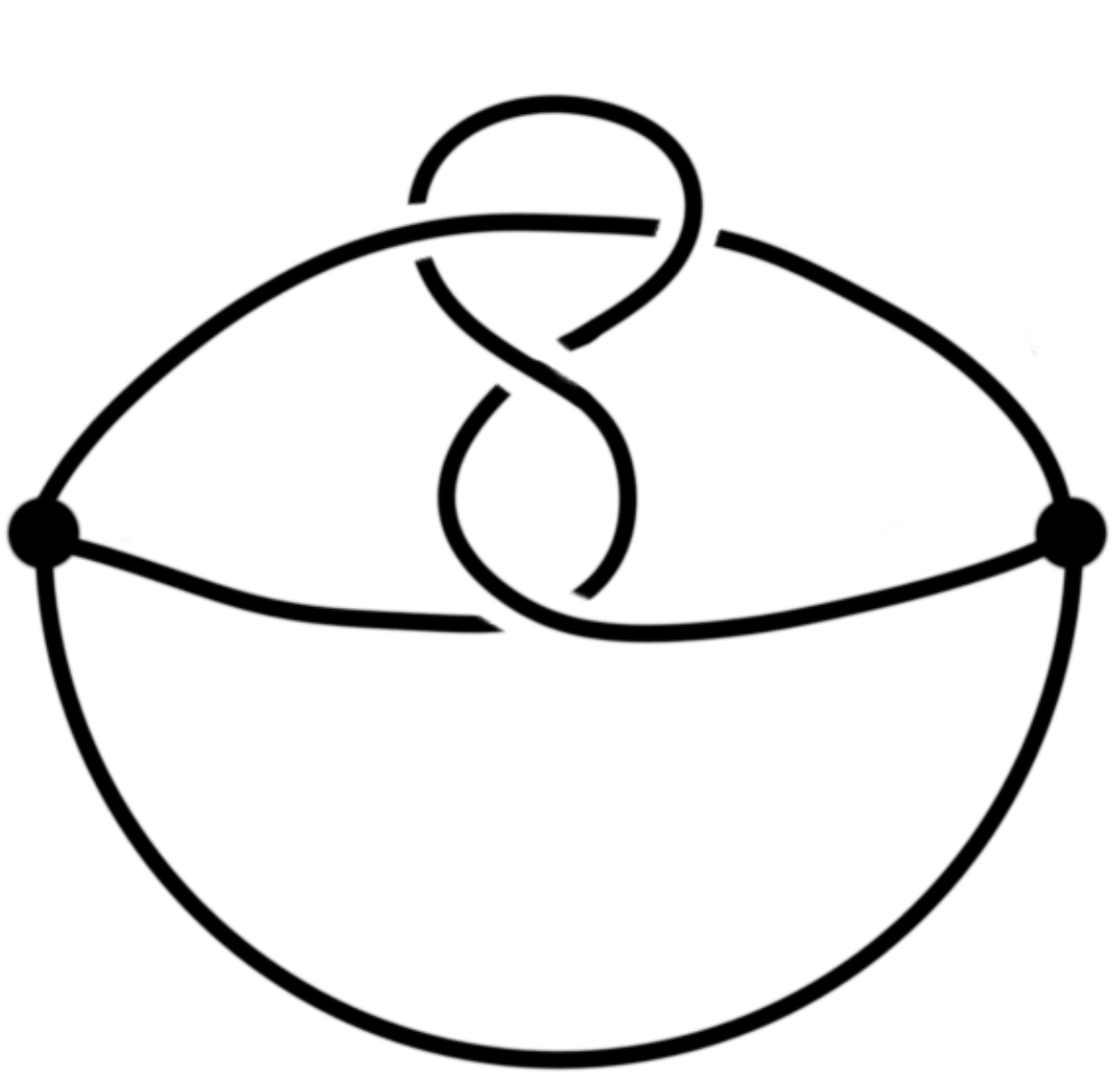}
        \caption{$4_1$}
        \label{fig:4_1}
    \end{subfigure}
    \begin{subfigure}[b]{0.22\textwidth}
    \centering
        \includegraphics[width=0.55\textwidth]{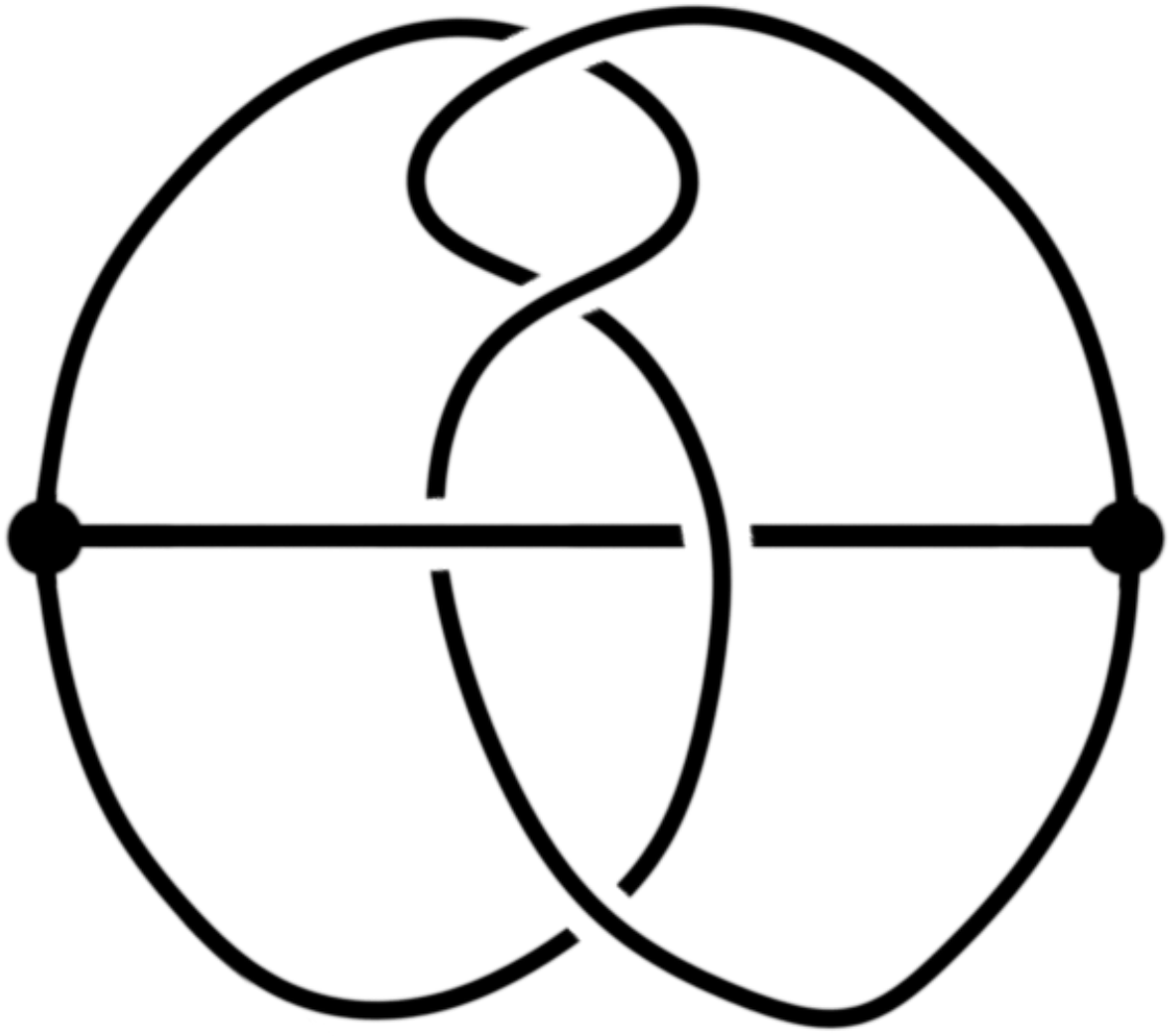}
        \caption{$5_1$}
        \label{fig:5_1}
    \end{subfigure}
    \caption{Diagrams of various embeddings of the $\theta$ graph: The names are taken from the Litherland-Moriuchi enumeration \cite{theta_curves}. The $5_1$ embedding is a ravel, which means that its cycles are neither knotted nor linked, yet the embedding is tangled.}
    \label{fig:theta_graph}
\end{figure}

A finite graph is a collection of finitely many vertices connected by finitely many edges. Entanglement may also arise in embeddings of finite graphs in space, the characterisation of which, although different, is closely related to that of knots and links \cite{R_moves_graph,theta_curves,KAUR2017586}. One useful method is that of a cycle analysis, which involves examining the knotting and linking of disjoint cycles, or loops, formed by the edges of a given graph embedding \cite{R_moves_graph}. An example is given in Fig. \ref{fig:K6} where two cycles of an embedding of the graph called $K_6$ form a Hopf link (compare with the diagram shown in Fig. \ref{fig:hopf_link_a}). There exist, however, some entanglement types that are not captured through a cycle analysis. For example, the embedding $5_1$ of the $\theta$ graph, whose diagram is shown in Fig. \ref{fig:5_1}, is a \textit{ravel} \cite{ravels}, which is a graph embedding that does not possess any knotted or linked cycles, yet is tangled. Ravels are not just theoretical constructions illustrating the limits of a cycle analysis; some of them have in fact been synthesised \cite{synthesised_ravel}.
Furthermore, unlike links, some graphs, such as the $K_6$ graph for example, do not possess an embedding that is unlinked \cite{k6_linked}. This shows that, first, typical methods used to quantify the complexity of knots and links do not necessarily extend automatically to finite graph embeddings, and second, even for finite graphs, the notion of a least tangled embedding is not quite well determined.

\begin{figure*}[hbtp]
    \centering
    \begin{subfigure}[b]{0.3\textwidth}
    \centering
        \includegraphics[width=\textwidth]{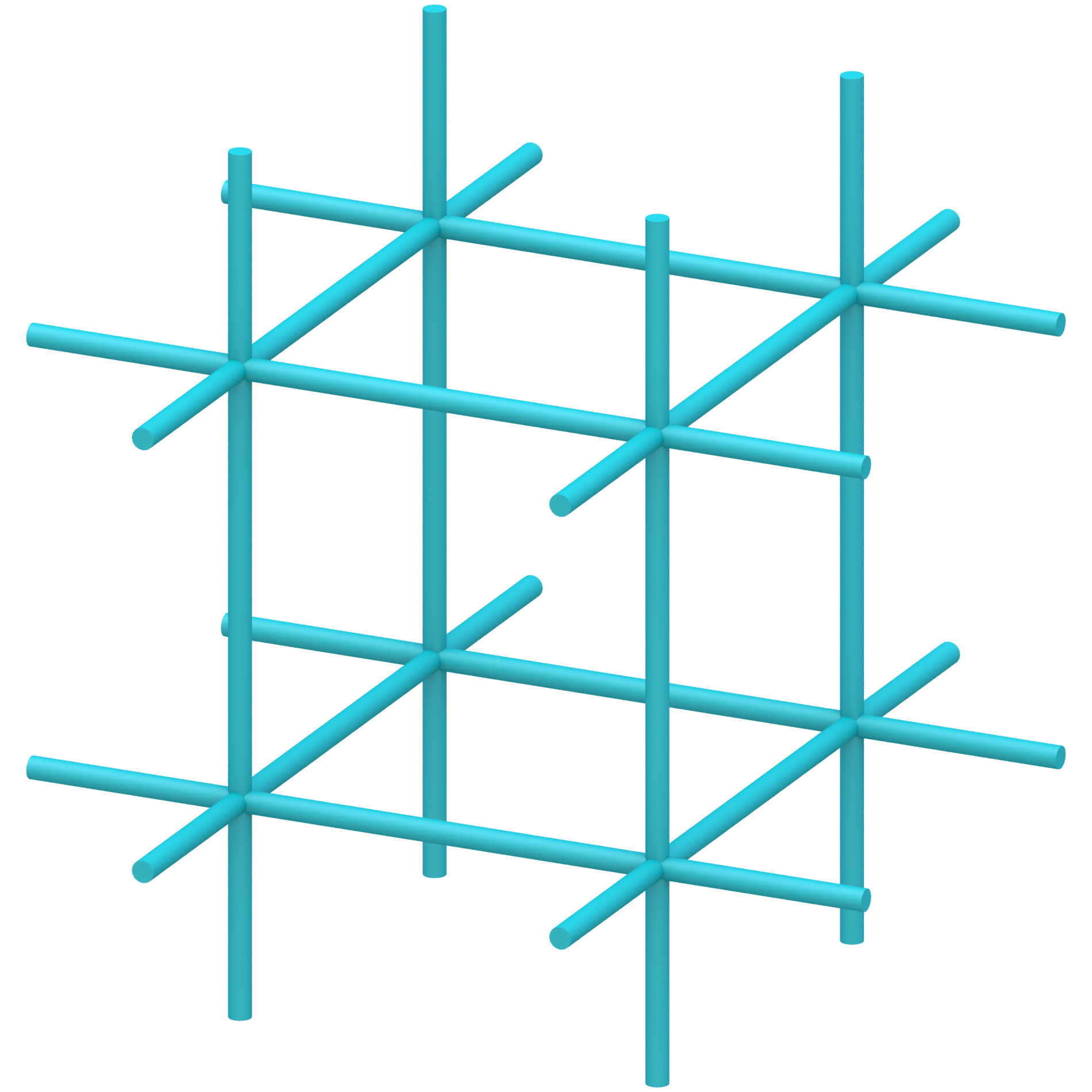}
        \caption{}
        \label{fig:pcu_extended}
    \end{subfigure}
    \begin{subfigure}[b]{0.3\textwidth}
    \centering
        \includegraphics[width=\textwidth]{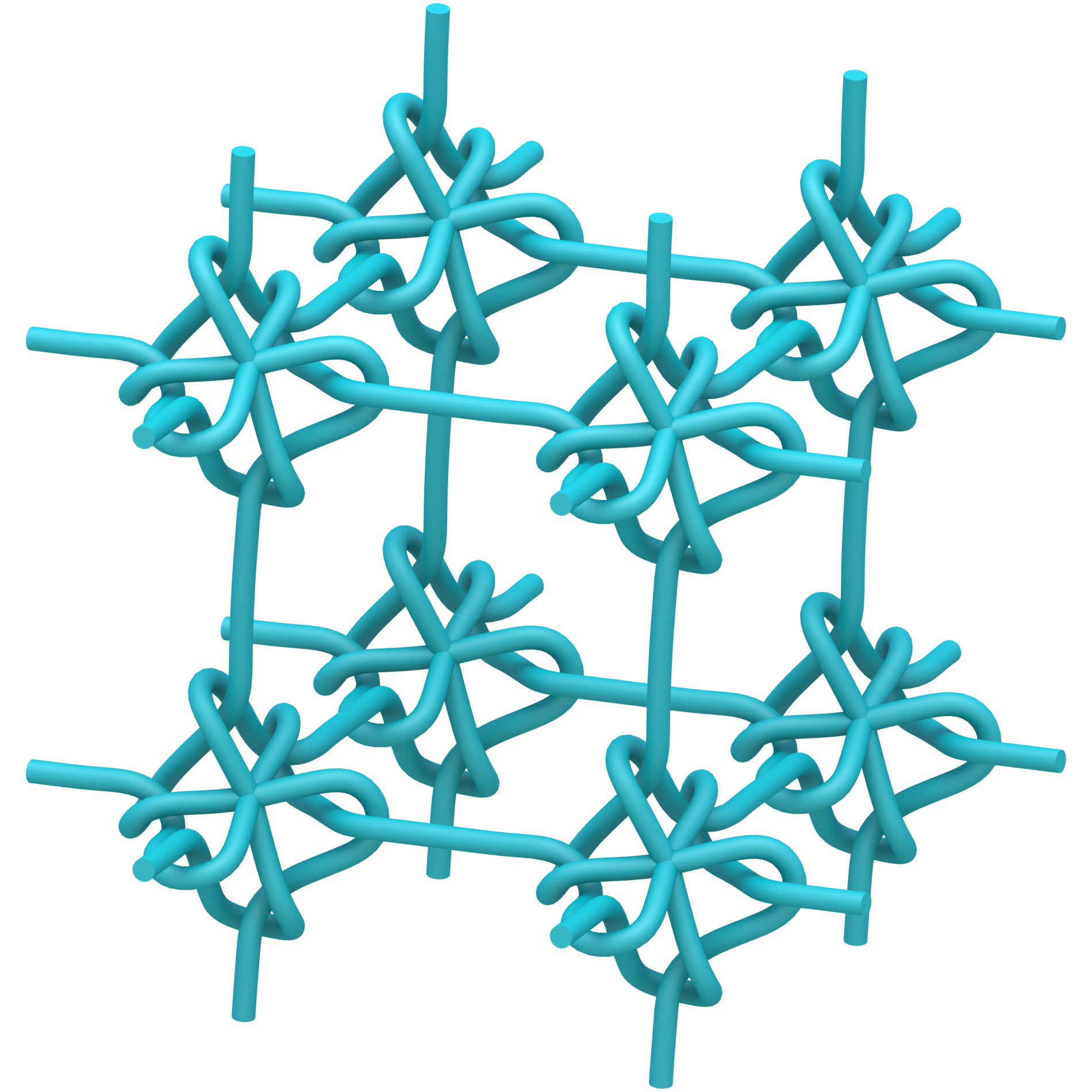}
        \caption{}
        \label{fig:ravelled_pcu_extended}
    \end{subfigure}
    \caption{Two embeddings of the \textbf{pcu} network whose strong rings are neither knotted nor linked, and whose HRNs are therefore trivial: (a) The standard embedding of \textbf{pcu}. (b) A pure 3-periodic ravel, demonstrating that the absence of knotted or linked cycles is not sufficient for an embedding to be the least tangled one.}
    \label{fig:pcu_and_ravelled_pcu_extended}
\end{figure*}

\begin{figure}[hbtp]
    \centering
        \includegraphics[width=0.3\textwidth]{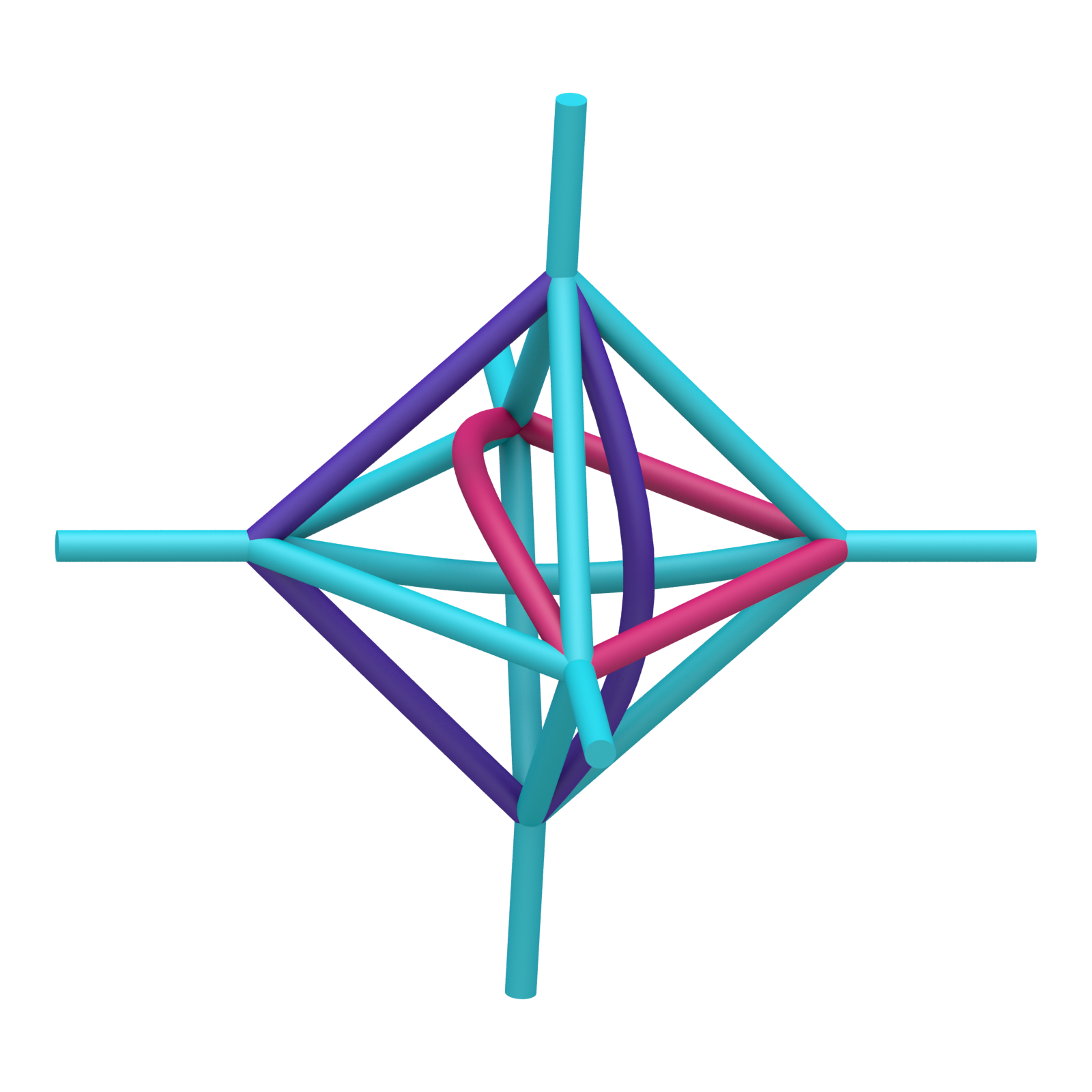}
    \caption{The likely least tangled embedding of a 3-periodic network obtained by replacing the vertices of \textbf{pcu} with $K_6$ graphs. Due the presence of the $K_6$ graphs, any embedding of the network must have linked cycles.}
    \label{fig:pcu-K6}
\end{figure}

Cycle analysis has been extended to periodic graphs as a measure of entanglement complexity in these structures \cite{CARLUCCI2003247,sym14040822,tesselate_decussate}. In this setting, it consists of analysing the knotting and linking of particular cycles in the graph embedding called \textit{strong rings}, which are, essentially, cycles that cannot be subdivided into shorter ones \cite{DELGADOFRIEDRICHS20052480}. This can be done both within a single component and between disjoint components in the graph, with a difference in this signature indicating different entanglements of the embeddings. Based on cycle analysis, another method for quantifying entanglement complexity was introduced in \cite{Alexandrov:eo5016}: the computation of the \textit{Hopf Ring Net} (HRN) of a given graph embedding and the analysis of its complexity. The HRN of a given embedding can be obtained by representing pairwise linked strong rings by their barycentres and connecting them. These two methods, both computable using the software ToposPro \cite{topospro}, can in fact be used to distinguish and establish a hierarchy between the three embeddings of two left-handed \textbf{srs} networks displayed in Fig. \ref{fig:srs_c_star_srs_c_star_star_0p6on2_to_the_6}. Indeed, the links formed by the strong rings of these embeddings are different. Those of \textbf{srs-c*} form only Hopf links. Those of \textbf{srs-c**} form Hopf links and $4^2_1$ links (we use the Alexander-Briggs notation for knots and links \cite{Adams.book}). Finally, the strong rings of the embedding in Fig. \ref{fig:0p6on2_to_the_6_extended} form Hopf links, $4^2_1$ links, and $6^2_1$ links. These links have increasing crossing numbers, 2, 4, and 6 for the Hopf, $4^2_1$, and $6^2_1$ links, respectively, which agrees with the intuition that the complexity of the graph embeddings increases accordingly. However, cycle analysis and HRN computation both rely on the presence of knots and links in cycles, which is not always the case, even for finite graphs, as mentioned earlier. For example, the structure shown in Fig. \ref{fig:ravelled_pcu_extended} is a ravelled embedding of the \textbf{pcu} network, that is, its strong rings are neither knotted nor linked. Consequently, its HRN is trivial. In particular, this example shows that the absence of knotted or linked cycles is not sufficient for an embedding to be the least tangled one. Moreover, it is not a necessary condition either. An exotic example, very unlikely to occur in an actual chemical structure, can be built by  replacing the vertices of the usual embedding of \textbf{pcu} with standard embeddings of the $K_6$ graph. The resulting embedding, shown in Fig. \ref{fig:pcu-K6}, is likely the least tangled embedding of the network, while still having linked cycles. These two examples highlight the need for new methods to capture entanglement complexity in 3-periodic networks, in addition to existing approaches such as cycle analysis and HRN computation.

Various studies \cite{Wells:a01232, Wells:a01308, tesselate_decussate, Delgado-Friedrichs:au5000, Power:ib5087, myf2011_entanglement_graphs, myf_ideal_geo} address the problem of finding an embedding that may be considered the canonical least tangled form of networks.
For example, in a series of papers on the geometrical basis of crystal chemistry, Wells noted instances of the embedding of two enantiomorphic \textbf{srs} networks \cite{Wells:a01232, Wells:a01308}, namely \textbf{srs-c}, shown in Fig. \ref{fig:interpenetrating_enantiomorphic_srs_nets_extended}. In fact, this structure can be derived from a single \textbf{srs} network in its \textit{barycentric embedding}, shown in Fig. \ref{fig:srs_extended}, along with its \textit{dual network} \cite{Hyde2008ASH}. A \textit{barycentric embedding} is an embedding obtained by placing each vertex of the network in the barycentre of its neighbours \cite{equilibrium_placement}. It is computable with the Systre algorithm \cite{Delgado-Friedrichs:au5000}, and is intuitively known to avoid unnecessary entanglement. The \textit{dual} of a given network is the network carried by the dual of the \textit{tiling} of the network. A \textit{tiling} is a covering of space by generalised polyhedra. A simple example is the tiling of space by cubes sharing faces. A 3-dimensional tiling naturally carries a network formed by the edges and vertices of the tiles. For example, the network carried by the tiling with cubes is \textbf{pcu}, whose embedding is shown in Fig. \ref{fig:pcu_extended}. The dual of a tiling is obtained by placing new vertices in the centre of old tiles and by joining them with new edges passing through the faces of the old tiles. The new tiling carries a new network, which is the \textit{dual} of the original network \cite{DELGADOFRIEDRICHS20052533}. As already noted by Wells \cite{Wells:a01308}, embeddings of dual networks obtained via the aforementioned construction tend to exhibit elegant intergrowth between their two components \cite{Hyde2008ASH}, and are therefore likely to represent least tangled embeddings. In the case of \textbf{srs}, its dual is itself but of opposite handedness, which gives the embedding \textbf{srs-c}. Furthermore, the barycentric embedding of \textbf{srs} and the embedding \textbf{srs-c} of two dual \textbf{srs} networks, respectively shown in Fig. \ref{fig:srs_extended} and Fig. \ref{fig:interpenetrating_enantiomorphic_srs_nets_extended}, are also the highest-symmetry embeddings \cite{tesselate_decussate,Hyde2008ASH}. Moreover, they minimise the \textit{ropelength energy}, defined as the ratio of the total length of edges in a unit cell of the embedding to the diameter of the edges, which is computable via the PB-SONO algorithm \cite{myf_ideal_geo}. It is suggested in \cite{myf2011_entanglement_graphs} that an embedding minimising the ropelength energy is untangled. These observations regarding \textbf{srs} and \textbf{srs-c} may serve as a basis for characterising canonical least tangled embeddings of other networks.

In this paper, we formally define the least tangled embeddings of 3-periodic networks, which we refer to as \textit{ground states}, as well as a measure of entanglement complexity called the \textit{untangling number} of 3-periodic networks. The techniques employed here share some similarities with those described in \cite{andriamanalina2025untanglingnumber3periodictangles}, where the least tangled embeddings and the untangling number of 3-periodic filament entanglements (\textit{3-periodic tangles}) are defined. The paper is organised as follows. In Sect. \ref{sec:2}, we describe 3-periodic graphs and diagrams of their embeddings along the lines of the work in \cite{ANDRIAMANALINA2025109346} on 3-periodic tangles. In Sect. \ref{sec:least_tangled_embeddings}, we define the \textit{ground states}, which are the least tangled embeddings of 3-periodic graphs and illustrate the definition with examples and discussions. In Sect. \ref{sec:untangling_number}, we define the \textit{untangling number}, a measure of entanglement complexity. In Sect. \ref{sec:computability}, we discuss the computability of the untangling number. Finally, Sect. \ref{sec:conclusion} is dedicated to the conclusions of this work.

\section{Graph embeddings and their diagrams}\label{sec:2}
We recall that a \textit{graph} is an abstract collection of vertices connected by edges. An \textit{embedding} of a graph is a geometric realisation of the graph in a given space. In this paper, unless otherwise mentioned, we only consider graphs with 3-periodic embeddings in Euclidean space $\mathbb{R}^3$, that is, the embeddings are periodic along the three directions of space. In this section, we introduce the notions necessary for defining the concepts of a least tangled embedding and the untangling number.

\subsection{Unit cells, unit-diagrams and tridiagrams}
\begin{figure*}[hbtp]
    \centering
    \begin{subfigure}[b]{0.9\textwidth}
    \centering
        \includegraphics[width=0.585\textwidth]{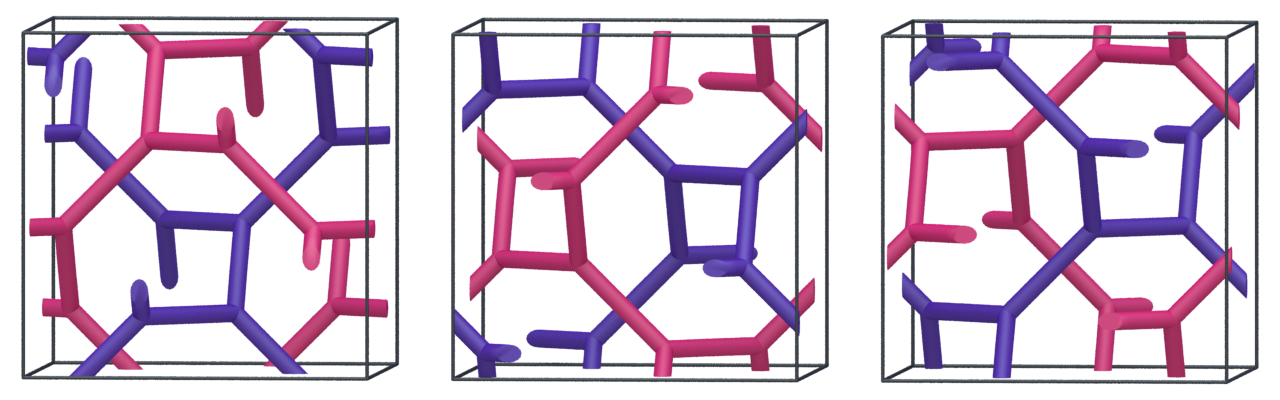}
    \caption{}
    \label{fig:interpenetrating_enantiomorphic_srs_nets_uc_xyz}
    \end{subfigure}
    \vskip\baselineskip
    \begin{subfigure}[b]{0.9\textwidth}
        \centering
    \includegraphics[width=0.18\textwidth]{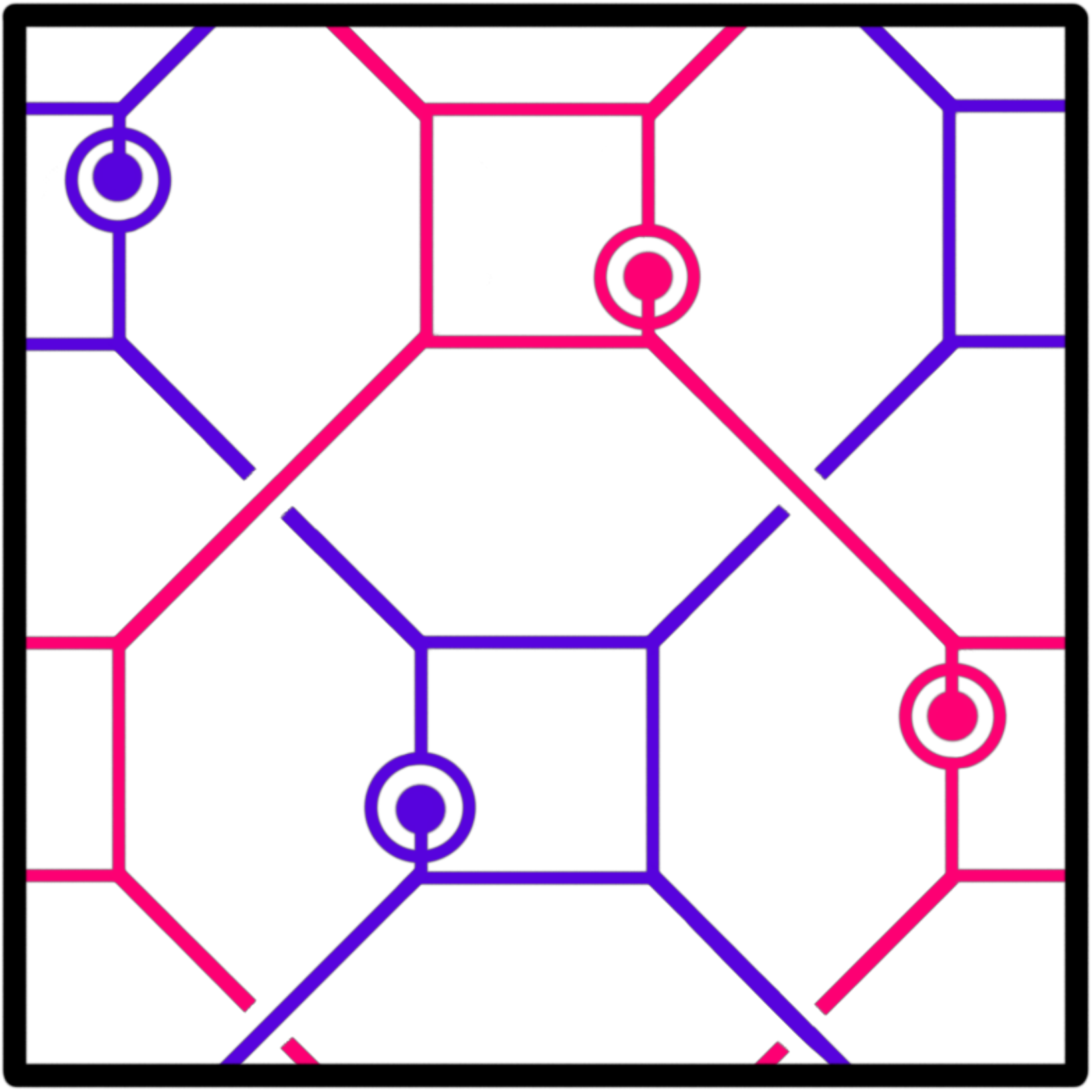}
        \hspace{0.1cm}
        \includegraphics[width=0.18\textwidth]{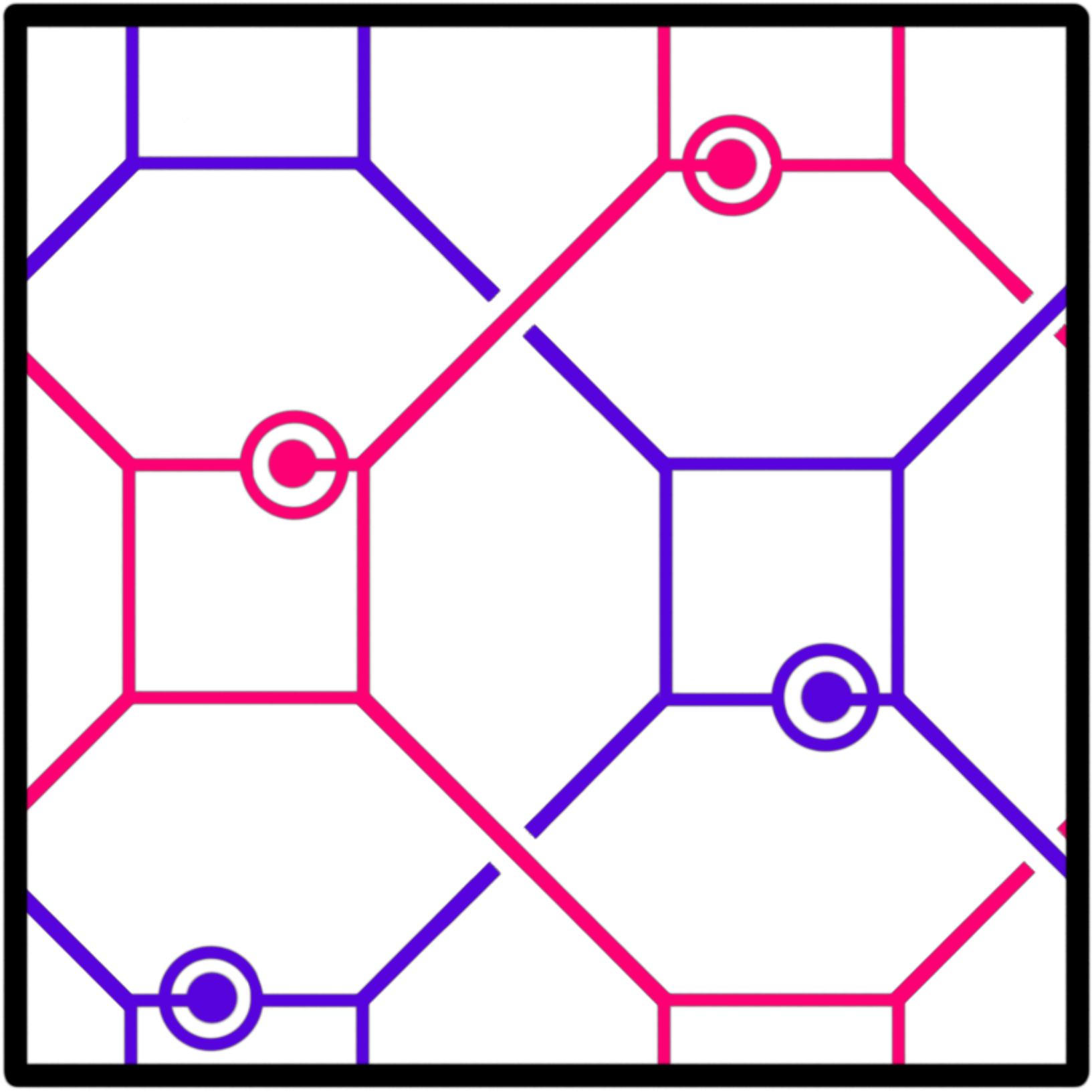}
        \hspace{0.1cm}
        \includegraphics[width=0.18\textwidth]{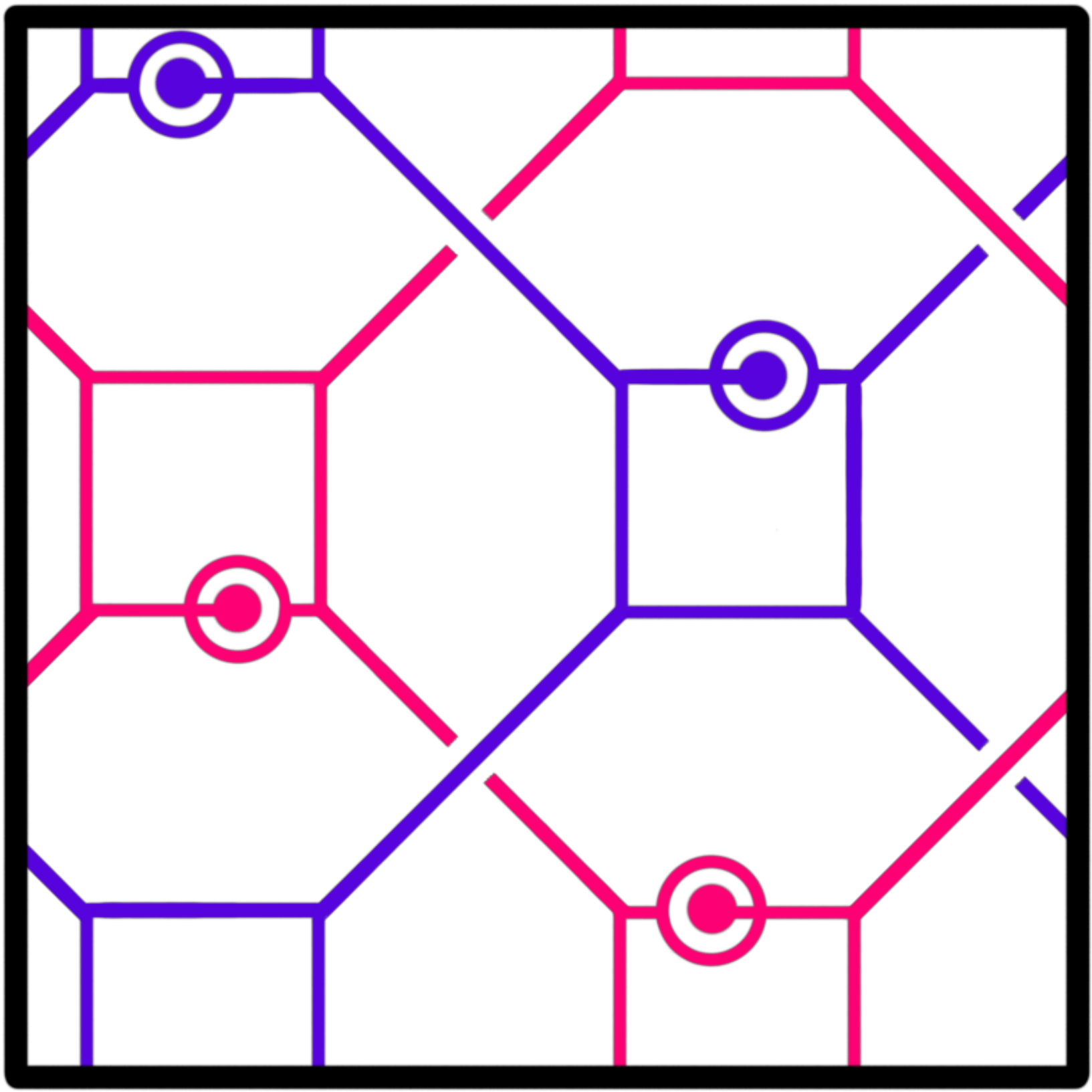}
        \caption{}
        \label{fig:interpenetrating_enantiomorphic_srs_nets_tridia}
    \end{subfigure}
    \caption{A unit cell of a graph embedding and a tridiagram obtained from it: (a)
    A unit cell of \textbf{srs-c} rotated in space so as to be viewed along the front, the top and the right sides. (b) A tridiagram, that is, a set of three diagrams obtained from projections from three non-coplanar axes of the unit cell shown in (a). An intersection with the front face is encoded using a thick dot and an intersection with the back face is encoded with an open circle.}
    \label{fig:interpenetrating_enantiomorphic_srs_nets_uc_and_tridia}
\end{figure*}

Consider a graph $\mathcal{K}$ and an embedding $K$ of $\mathcal{K}$ in Euclidean space $\mathbb{R}^3$.
A \textit{unit cell} of $K$ is a domain delimited by a parallelepiped with identified opposite faces that generates the whole embedding $K$ under translations that preserve the periodicity of the structure. An example is shown in Fig. \ref{fig:interpenetrating_enantiomorphic_srs_nets_uc_xyz}, which is a unit cell of \textbf{srs-c}, shown in Fig. \ref{fig:interpenetrating_enantiomorphic_srs_nets_extended}.
A parallelepiped delimiting a unit cell can always be rectified to a cube, which we do for simplicity in this paper.

A \textit{unit-diagram} or \textit{diagram} can be obtained from a unit cell by projecting it onto a square with identified edges. Each intersection in the projection is assigned a \textit{crossing type} indicating which strand passes over the other. In addition, a thick dot is used to represent an intersection with the front face of the unit cell, and an open circle is used to represent an intersection with the back face. Furthermore, a set of three diagrams obtained from the projections along the three vectors delimiting a unit cell constitutes a \textit{tridiagram}. Fig. \ref{fig:interpenetrating_enantiomorphic_srs_nets_tridia} shows a tridiagram obtained from the unit cell of Fig. \ref{fig:interpenetrating_enantiomorphic_srs_nets_uc_xyz}. Tridiagrams are necessary for converting the full spectrum of transformations in $3$-space into transformations in two-dimensional diagrams, where the crossing information of each diagram is different and equally important.
Unless otherwise mentioned, every tridiagram that we consider in this paper is obtained from the projections from the front, top, and right faces of a unit cell.

\subsection{Equivalence of embeddings}
Essentially, deforming graph embeddings in space does not change their entanglement complexity, provided that edges are not allowed to pass through one another or through vertices. This corresponds to the mathematical notion of \textit{ambient isotopy}. It is therefore natural, when considering entanglement, to regard embeddings up to their ambient isotopy classes. In this paper, we restrict ourselves to ambient isotopies realised within unit cells, preserving the periodicity of the networks. In a diagram, an ambient isotopy is reflected by a set of moves called $R$-moves, whose details are given in the Supplementary Information. From now on, the embeddings and the unit cells considered in this paper represent their ambient isotopy classes.

\subsection{Crossing number}
Diagrams are used to distinguish and classify embeddings via the use of \textit{invariants}, which are quantities that do not change under isotopy deformations made to a diagram. One of the simplest diagrammatic invariants is the \textit{crossing number}, which we present in the following in the context of 3-periodic graphs. The definition aligns with that given in \cite{ANDRIAMANALINA2025109346} for 3-periodic tangles.

\begin{figure*}[hbtp]
    \centering
    \begin{subfigure}[b]{0.9\textwidth}
    \centering
        \includegraphics[width=0.585\textwidth]{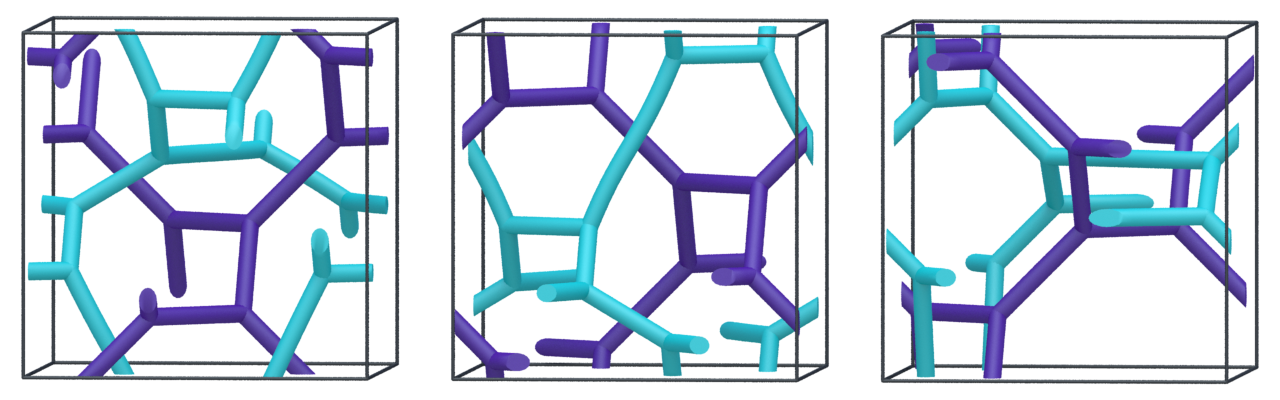}
    \caption{}
    \label{fig:interpenetrating_LH_srs_nets_uc_xyz}
    \end{subfigure}
    \vskip\baselineskip
    \begin{subfigure}[b]{0.9\textwidth}
        \centering
    \includegraphics[width=0.18\textwidth]{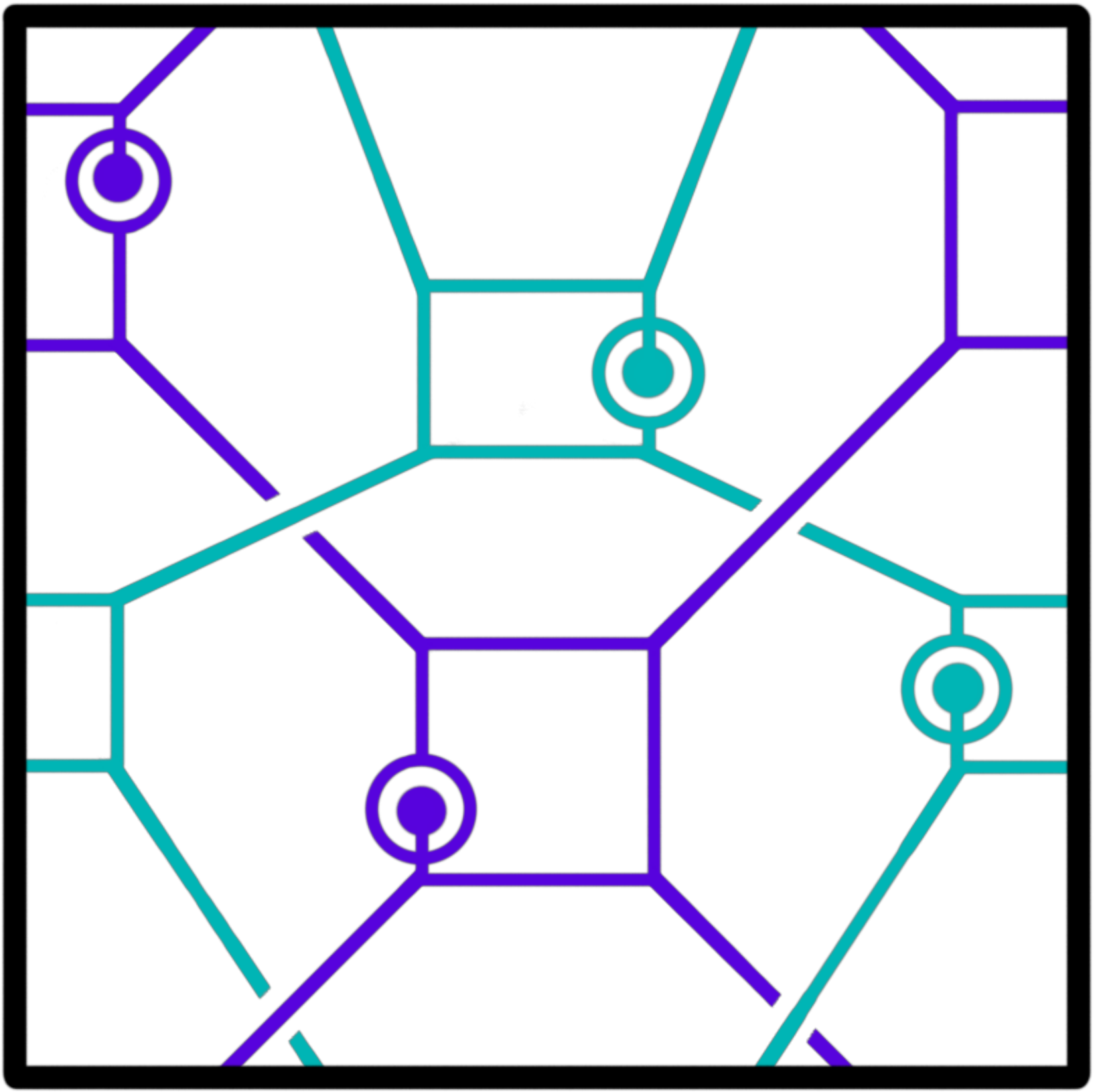}
        \hspace{0.1cm}
        \includegraphics[width=0.18\textwidth]{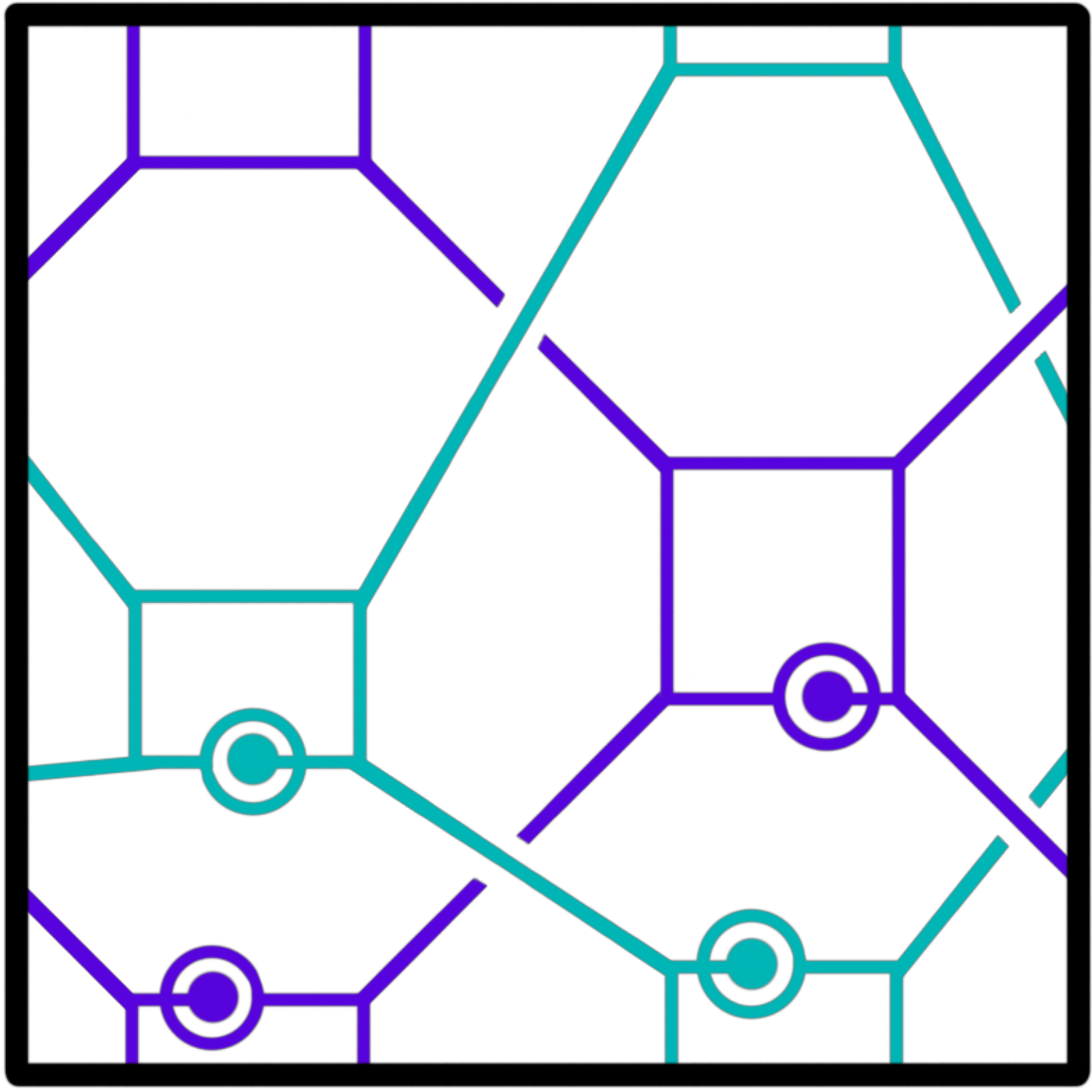}
        \hspace{0.1cm}
        \includegraphics[width=0.18\textwidth]{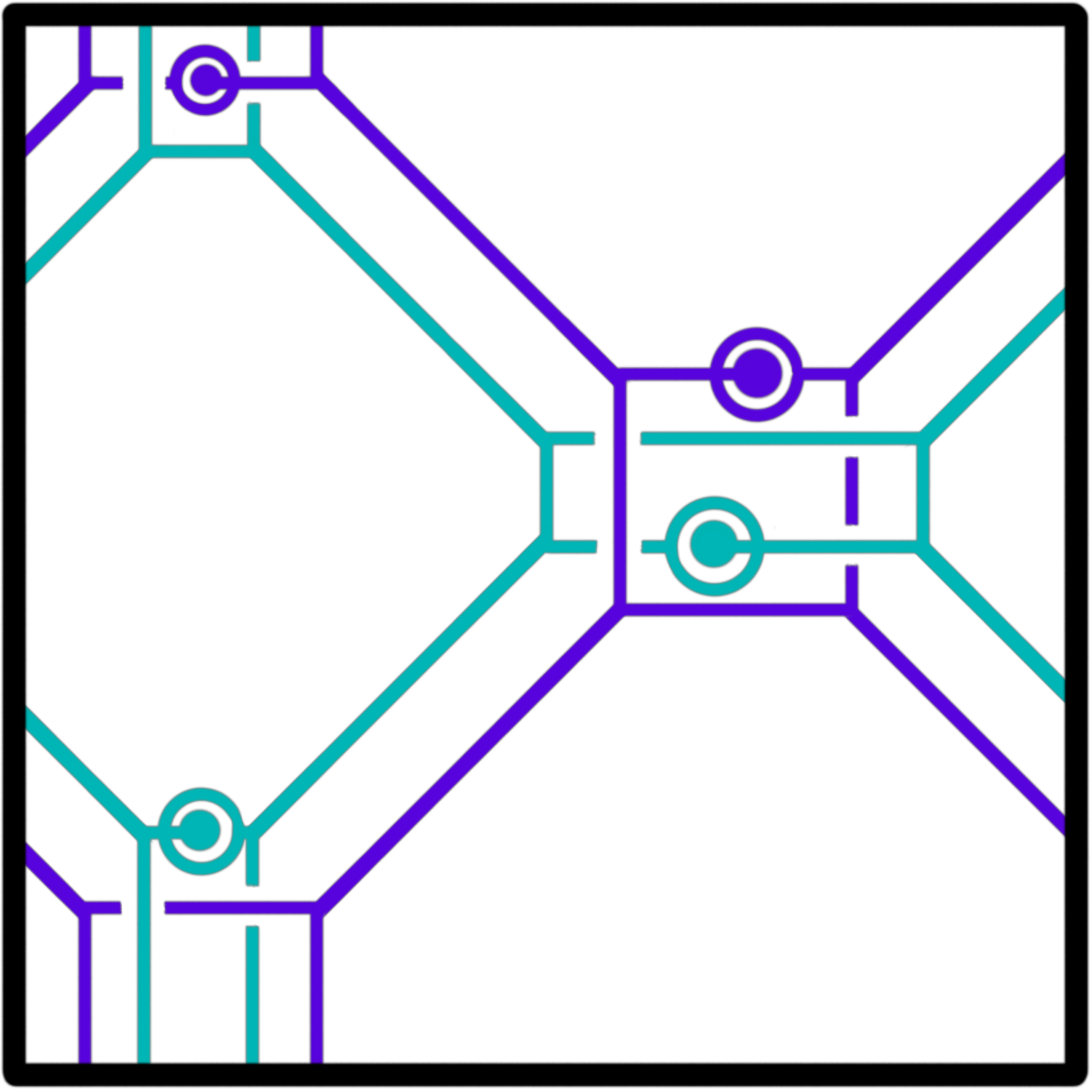}
        \caption{}
        \label{fig:interpenetrating_LH_srs_nets_tridia}
    \end{subfigure}
    \caption{An example of a minimal tridiagram: (a) A unit cell of \textbf{srs-c*}, shown in Fig. \ref{fig:interpenetrating_left_handed_srs_nets_extended}, viewed from the front, the top and the right sides. The embedding has been isotopically deformed to minimise the number of crossings. (b) A minimal tridiagram obtained from the unit cell in (a), showing that the minimum crossing number triplet is $(4,4,8)$.}
    \label{fig:interpenetrating_LH_srs_nets_uc_and_tridia}
\end{figure*}

Consider an embedding $K$ of a graph $\mathcal{K}$ and a unit cell $U$ of $K$. To a tridiagram $T = \lbrace D_1,D_2,D_3 \rbrace$ associated to $U$, one can associate a triplet $(a,b,c)$, where $a$, $b$ and $c$ are the numbers of crossings of the diagrams $D_1$, $D_2$ and $D_3$, respectively. Such a triplet is referred to as a \textit{triplet of crossings}. The \textit{crossing number of $K$ with respect to $U$}, denoted by $c(K,U)$, is the minimum of $c(T) = a^2 + b^2 +c^2$ among all tridiagrams $T$ representing the isotopy class of the unit cell $U$.
A tridiagram and its triplet of crossings realising $c(K,U)$ are respectively called a \textit{minimal tridiagram of $K$ with respect to $U$}, and a \textit{minimum crossing number triplet with respect to $U$}.
The \textit{crossing number} of $K$ is defined as the minimum of all $c(K,U)$ among all unit cells $U$ of $K$, and it is an invariant.
To illustrate this, see the example given in Fig. \ref{fig:interpenetrating_enantiomorphic_srs_nets_uc_and_tridia} where the minimum crossing number triplet of \textbf{srs-c} with respect to the unit cell shown in Fig. \ref{fig:interpenetrating_enantiomorphic_srs_nets_uc_xyz}, is $(4,4,4)$ and its crossing number is $48$, as given by the minimal tridiagram shown in Fig. \ref{fig:interpenetrating_enantiomorphic_srs_nets_tridia}. Another example is given in Fig. \ref{fig:interpenetrating_LH_srs_nets_uc_and_tridia}, where we display a unit cell and the associated minimal tridiagram of \textbf{srs-c*}, the embedding displayed in Fig. \ref{fig:interpenetrating_left_handed_srs_nets_extended}. The minimum crossing number triplet is $(4,4,8)$ and the crossing number is $96$. Notice how the structure has been isotopically deformed to obtain a minimal tridiagram.

\section{Least tangled embeddings}\label{sec:least_tangled_embeddings}
Using the framework established in Sect. \ref{sec:2}, namely the crossing diagrams, we are now able to characterise the least tangled embeddings of 3-periodic graphs that we call ground states.\\

\begin{figure*}[hbtp]
    \centering
    \begin{subfigure}[b]{0.18\textwidth}
        \centering
        \includegraphics[width=0.9\textwidth]{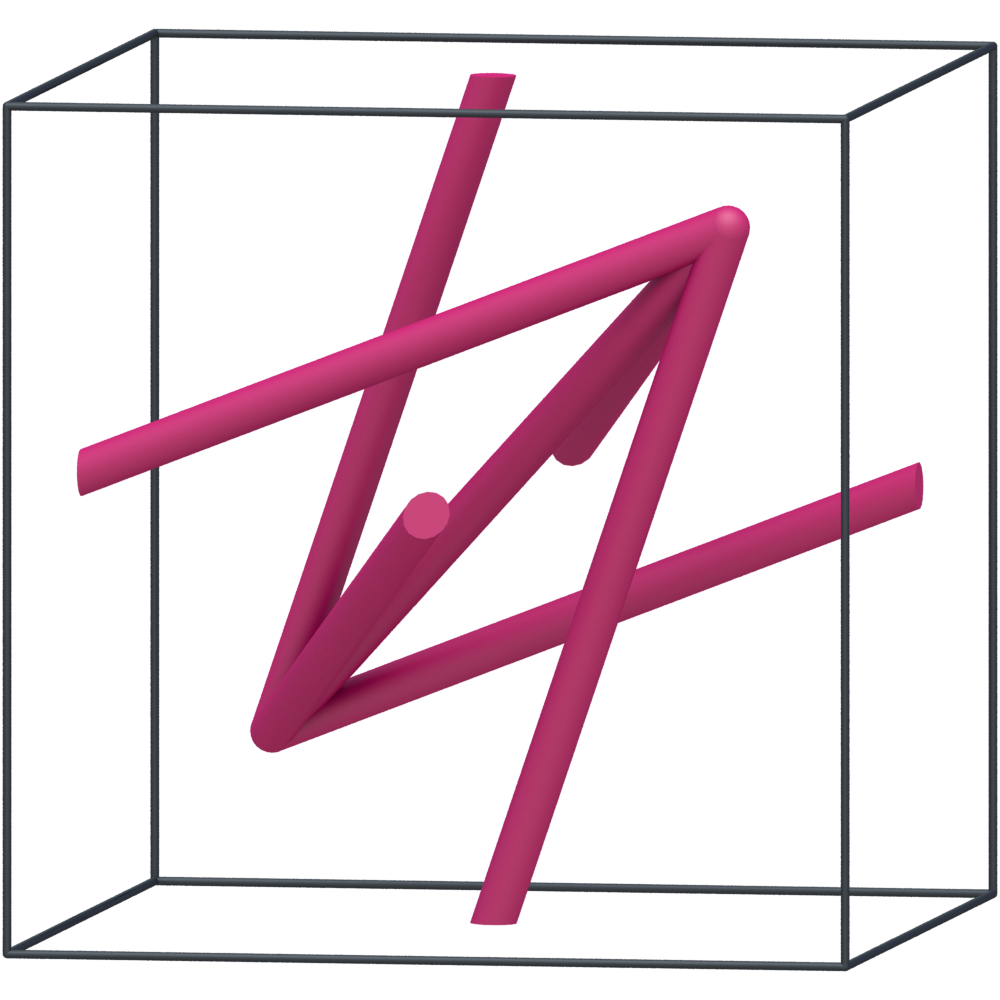}
        \caption{}
        \label{fig:dia-z_uc}
    \end{subfigure}
    \begin{subfigure}[b]{0.18\textwidth}
        \centering
        \includegraphics[width=0.9\textwidth]{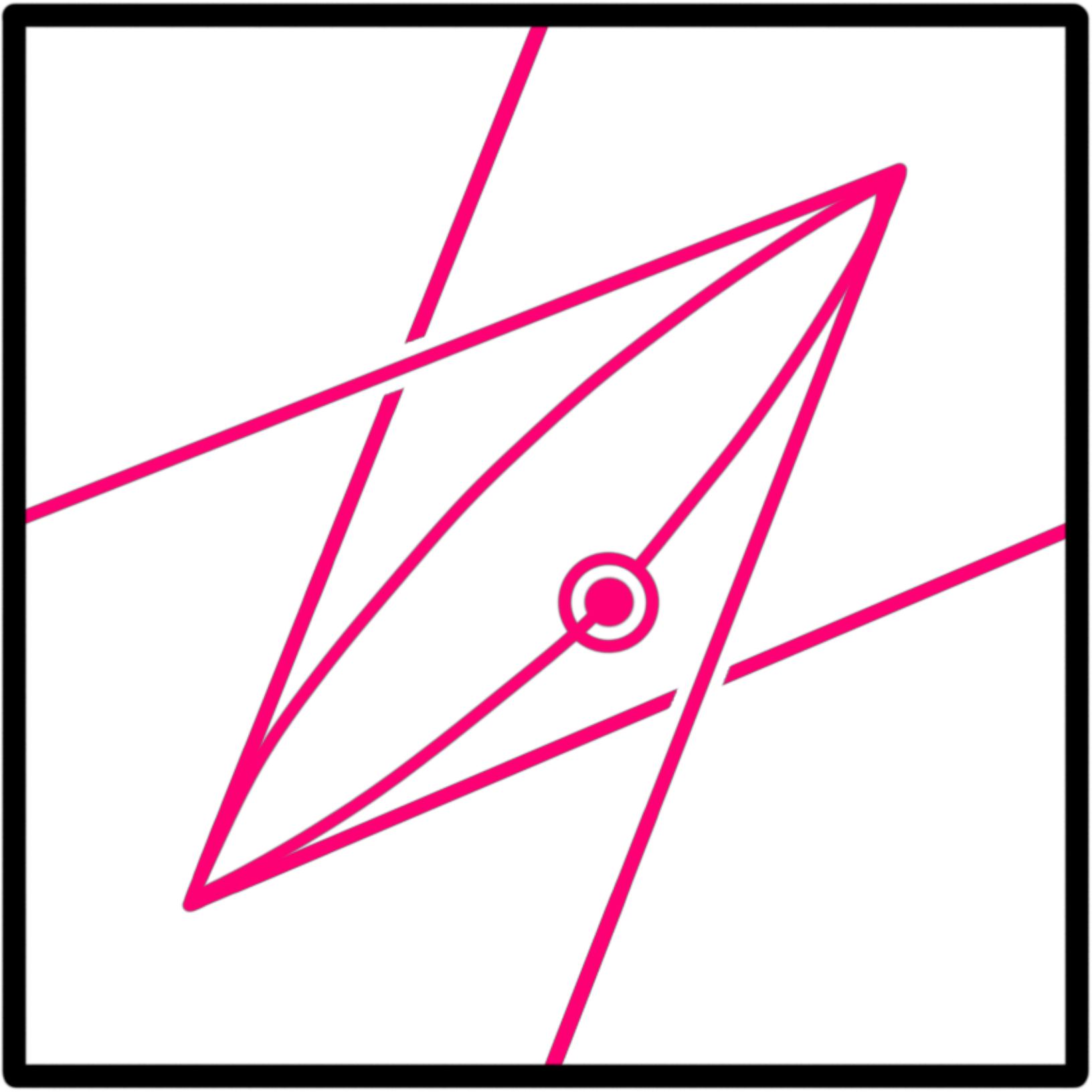}
        \caption{}
        \label{fig:dia-z_untangling_dia_1}
    \end{subfigure}
    \begin{subfigure}[b]{0.18\textwidth}
        \centering
        \includegraphics[width=0.9\textwidth]{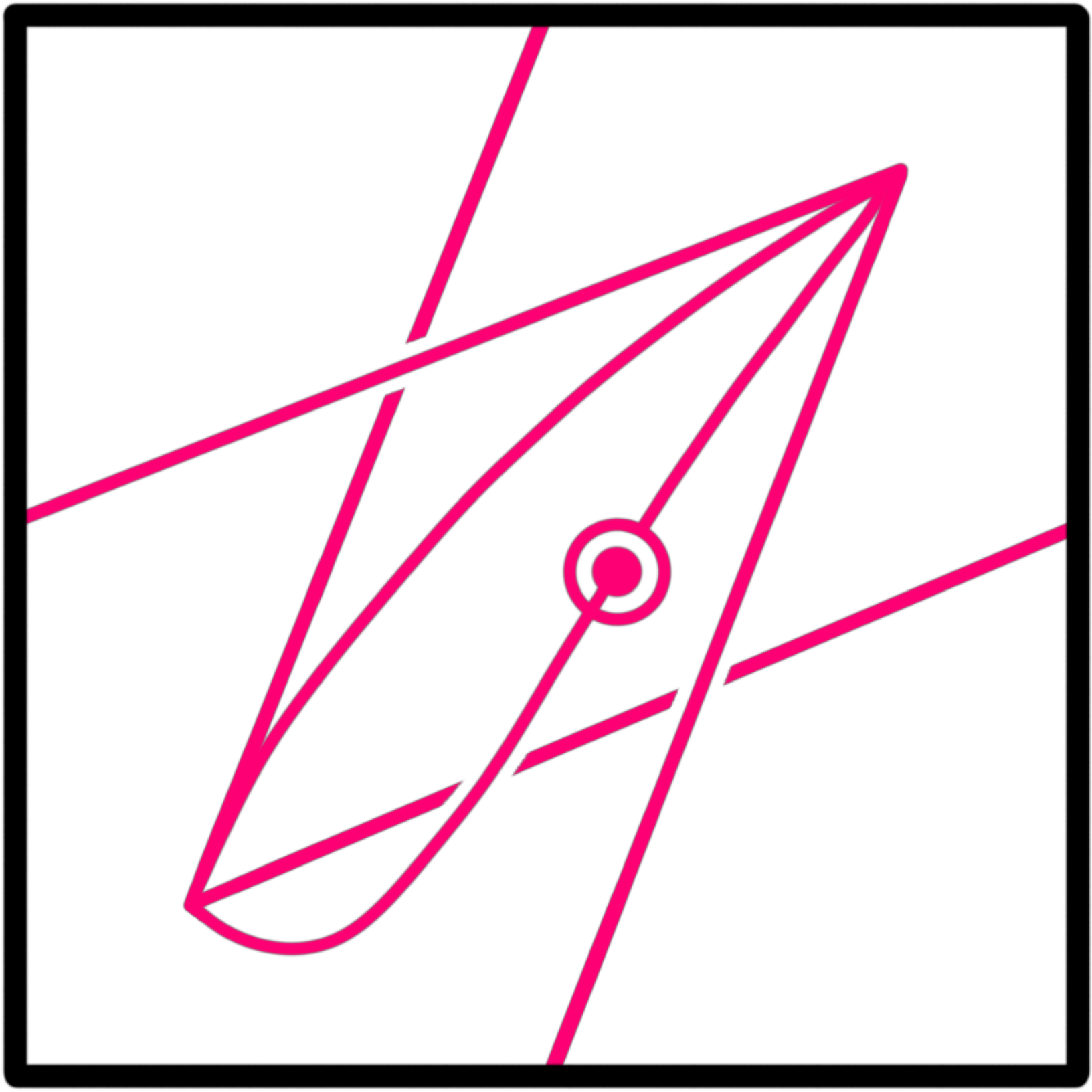}
        \caption{}
        \label{fig:dia-z_untangling_dia_2}
    \end{subfigure}
    \begin{subfigure}[b]{0.18\textwidth}
        \centering
        \includegraphics[width=0.9\textwidth]{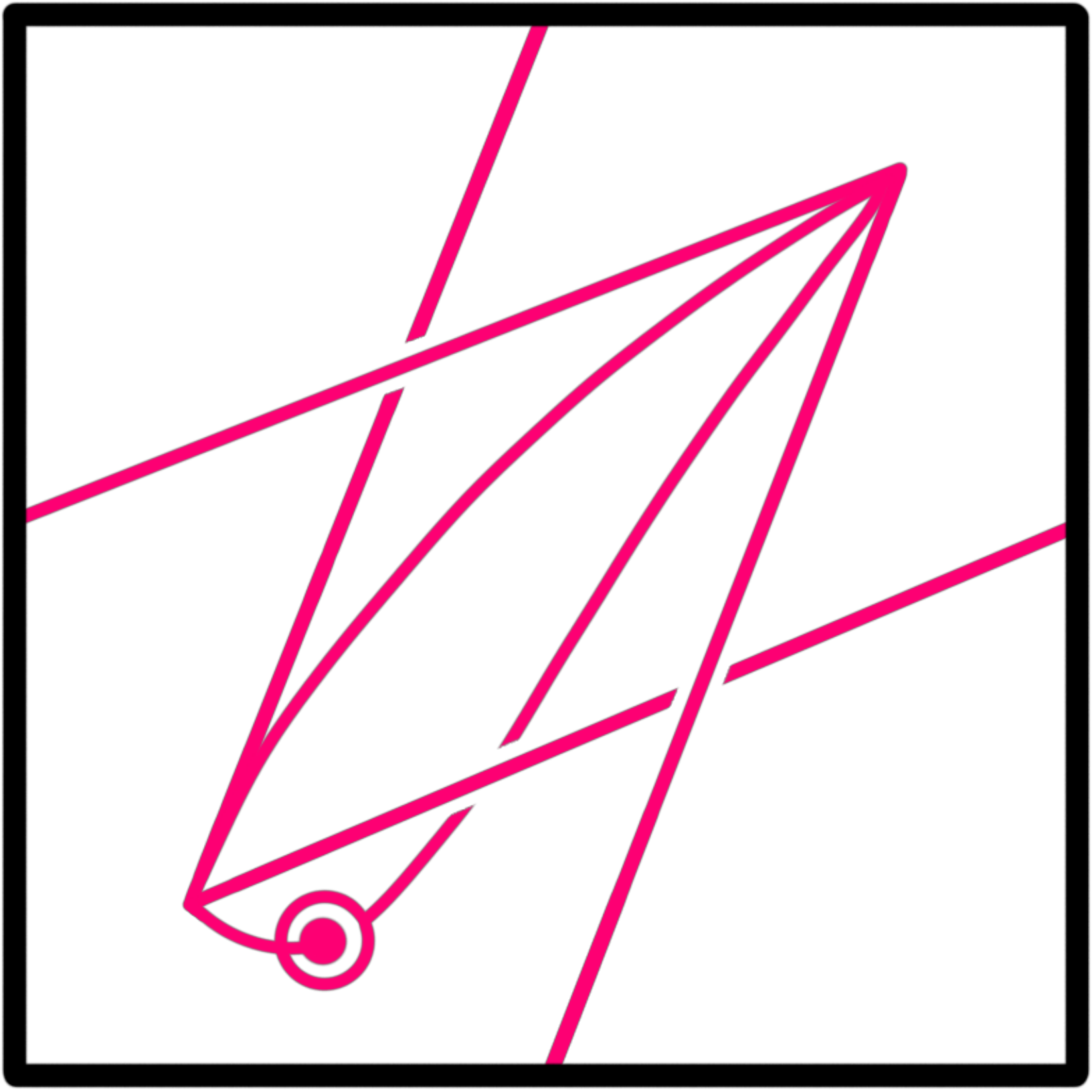}
        \caption{}
        \label{fig:dia-z_untangling_dia_3}
    \end{subfigure}
    \begin{subfigure}[b]{0.18\textwidth}
        \centering
        \includegraphics[width=0.9\textwidth]{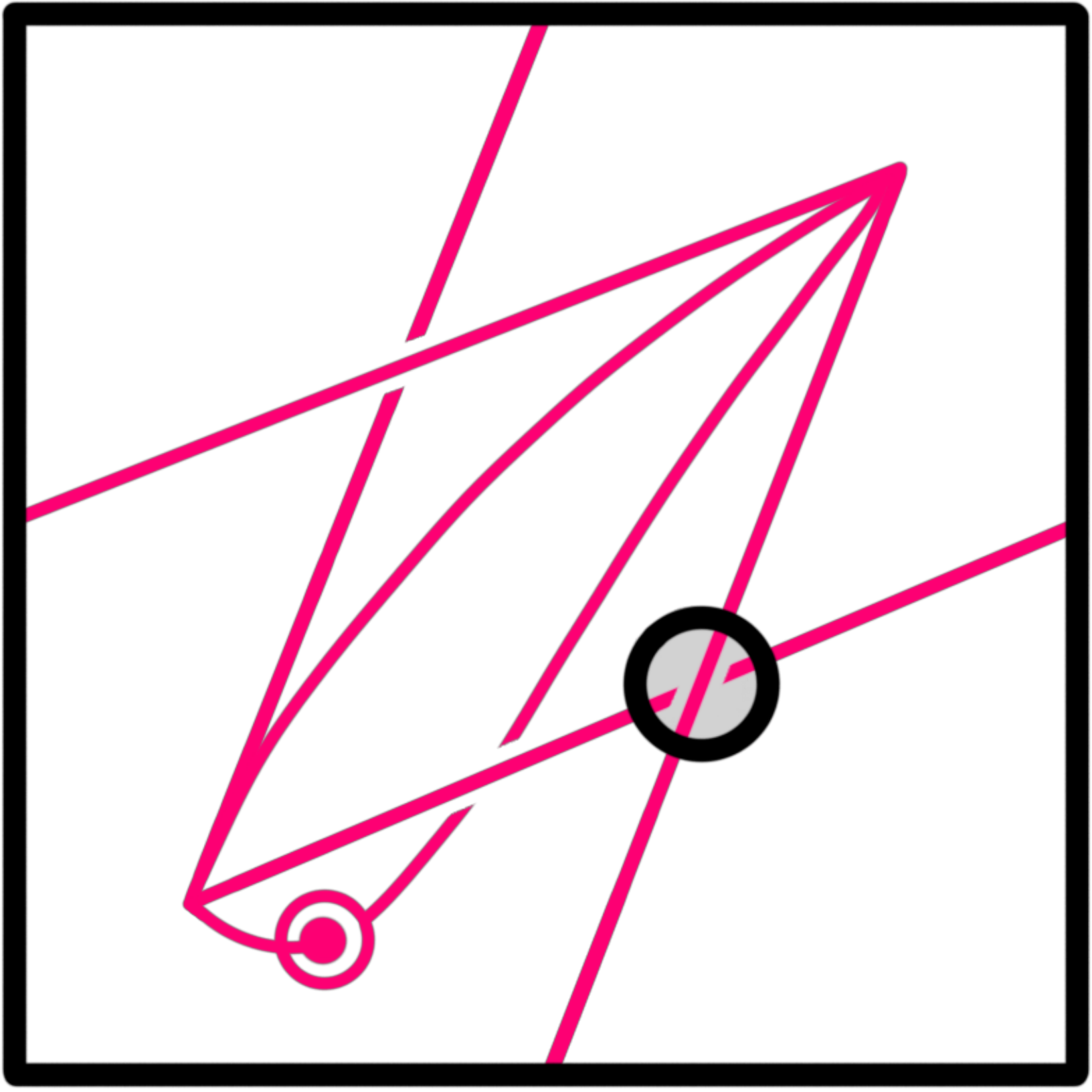}
        \caption{}
        \label{fig:dia-z_untangling_dia_4}
    \end{subfigure}
    
    \vskip\baselineskip
    
    \begin{subfigure}[b]{0.18\textwidth}
        \centering
        \includegraphics[width=0.9\textwidth]{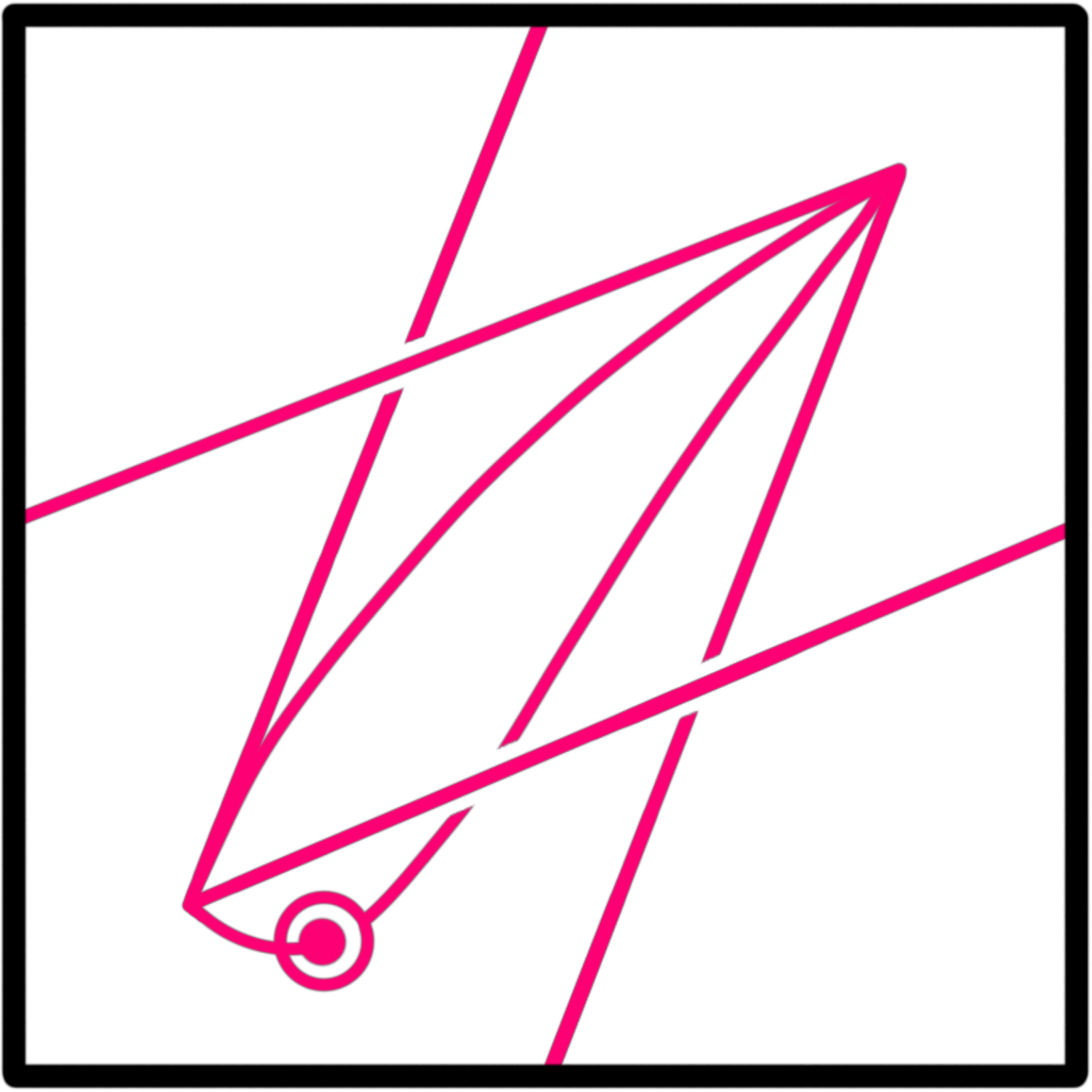}
        \caption{}
        \label{fig:dia-z_untangling_dia_5}
    \end{subfigure}
    \begin{subfigure}[b]{0.18\textwidth}
        \centering
        \includegraphics[width=0.9\textwidth]{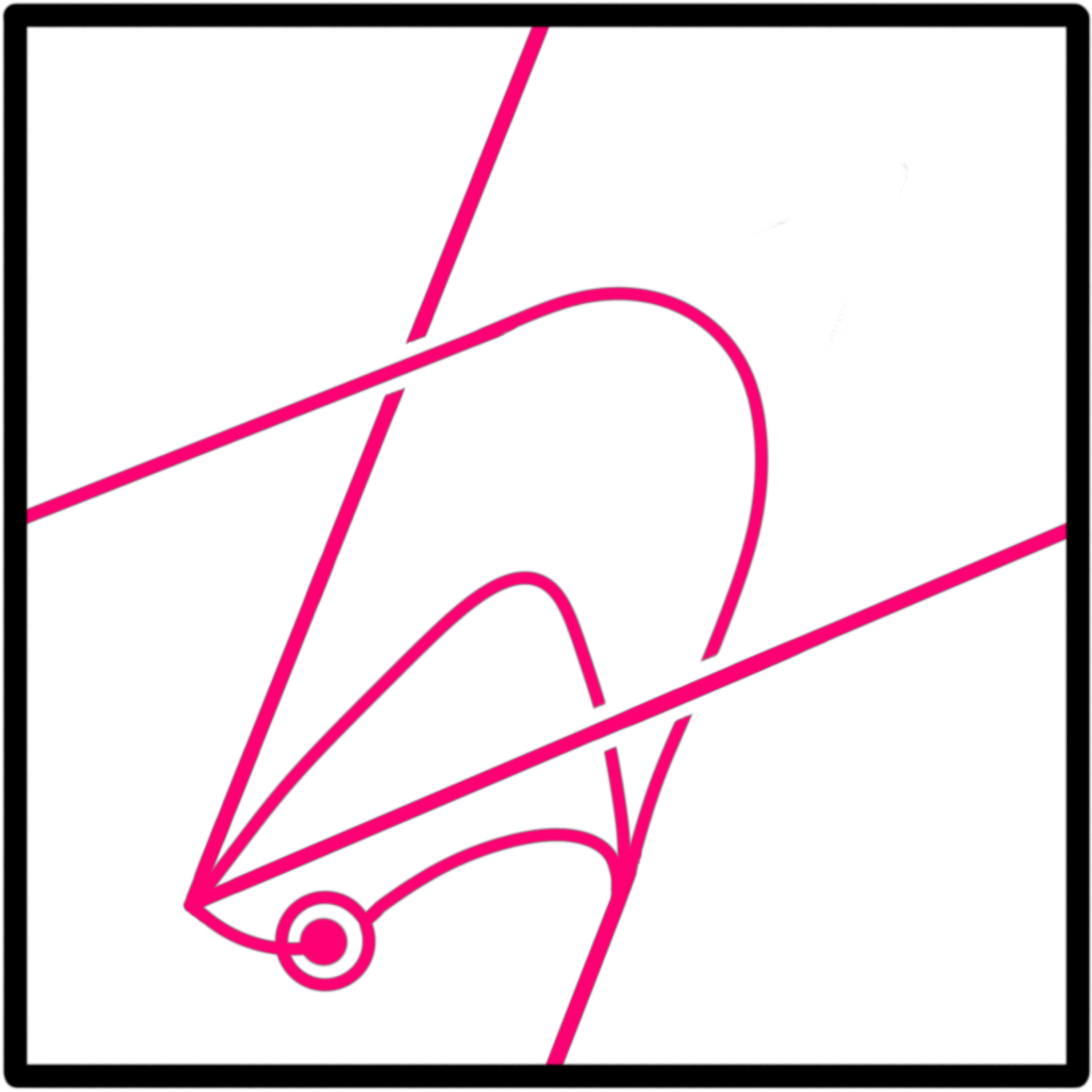}
        \caption{}
        \label{fig:dia-z_untangling_dia_6}
    \end{subfigure}
    \begin{subfigure}[b]{0.18\textwidth}
        \centering
        \includegraphics[width=0.9\textwidth]{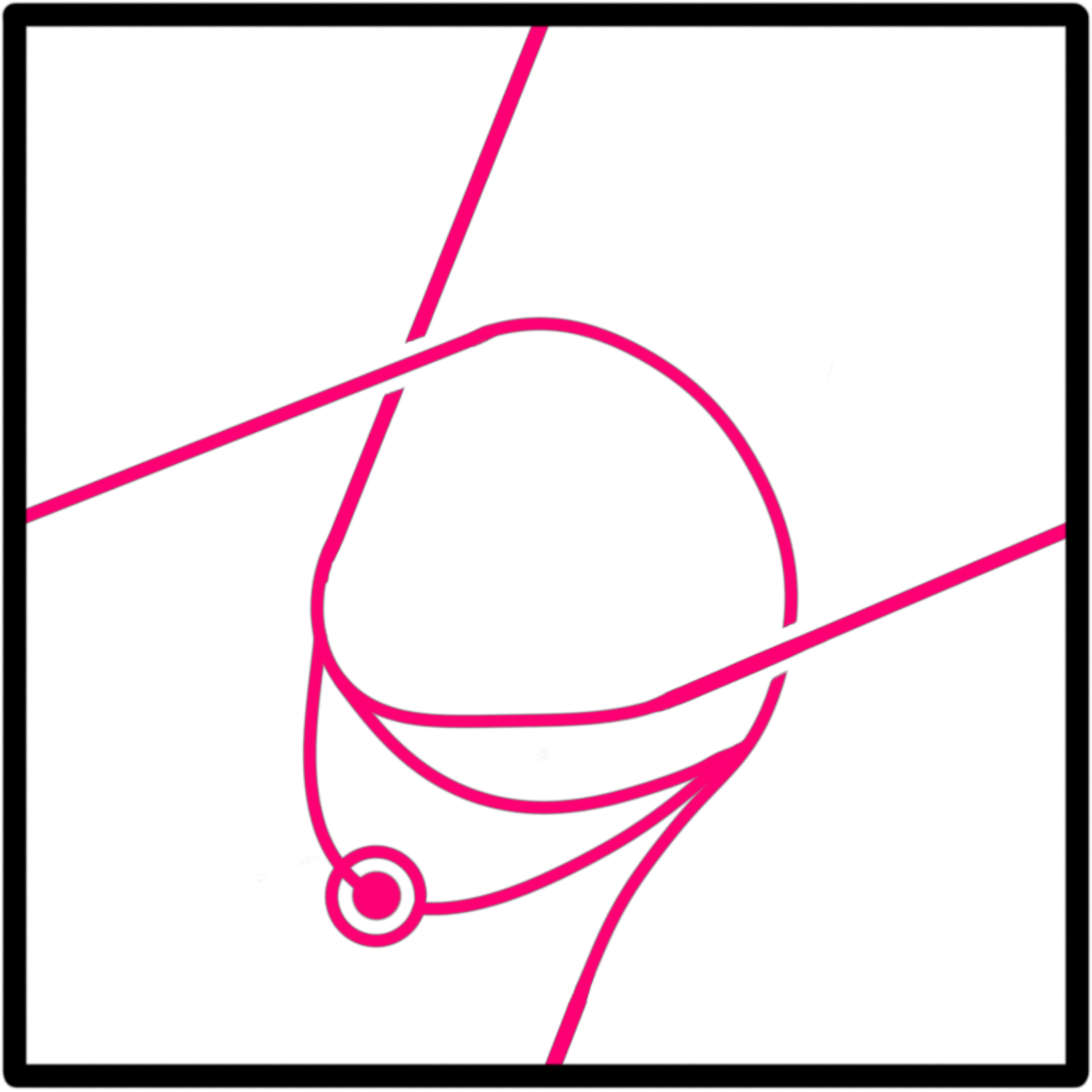}
        \caption{}
        \label{fig:dia-z_untangling_dia_7}
    \end{subfigure}
    \begin{subfigure}[b]{0.18\textwidth}
        \centering
        \includegraphics[width=0.9\textwidth]{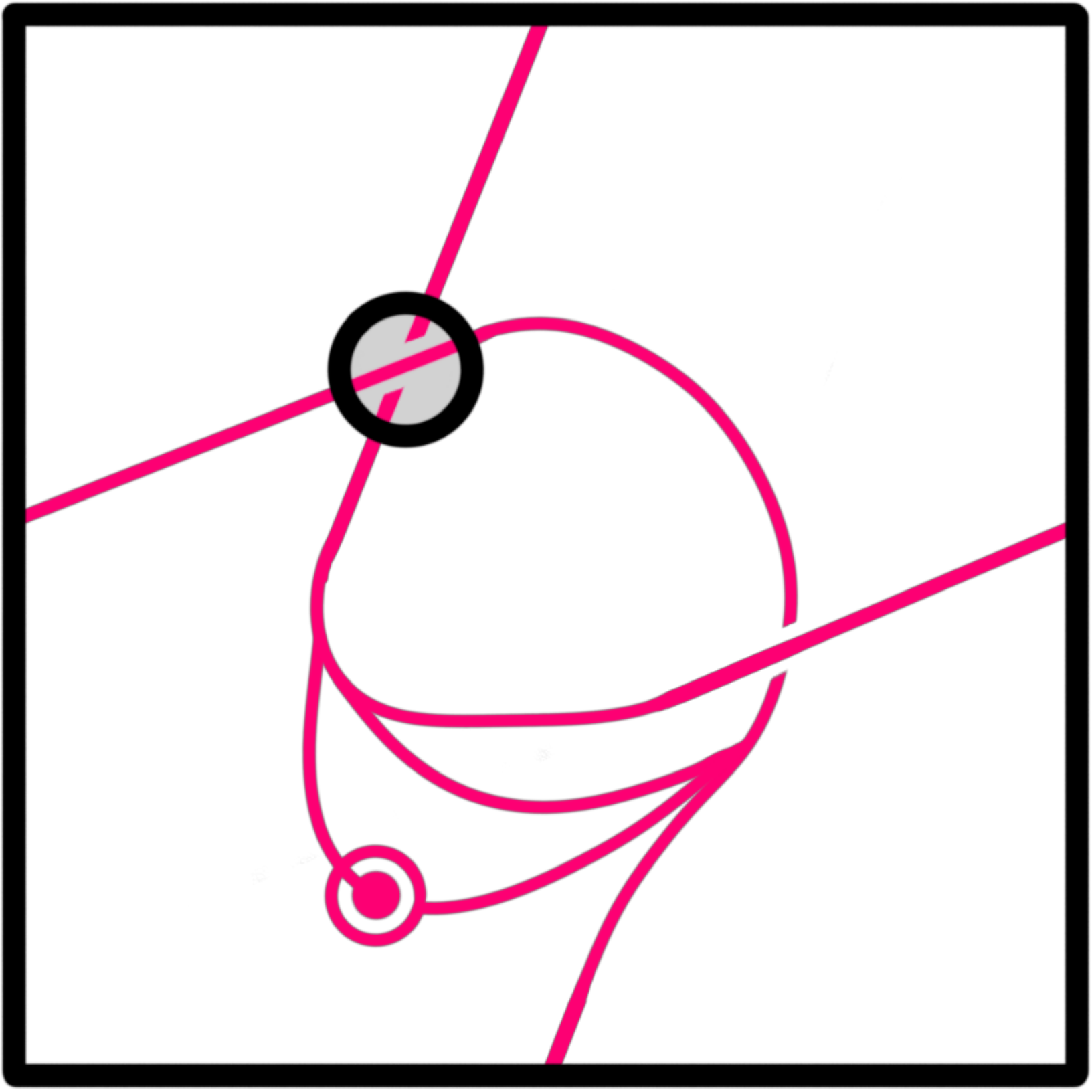}
        \caption{}
        \label{fig:dia-z_untangling_dia_8}
    \end{subfigure}
    \begin{subfigure}[b]{0.18\textwidth}
        \centering
        \includegraphics[width=0.9\textwidth]{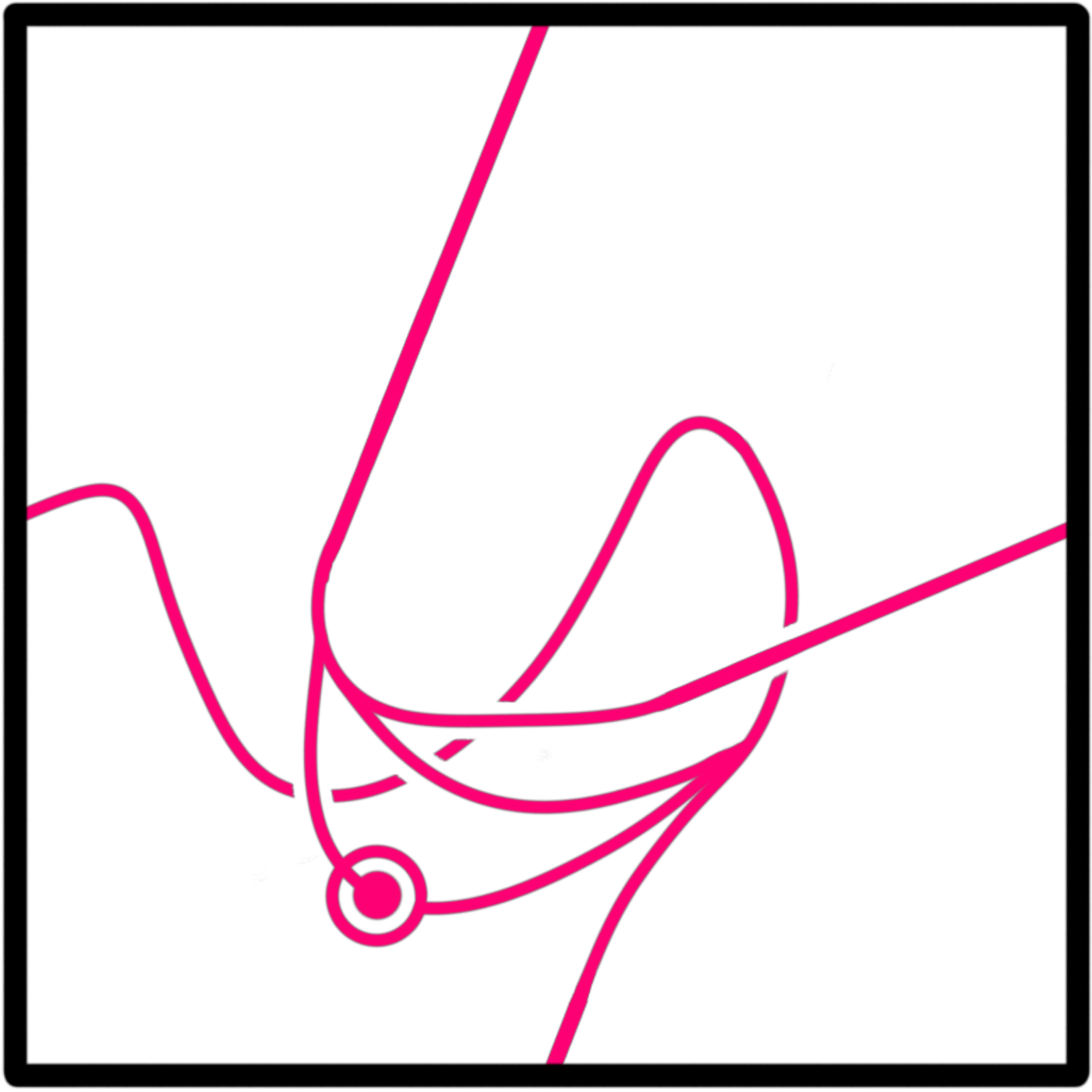}
        \caption{}
        \label{fig:dia-z_untangling_dia_9}
    \end{subfigure}

    \vskip\baselineskip
    
    \begin{subfigure}[b]{0.18\textwidth}
        \centering
        \includegraphics[width=0.9\textwidth]{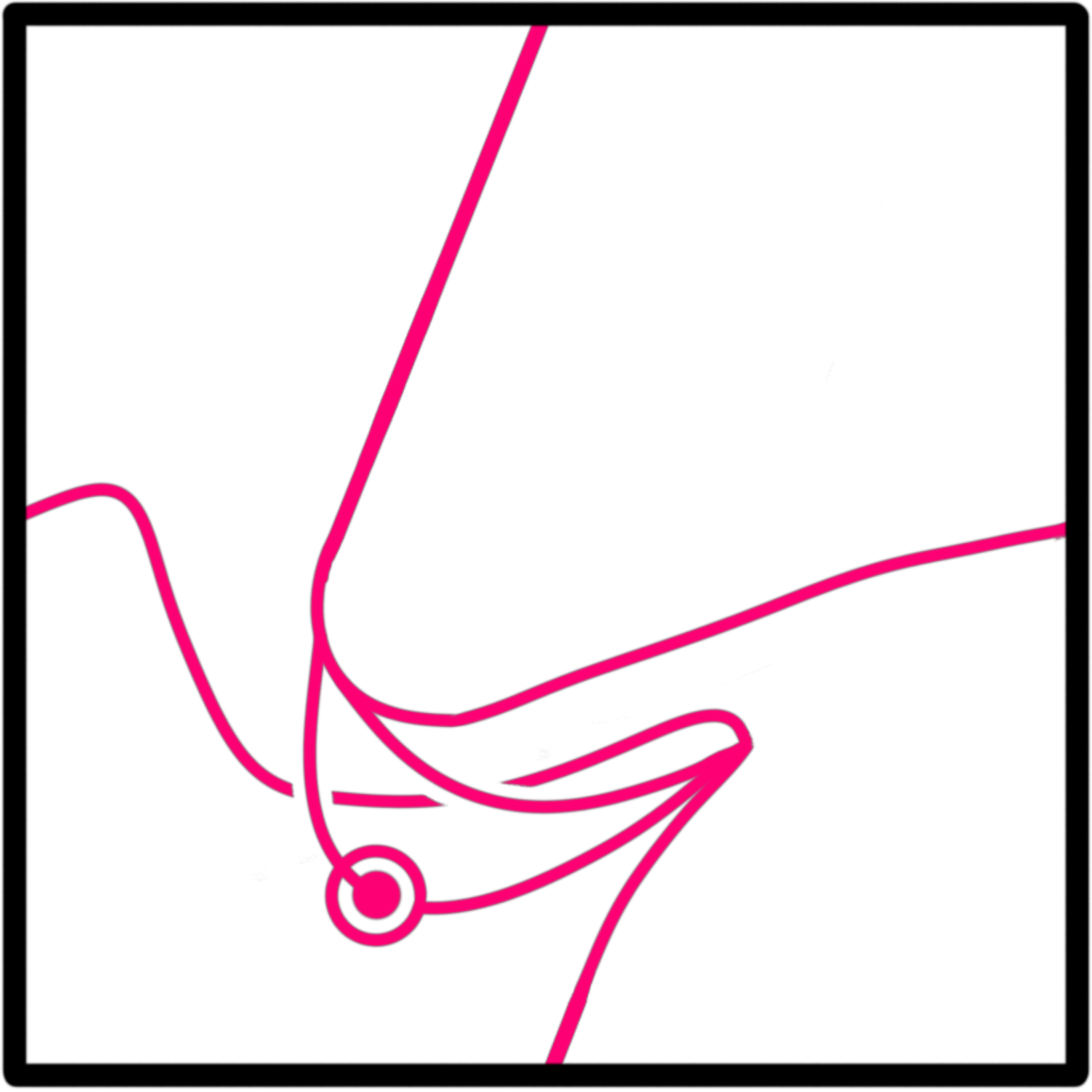}
        \caption{}
        \label{fig:dia-z_untangling_dia_10}
    \end{subfigure}
    \begin{subfigure}[b]{0.18\textwidth}
        \centering
        \includegraphics[width=0.9\textwidth]{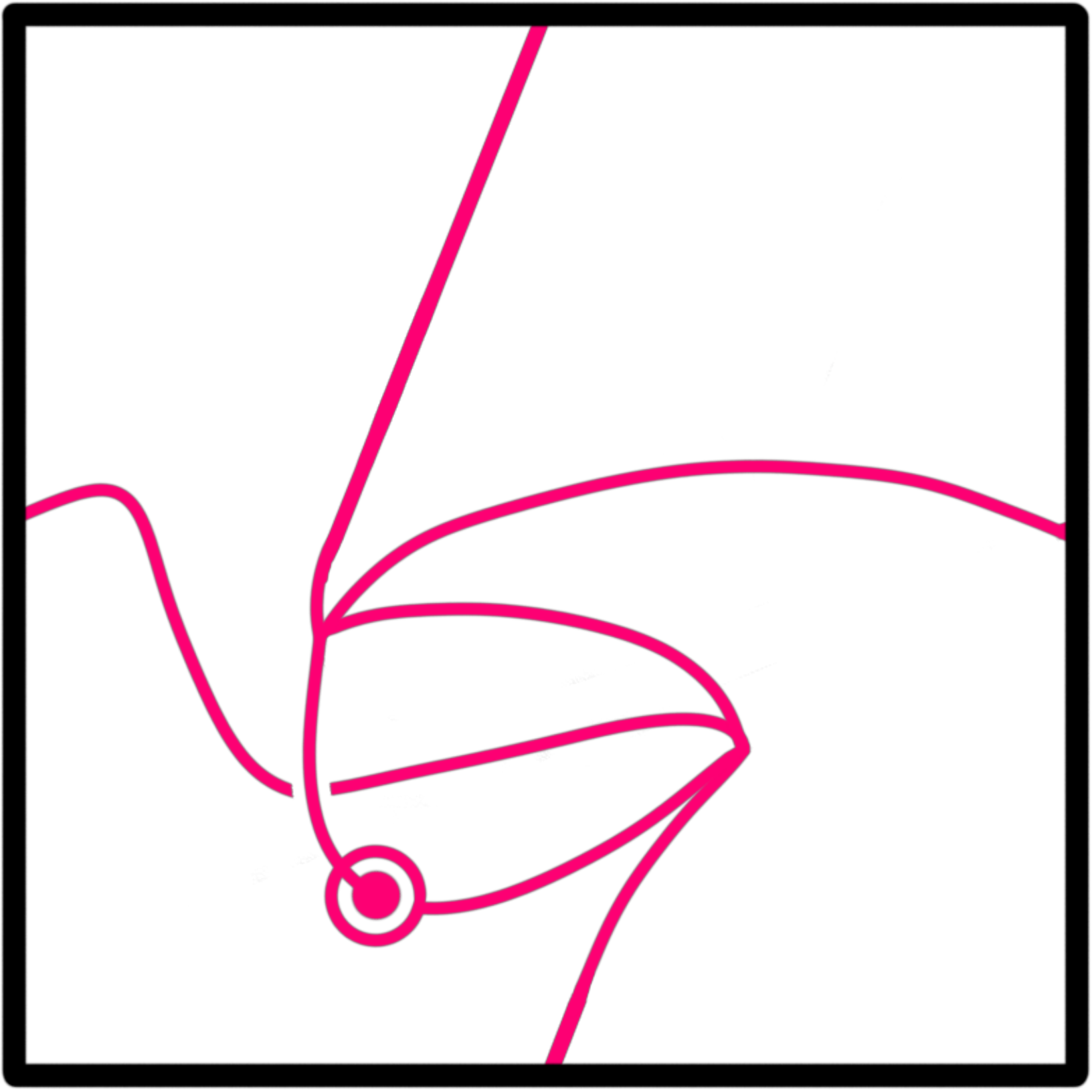}
        \caption{}
        \label{fig:dia-z_untangling_dia_11}
    \end{subfigure}
    \begin{subfigure}[b]{0.18\textwidth}
        \centering
        \includegraphics[width=0.9\textwidth]{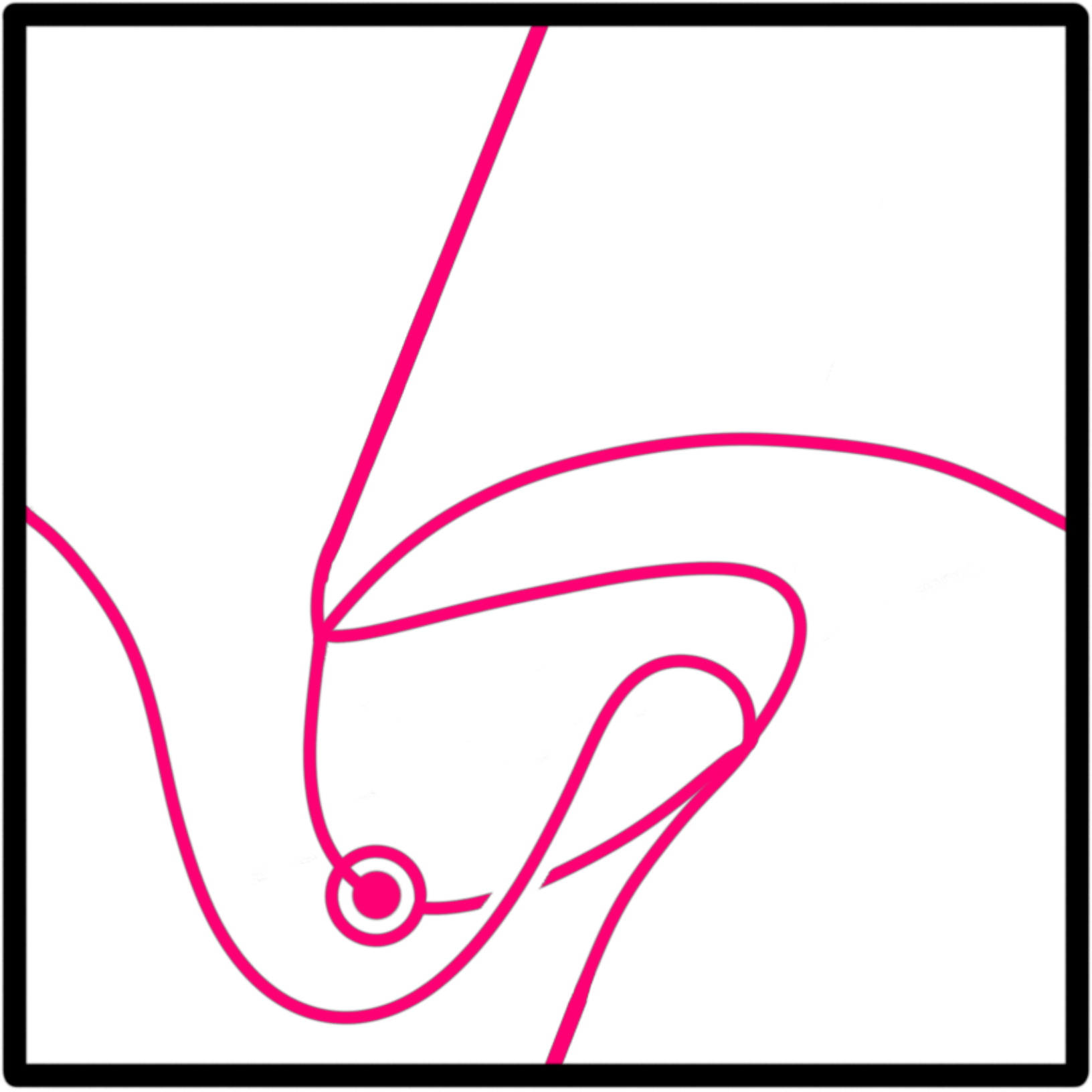}
        \caption{}
        \label{fig:dia-z_untangling_dia_12}
    \end{subfigure}
    \begin{subfigure}[b]{0.18\textwidth}
        \centering
        \includegraphics[width=0.9\textwidth]{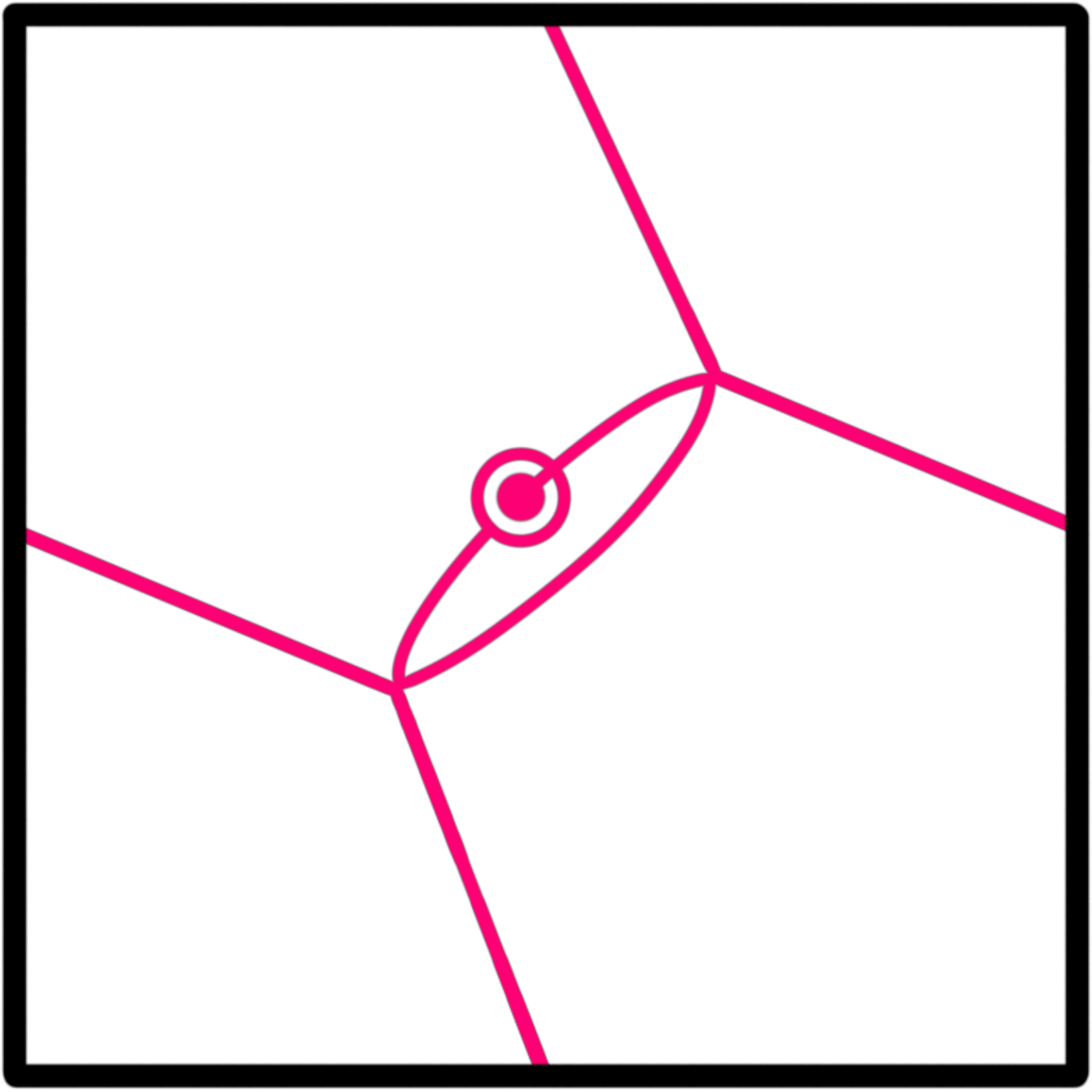}
        \caption{}
        \label{fig:dia-z_untangling_dia_13}
    \end{subfigure}
    \begin{subfigure}[b]{0.18\textwidth}
        \centering
        \includegraphics[width=0.9\textwidth]{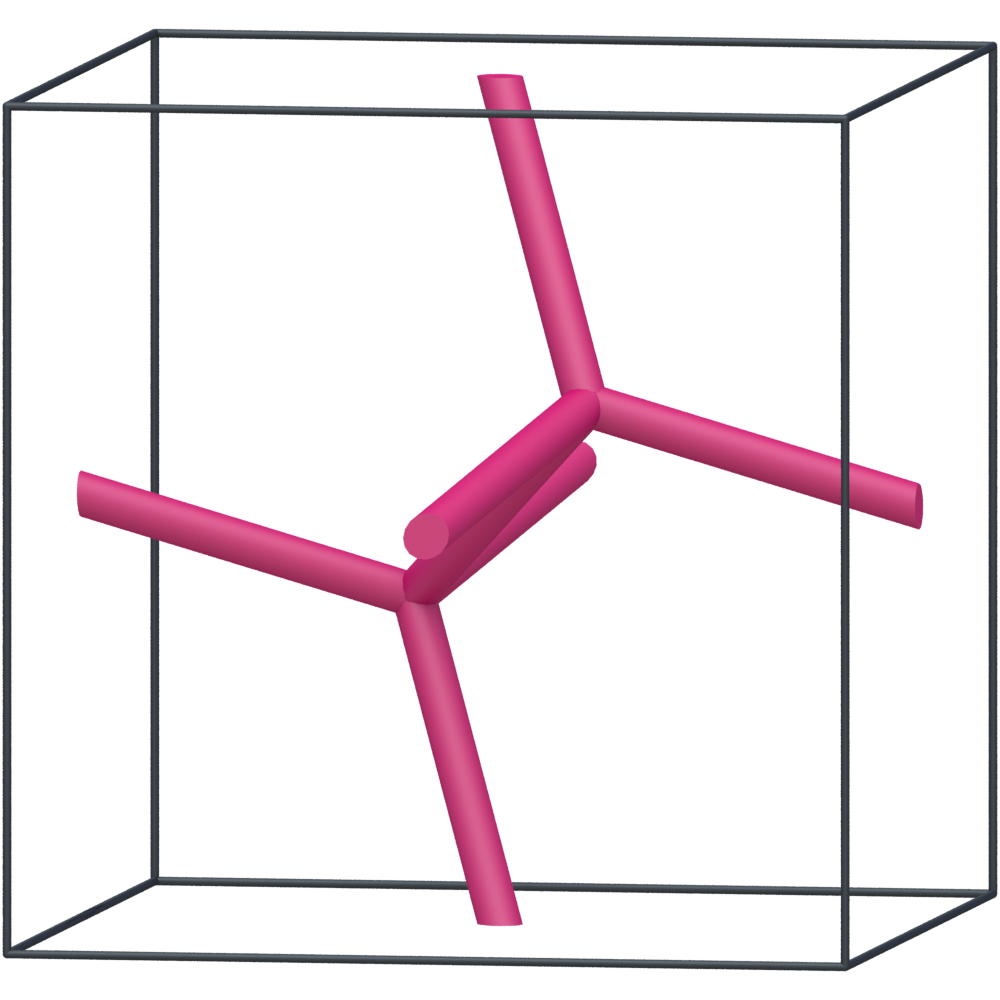}
        \caption{}
        \label{fig:dia_primitive_uc}
    \end{subfigure}
    \caption{Unit cells connected via finitely many crossing changes: (a) The primitive unit cell of \textbf{dia-z}. (o) The primitive unit cell of \textbf{dia}. The two unit cells are connected by two crossing changes applied on the crossings highlighted by the black circles in the diagrams shown in (e) and (i). Note that the unit cells have been sheared into cubes, which breaks the symmetries of the structures.}
    \label{fig:untangling_dia-z}
\end{figure*}

Classical links are usually distinguished by their number of components. For a given number $m$, links with $m$ components are connected via crossing changes to the simplest such link, namely the disjoint union of $m$ circles, which in particular minimises the crossing number. Recall the example shown in Fig. \ref{fig:hopf_link}, which illustrates links with two components. The concept of a crossing change can naturally be applied to diagrams of 3-periodic graphs, as shown in the example given in Fig. \ref{fig:untangling_dia-z}. Therefore, by analogy with classical links, we shall regroup all embeddings of a 3-periodic graph connected by crossing changes into a family, within which the least tangled states would be those that minimise the crossing number.

Consider a graph $\mathcal{K}$, an embedding $K$ of $\mathcal{K}$ and a unit cell $U$ of $K$. Let us call $\mathcal{U}$-family of $K$ with respect to $U$ and denote it by $\mathcal{U}(K,U)$ the family of all graph embeddings that possess diagrams connected to a diagram of $U$ by finitely many crossing changes. For example, if $U$ is the unit cell shown in Fig. \ref{fig:dia-z_uc}, which is the primitive unit cell of \textbf{dia-z}, then the primitive unit cell of \textbf{dia}, displayed in Fig. \ref{fig:dia_primitive_uc} belongs to the family $\mathcal{U}(K,U)$. The converse is also true, that is, if $U$ is the primitive unit cell of \textbf{dia}, then the primitive unit cell of \textbf{dia-z} belongs to its family $\mathcal{U}(K,U)$. By contrast, \textbf{srs-c} and \textbf{srs-c*} are examples of structures that do not belong to the same $\mathcal{U}$-family. Indeed, even though they both are embeddings of the same graph, that is, a graph comprising two \textbf{srs} networks, their unit cells are not connected by finitely many crossing changes. In fact, there is no transformation in the ambient space $\mathbb{R}^3$ preserving the underlying abstract graph that transforms one into the other, even when allowing the edges to pass through each other. Within the family $\mathcal{U}(K,U)$, let us call $\mathcal{G}$-family of $K$ with respect to $U$ and denote it by $\mathcal{G}(K,U)$ the subfamily consisting of all graph embeddings that share a common crossing number that is the least among all the crossing numbers of the elements of $\mathcal{U}(K,U)$, with respect to the unit cells connected to $U$ by crossing changes. We define as \textit{ground states} the elements of the family $\mathcal{G}(K,U)$, which we claim to be the least tangled embeddings of the graph $\mathcal{K}$. For example, \textbf{dia}, with respect to its primitive unit cell, is a ground state since its minimum crossing number triplet is $(0,0,0)$, which is associated to the least possible crossing number, that is, $0$. Further mathematical details on the families $\mathcal{U}(K,U)$ and $\mathcal{G}(K,U)$ are given in the Supplementary Information.

\begin{figure*}[hbtp]
    \centering
    \begin{subfigure}[b]{0.18\textwidth}
        \centering
        \includegraphics[width=0.9\textwidth]{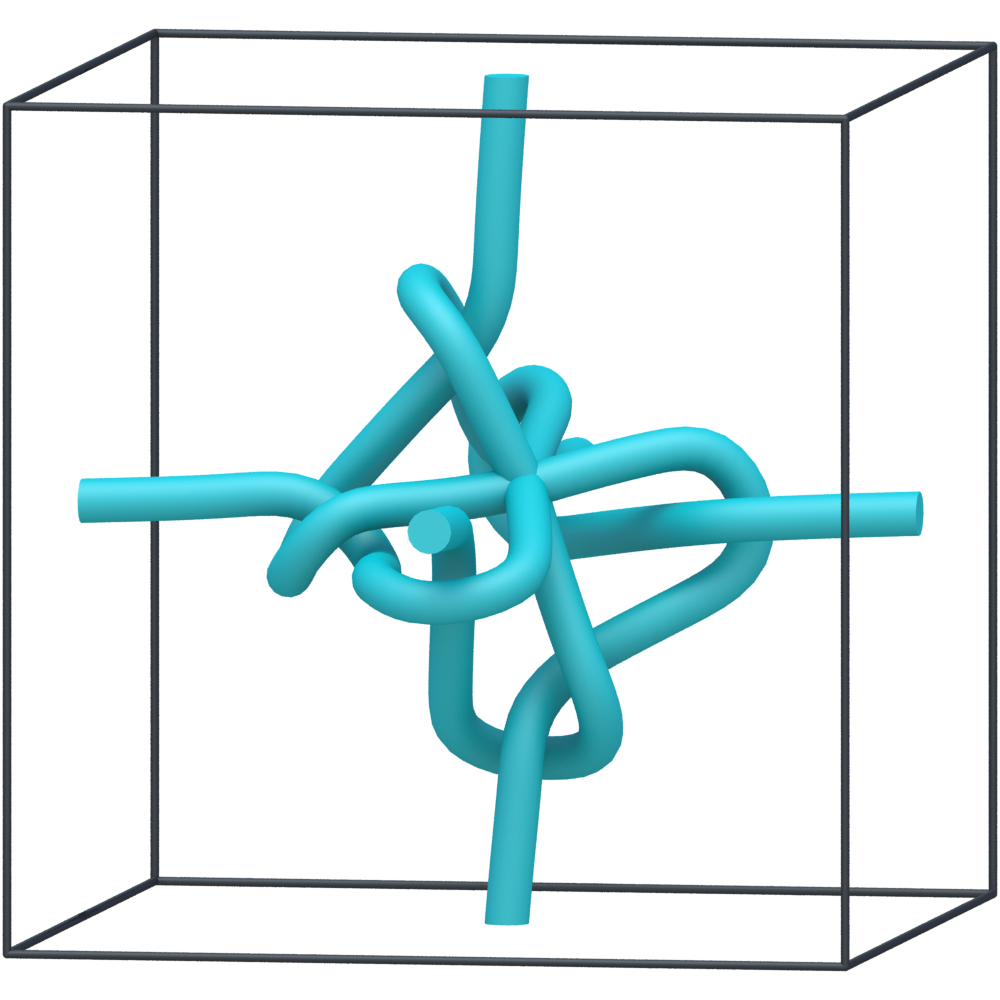}
        \caption{}
        \label{fig:ravelled_pcu_uc}
    \end{subfigure}
    \begin{subfigure}[b]{0.18\textwidth}
        \centering
        \includegraphics[width=0.9\textwidth]{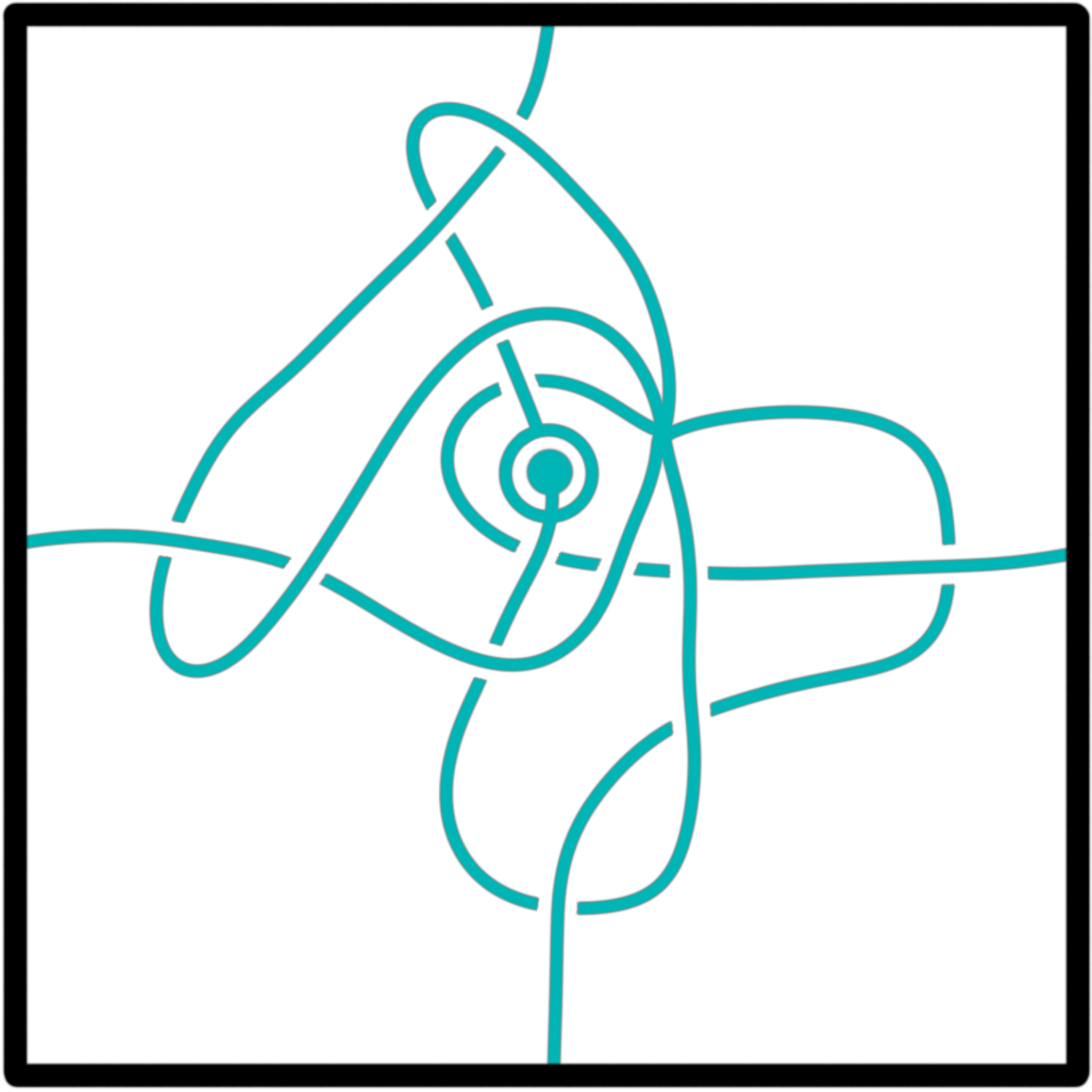}
        \caption{}
        \label{fig:untangling_ravelled_pcu_dia_1}
    \end{subfigure}
    \begin{subfigure}[b]{0.18\textwidth}
        \centering
        \includegraphics[width=0.9\textwidth]{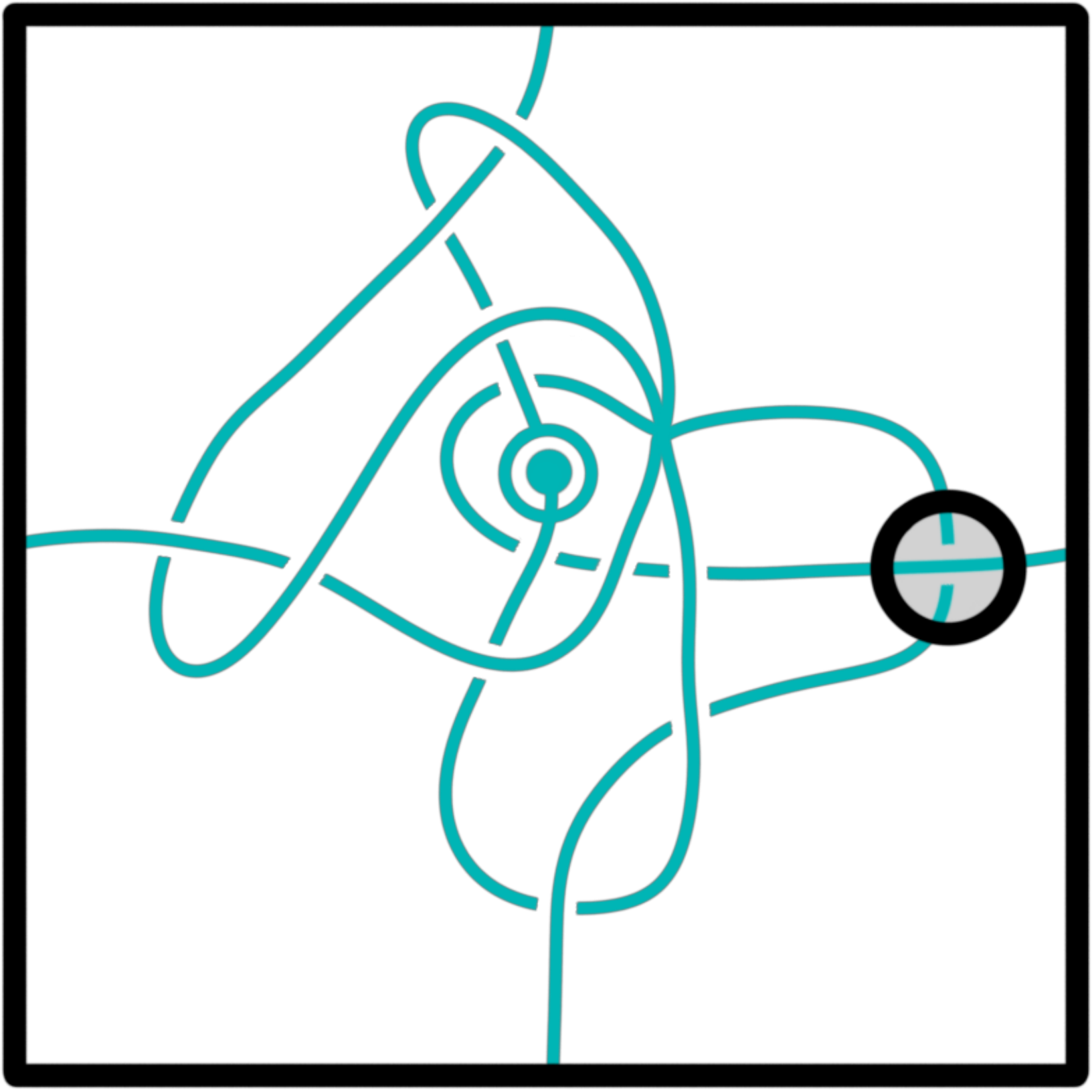}
        \caption{}
        \label{fig:untangling_ravelled_pcu_dia_2}
    \end{subfigure}
    \begin{subfigure}[b]{0.18\textwidth}
        \centering
        \includegraphics[width=0.9\textwidth]{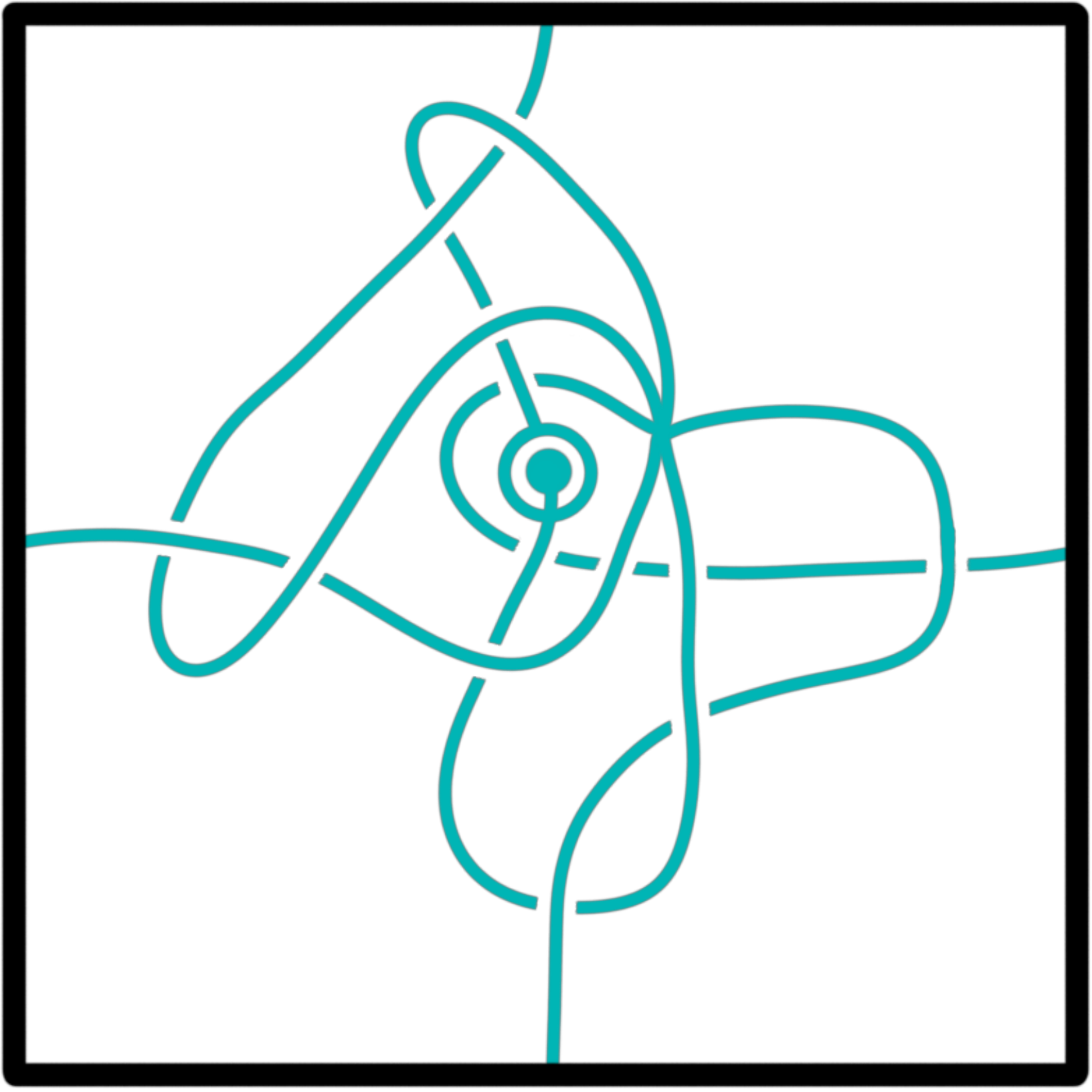}
        \caption{}
        \label{fig:untangling_ravelled_pcu_dia_3}
    \end{subfigure}
    \begin{subfigure}[b]{0.18\textwidth}
        \centering
        \includegraphics[width=0.9\textwidth]{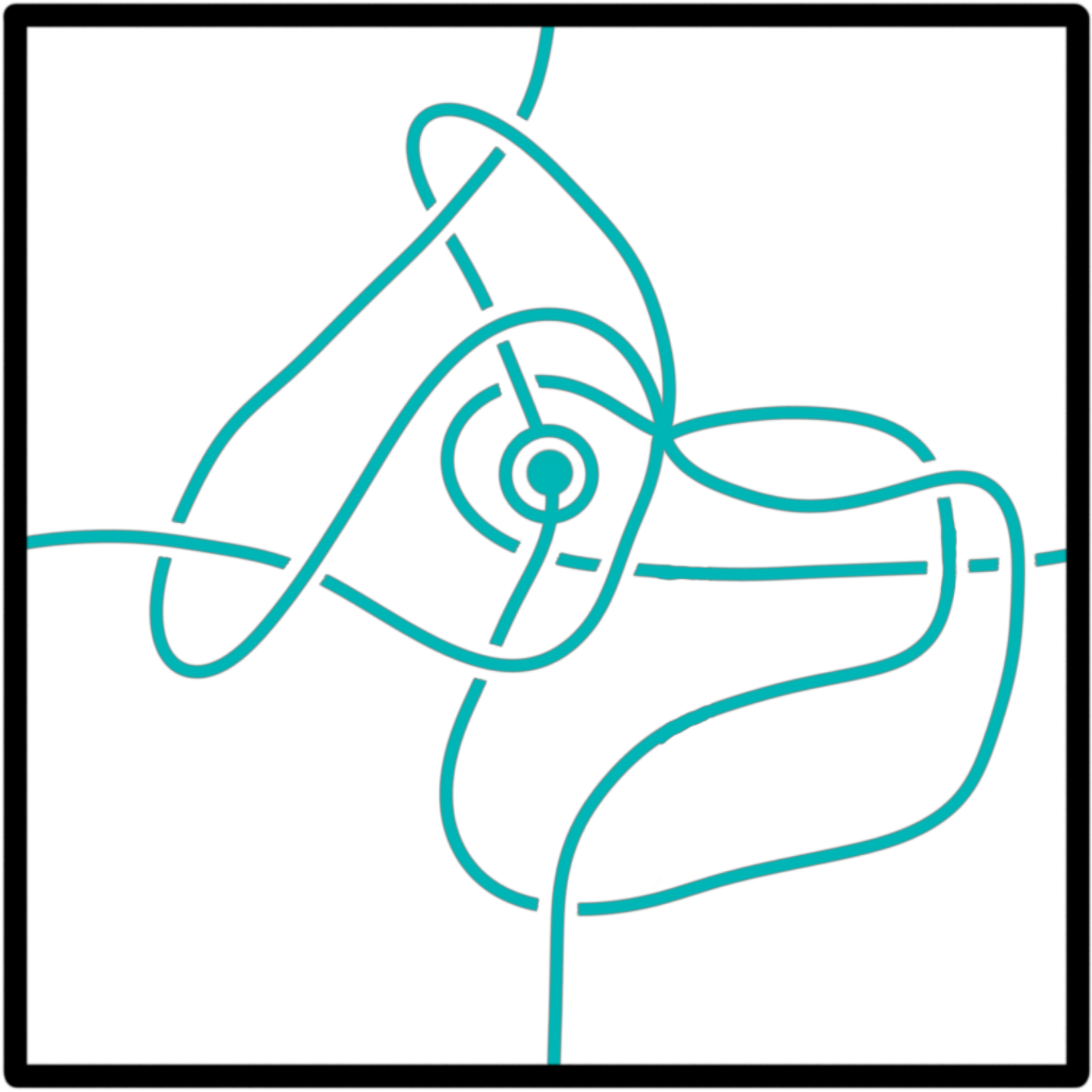}
        \caption{}
        \label{fig:untangling_ravelled_pcu_dia_4}
    \end{subfigure}
    
    \vskip\baselineskip
    
    \begin{subfigure}[b]{0.18\textwidth}
        \centering
        \includegraphics[width=0.9\textwidth]{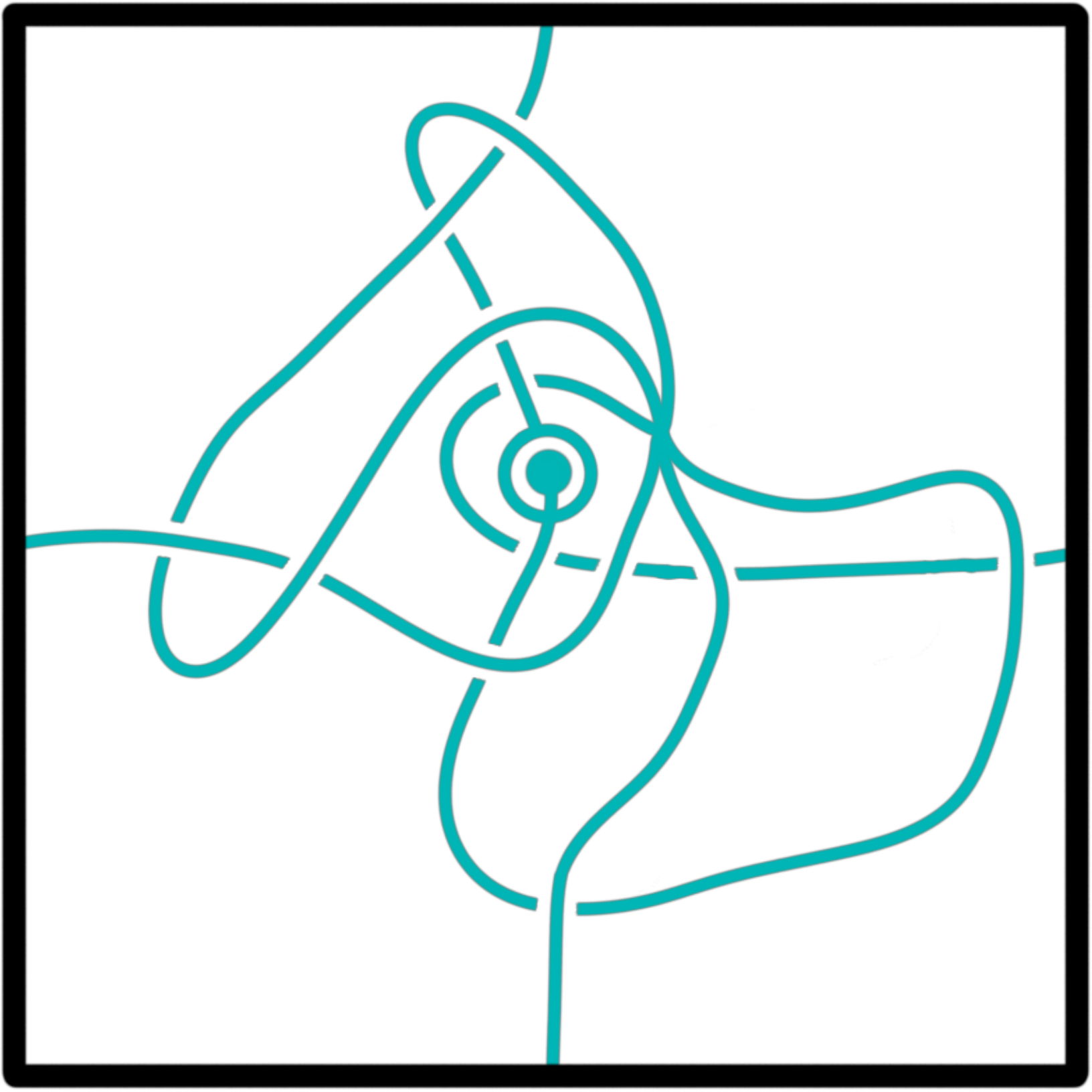}
        \caption{}
        \label{fig:untangling_ravelled_pcu_dia_5}
    \end{subfigure}
    \begin{subfigure}[b]{0.18\textwidth}
        \centering
        \includegraphics[width=0.9\textwidth]{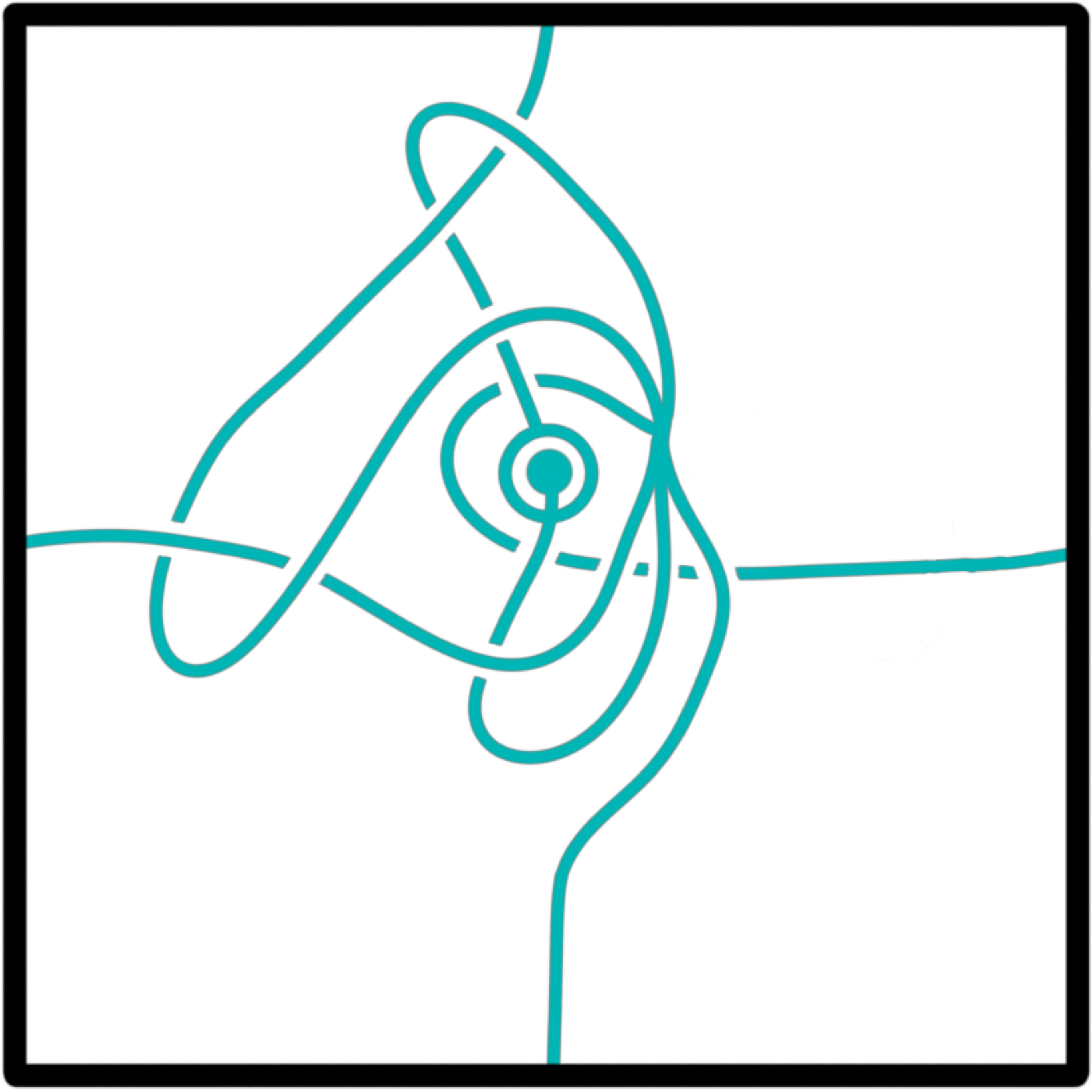}
        \caption{}
        \label{fig:untangling_ravelled_pcu_dia_6}
    \end{subfigure}
    \begin{subfigure}[b]{0.18\textwidth}
        \centering
        \includegraphics[width=0.9\textwidth]{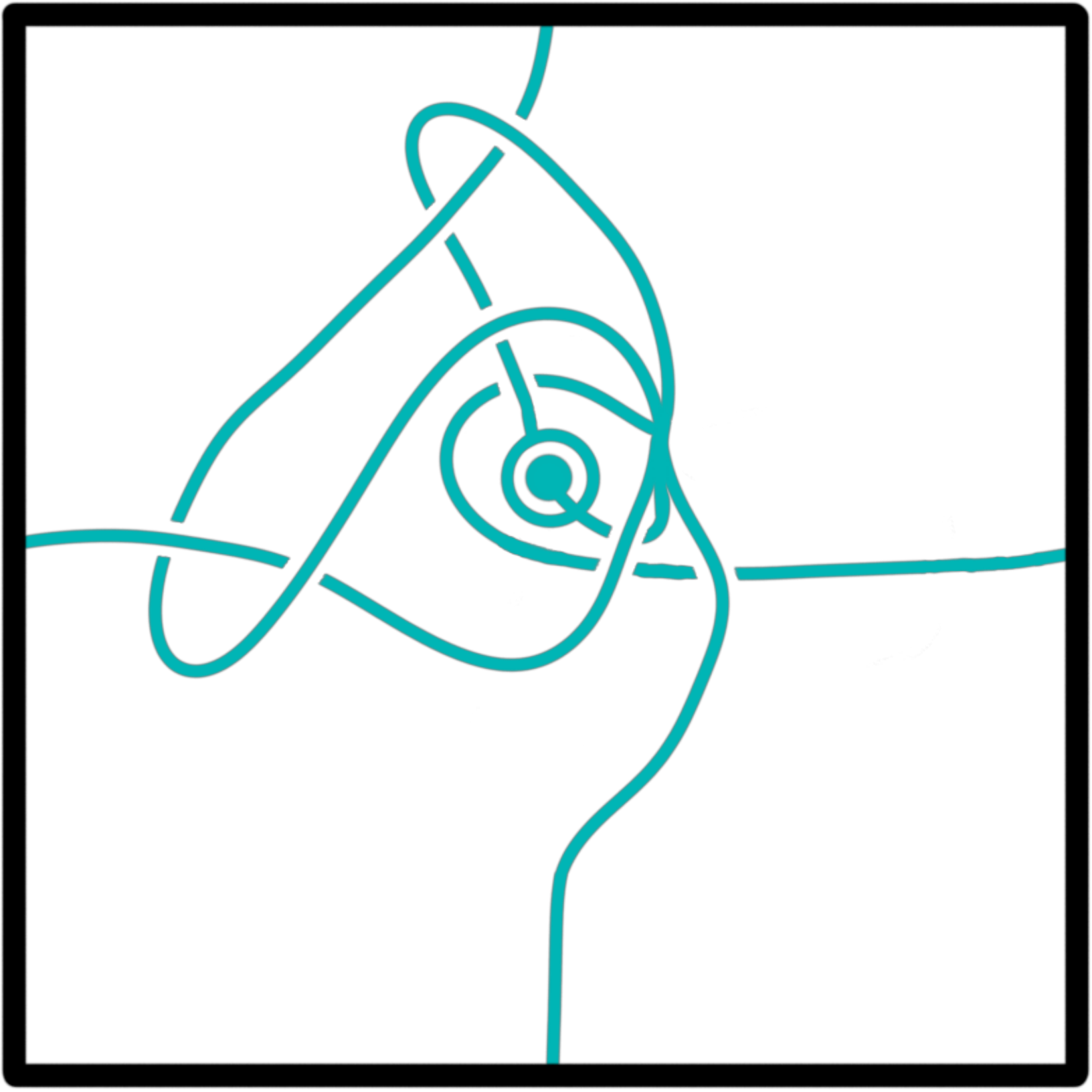}
        \caption{}
        \label{fig:untangling_ravelled_pcu_dia_7}
    \end{subfigure}
    \begin{subfigure}[b]{0.18\textwidth}
        \centering
        \includegraphics[width=0.9\textwidth]{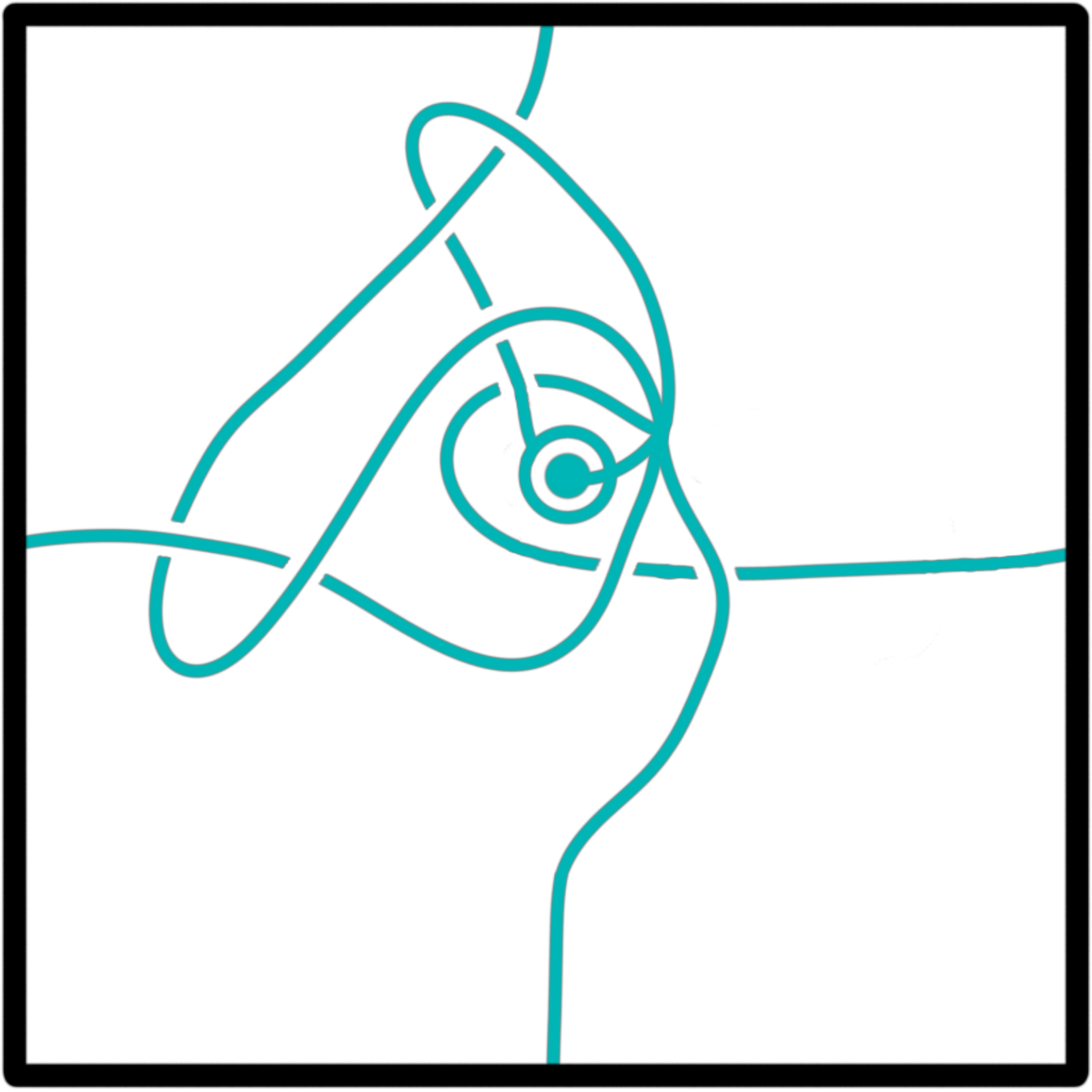}
        \caption{}
        \label{fig:untangling_ravelled_pcu_dia_8}
    \end{subfigure}
    \begin{subfigure}[b]{0.18\textwidth}
        \centering
        \includegraphics[width=0.9\textwidth]{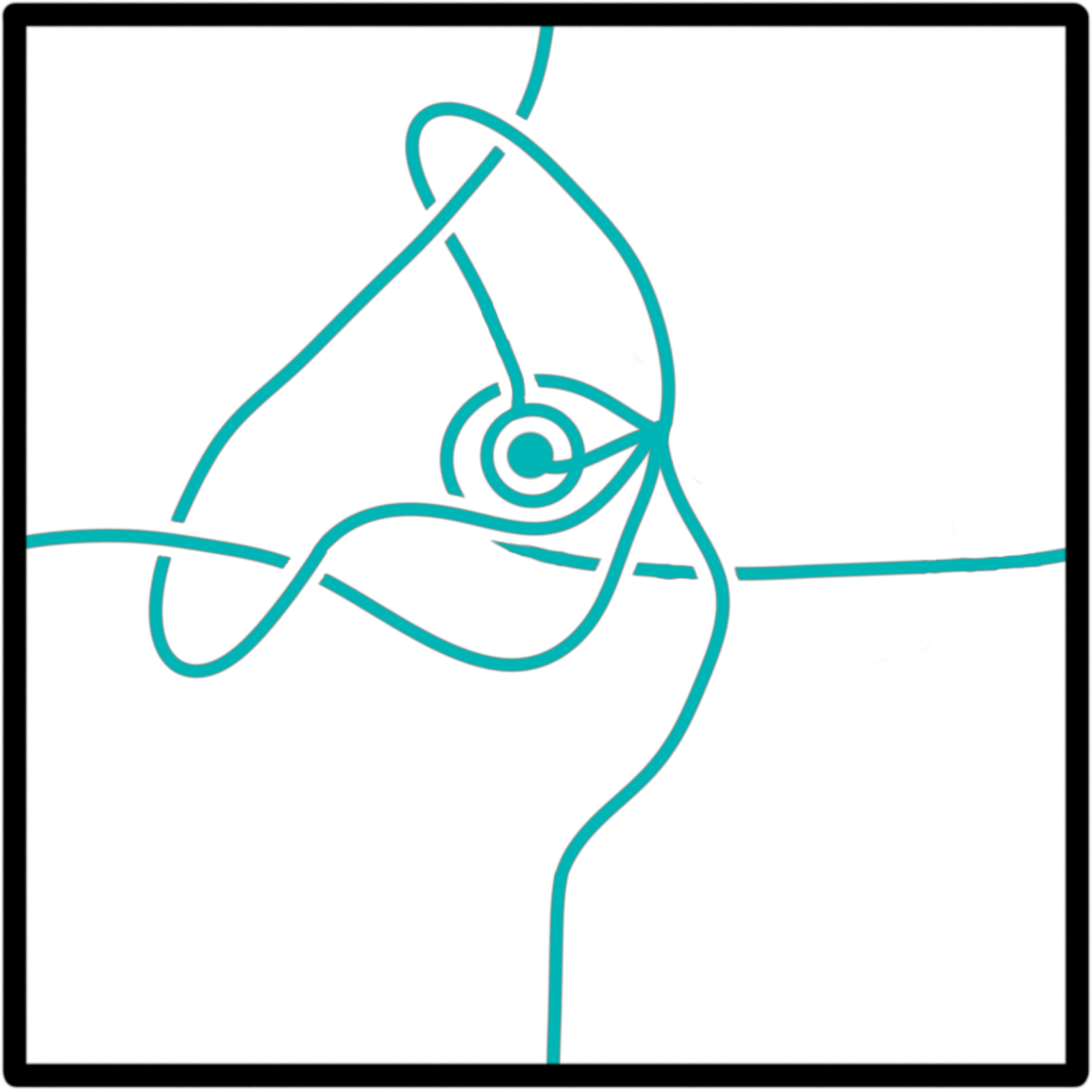}
        \caption{}
        \label{fig:untangling_ravelled_pcu_dia_9}
    \end{subfigure}

    \vskip\baselineskip

    \begin{subfigure}[b]{0.18\textwidth}
        \centering
        \includegraphics[width=0.9\textwidth]{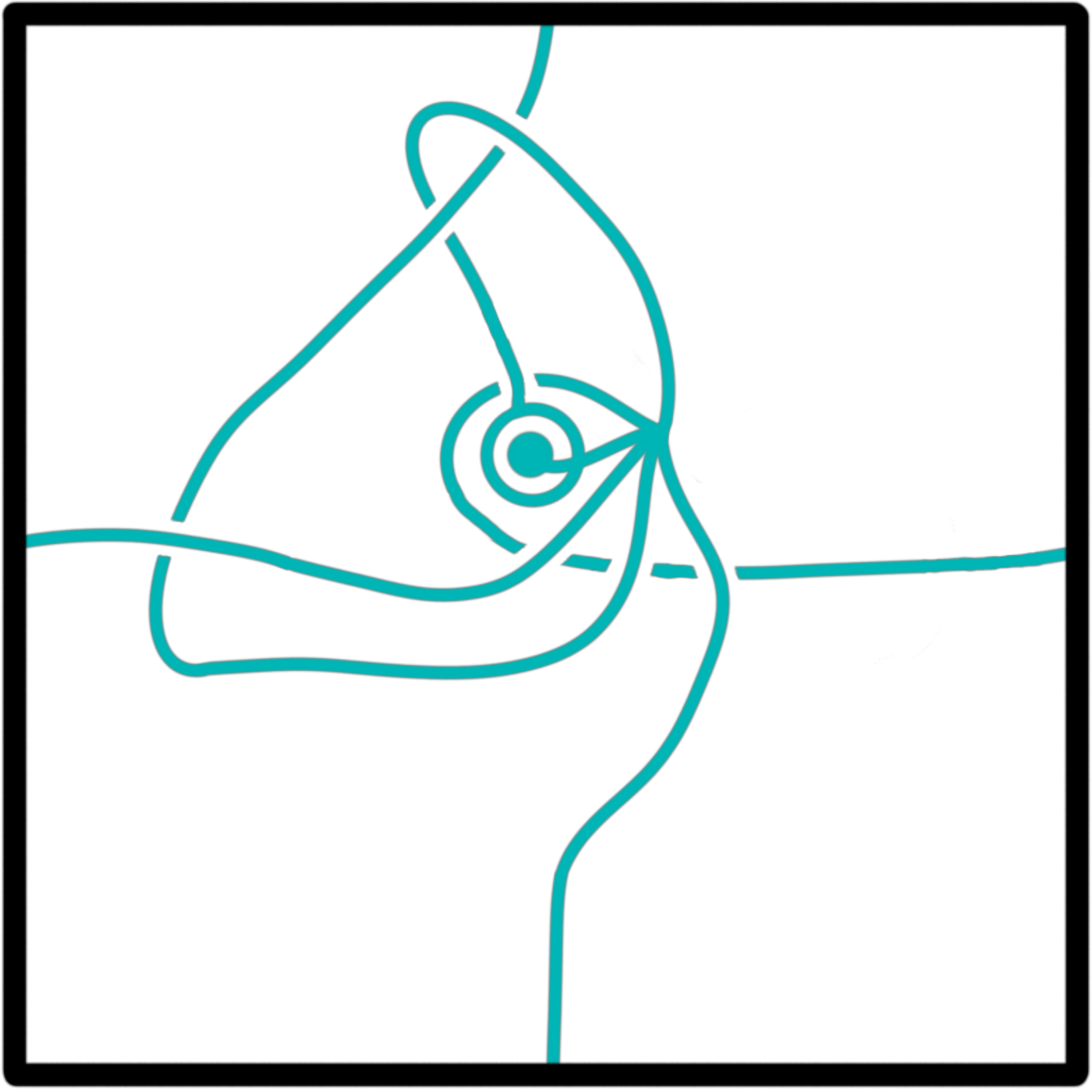}
        \caption{}
        \label{fig:untangling_ravelled_pcu_dia_10}
    \end{subfigure}
    \begin{subfigure}[b]{0.18\textwidth}
        \centering
        \includegraphics[width=0.9\textwidth]{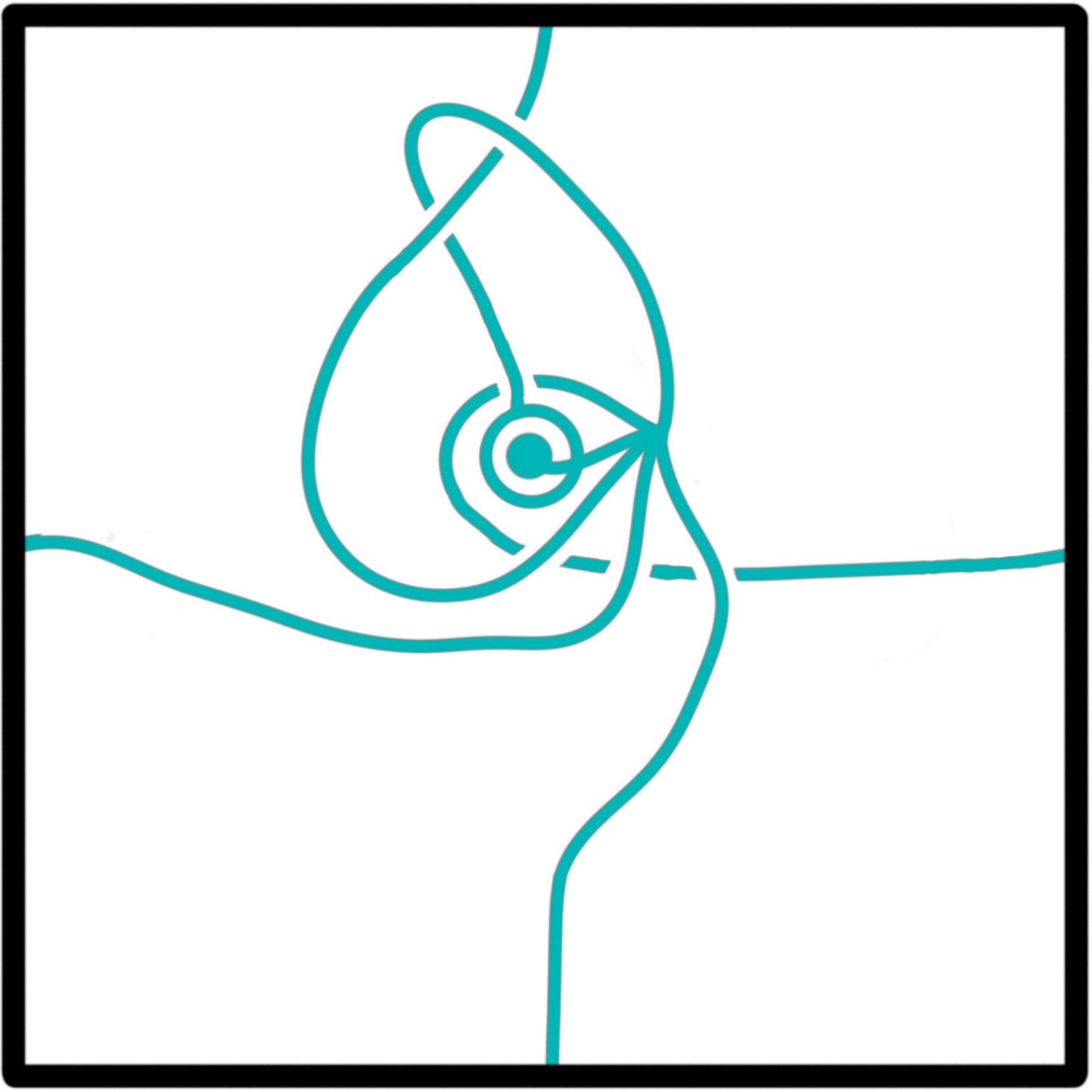}
        \caption{}
        \label{fig:untangling_ravelled_pcu_dia_11}
    \end{subfigure}
    \begin{subfigure}[b]{0.18\textwidth}
        \centering
        \includegraphics[width=0.9\textwidth]{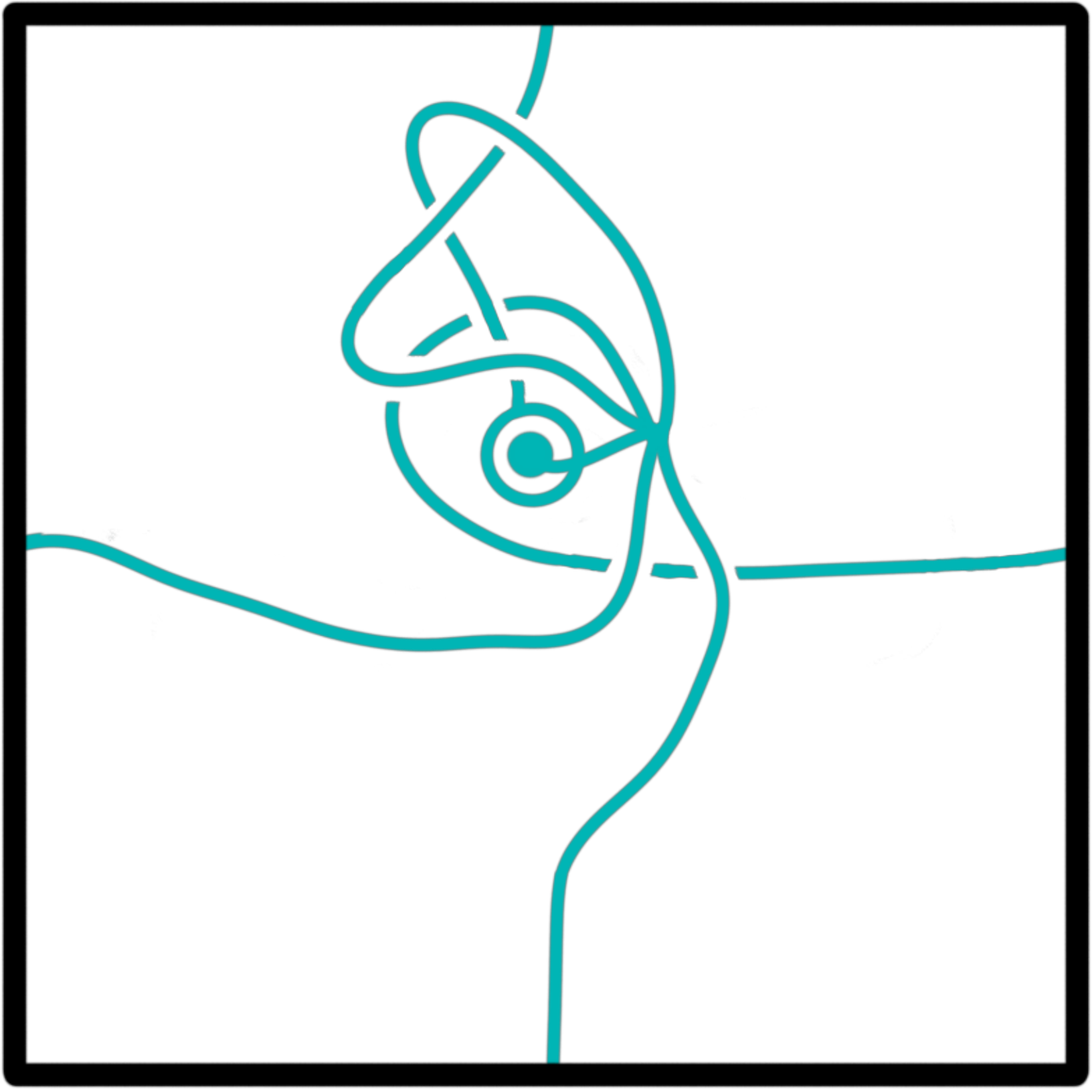}
        \caption{}
        \label{fig:untangling_ravelled_pcu_dia_12}
    \end{subfigure}
    \begin{subfigure}[b]{0.18\textwidth}
        \centering
        \includegraphics[width=0.9\textwidth]{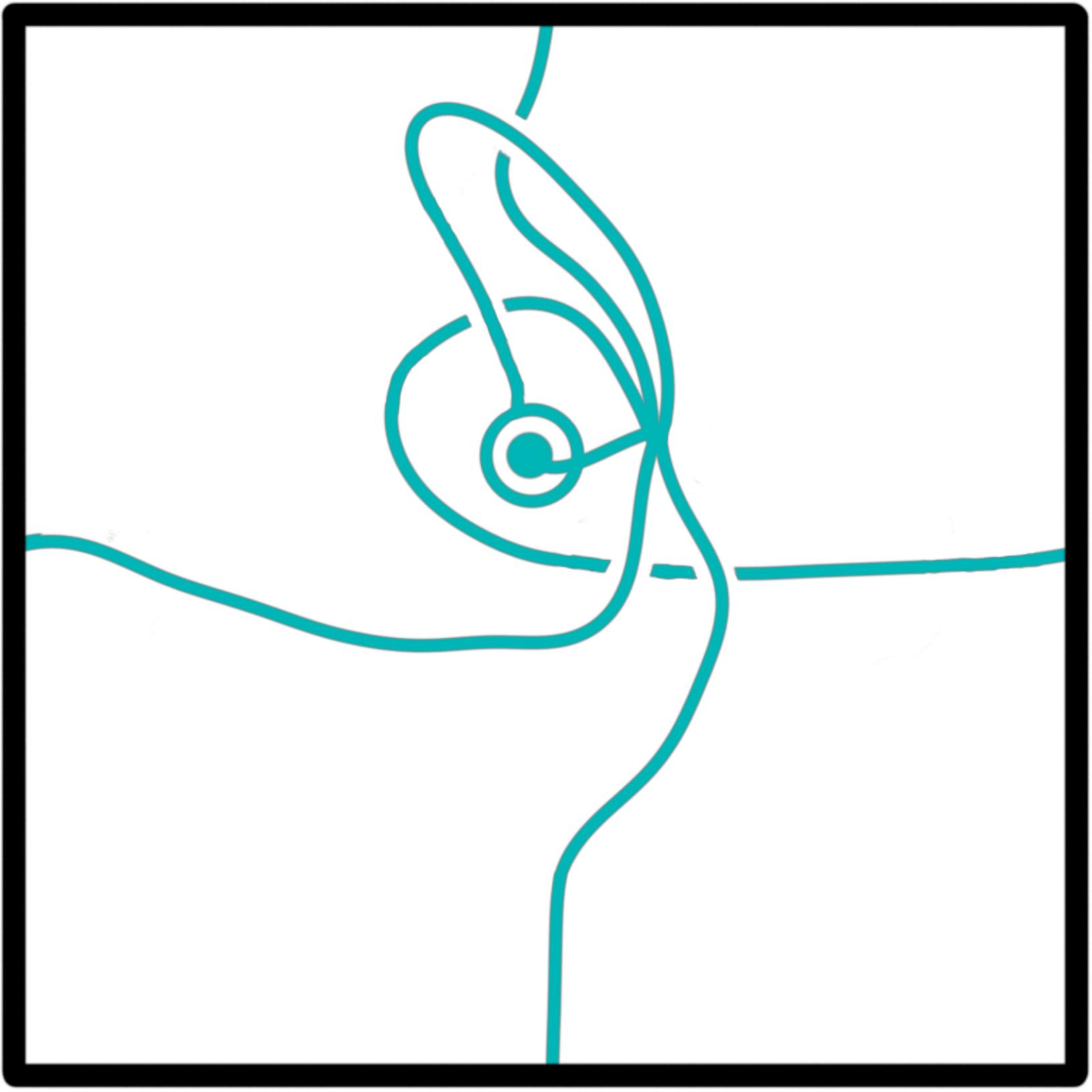}
        \caption{}
        \label{fig:untangling_ravelled_pcu_dia_13}
    \end{subfigure}
    \begin{subfigure}[b]{0.18\textwidth}
        \centering
        \includegraphics[width=0.9\textwidth]{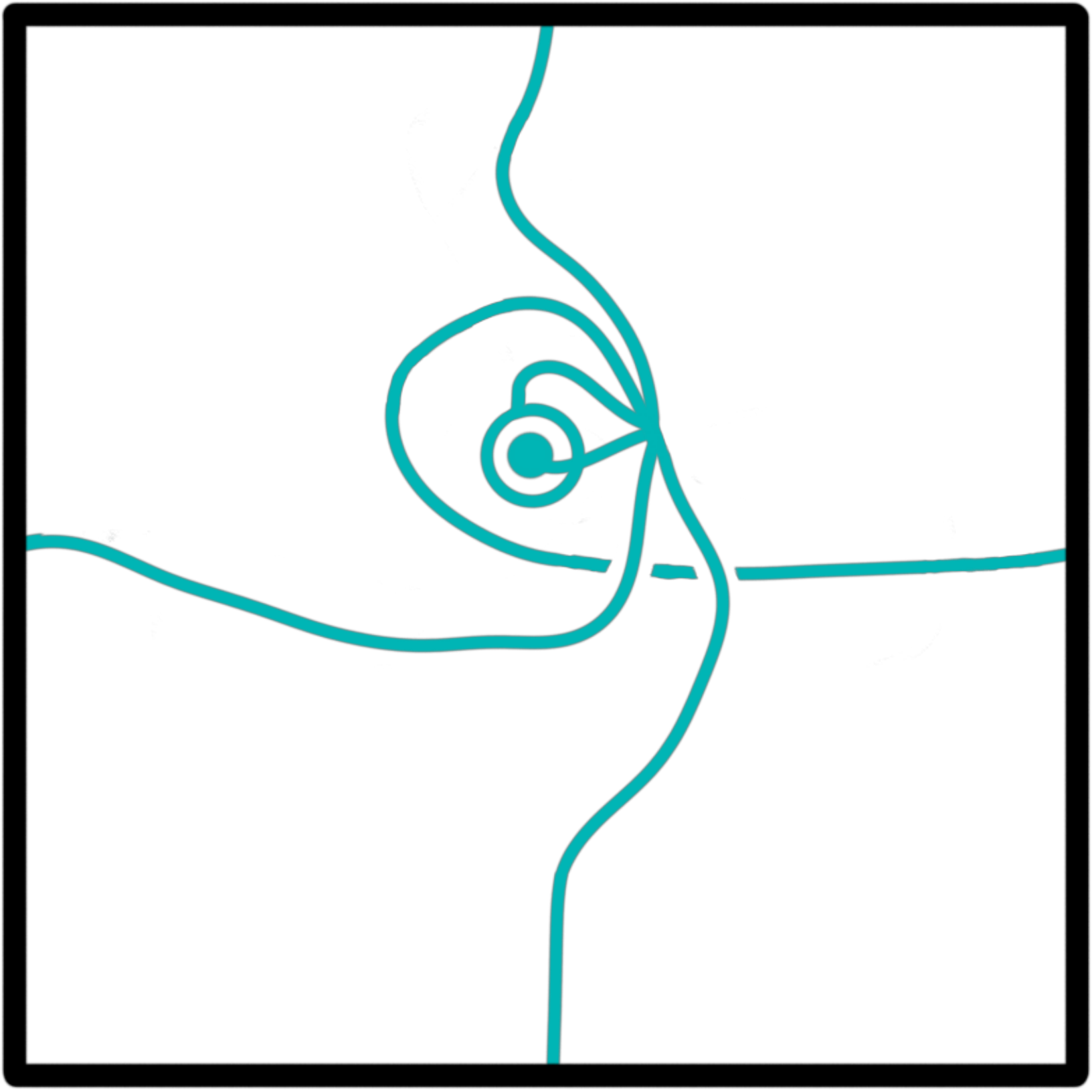}
        \caption{}
        \label{fig:untangling_ravelled_pcu_dia_14}
    \end{subfigure}

    \vskip\baselineskip
    
    \begin{subfigure}[b]{0.18\textwidth}
        \centering
        \includegraphics[width=0.9\textwidth]{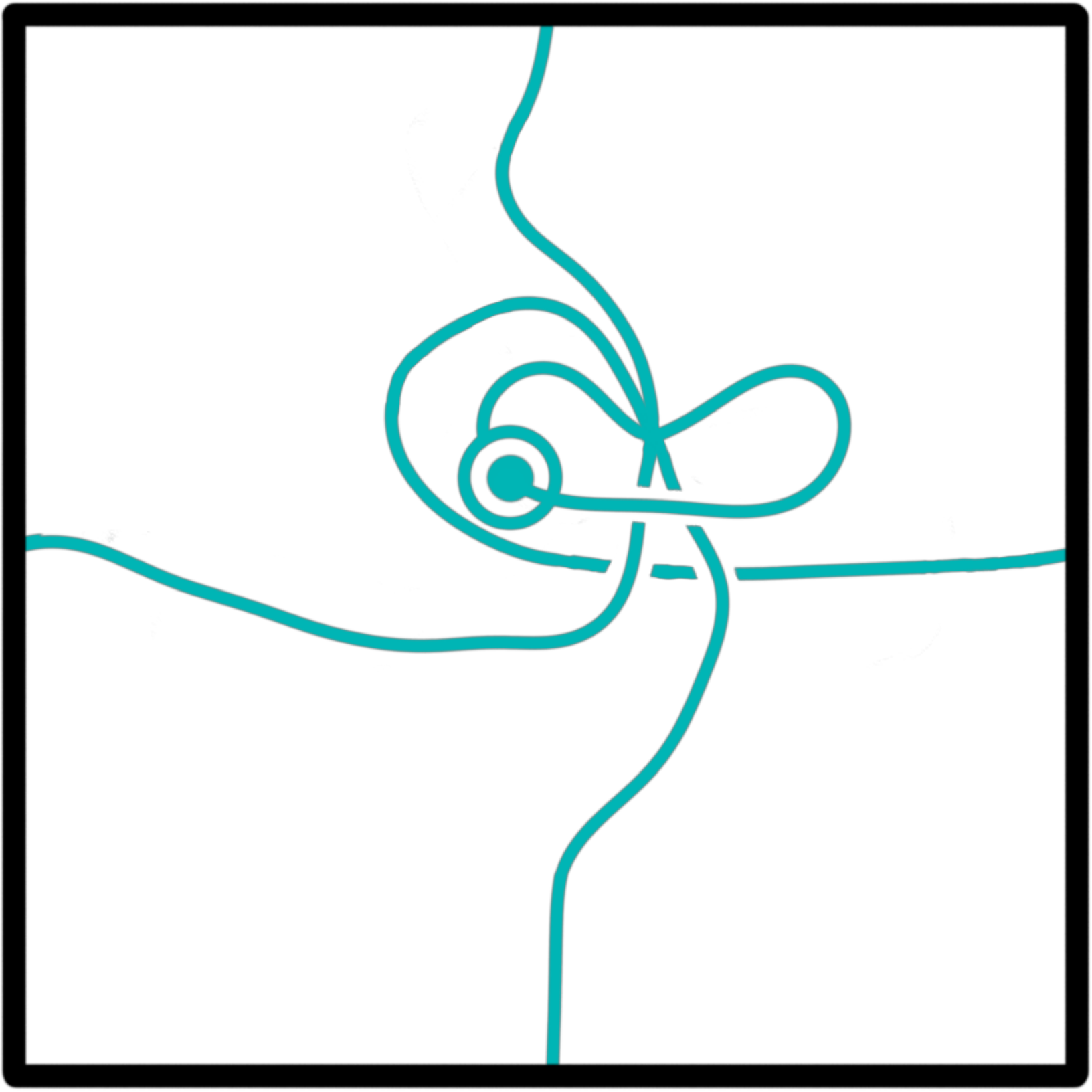}
        \caption{}
        \label{fig:untangling_ravelled_pcu_dia_15}
    \end{subfigure}
    \begin{subfigure}[b]{0.18\textwidth}
        \centering
        \includegraphics[width=0.9\textwidth]{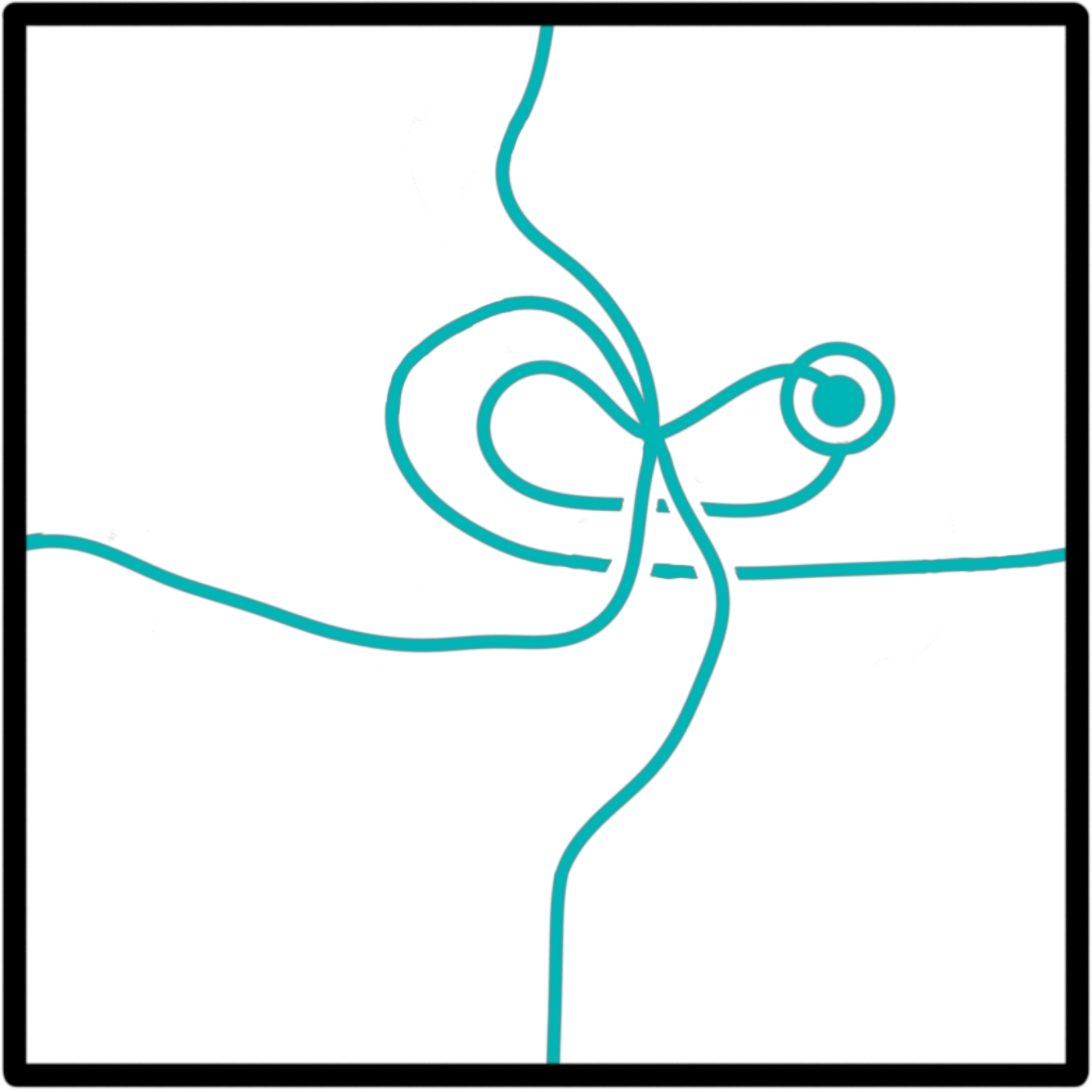}
        \caption{}
        \label{fig:untangling_ravelled_pcu_dia_16}
    \end{subfigure}
    \begin{subfigure}[b]{0.18\textwidth}
        \centering
        \includegraphics[width=0.9\textwidth]{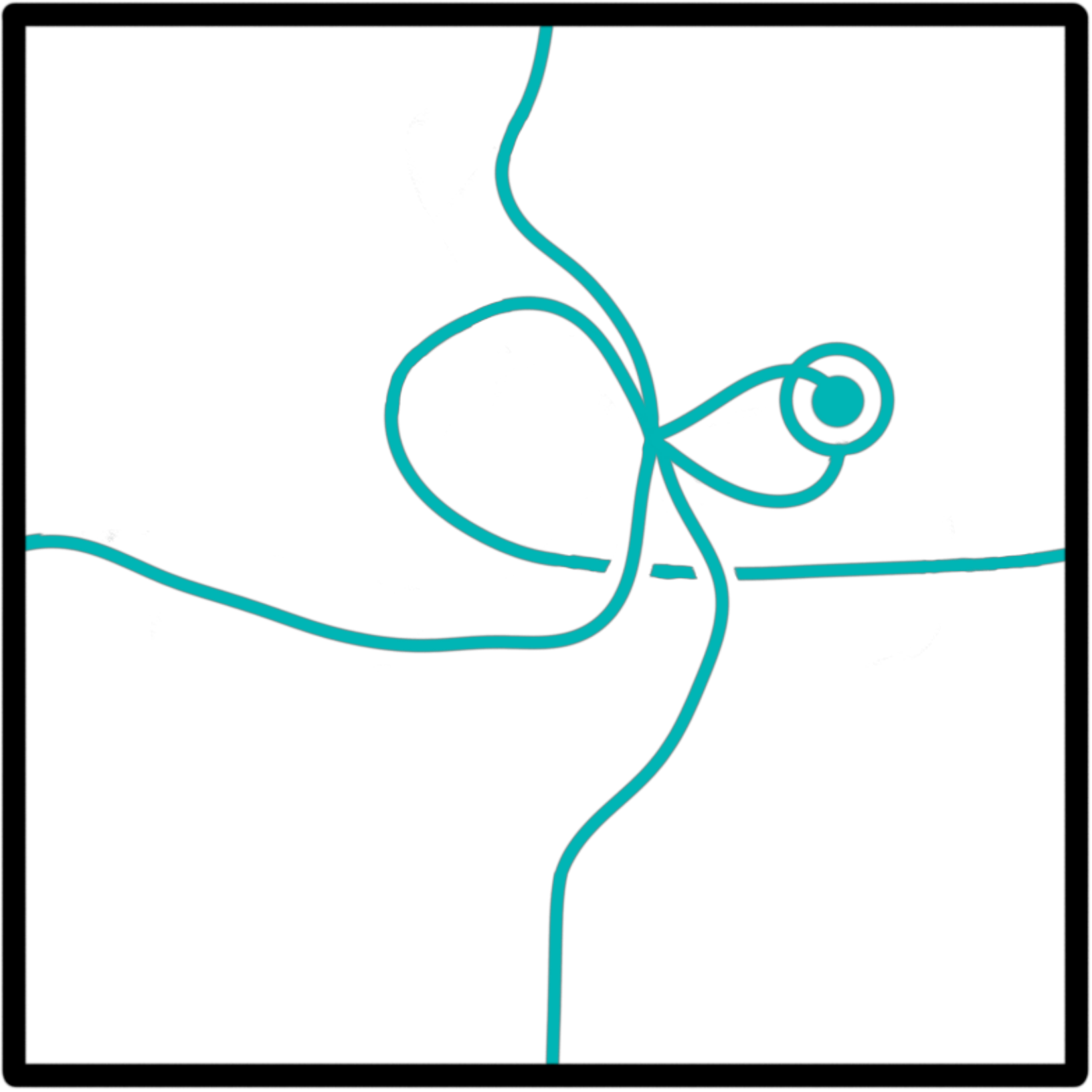}
        \caption{}
        \label{fig:untangling_ravelled_pcu_dia_17}
    \end{subfigure}
    \begin{subfigure}[b]{0.18\textwidth}
        \centering
        \includegraphics[width=0.9\textwidth]{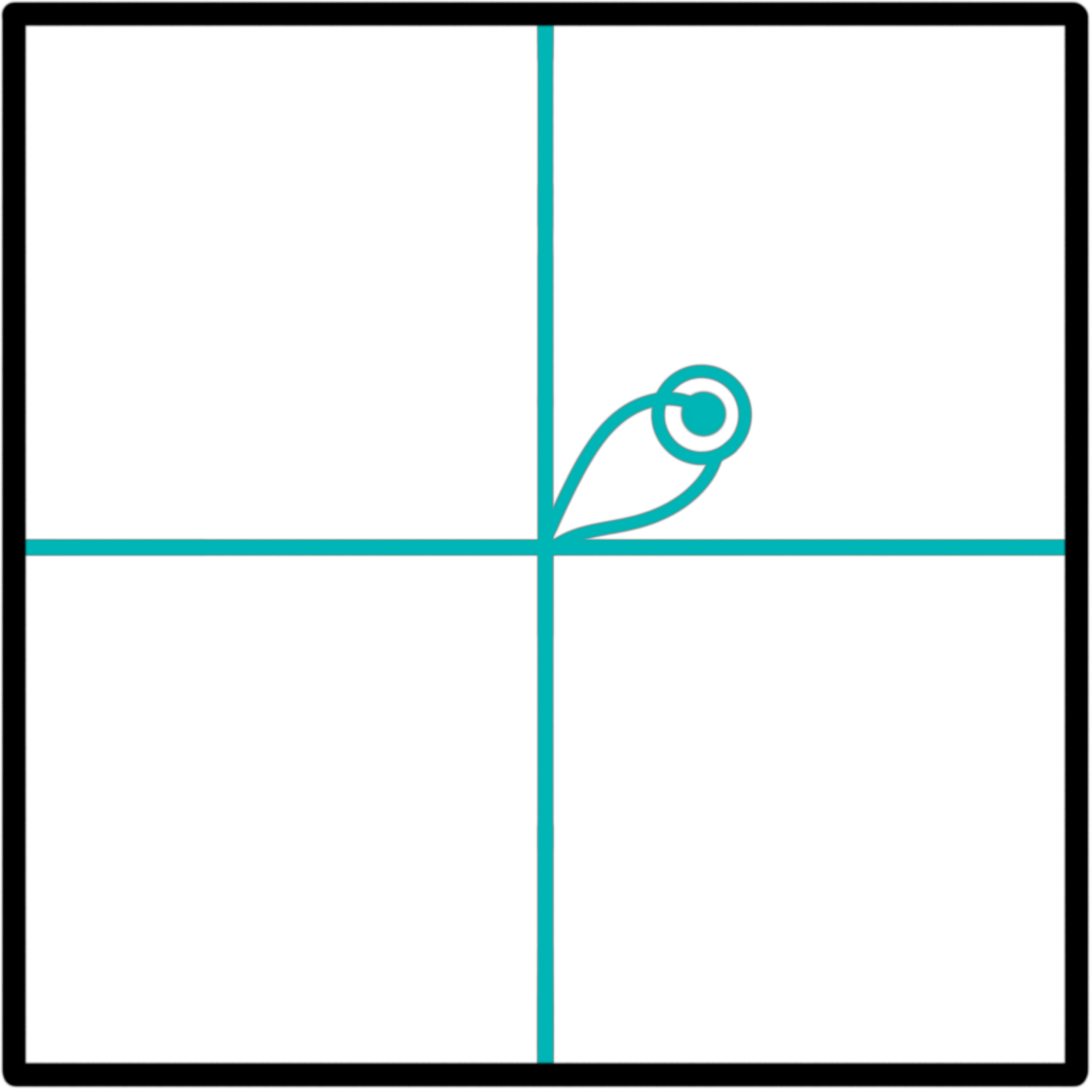}
        \caption{}
        \label{fig:untangling_ravelled_pcu_dia_18}
    \end{subfigure}
    \begin{subfigure}[b]{0.18\textwidth}
        \centering
        \includegraphics[width=0.9\textwidth]{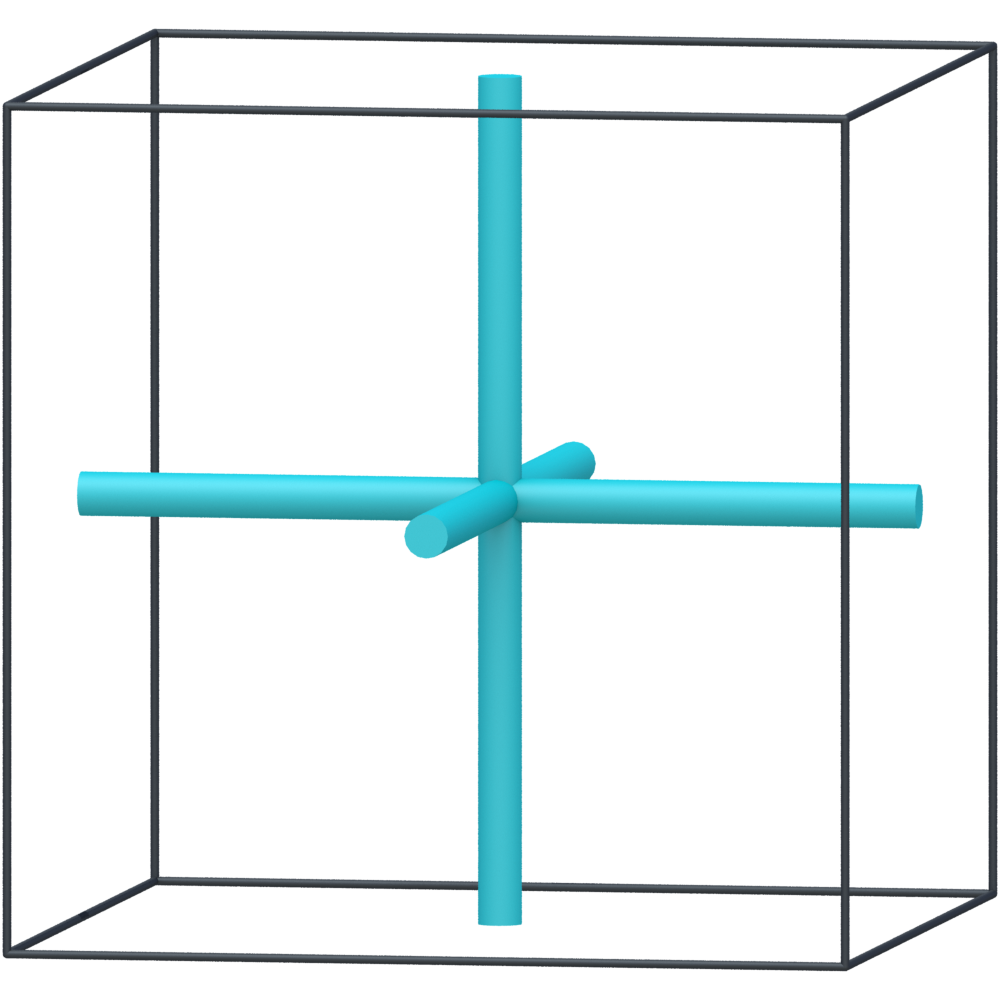}
        \caption{}
        \label{fig:pcu_uc_untangling_ravelled_pcu}
    \end{subfigure}
    \caption{Untangling the ravelled embedding of \textbf{pcu} shown in Fig. \ref{fig:ravelled_pcu_extended}: Although the strong rings of the embedding are neither knotted nor linked, and although the HRN is trivial, it is not a ground state and it can be untangled to the barycentric embedding of \textbf{pcu}. This can be done by applying a single crossing change on the crossing highlighted by a black circle.}
    \label{fig:untangling_ravelled_pcu}
\end{figure*}

We present in Fig. \ref{fig:examples_of_ground_states} some examples of ground states. These structures are, respectively from top to bottom, \textbf{sql}, \textbf{srs}, \textbf{dia}, \textbf{pcu} and \textbf{dia-c}. The first four have minimum crossing number triplet $(0,0,0)$ with respect to their usual unit cells, and so does \textbf{dia-c} with respect to its primitive unit cell. We note that the \textbf{sql} network is a 2-periodic network, but a 3-periodic structure can be obtained from it by layering infinitely many copies of the network in the third direction of space. In particular, this example shows that the ground states are well defined for lower-periodic structures, by considering the crossing information obtained from three diagrams projected from three non-coplanar axes. We also note that \textbf{dia} thus possesses two distinct unit cells with respect to which it has no crossings; the primitive unit cell in Fig. \ref{fig:dia_primitive_uc}, and the unit cell in Fig. \ref{fig:dia_uc}. The example of \textbf{dia-c} constitutes an even more interesting case than what was stated in Sect. \ref{sec:1} regarding ground states that might have linked cycles and non-trivial HRNs. Indeed, \textbf{dia-c}'s cycles are linked and its HRN is \textbf{hxg} \cite{Alexandrov:eo5016}, yet it possesses a unit cell with no crossings.  Conversely, we also mentioned in Sect. \ref{sec:1}, via the example of the ravelled embedding of \textbf{pcu} shown in Fig. \ref{fig:ravelled_pcu_extended}, that the absence of knots and links in cycles and HRN information is not sufficient for an embedding to be a ground state. We untangle the said structure in Fig. \ref{fig:untangling_ravelled_pcu} into the barycentric embedding of \textbf{pcu}. In particular, these two examples show that the knots and links in the cycles and the HRN appear not to provide any information about the crossings in an embedding of a 3-periodic graph.

\begin{figure*}

    \centering
    
    \begin{subfigure}[b]{0.27\textwidth}
        \includegraphics[width=0.6\textwidth]{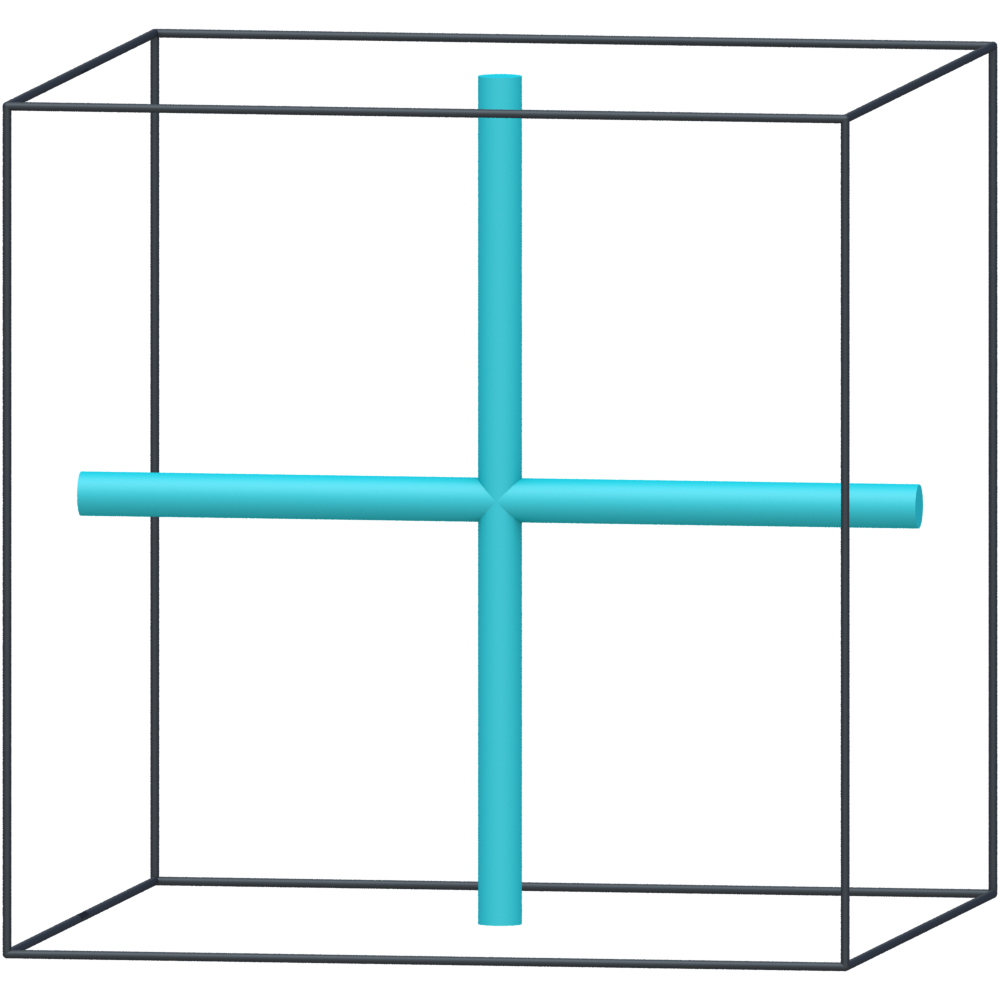}
        \caption{A unit cell of \textbf{sql}.}
        \label{fig:sql_uc}
    \end{subfigure}
    \begin{subfigure}[b]{0.54\textwidth}
        \includegraphics[width=0.29\textwidth]{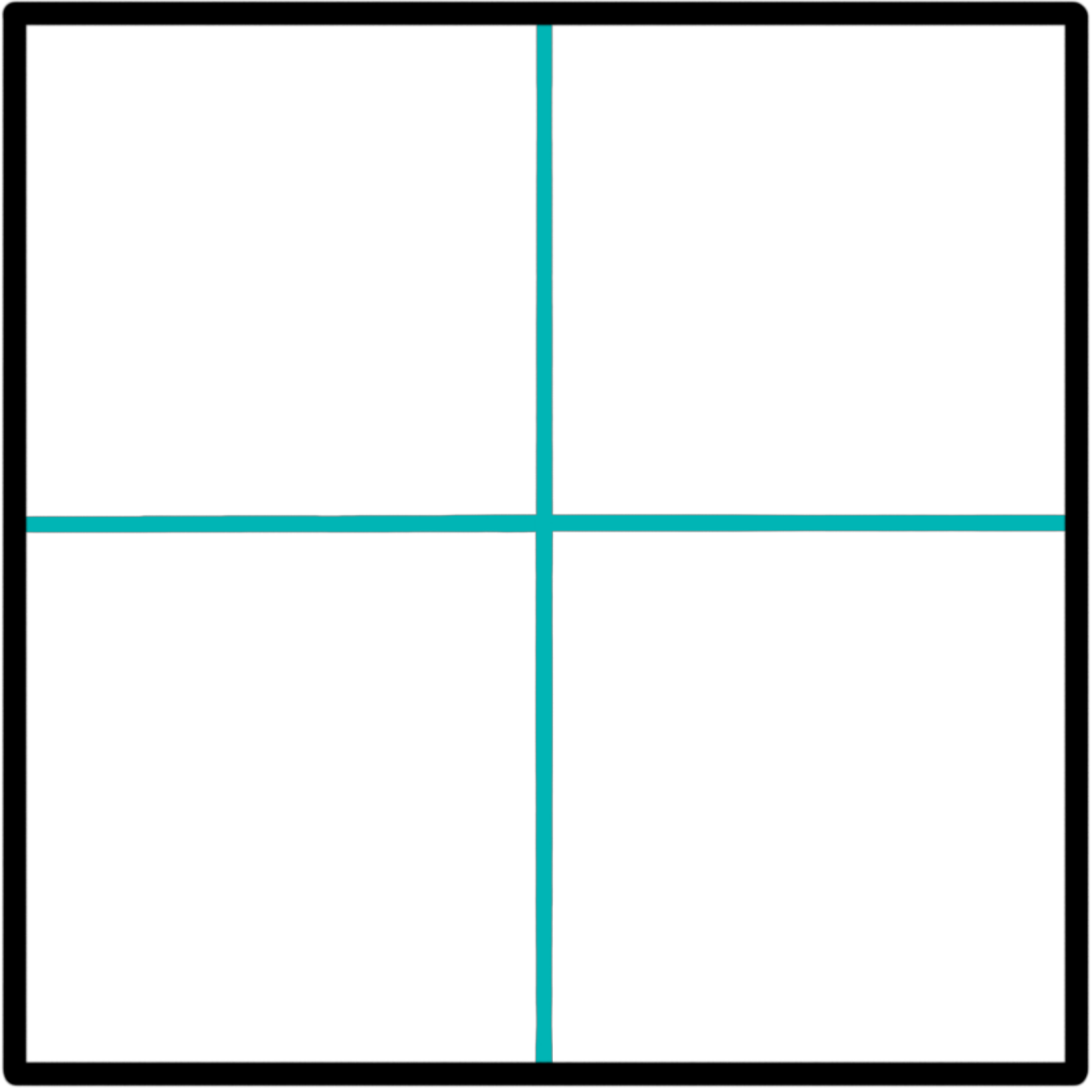}
        \hspace{0.2cm}
        \includegraphics[width=0.29\textwidth]{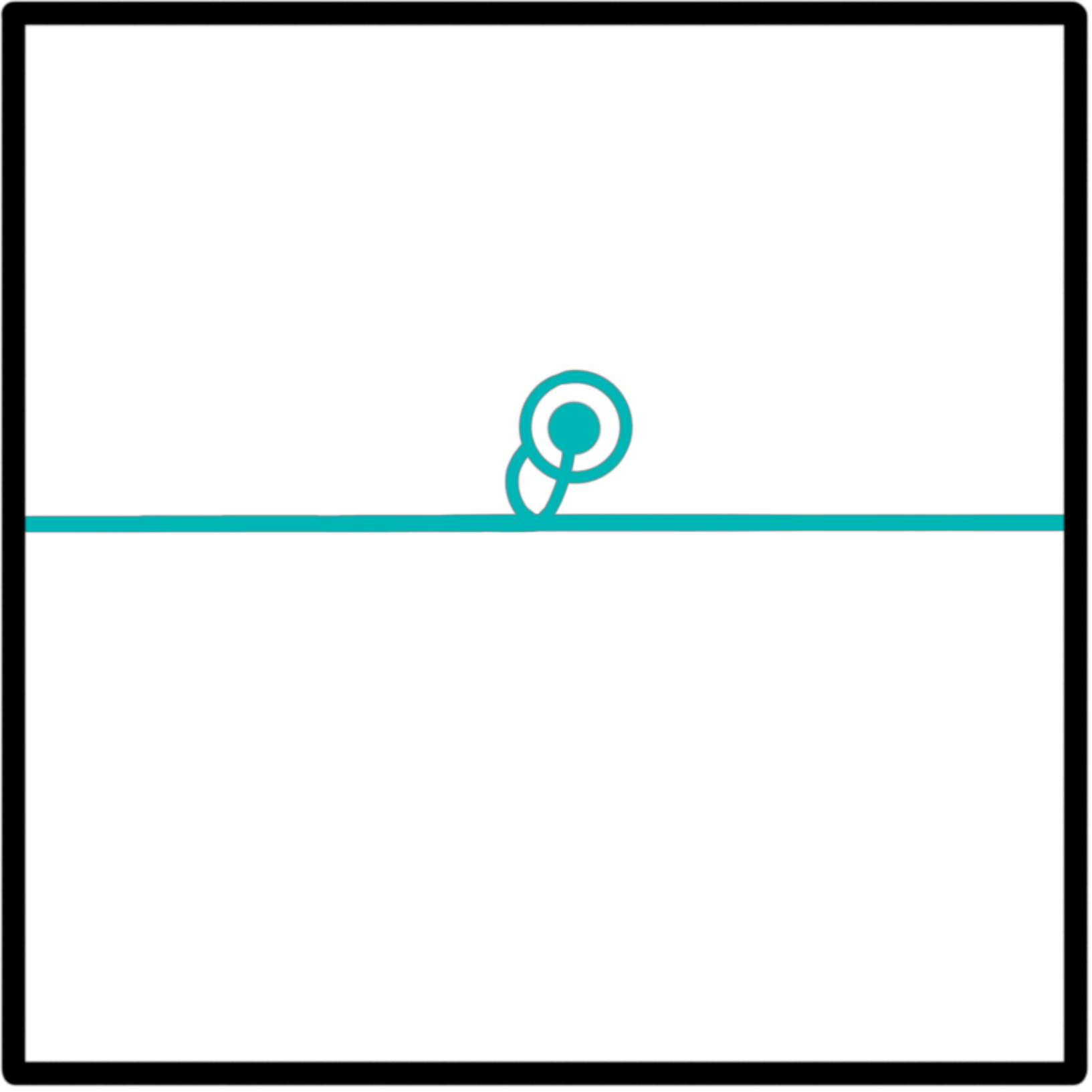}
        \hspace{0.2cm}
        \includegraphics[width=0.29\textwidth]{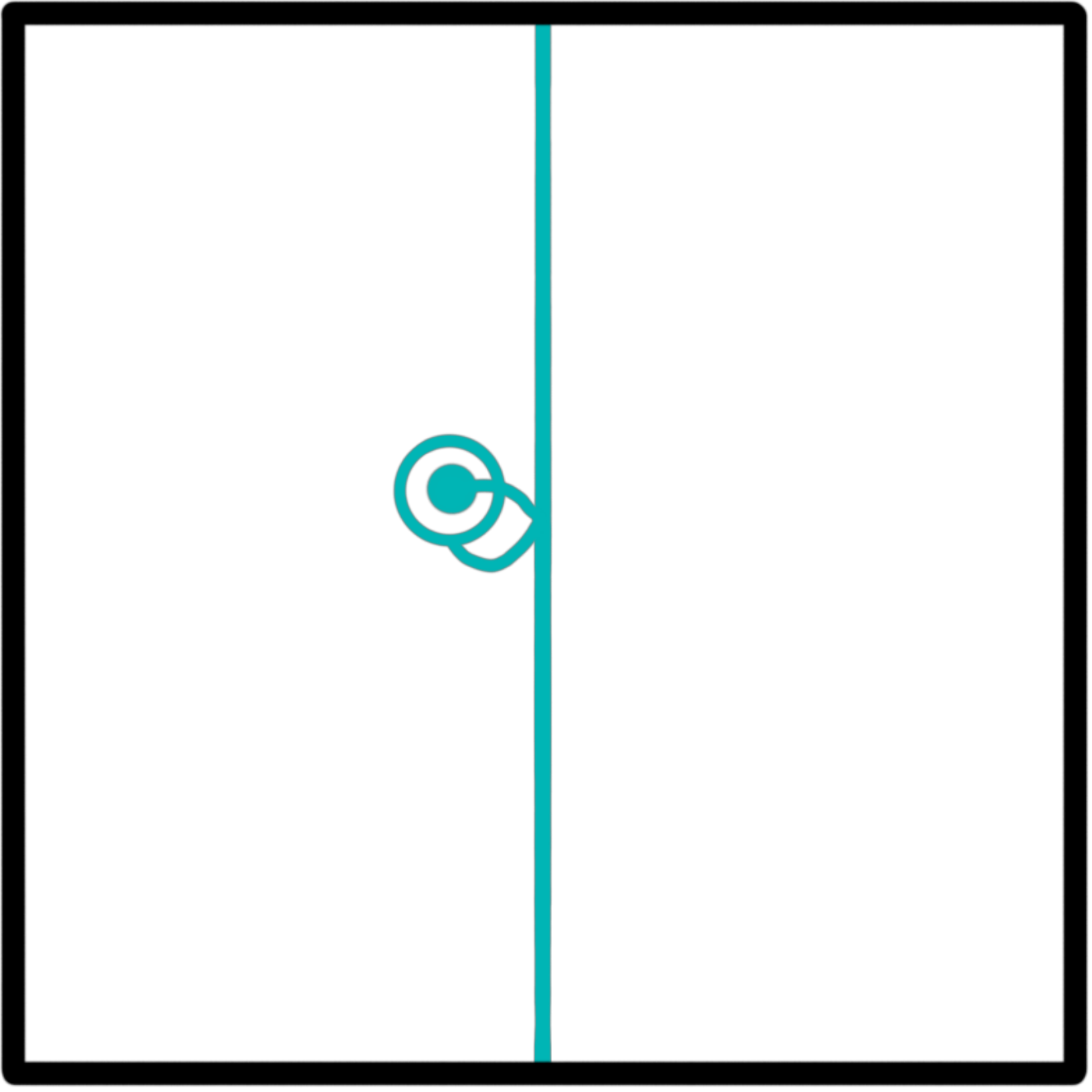}
        \caption{A tridiagram of \textbf{sql}.}
        \label{fig:sql_tridia}
    \end{subfigure}
    
    \vskip\baselineskip
    
    \begin{subfigure}[b]{0.27\textwidth}
        \includegraphics[width=0.6\textwidth]{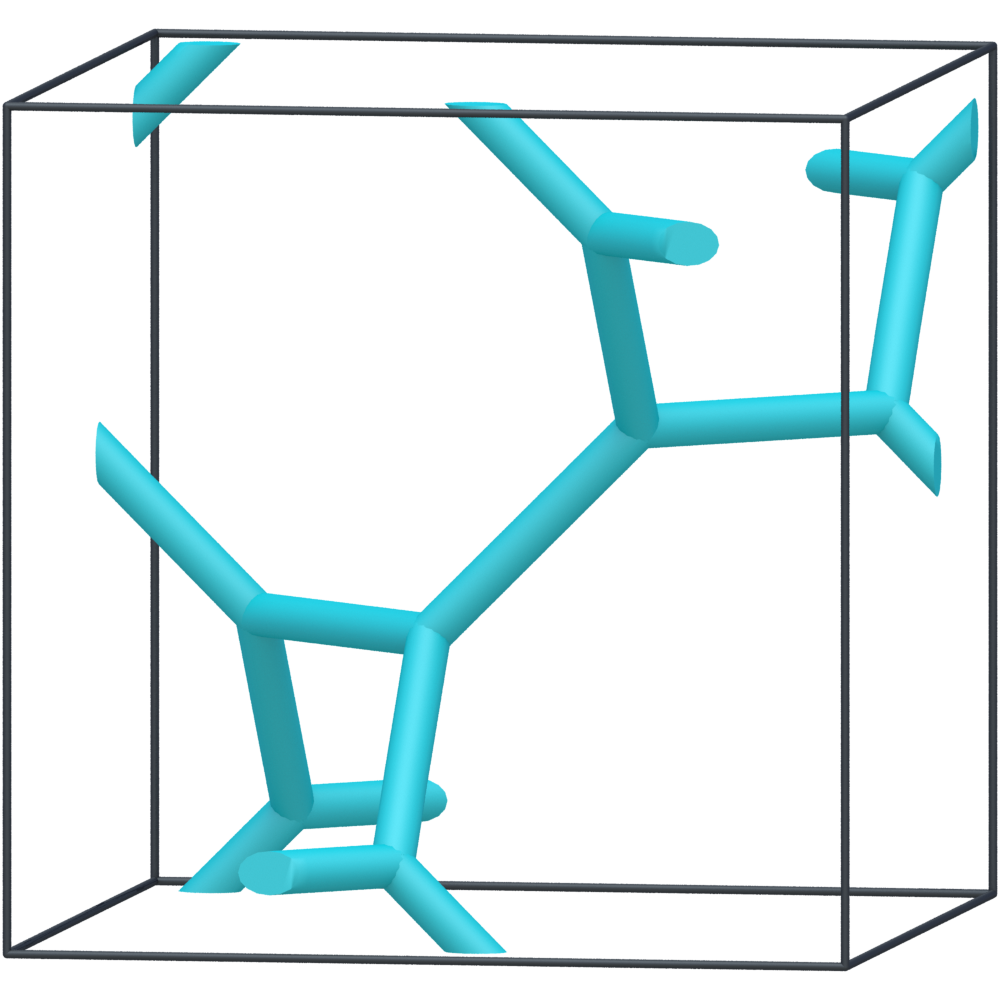}
        \caption{A unit cell of \textbf{srs}.}
        \label{fig:srs_uc}
    \end{subfigure}
    \begin{subfigure}[b]{0.54\textwidth}
        \includegraphics[width=0.29\textwidth]{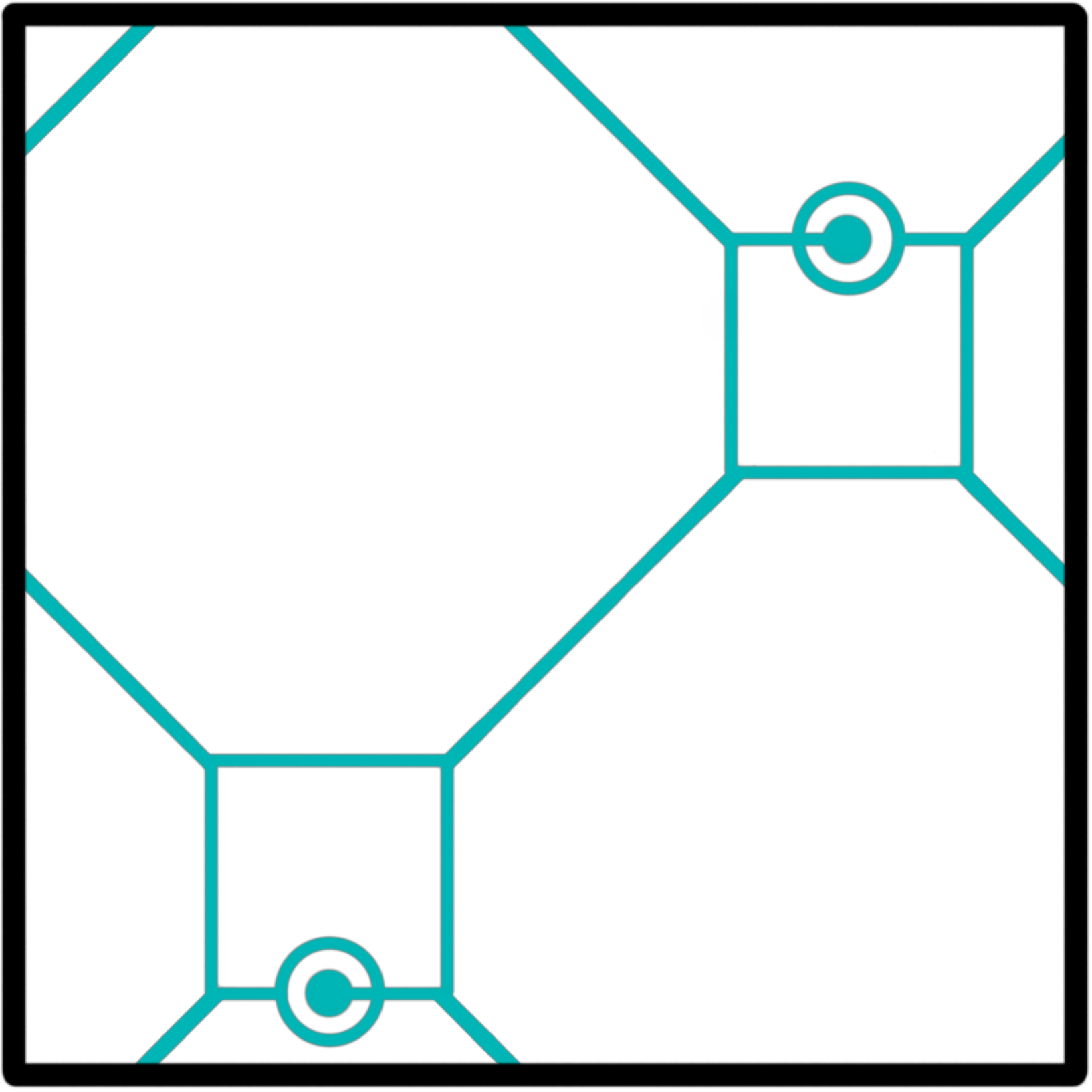}
        \hspace{0.2cm}
        \includegraphics[width=0.29\textwidth]{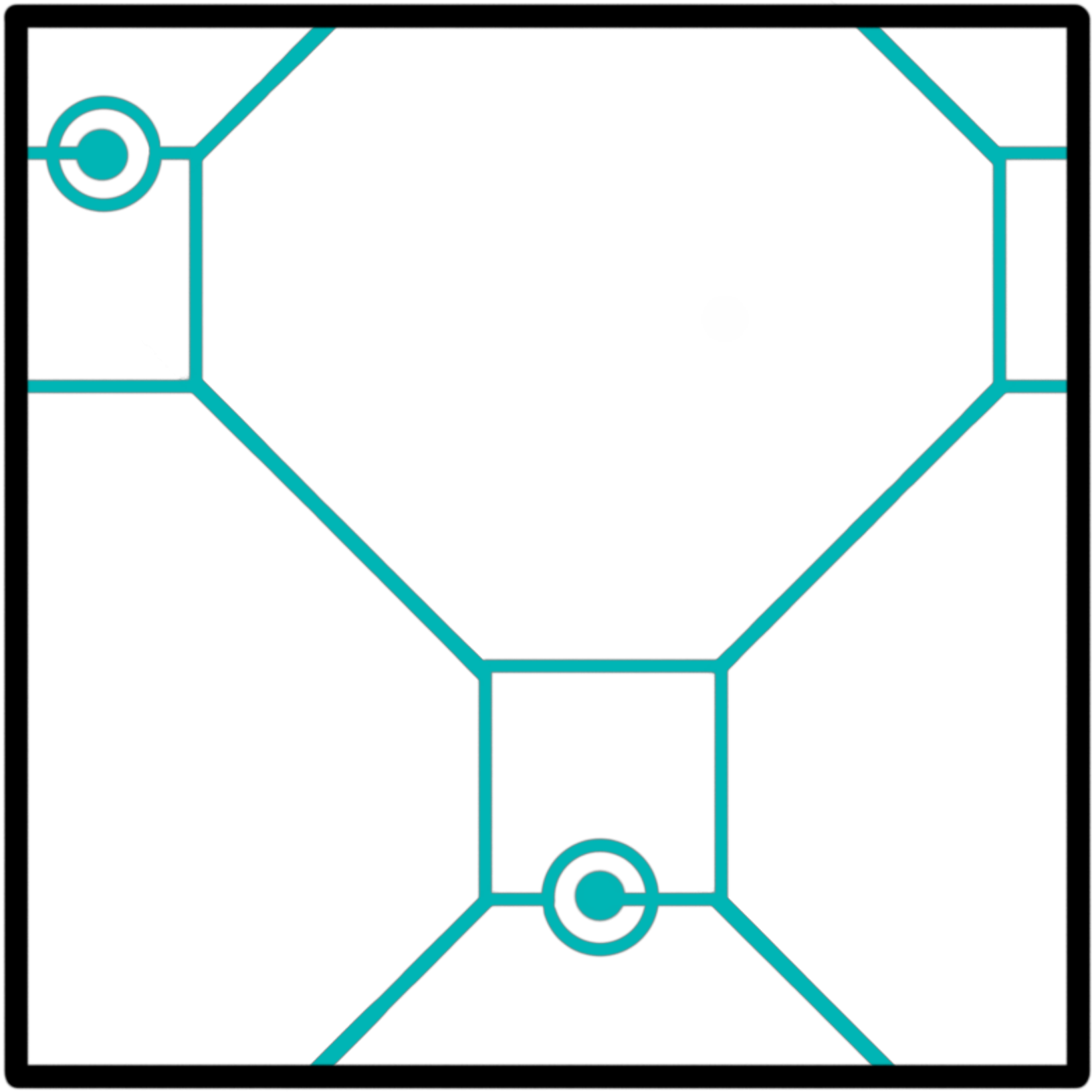}
        \hspace{0.2cm}
        \includegraphics[width=0.29\textwidth]{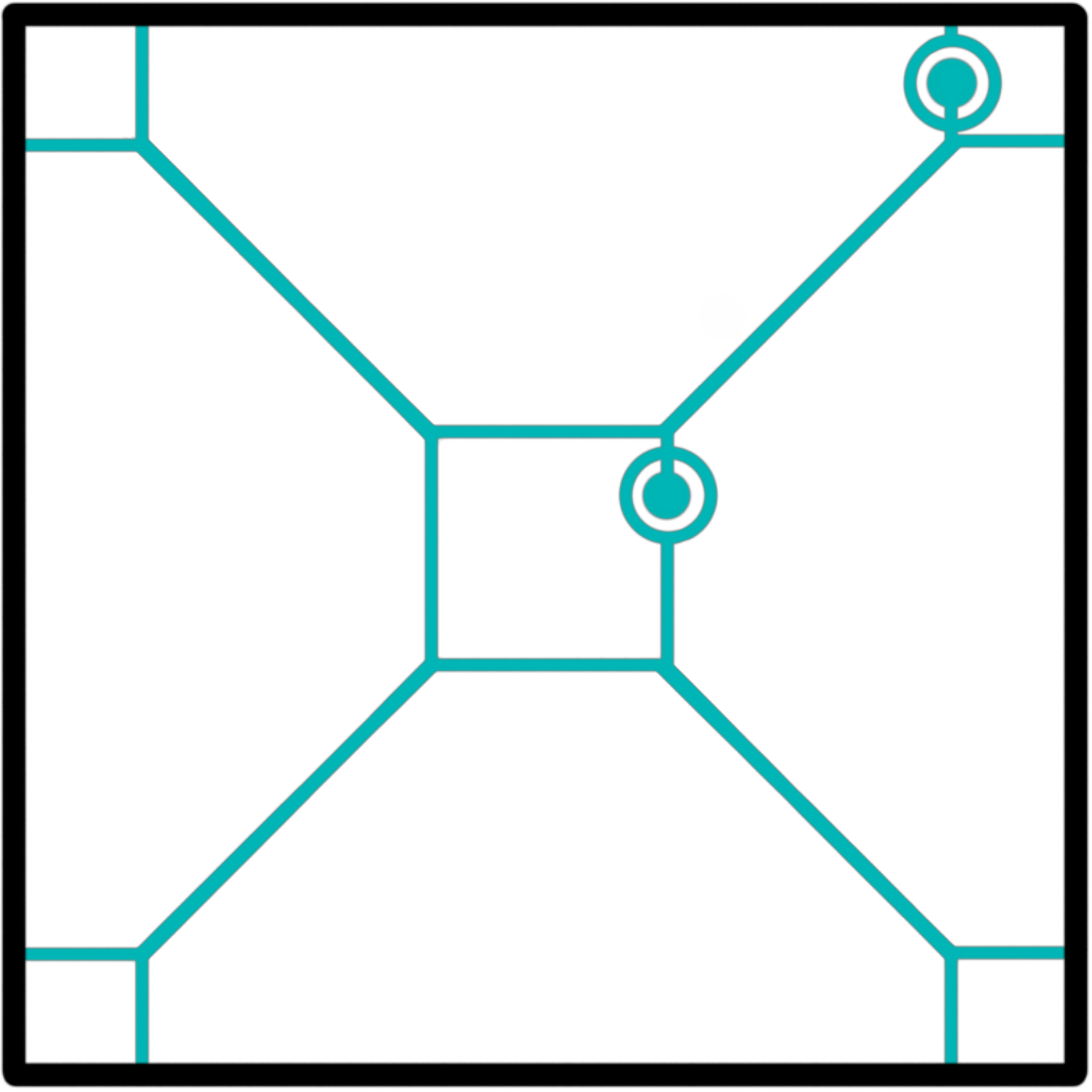}
        \caption{A tridiagram of \textbf{srs}.}
        \label{fig:srs_tridia}
    \end{subfigure}
    
    \vskip\baselineskip
    
    \begin{subfigure}[b]{0.27\textwidth}
        \includegraphics[width=0.6\textwidth]{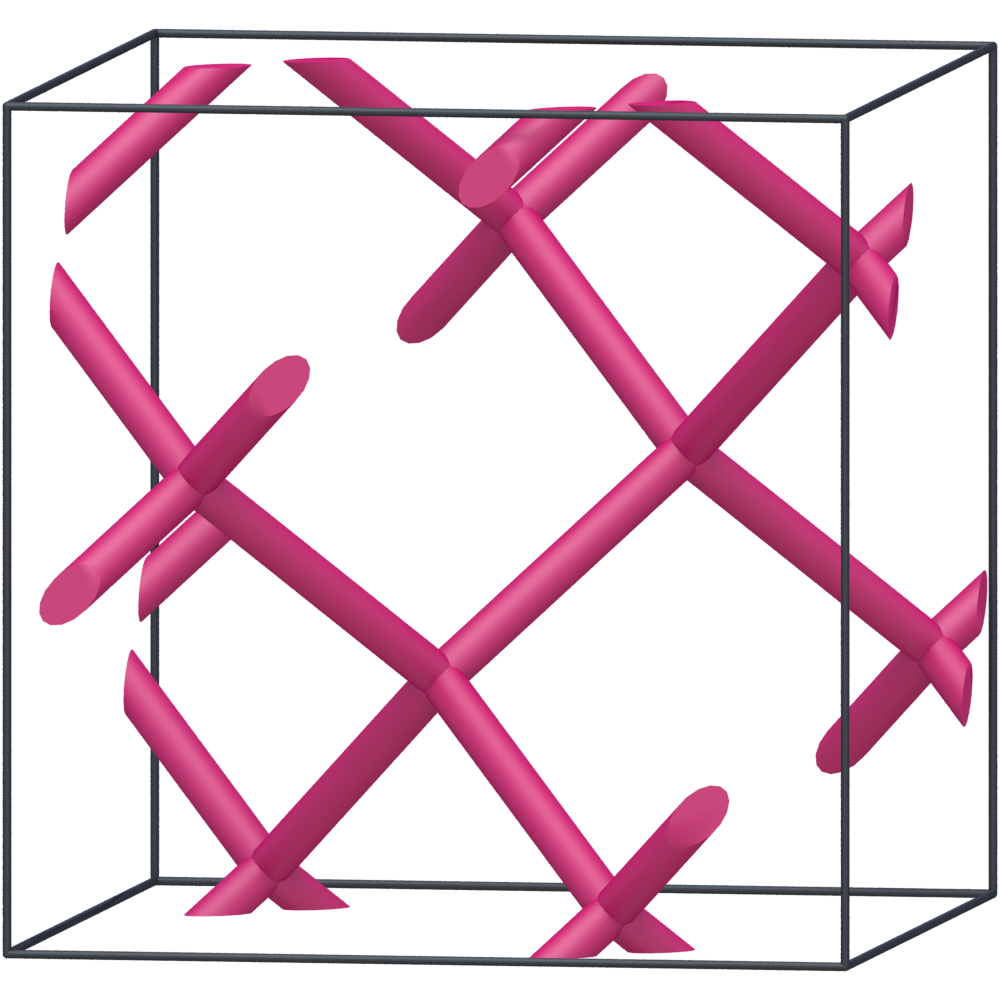}
        \caption{A unit cell of \textbf{dia}.}
        \label{fig:dia_uc}
    \end{subfigure}
    \begin{subfigure}[b]{0.54\textwidth}
        \includegraphics[width=0.29\textwidth]{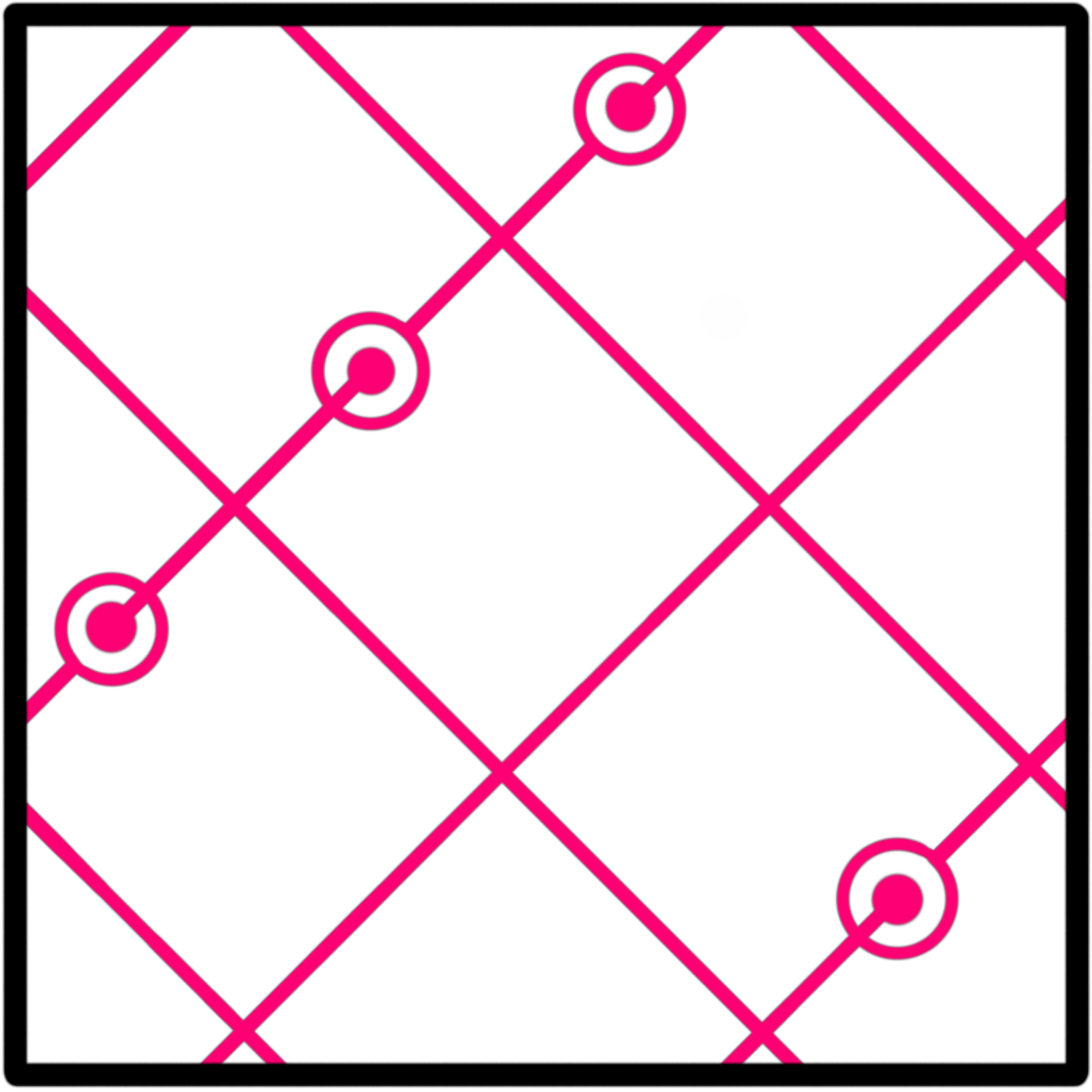}
        \hspace{0.2cm}
        \includegraphics[width=0.29\textwidth]{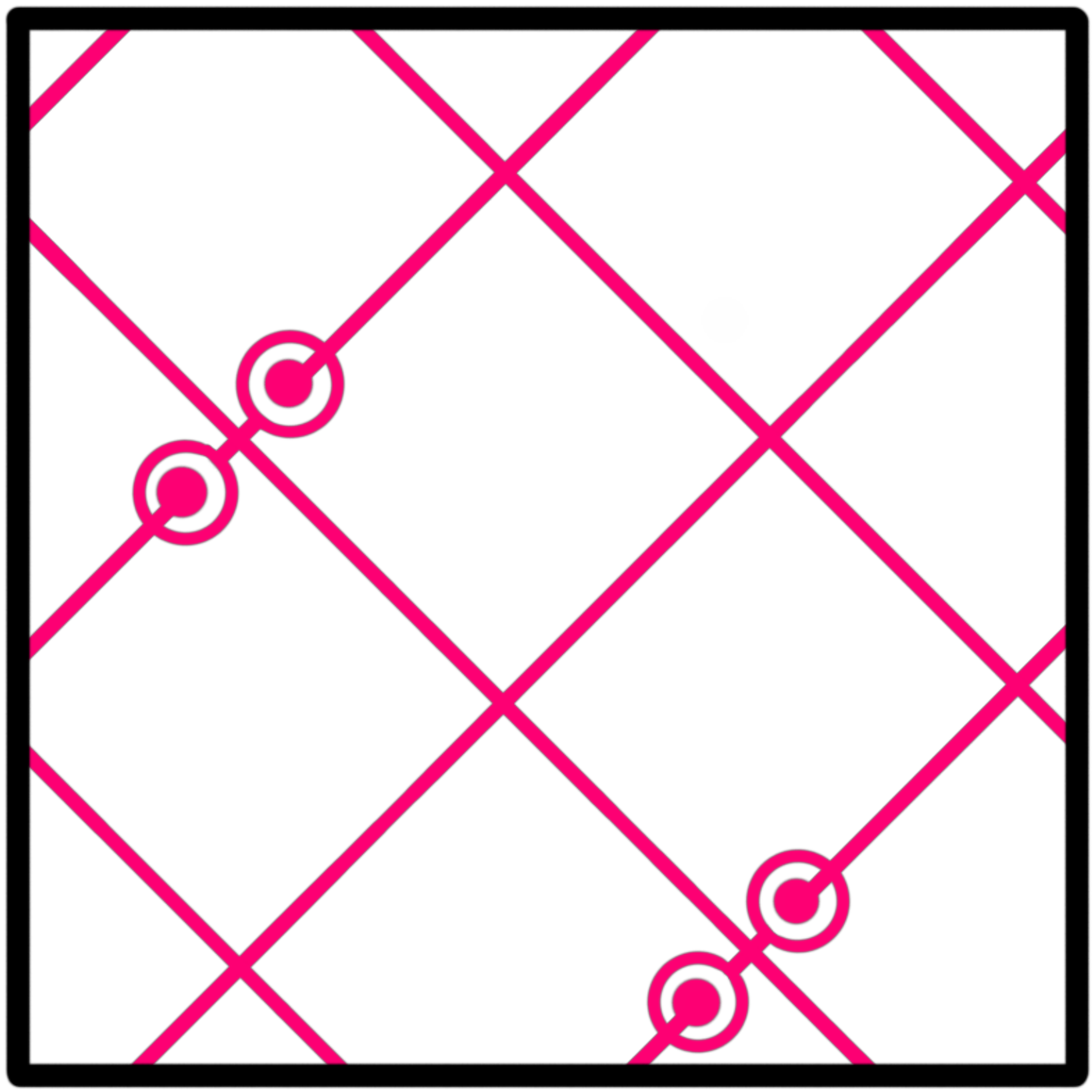}
        \hspace{0.2cm}
        \includegraphics[width=0.29\textwidth]{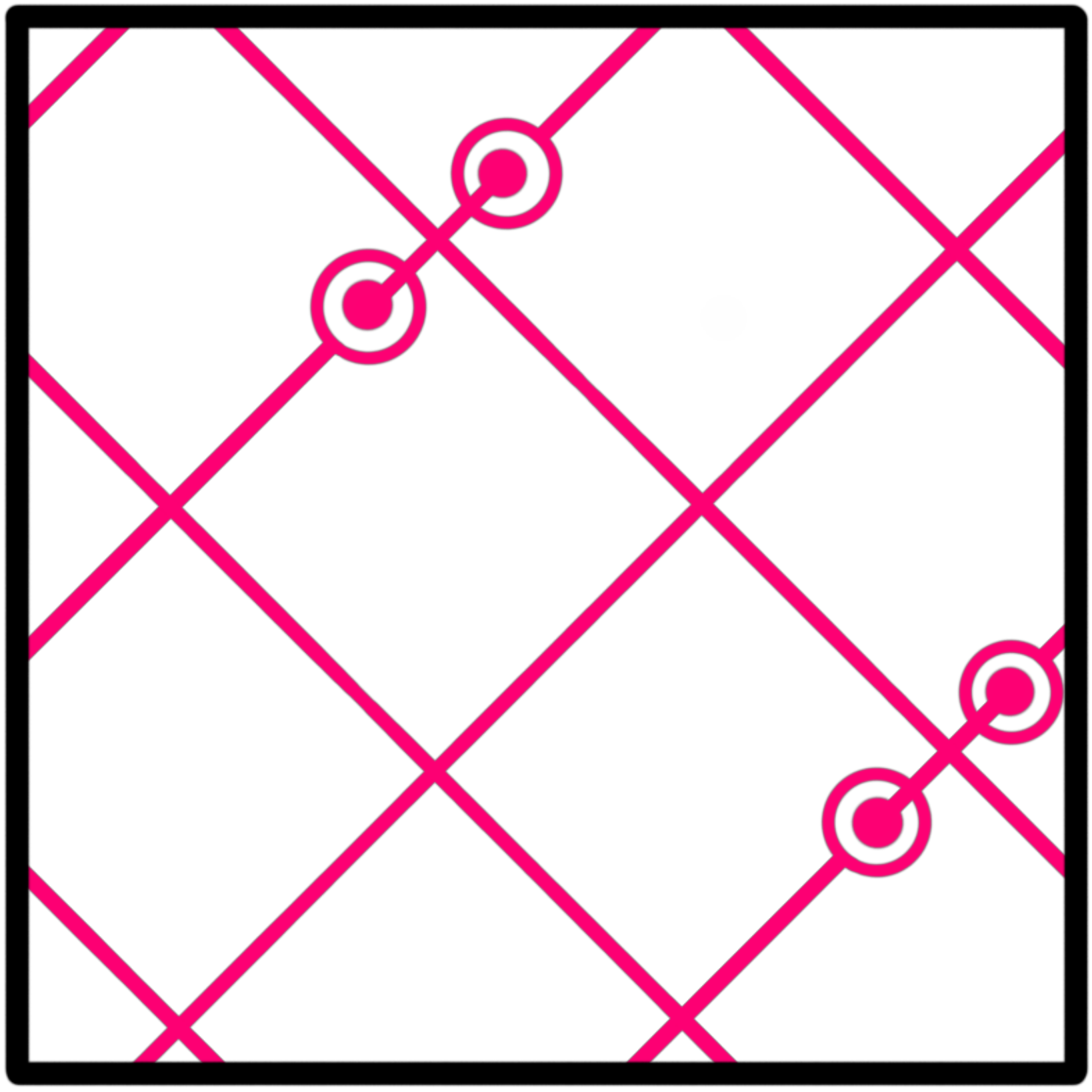}
        \caption{A tridiagram of \textbf{dia}.}
        \label{fig:dia_tridia}
    \end{subfigure}
    
    \vskip\baselineskip
    
    \begin{subfigure}[b]{0.27\textwidth}
        \includegraphics[width=0.6\textwidth]{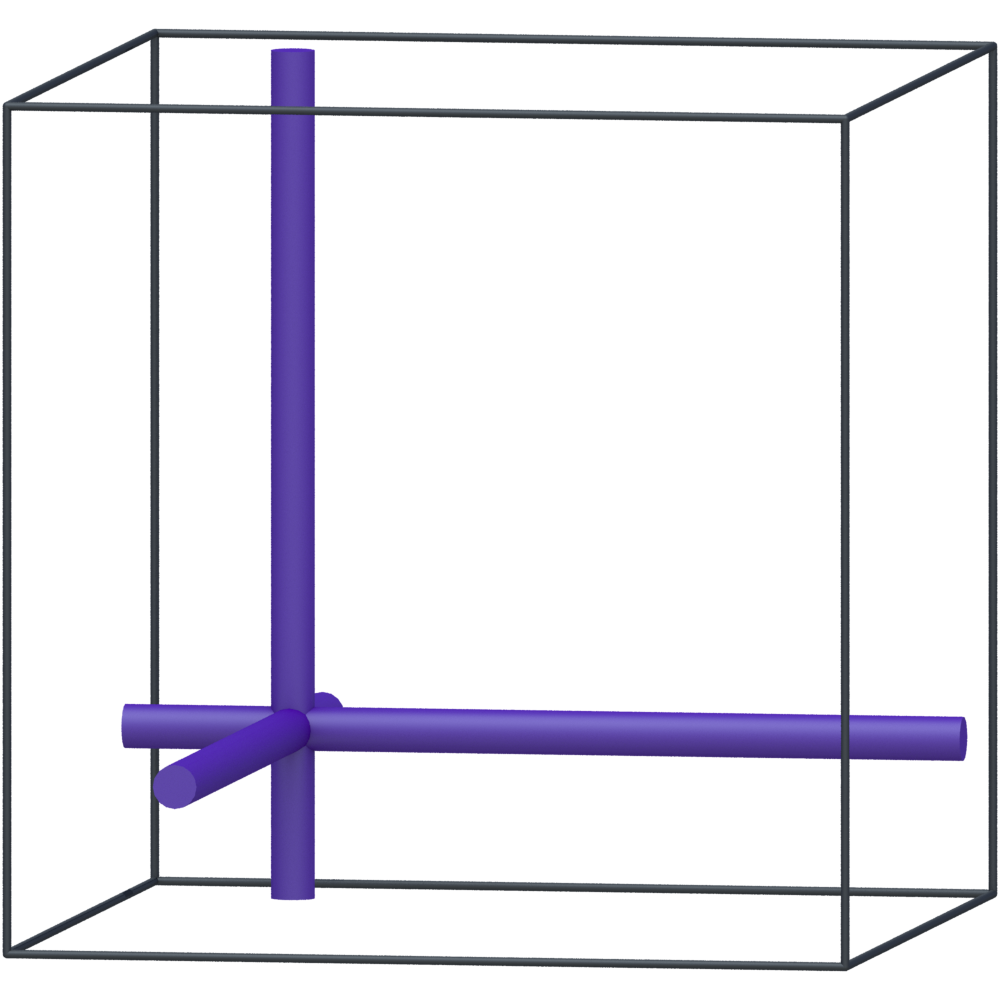}
        \caption{A unit cell of \textbf{pcu}.}
        \label{fig:pcu_uc}
    \end{subfigure}
    \begin{subfigure}[b]{0.54\textwidth}
        \includegraphics[width=0.29\textwidth]{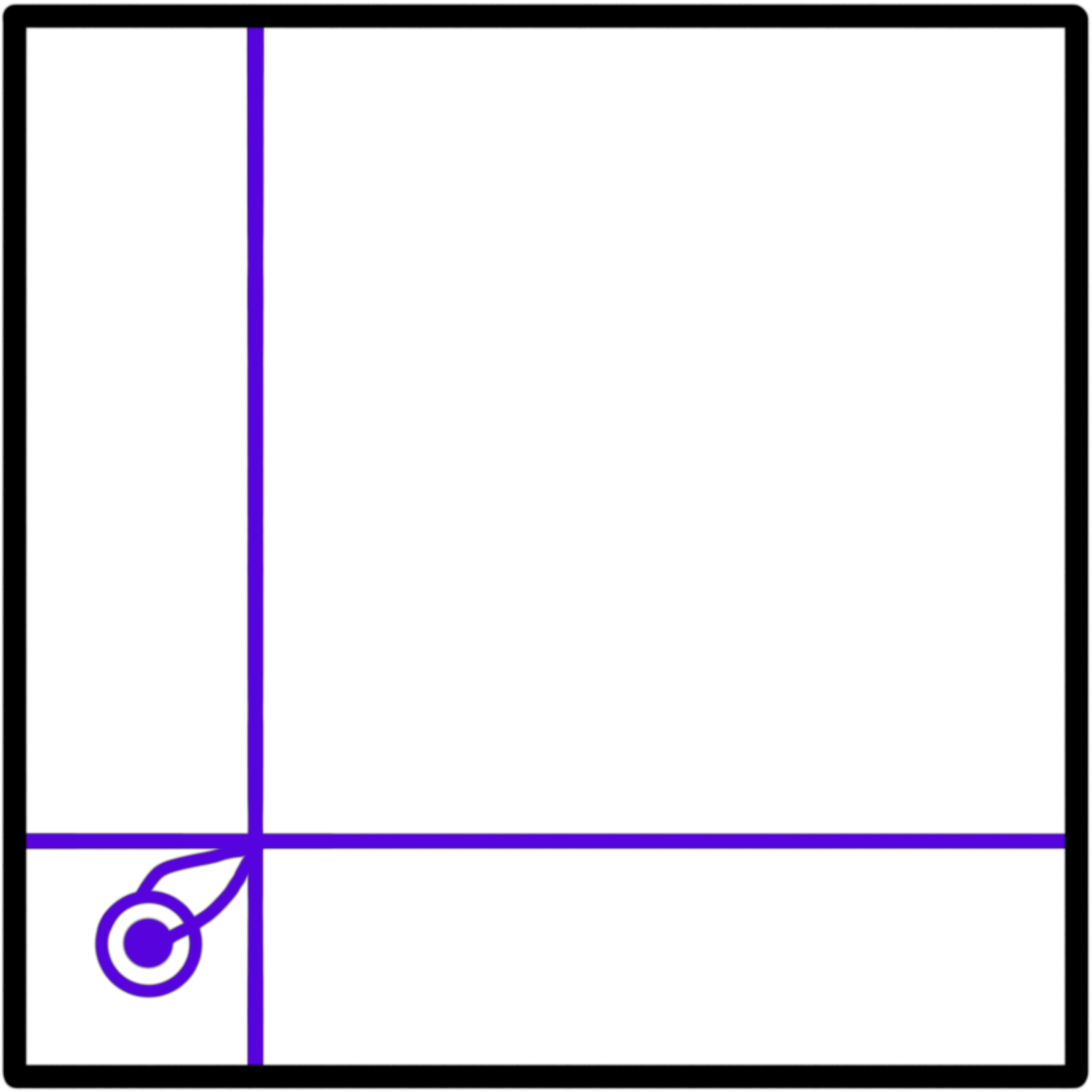}
        \hspace{0.2cm}
        \includegraphics[width=0.29\textwidth]{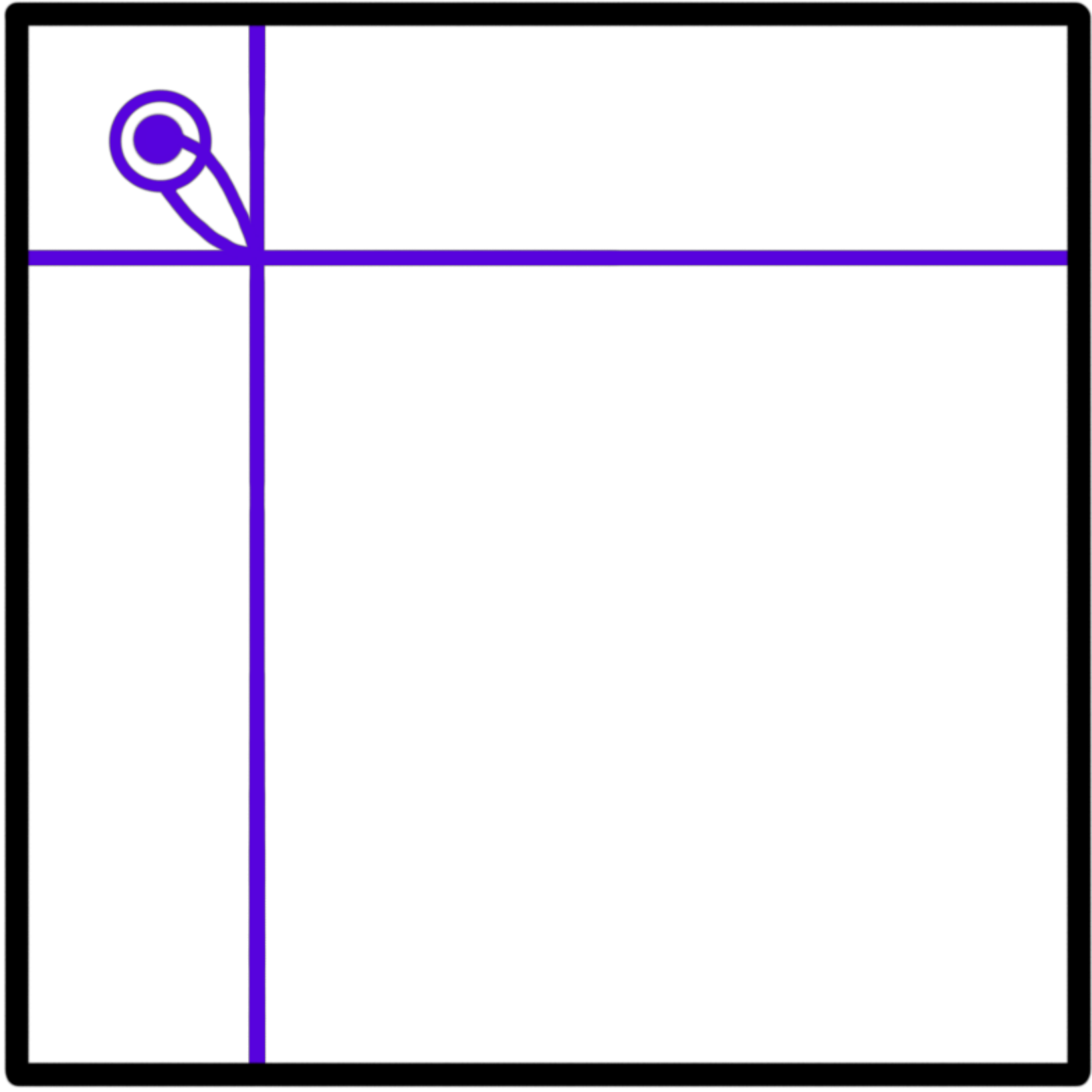}
        \hspace{0.2cm}
        \includegraphics[width=0.29\textwidth]{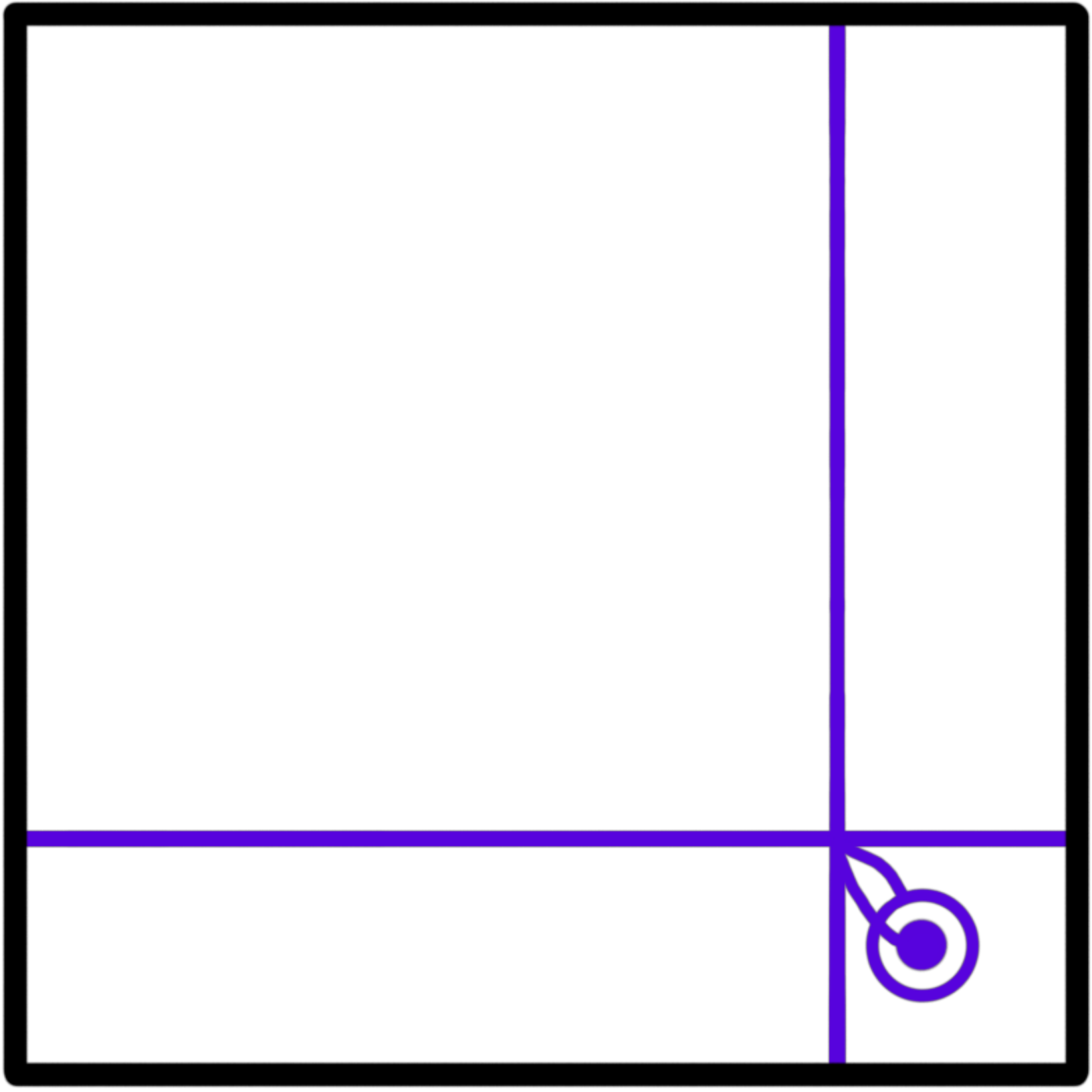}
        \caption{A tridiagram of \textbf{pcu}.}
        \label{fig:pcu_tridia}
    \end{subfigure}

    \vskip\baselineskip
    
    \begin{subfigure}[b]{0.27\textwidth}
        \includegraphics[width=0.6\textwidth]{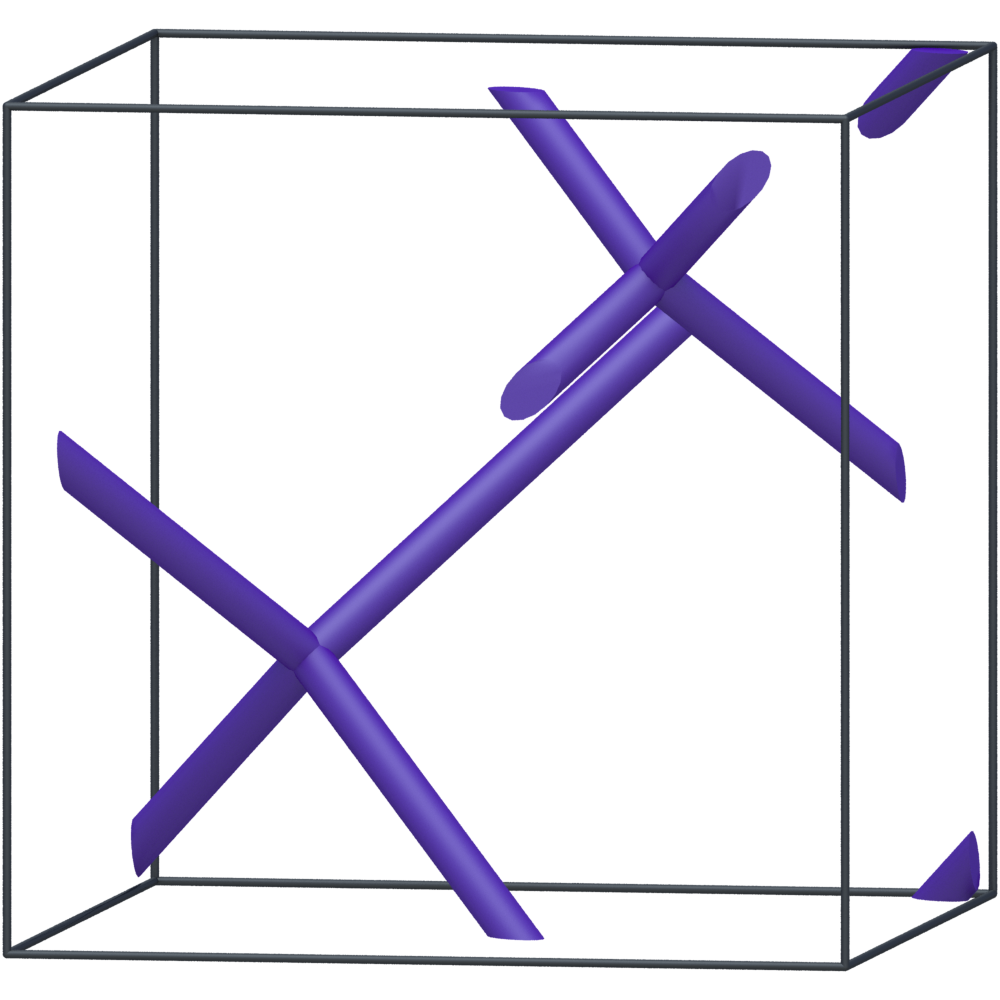}
        \caption{A unit cell of \textbf{dia-c}.}
        \label{fig:dia-c_uc}
    \end{subfigure}
    \begin{subfigure}[b]{0.54\textwidth}
        \includegraphics[width=0.29\textwidth]{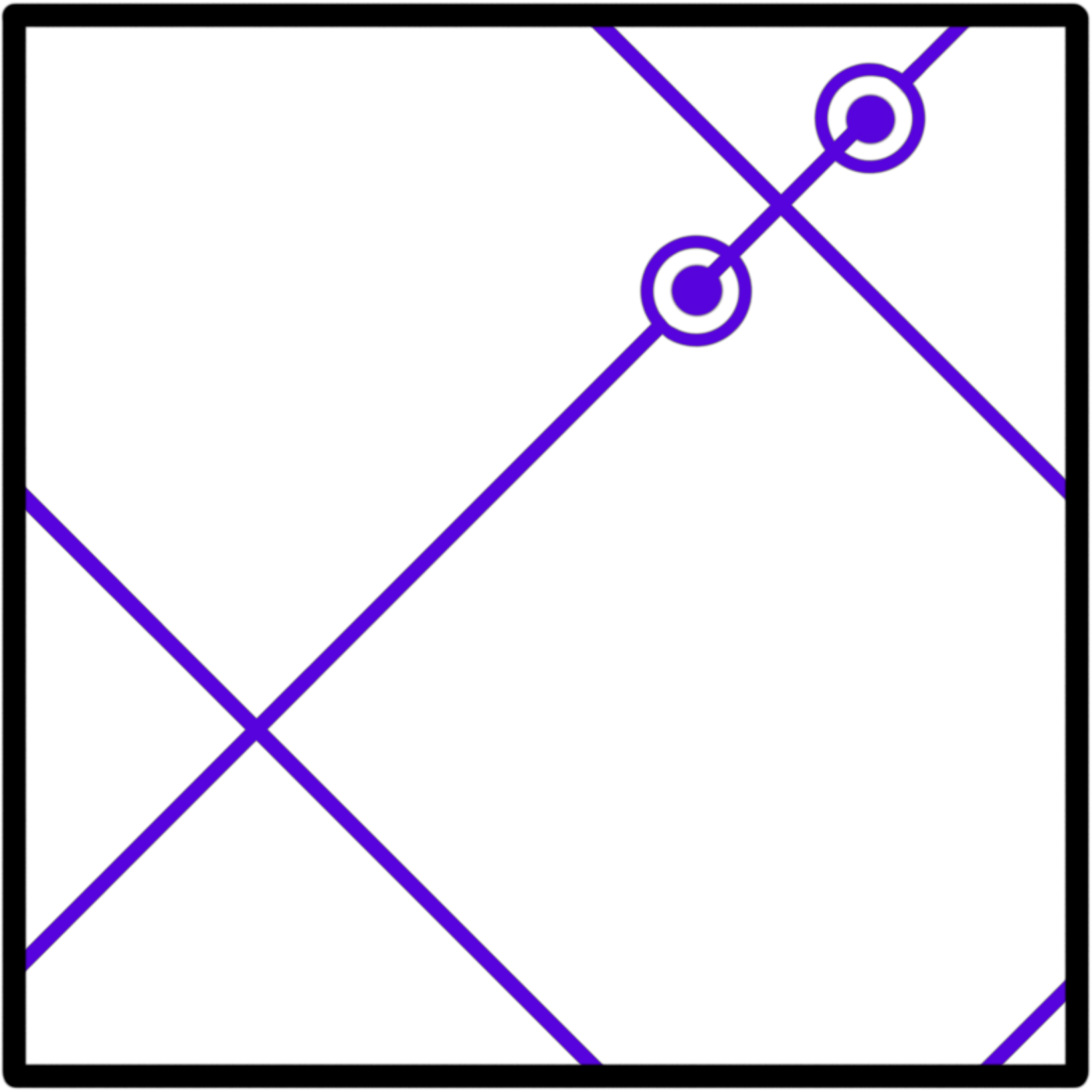}
        \hspace{0.2cm}
        \includegraphics[width=0.29\textwidth]{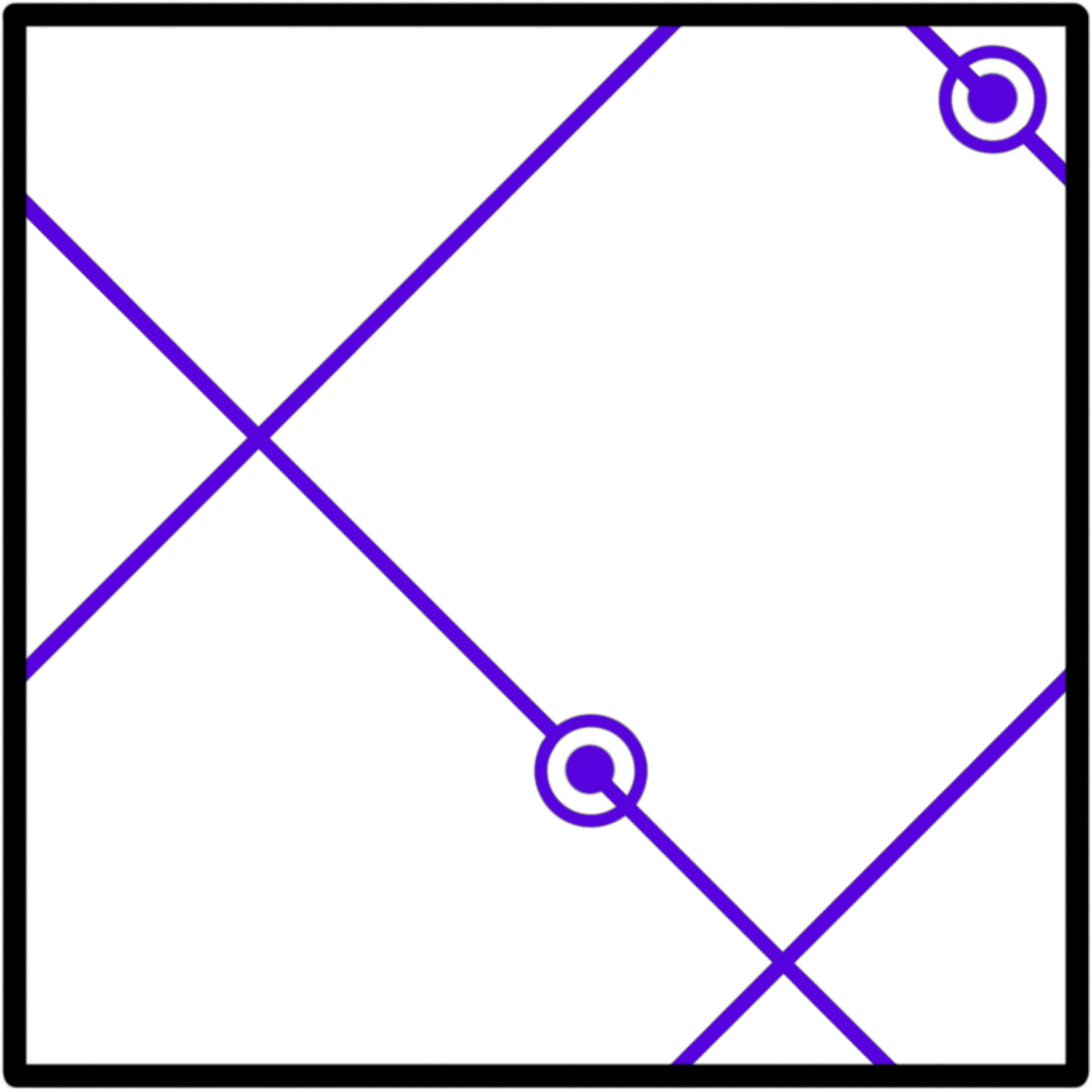}
        \hspace{0.2cm}
        \includegraphics[width=0.29\textwidth]{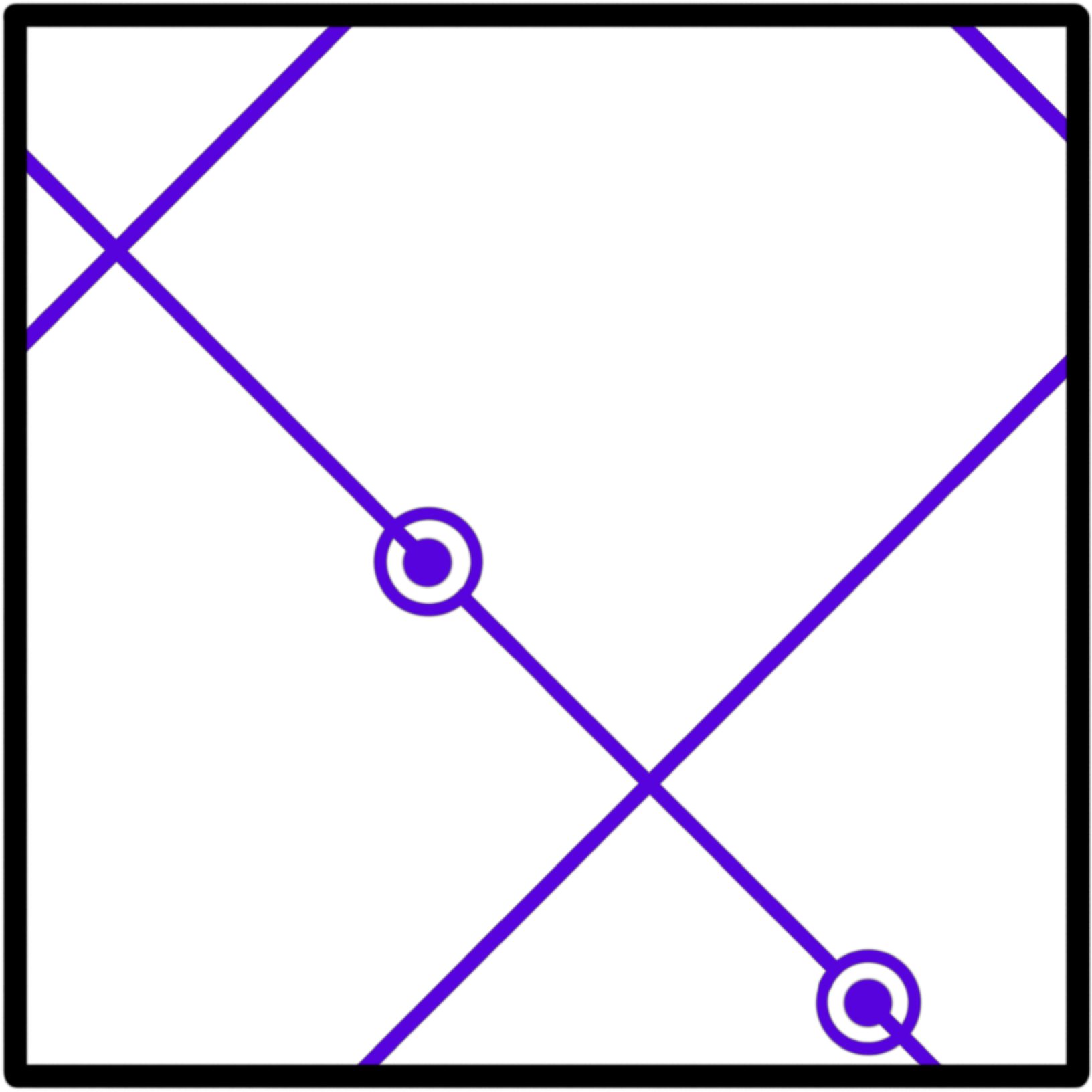}
        \caption{A tridiagram of \textbf{dia-c}.}
        \label{fig:dia-c_tridia}
    \end{subfigure}
    
    \caption{Examples of ground states: From top to bottom, the unit cells and tridiagrams of \textbf{sql}, \textbf{srs}, \textbf{dia}, \textbf{pcu} and \textbf{dia-c}. The tridiagrams are obtained from projections from the front, top and right sides of the unit cells. These structures are ground states with respect to these unit cells because their minimum crossing number triplet is $(0,0,0)$. The \textbf{sql} network shows that the concept of a ground state is readily defined for lower-periodic structures. The example of \textbf{dia-c}, whose strong rings are linked and whose HRN is \textbf{hxg}, shows that neither the knots and links present in the cycles of an embedding nor its HRN seems to provide information about the crossings of the embedding.}
    \label{fig:examples_of_ground_states}
\end{figure*}

Contrary to the case of 3-periodic tangles however \cite{andriamanalina2025untanglingnumber3periodictangles}, it is in general difficult to characterise ground states of 3-periodic graphs. For example, it is reasonable to believe that \textbf{srs-c} and \textbf{srs-c*} are ground states in their respective $\mathcal{U}$-families. However, since their respective minimum crossing number triplets are $(4,4,4)$ and $(4,4,8)$ with respect to the unit cells displayed in Fig. \ref{fig:interpenetrating_enantiomorphic_srs_nets_uc_xyz} and Fig. \ref{fig:interpenetrating_LH_srs_nets_uc_xyz}, one cannot definitively claim that no other embeddings with fewer crossings exist. Nevertheless, changing the over-under information of any crossings in any diagrams of these structures would only increase the total number of crossings, suggesting that these are indeed ground states. Further discussion on this matter is given in Sect. \ref{sec:computability}.

The fact that \textbf{dia-c}, \textbf{srs-c} and \textbf{srs-c*} are ground states, coincides with some other characterisations of least tangled embeddings. Indeed, \textbf{dia-c} and \textbf{srs-c} are both obtained from the intergrowth of the barycentric embeddings of, respectively, \textbf{dia} and \textbf{srs}, and their duals, and both are also of maximum symmetry. The three structures all also minimise the ropelength energy among all embeddings of their respective $\mathcal{U}$-families, which coincides with the proposition for least tangled embeddings given in \cite{myf2011_entanglement_graphs}.

\begin{figure*}[hbtp]

    \centering
    
    \begin{subfigure}[b]{0.27\textwidth}
        \includegraphics[width=0.6\textwidth]{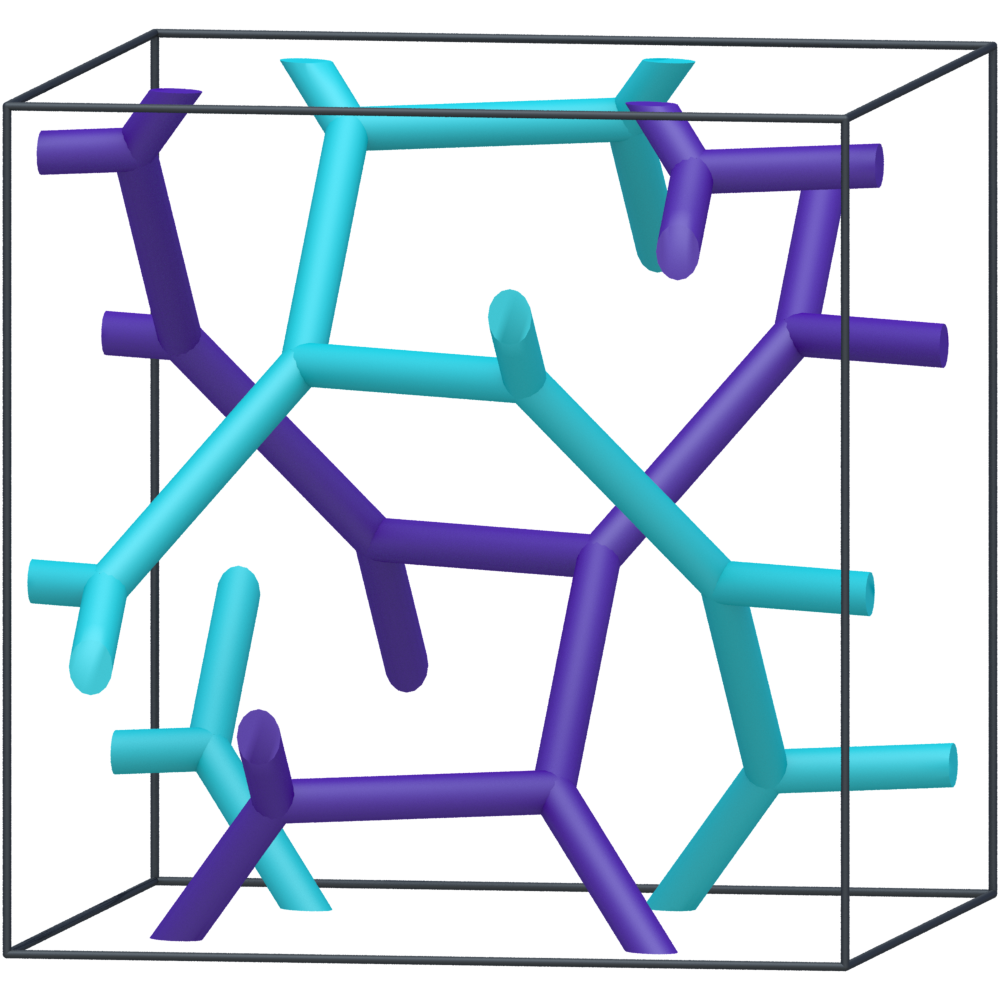}
        \caption{A unit cell of \textbf{utp-c*}.}
        \label{fig:utp-c_star_uc}
    \end{subfigure}
    \begin{subfigure}[b]{0.54\textwidth}
        \includegraphics[width=0.29\textwidth]{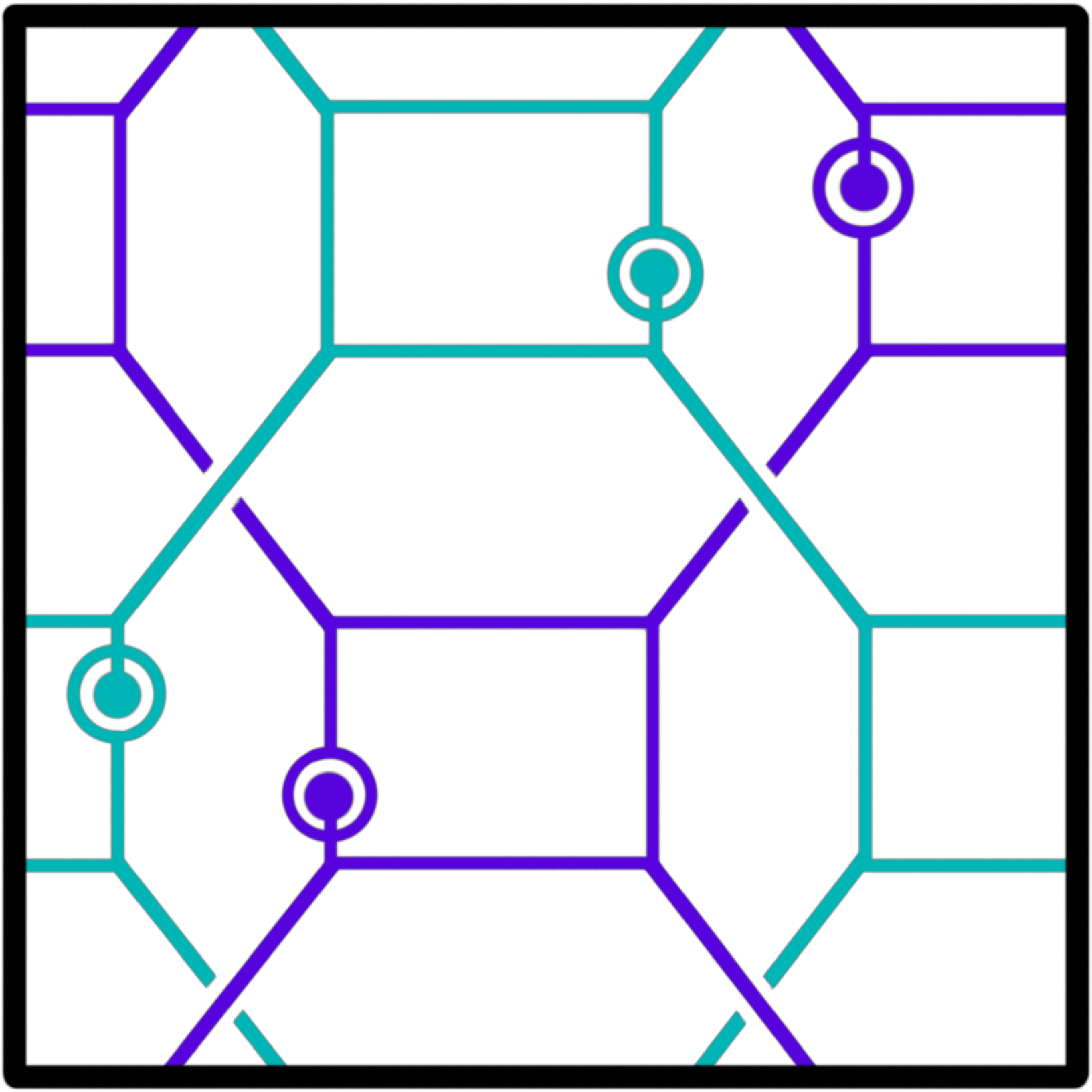}
        \hspace{0.2cm}
        \includegraphics[width=0.29\textwidth]{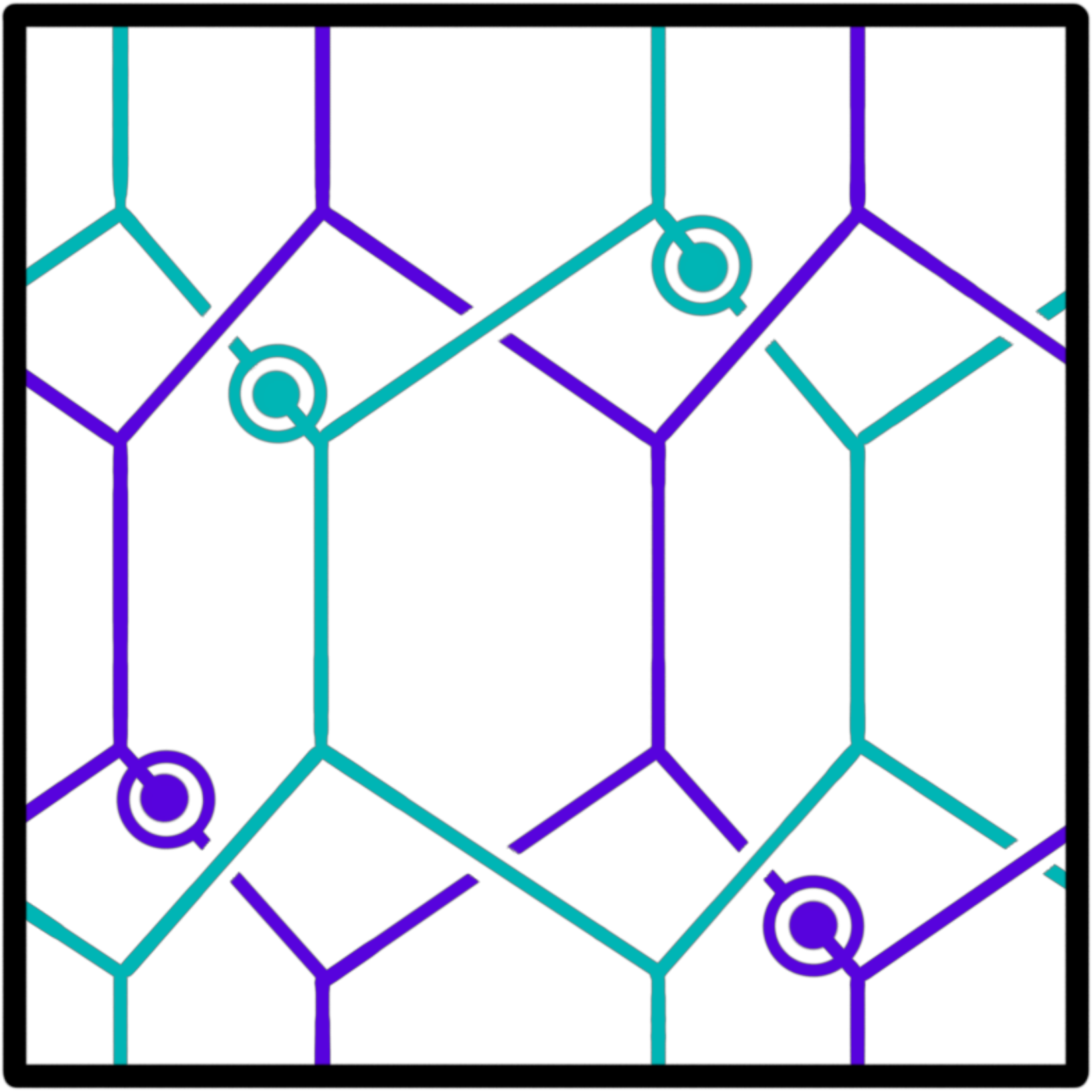}
        \hspace{0.2cm}
        \includegraphics[width=0.29\textwidth]{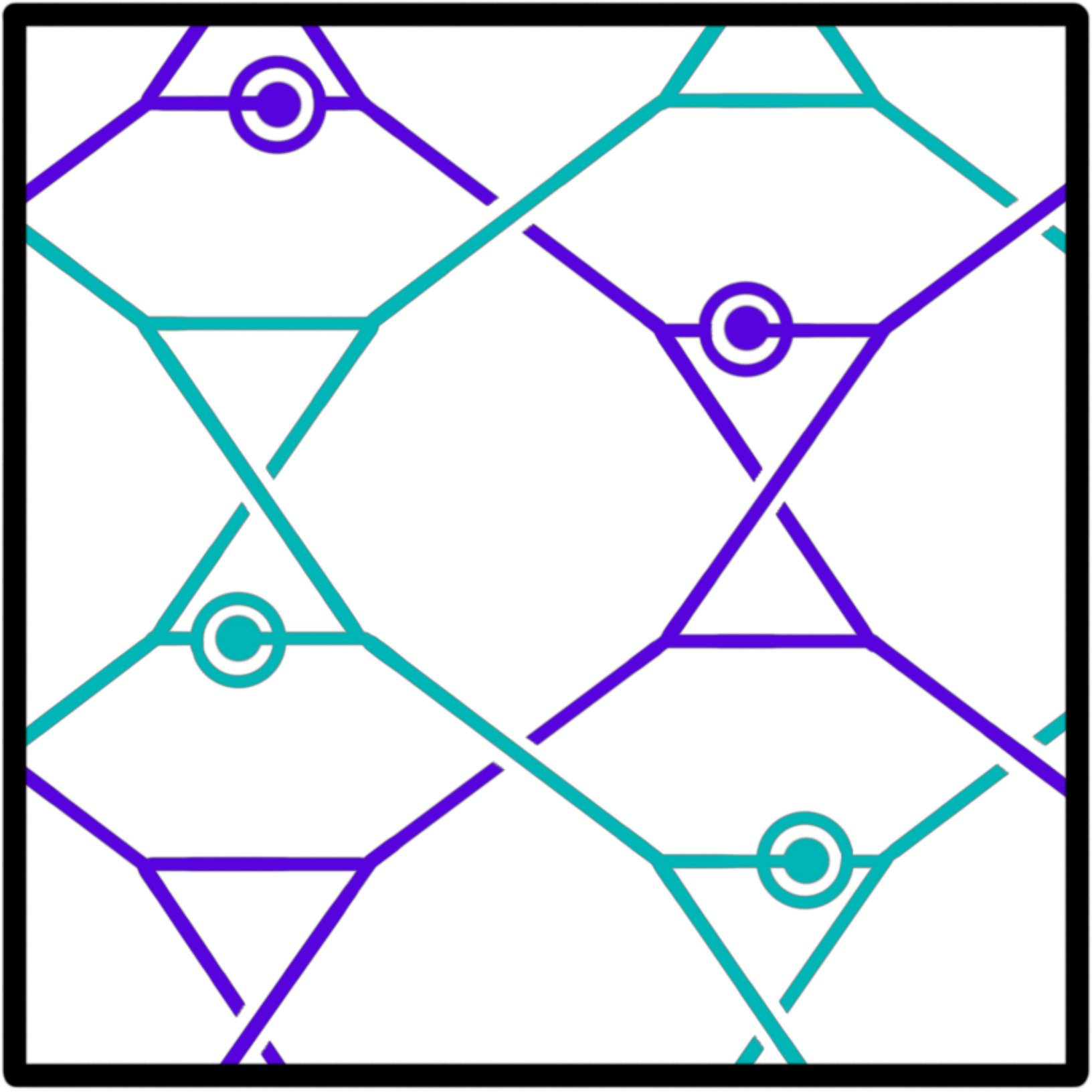}
        \caption{A tridiagram of \textbf{utp-c*}.}
        \label{fig:utp-c_star_tridia}
    \end{subfigure}

    \vskip\baselineskip
    
    \begin{subfigure}[b]{0.27\textwidth}
        \includegraphics[width=0.6\textwidth]{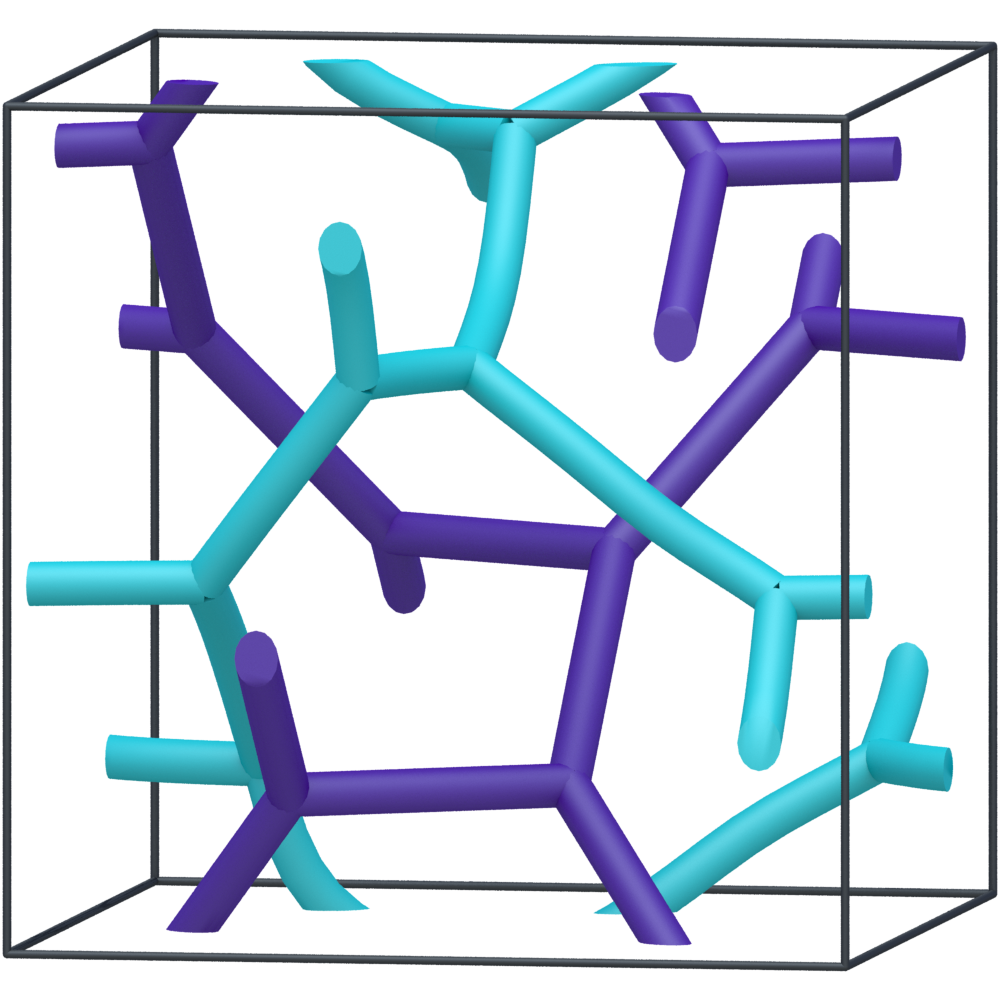}
        \caption{A unit cell of \textbf{utp-c**}.}
        \label{fig:utp-c_star_star_min_c_n_uc}
    \end{subfigure}
    \begin{subfigure}[b]{0.54\textwidth}
        \includegraphics[width=0.29\textwidth]{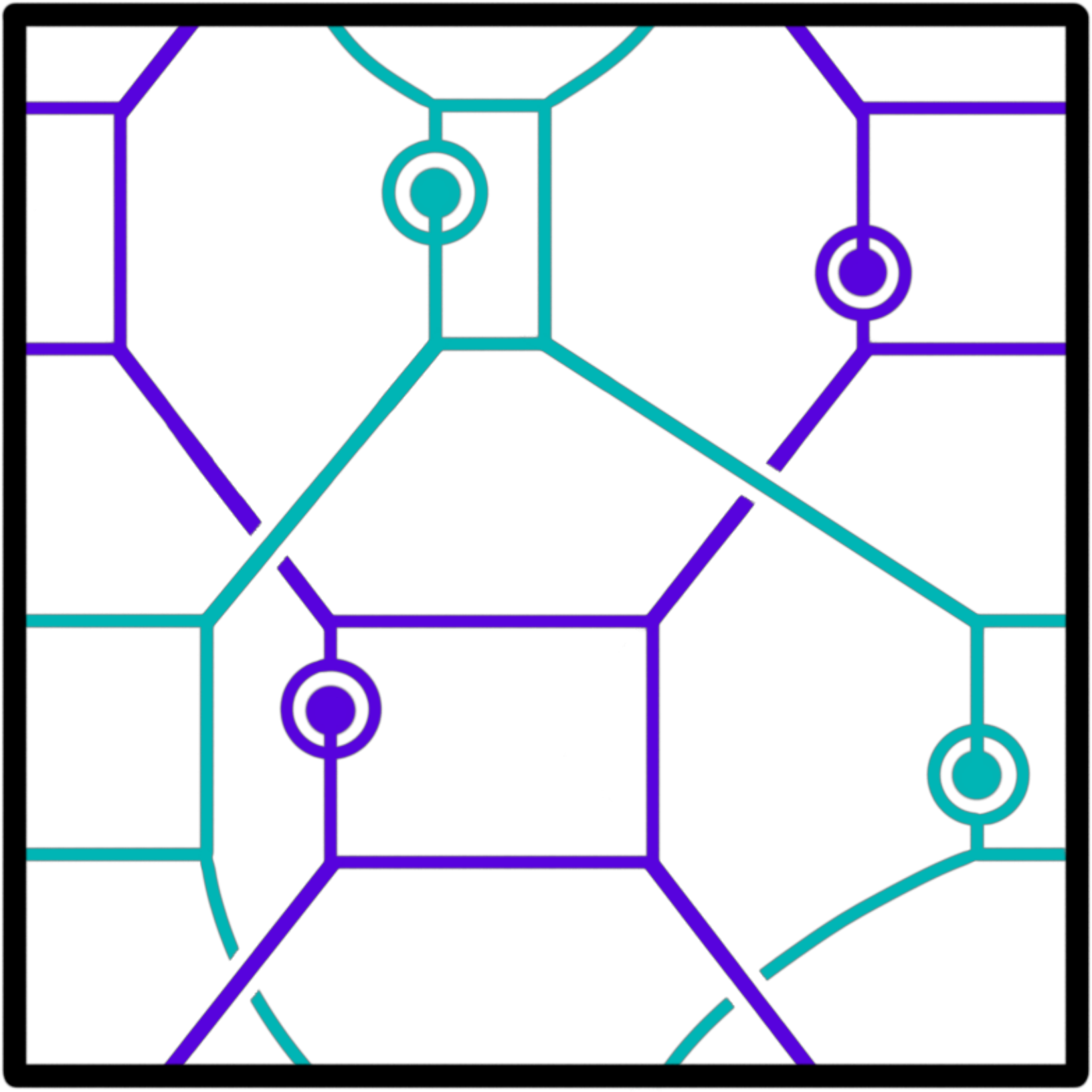}
        \hspace{0.2cm}
        \includegraphics[width=0.29\textwidth]{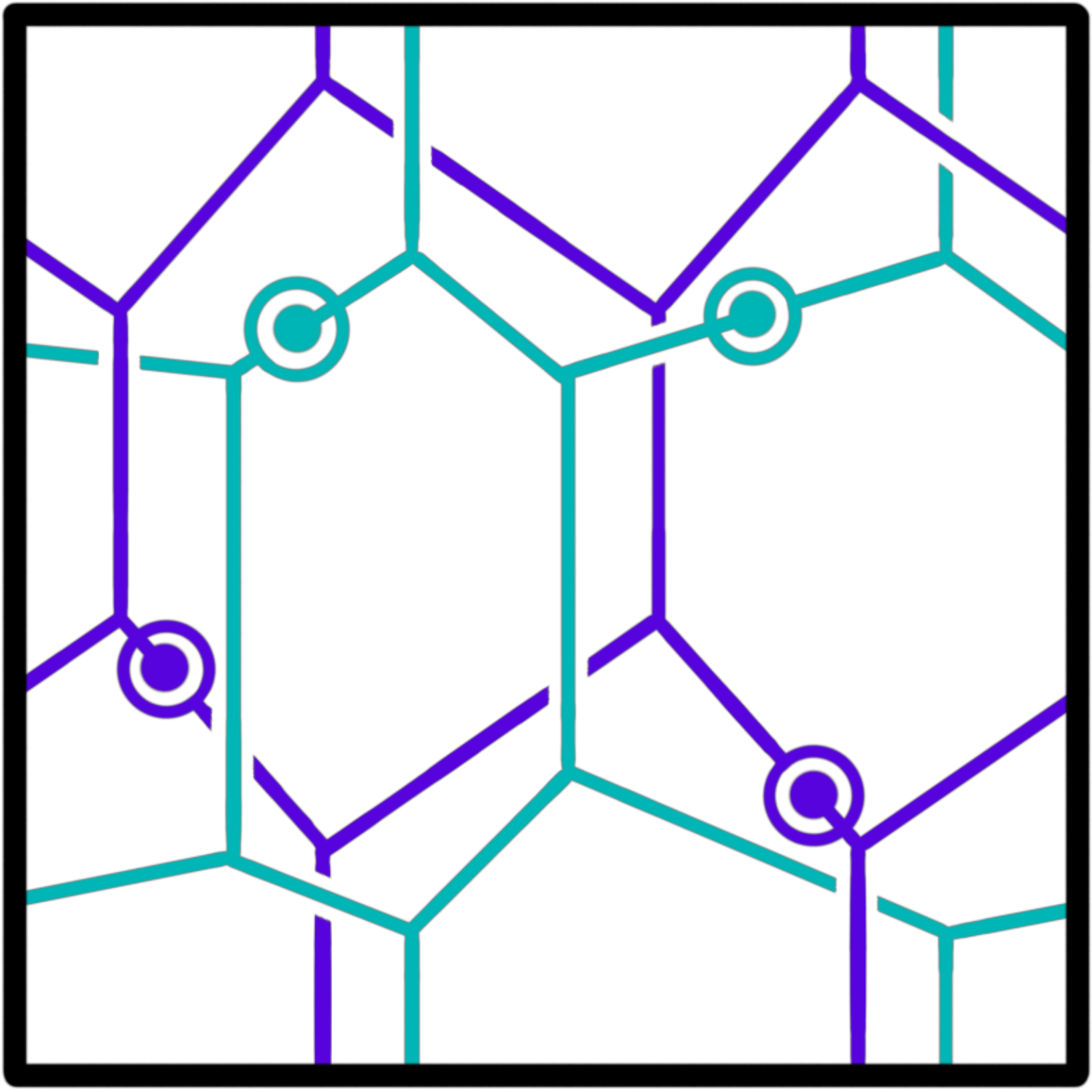}
        \hspace{0.2cm}
        \includegraphics[width=0.29\textwidth]{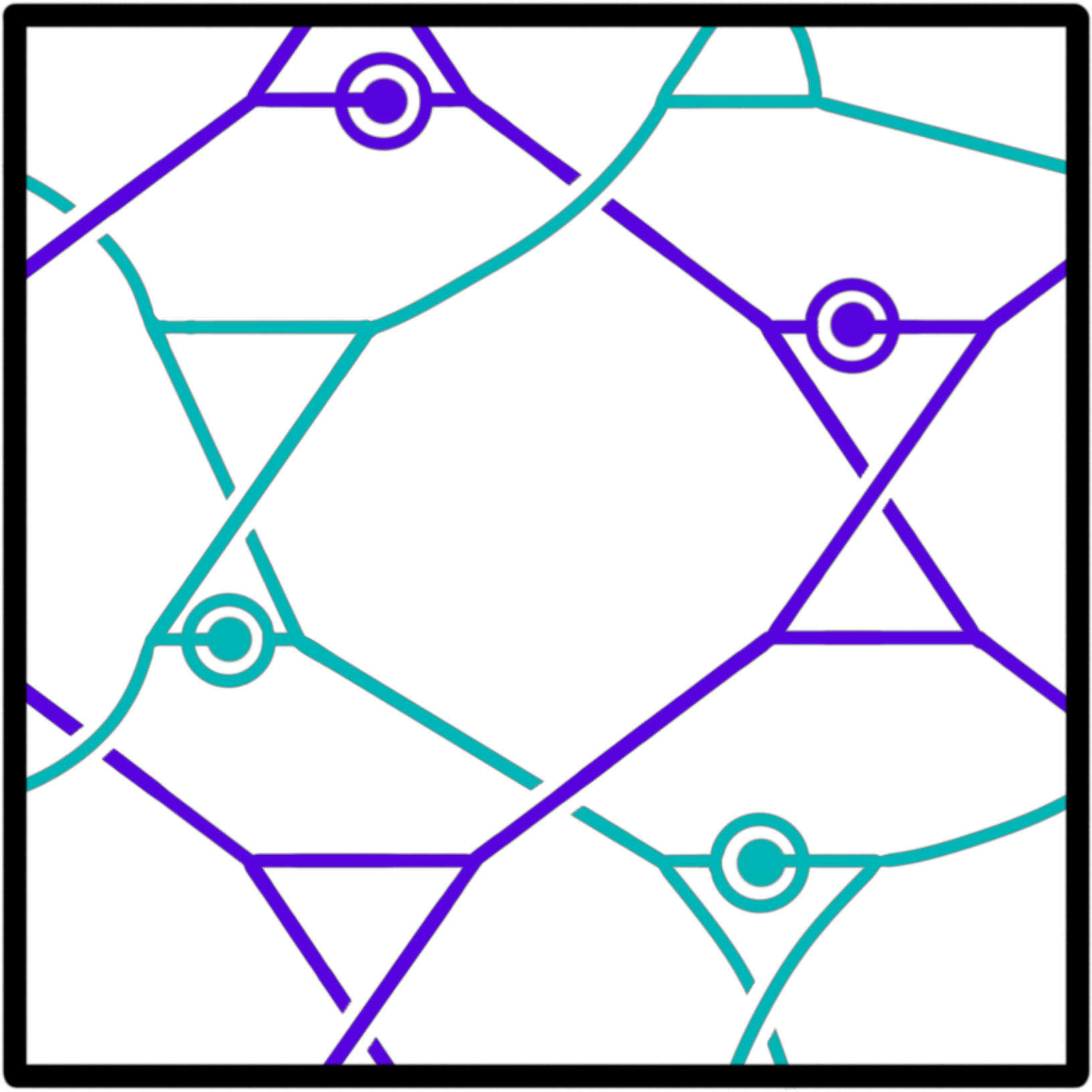}
        \caption{A tridiagram of \textbf{utp-c**}.}
        \label{fig:utp-c_star_star_tridia}
    \end{subfigure}
    \caption{Examples of ground states belonging to the same $\mathcal{G}$-family: With respect to the unit cell displayed in (a) and (c), the minimum crossing number triplet of both \textbf{utp-c*} and \textbf{utp-c**} is $(4,8,8)$, as can be seen with the tridiagrams shown in (b) and (d). These two structures likely have the least crossing number among all embeddings belonging to their common $\mathcal{U}$-family, suggesting that they are both ground states.}
    \label{fig:utp-c_star_and_utp-c_star_star_uc_and_tridia}
\end{figure*}

A given $\mathcal{G}$-family of ground states may indeed contain more than one element. An example is provided by the two embeddings of a graph comprising two \textbf{utp} networks described in \cite{Baburin:eo5056}, namely \textbf{utp-c*} and \textbf{utp-c**}. As noted in \cite{Baburin:eo5056}, these two structures have the interesting property of sharing the same HRN. More interestingly, they also have the same minimum crossing number triplet, $(4,8,8)$, with respect to the unit cells shown in Fig. \ref{fig:utp-c_star_and_utp-c_star_star_uc_and_tridia}. This minimum crossing number triplet, or more precisely the associated crossing number, is likely the least among all embeddings of the $\mathcal{U}$-family of the two structures, suggesting that they are both ground states. Besides \textbf{utp-c*} and \textbf{utp-c**}, networks that do not possess barycentric embeddings are other examples of structures that may possess several ground states. Indeed, we recall that for some networks, the computation of a barycentric embedding would result in vertices collapsing in the same position, or intersecting edges \cite{Delgado-Friedrichs:au5000}. In the latter case, an embedding can be obtained by moving one edge to either side of the other. The choice of side may result in potentially distinct embeddings that can all be ground states, provided that they have the same crossing number, the least in the $\mathcal{U}$-family.

\section{The untangling number, a measure of entanglement complexity}\label{sec:untangling_number}
The concept of a ground state described in Sect. \ref{sec:least_tangled_embeddings} allows us to define in the following a measure of entanglement complexity that we call the \textit{untangling number}, analogous to the unknotting number of classical knots. It is invariant under isotopy transformations applied to a chosen unit cell of a given embedding and can be minimised over all unit cells to give the \textit{minimum untangling number}, an invariant of the embedding. Further mathematical details on these notions are given in the Supplementary Information. \\

Consider an embedding $K$ of a graph $\mathcal{K}$ and a unit cell $U$ of $K$. Consider the family $\mathcal{U}(K,U)$ of $K$ as well as its subfamily $\mathcal{G}(K,U)$. We define \textit{untangling number of $K$ with respect to the unit cell $U$} to be the least number of crossing changes needed to transform $U$ to a unit cell of one element of $\mathcal{G}(K,U)$, that is, a ground state. It is invariant under isotopy transformations applied to $U$. Any ground state realising the untangling number is called a \textit{nearest ground state} of $K$.

We note that, although well defined, the computation of the untangling number is, in general, limited by the fact that only an upper bound can be given. This is similar to the limitation in computing the unknotting number of classical knots. Nevertheless, the computation provides an empirical understanding of the untangling number. This is discussed further in Sect. \ref{sec:computability}.

Fig. \ref{fig:untangling_dia-z} shows an example of the computation of an upper bound for the untangling number of a structure, that of \textbf{dia-z} with respect to its primitive unit cell. Since two crossing changes are needed to transform the unit cell to the primitive unit cell of \textbf{dia}, the untangling number is at most two. Another example is given in Fig. \ref{fig:untangling_ravelled_pcu}, where we untangle the ravelled embedding of \textbf{pcu} into its barycentric embedding. In this example, it is clear that one can apply several crossing changes to untangle the structure, but the computation shows that one single crossing change suffices to do so.

\begin{figure*}[hbtp]

    \centering
    
    \begin{subfigure}[b]{0.27\textwidth}
        \includegraphics[width=0.6\textwidth]{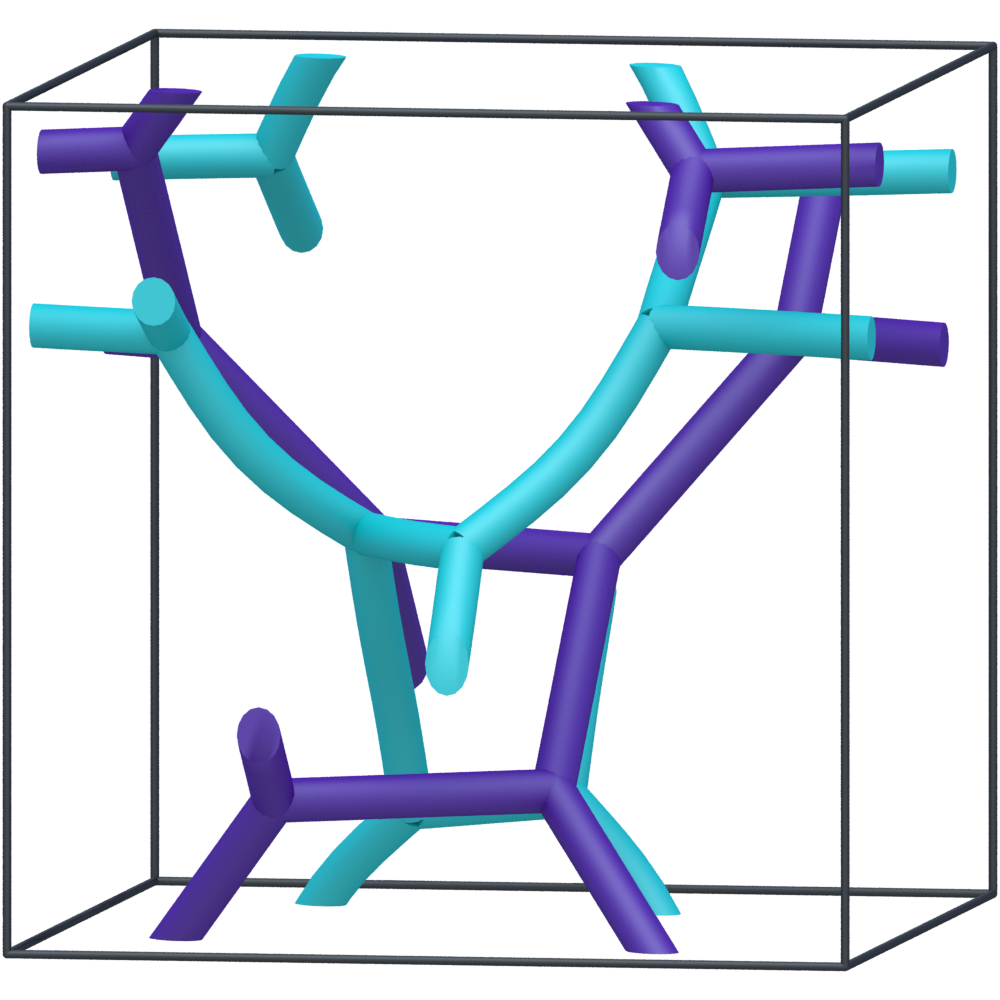}
        \caption{A unit cell of \textbf{utp-c***}.}
        \label{fig:utp-c_star_star_star_min_c_n_uc}
    \end{subfigure}
    \begin{subfigure}[b]{0.54\textwidth}
        \includegraphics[width=0.29\textwidth]{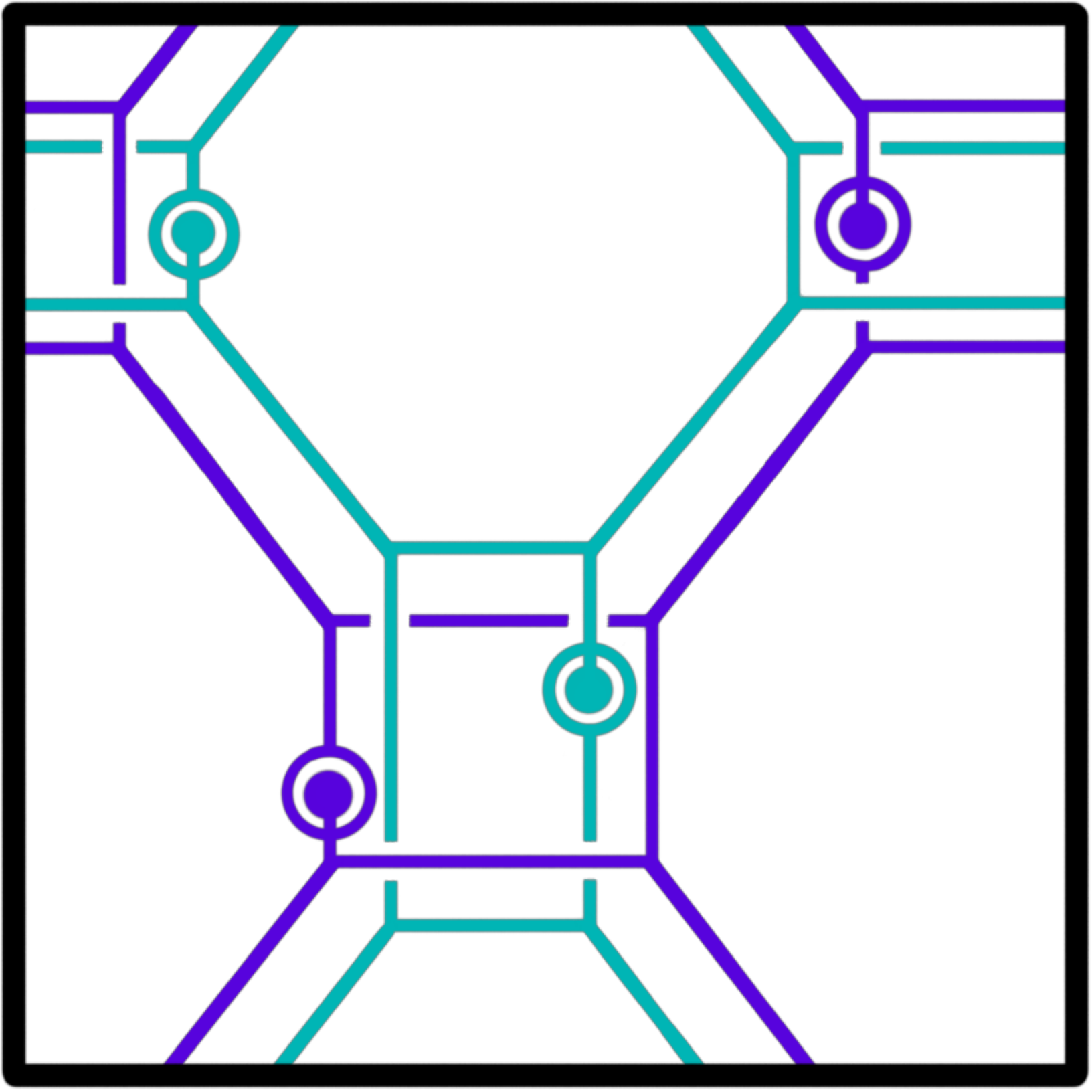}
        \hspace{0.2cm}
        \includegraphics[width=0.29\textwidth]{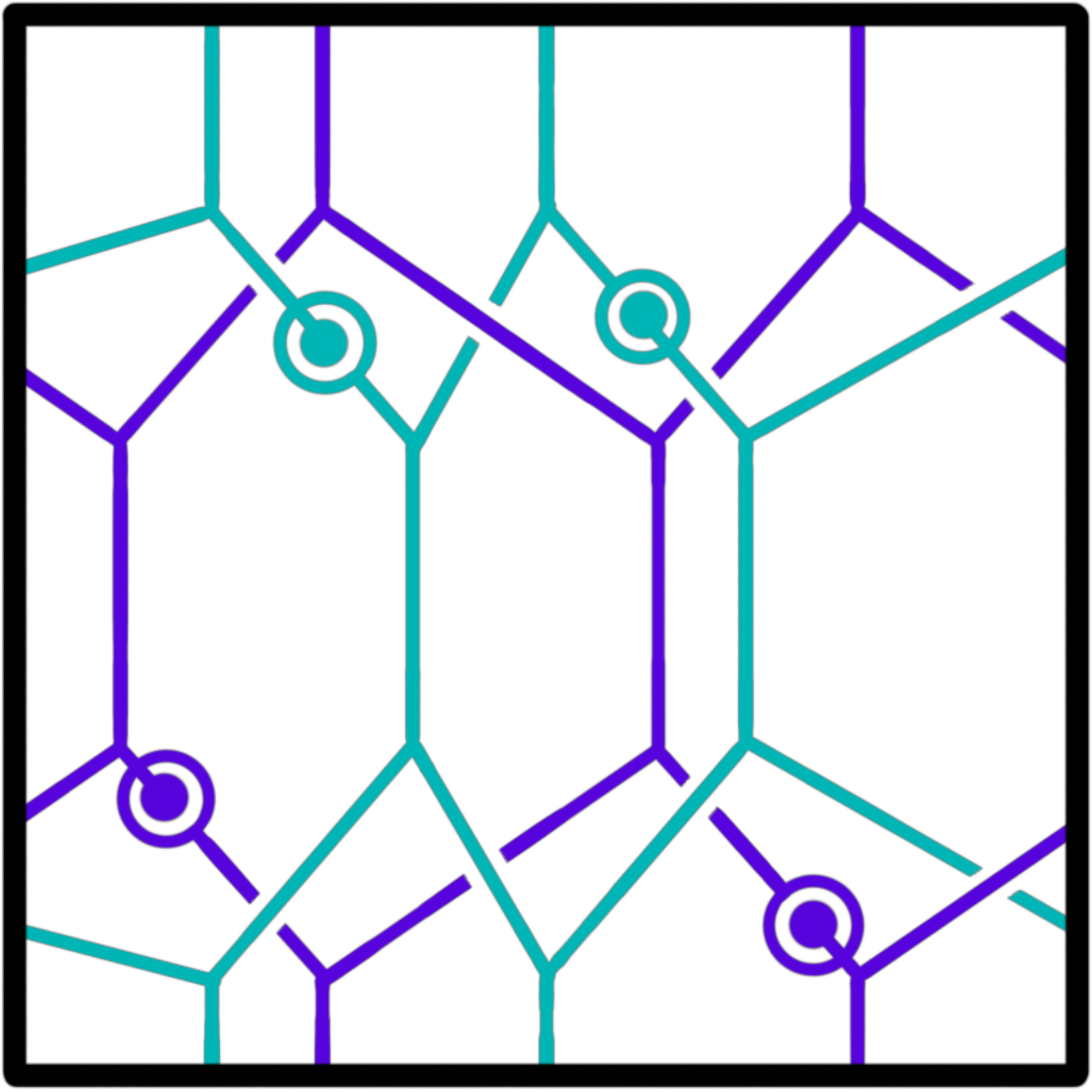}
        \hspace{0.2cm}
        \includegraphics[width=0.29\textwidth]{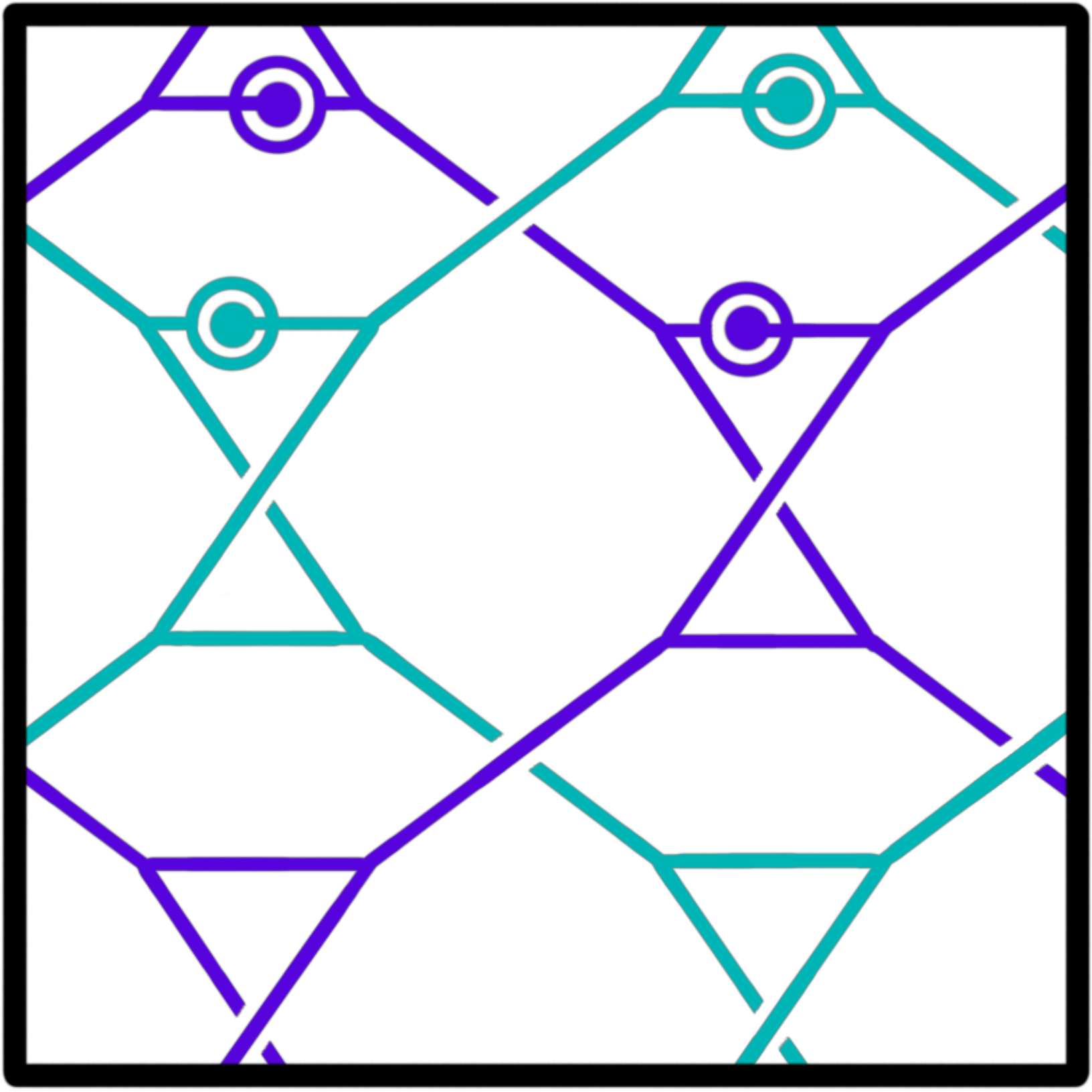}
        \caption{A tridiagram of \textbf{utp-c***}.}
        \label{fig:utp-c_star_star_star_tridia}
    \end{subfigure}
    \caption{The concept of a nearest ground state illustrated with \textbf{utp-c***}: As can be seen with the tridiagram in (b), the minimum crossing number triplet of \textbf{utp-c***} is $(8,8,8)$, which means that is has more crossings than \textbf{utp-c*} and \textbf{utp-c**}. Two crossing changes are needed to transform \textbf{utp-c***} into \textbf{utp-c*} while four are needed to transform it into \textbf{utp-c**}, suggesting that the former is the nearest ground state.}
    \label{fig:utp-c_star_star_star_uc_and_tridia}
\end{figure*}

Since a $\mathcal{G}$-family may contain more than one ground state, not all of them may realise the untangling number of a given embedding. For example, with respect to the unit cell shown in Fig. \ref{fig:utp-c_star_star_star_uc_and_tridia}, the minimum crossing number triplet of \textbf{utp-c***} is $(8,8,8)$. It has more crossings than \textbf{utp-c*} and \textbf{utp-c**}, both with $(4,8,8)$ minimum crossing number triplet, and is, therefore, not a ground state. However, although \textbf{utp-c***} can be transformed into either of \textbf{utp-c*} and \textbf{utp-c**}, only two crossings are needed to reach the former, while four are needed to attain the latter. Therefore, the nearest ground state of \textbf{utp-c***} is \textbf{utp-c*}. The computations of these upper bounds are given in the Supplementary Information.

\begin{figure*}[hbtp]
    \centering
    
    \begin{subfigure}[b]{0.27\textwidth}
        \includegraphics[width=0.6\textwidth]{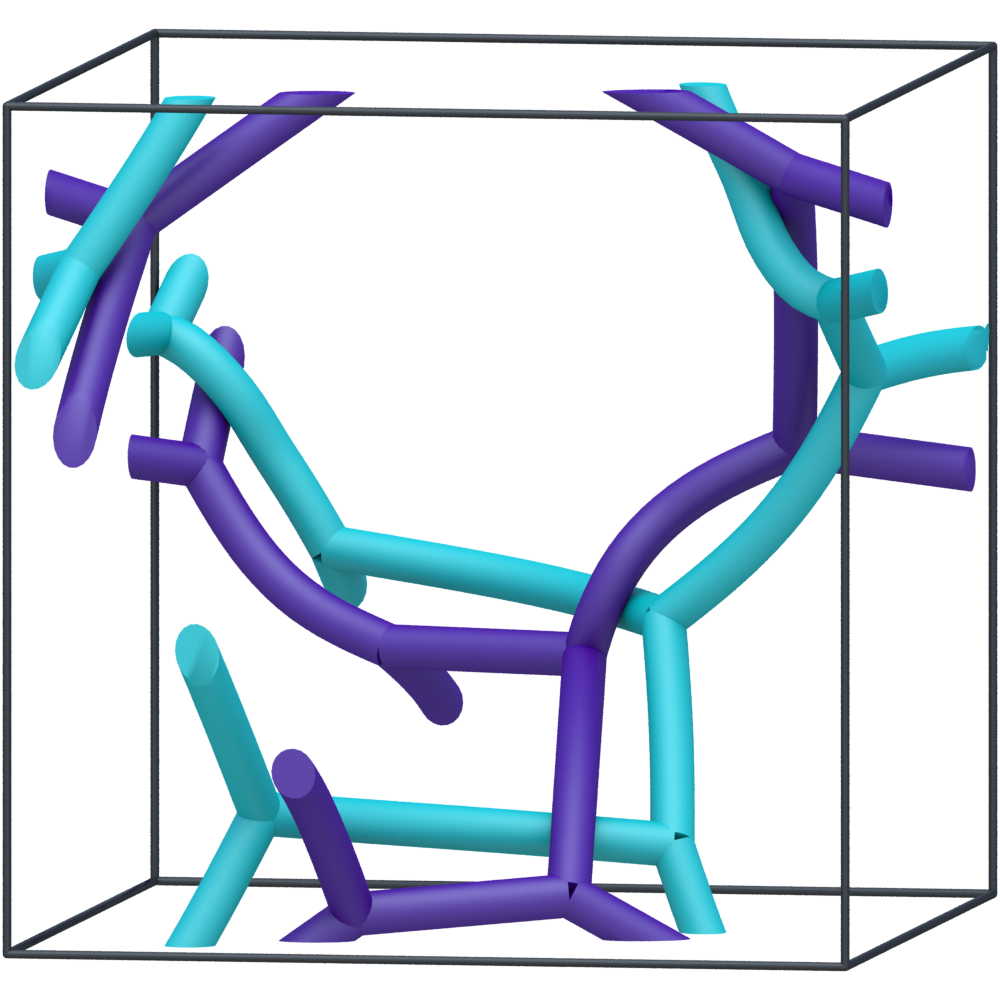}
        \caption{A unit cell of \textbf{srs-c**}.}
        \label{fig:n0p4on2_to_the_6_min_cn_t_uc}
    \end{subfigure}
    \begin{subfigure}[b]{0.54\textwidth}
        \includegraphics[width=0.29\textwidth]{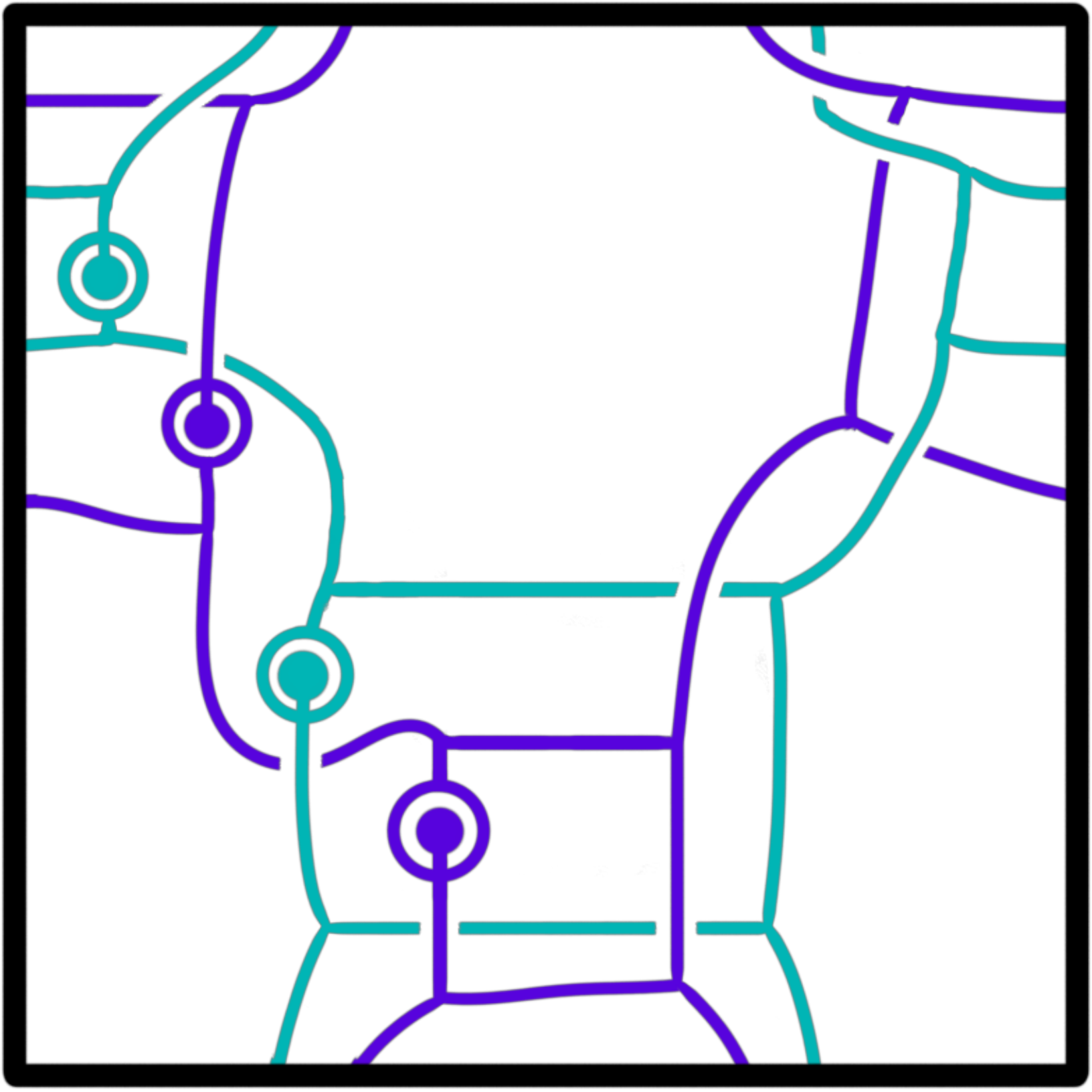}
        \hspace{0.2cm}
        \includegraphics[width=0.29\textwidth]{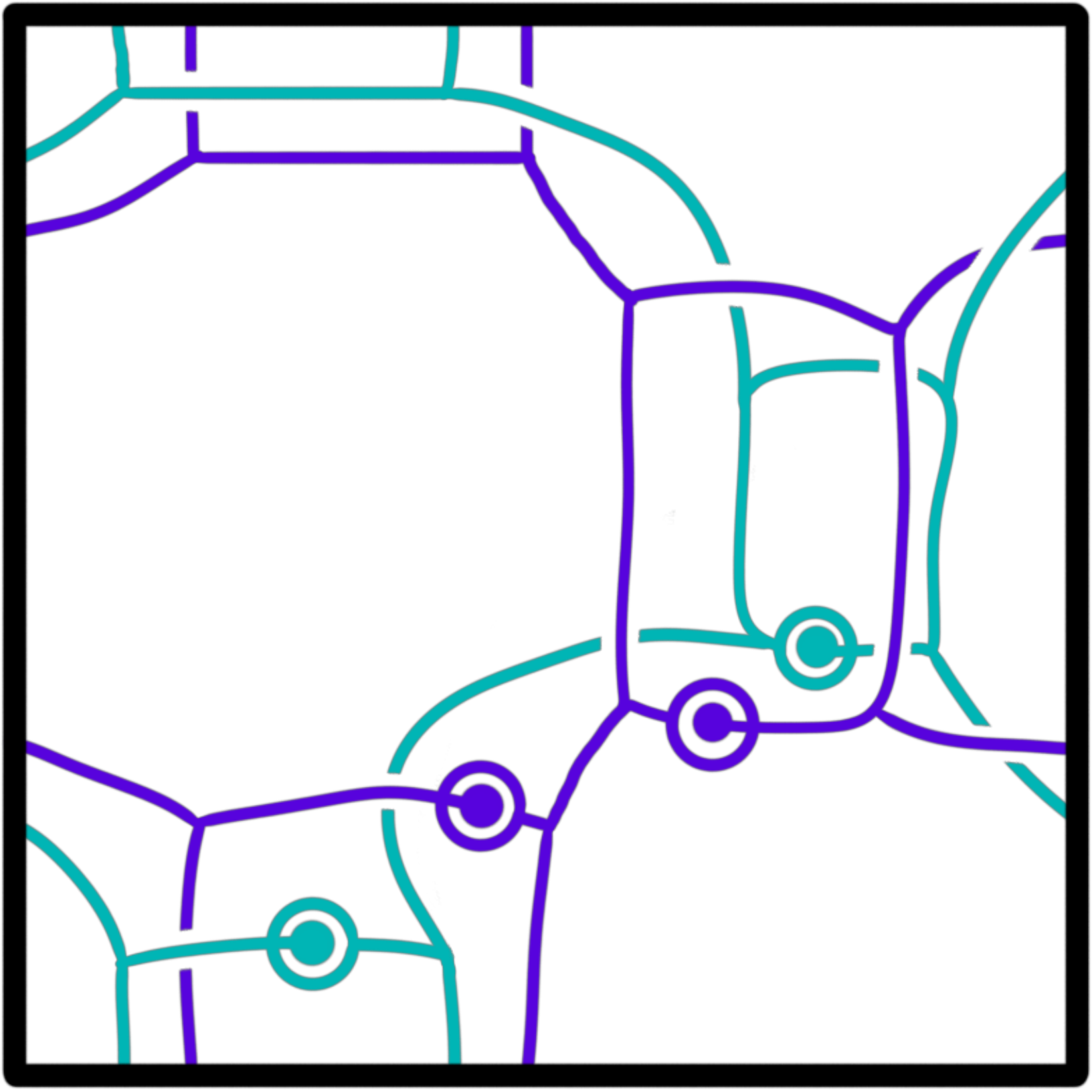}
        \hspace{0.2cm}
        \includegraphics[width=0.29\textwidth]{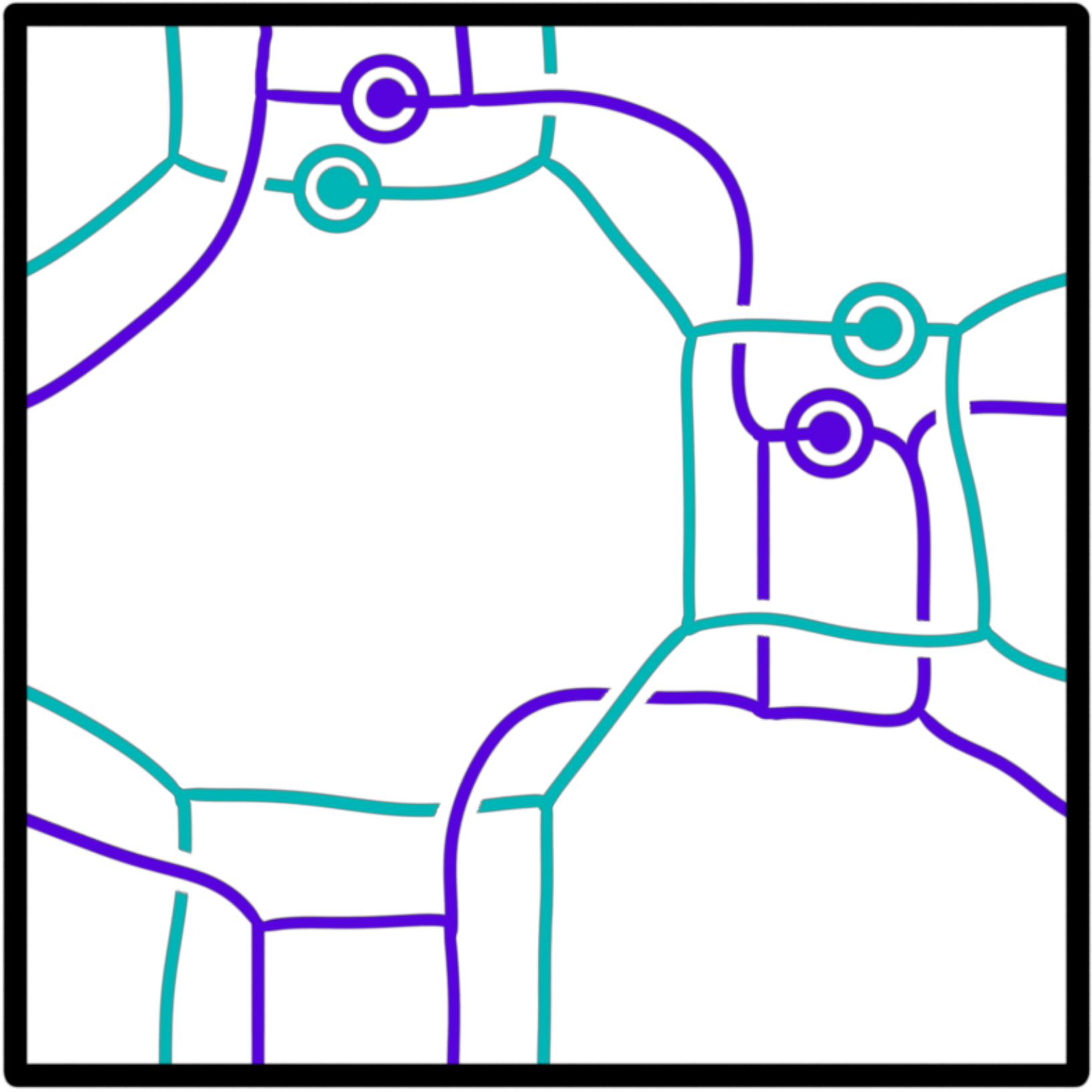}
        \caption{A tridiagram of \textbf{srs-c**}.}
        \label{fig:srs-c_star_star_tridia}
    \end{subfigure}

    \vskip\baselineskip

    \begin{subfigure}[b]{0.27\textwidth}
        \includegraphics[width=0.6\textwidth]{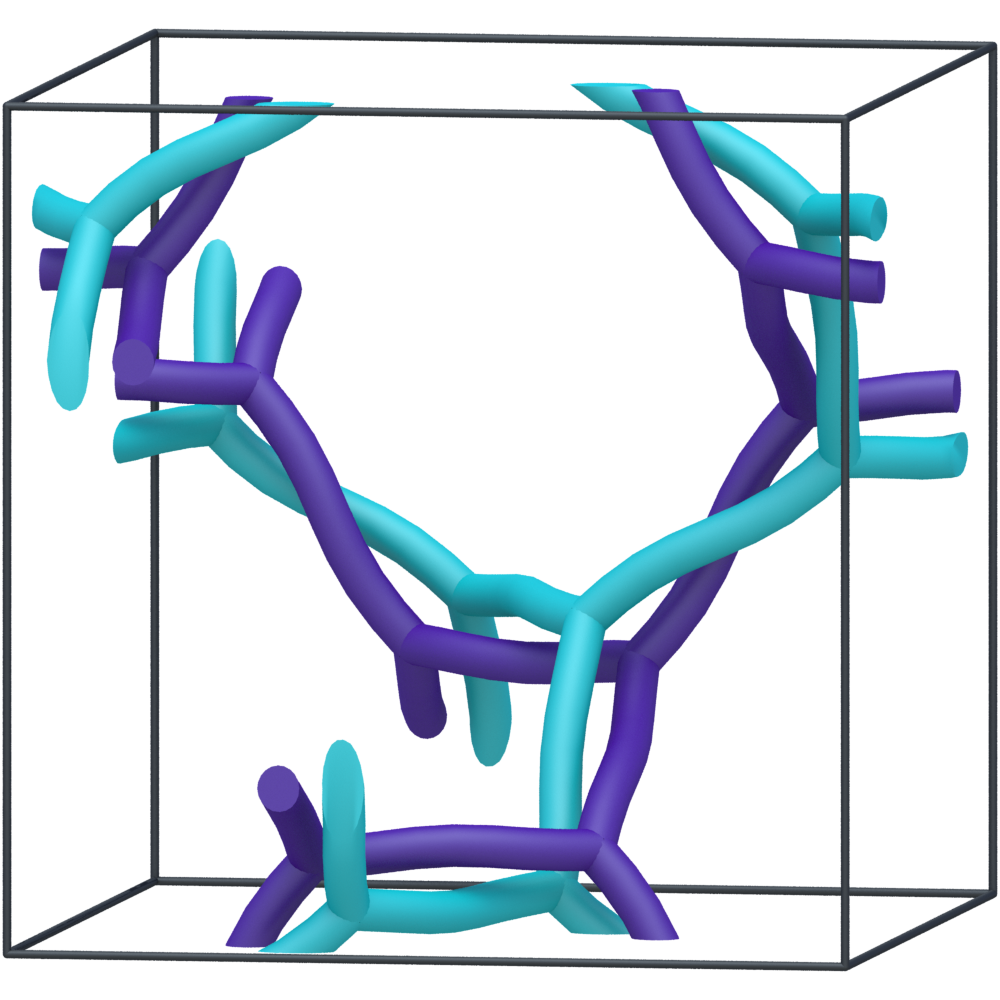}
        \caption{A unit cell of the embedding in Fig. \ref{fig:0p6on2_to_the_6_extended}.}
        \label{fig:0p6on2_to_the_6_uc}
    \end{subfigure}
    \begin{subfigure}[b]{0.54\textwidth}
        \includegraphics[width=0.29\textwidth]{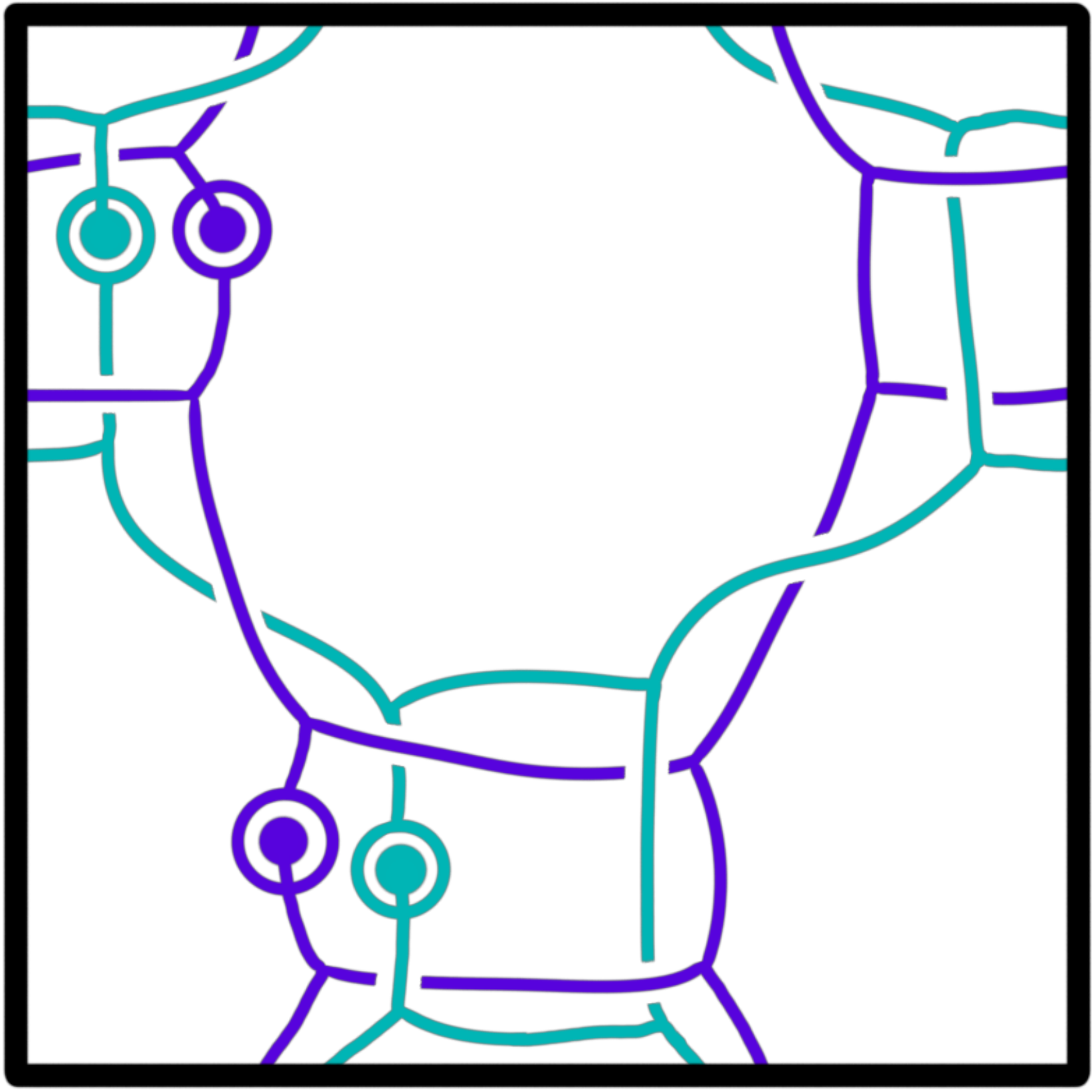}
        \hspace{0.2cm}
        \includegraphics[width=0.29\textwidth]{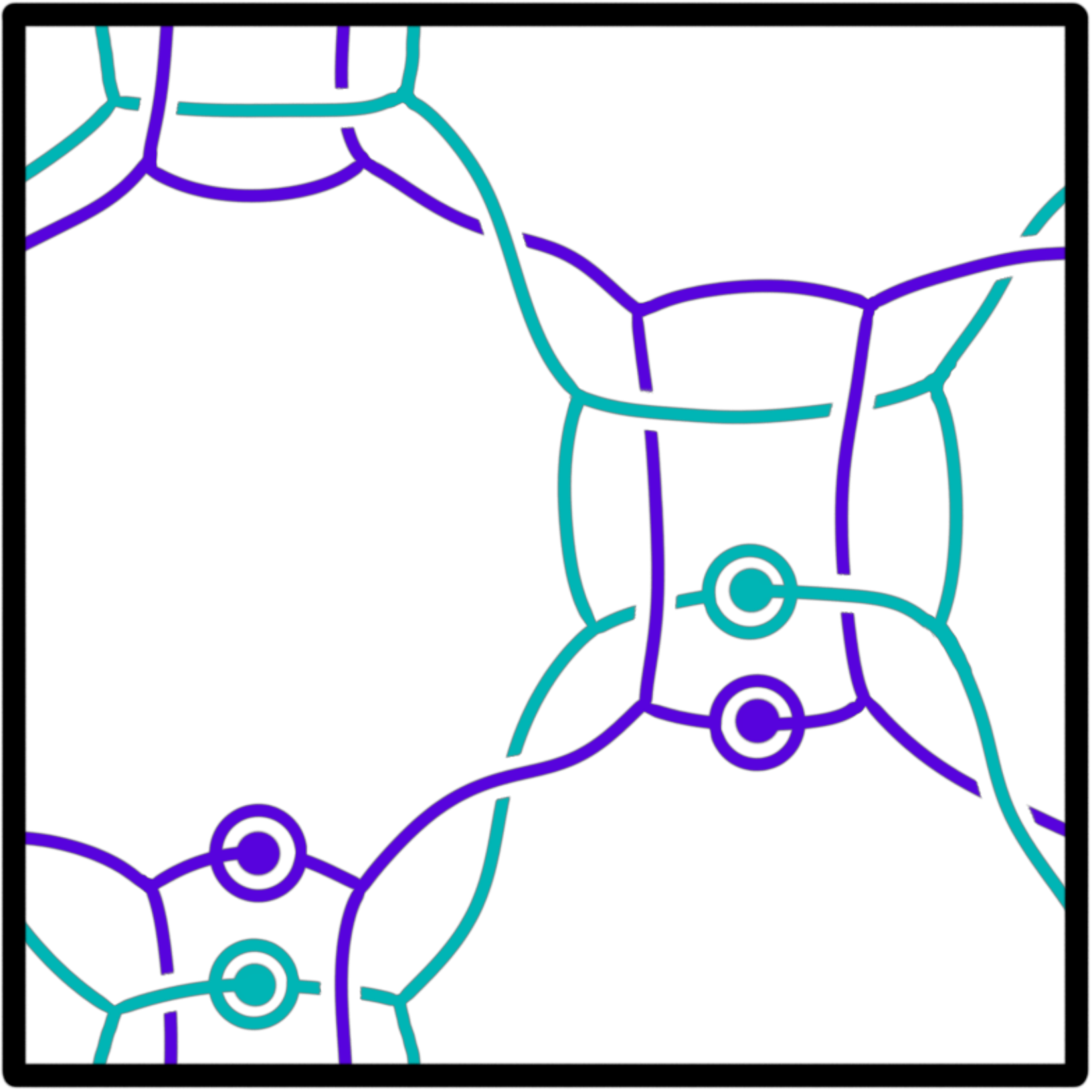}
        \hspace{0.2cm}
        \includegraphics[width=0.29\textwidth]{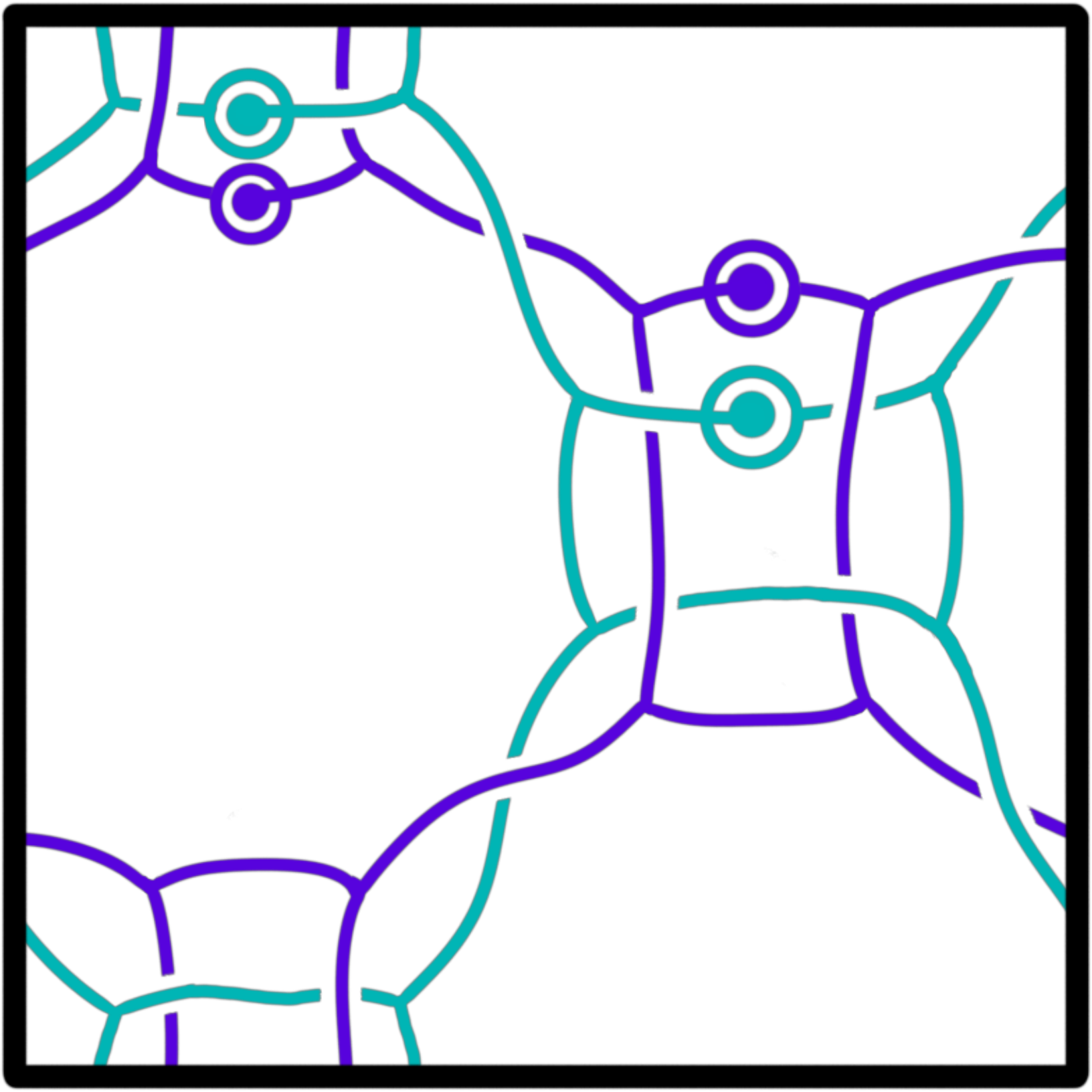}
        \caption{A tridiagram corresponding to the unit cell shown in (c).}
        \label{fig:0p6on2_to_the_6_tridia}
    \end{subfigure}
    \caption{Unit cells and minimal tridiagrams of \textbf{srs-c**} and the embedding of two same-handed \textbf{srs} shown in Fig. \ref{fig:0p6on2_to_the_6_extended}: With respect to the unit cells shown in (a) and (c), the minimum crossing number triplets are $(9,10,9)$ and $(12,12,12)$, respectively, and the untangling numbers likely are $4$ and $6$, the computations of which are given in the Supplementary Information.}
    \label{fig:0p6on2_to_the_6_and_n0p4on2_to_the_6_uc_and_tridia}
\end{figure*}

\begin{figure*}[htbp]
    \centering
    \begin{subfigure}[b]{0.3\textwidth}
    \centering
        \includegraphics[width=\textwidth]{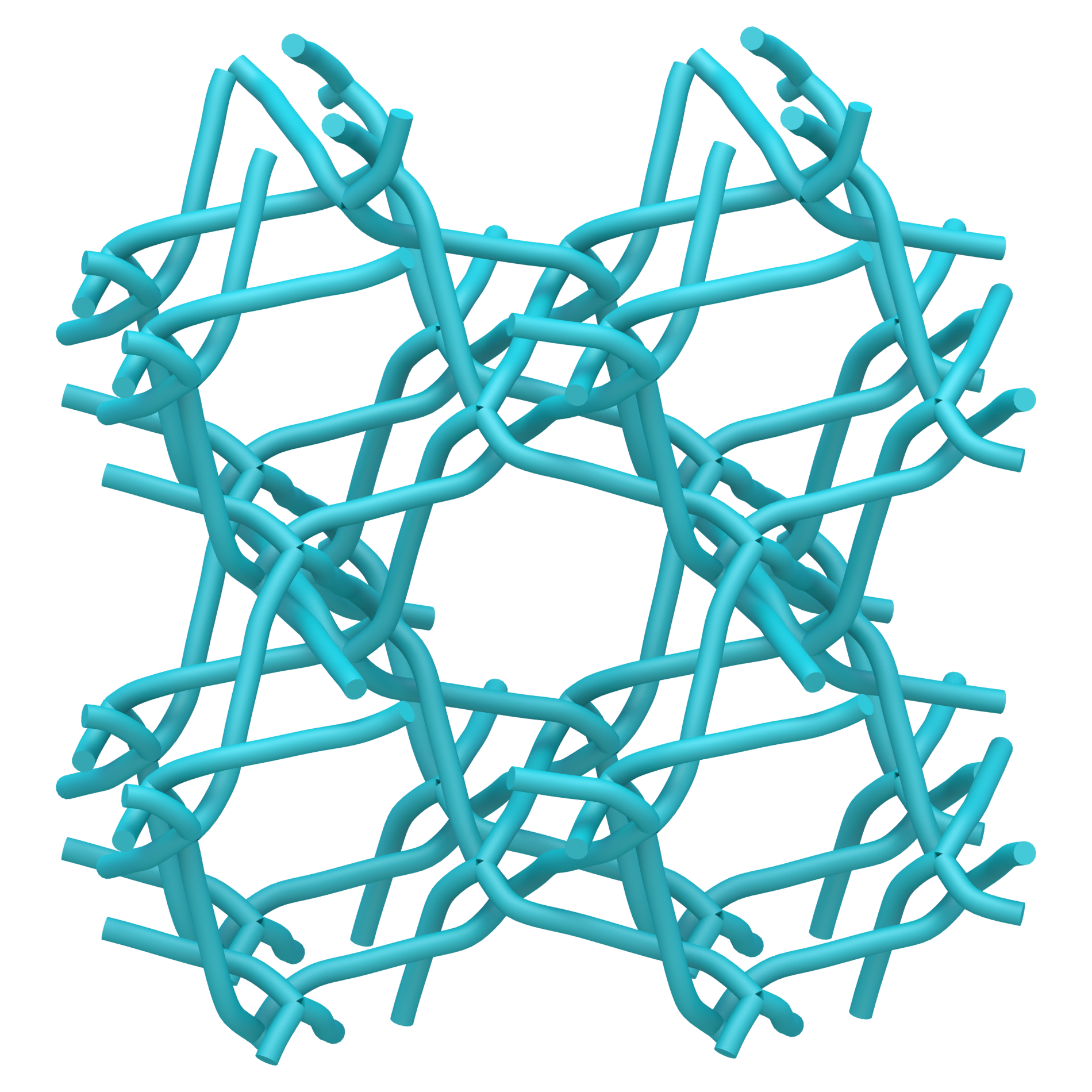}
        \caption{}
        \label{fig:srs-me_extended}
    \end{subfigure}
    \begin{subfigure}[b]{0.3\textwidth}
    \centering
        \includegraphics[width=\textwidth]{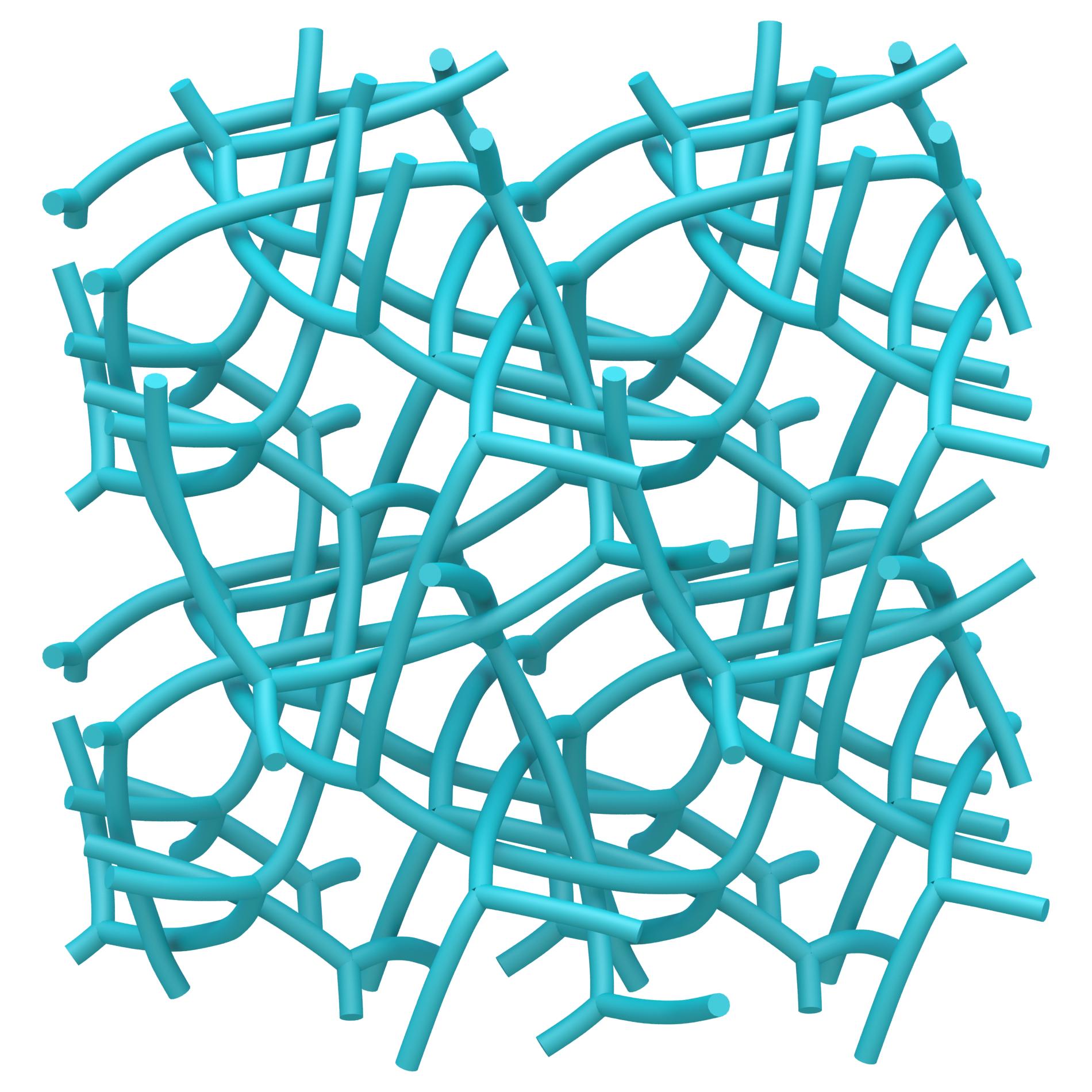}
        \caption{}
        \label{fig:srs-z_extended}
    \end{subfigure}
    \caption{Two tangled embeddings of the \textbf{srs} network: The strong rings of the embedding in (a) form Hopf links and $5^2_1$ links, whereas those of the embedding in (b), which is \textbf{srs-z}, form Hopf links only. In addition, both possess ravelled cages corresponding to the $5_1$ embedding of the $\theta$ graph, whose diagram is shown in Fig. \ref{fig:5_1}, undetected by cycle analysis.}
    \label{fig:srs-me_and_srs-z_extended}
\end{figure*}

\begin{figure*}[hbtp]

    \centering
    
    \begin{subfigure}[b]{0.27\textwidth}
        \includegraphics[width=0.6\textwidth]{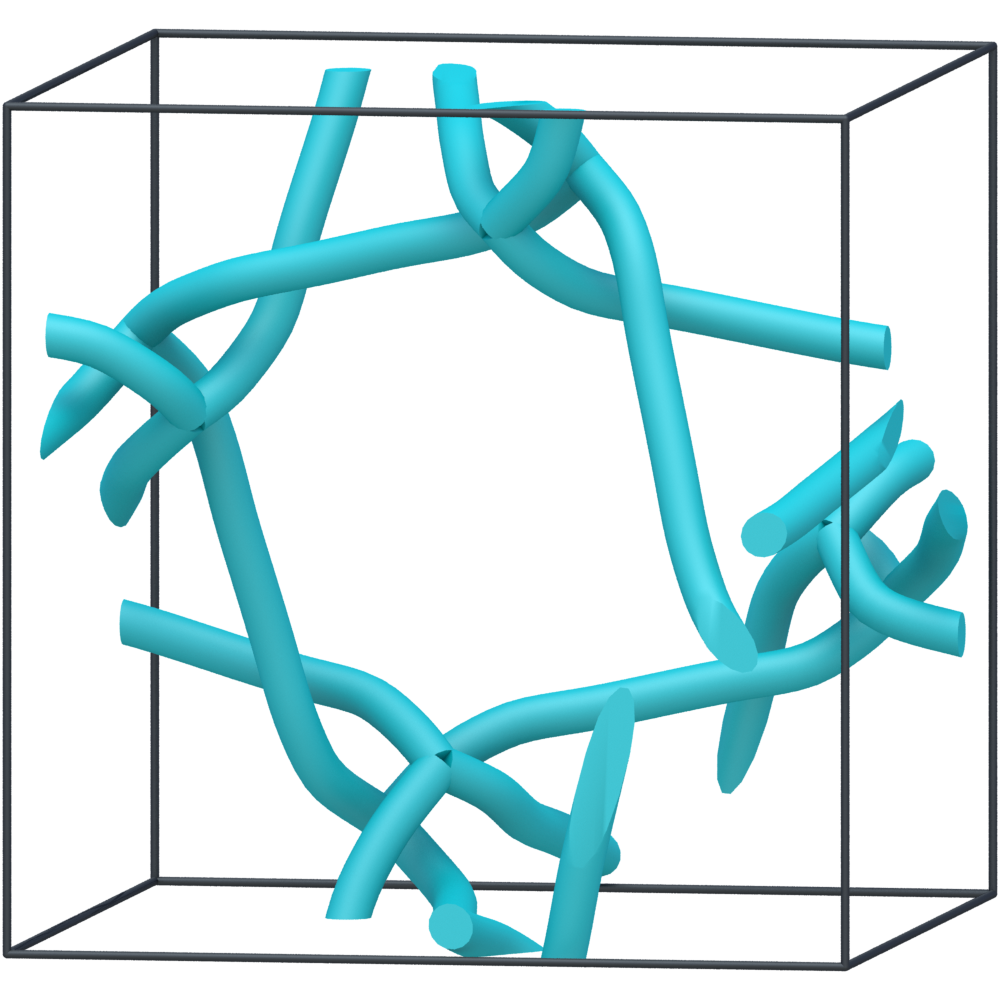}
        \caption{A unit cell of \textbf{srs} Fig. \ref{fig:srs-me_extended}.}
        \label{fig:srs-me_uc}
    \end{subfigure}
    \begin{subfigure}[b]{0.54\textwidth}
        \includegraphics[width=0.29\textwidth]{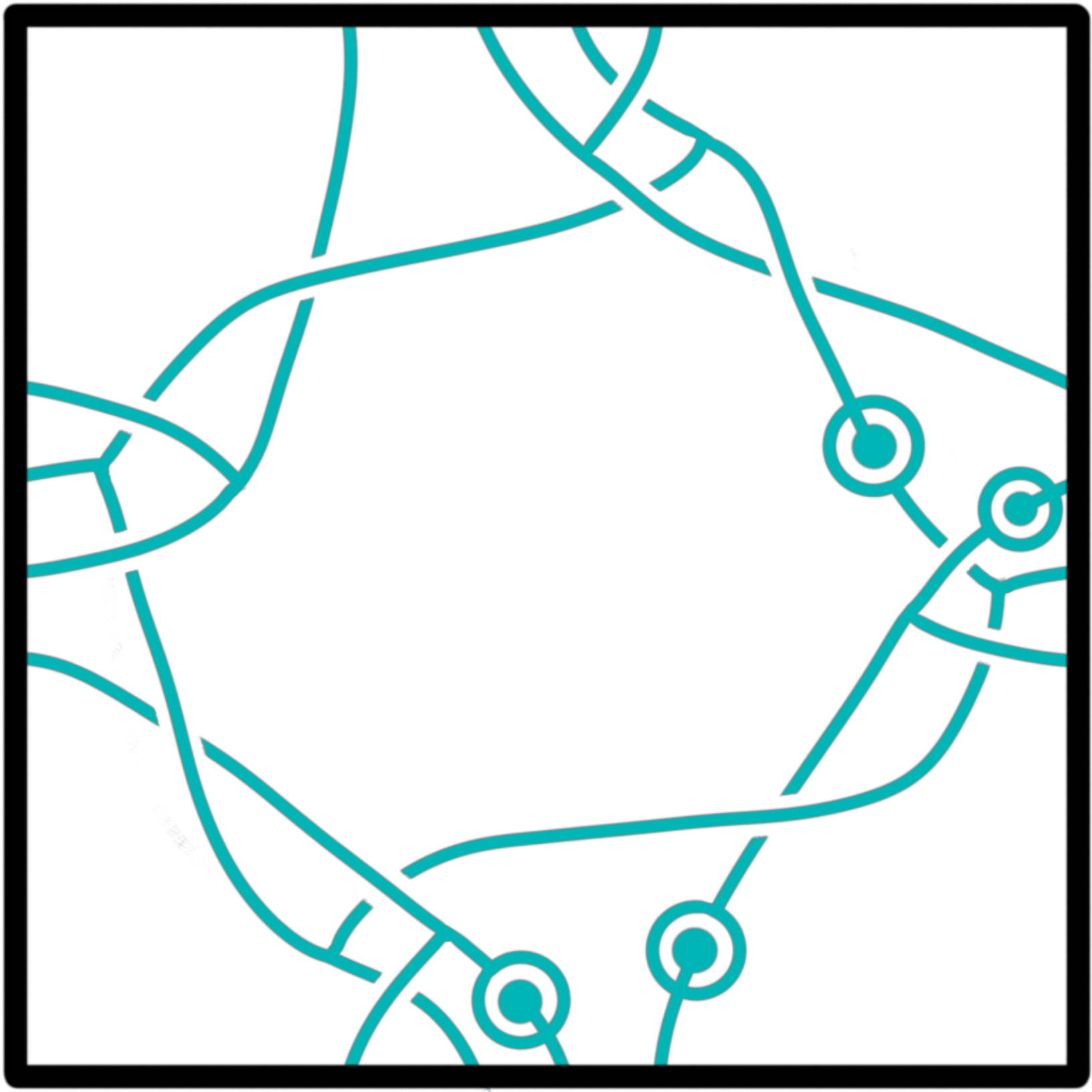}
        \hspace{0.2cm}
        \includegraphics[width=0.29\textwidth]{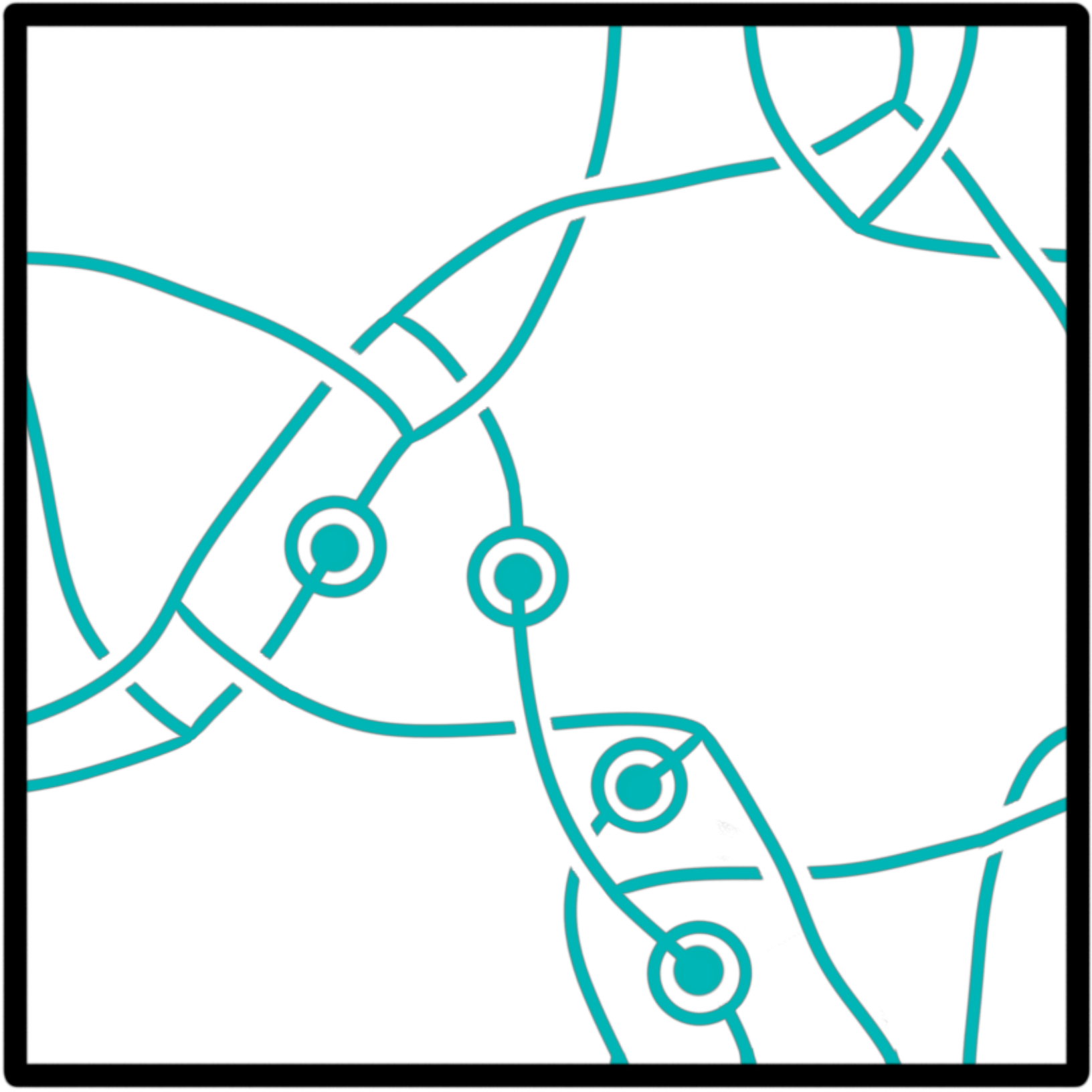}
        \hspace{0.2cm}
        \includegraphics[width=0.29\textwidth]{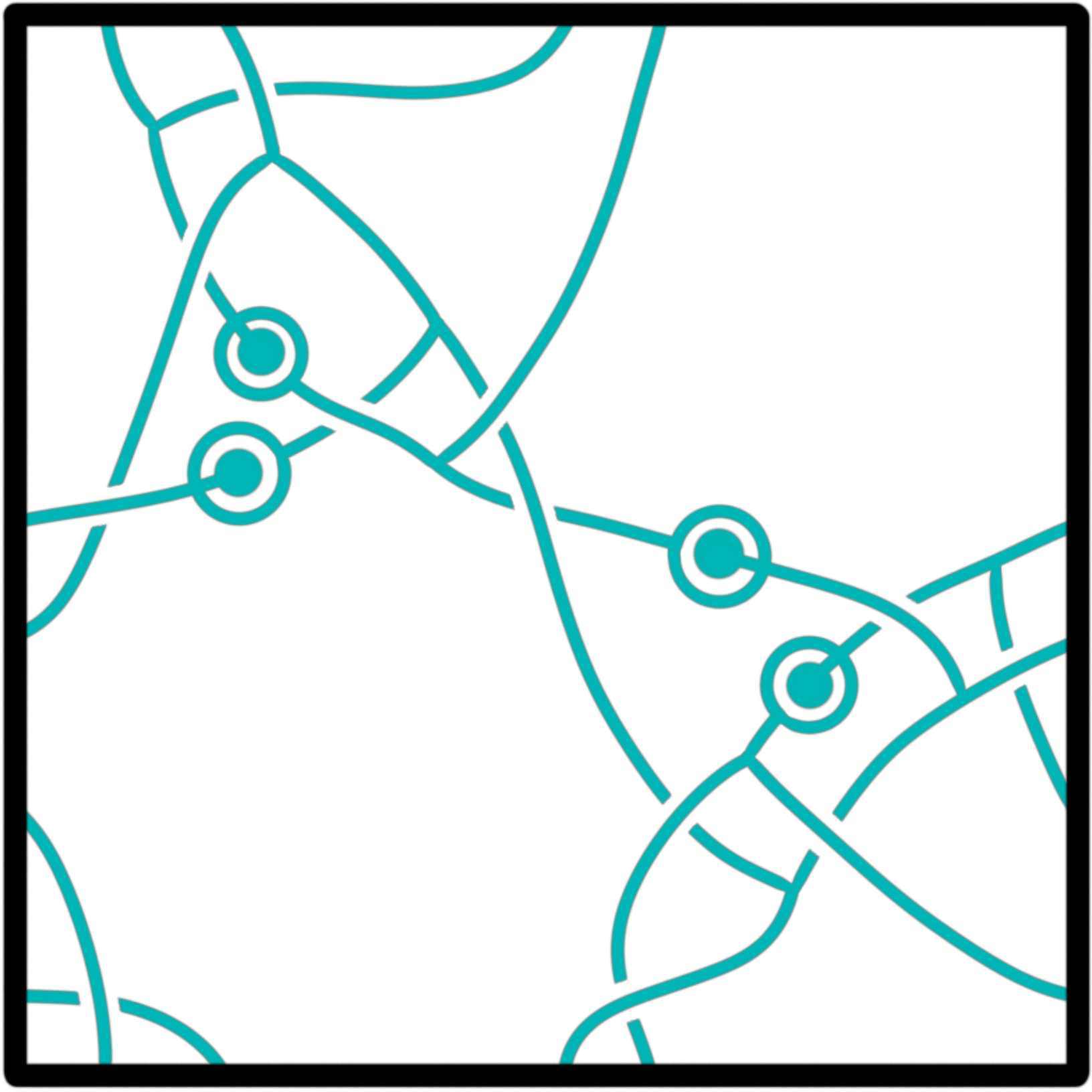}
        \caption{A tridiagram corresponding to the unit cell shown in (a).}
        \label{fig:srs-me_tridia}
    \end{subfigure}

    \vskip\baselineskip

    \begin{subfigure}[b]{0.27\textwidth}
        \includegraphics[width=0.6\textwidth]{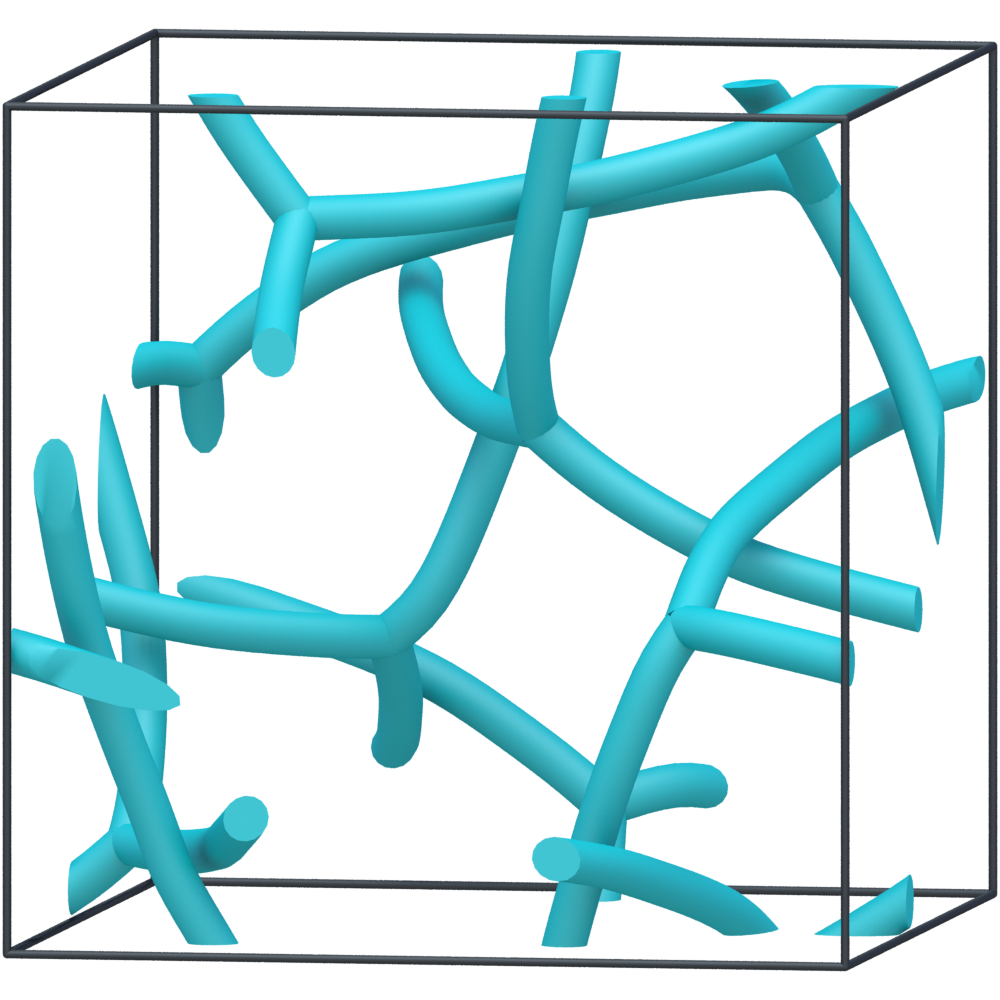}
        \caption{A unit cell of \textbf{srs-z}.}
        \label{fig:srs-z_uc}
    \end{subfigure}
    \begin{subfigure}[b]{0.54\textwidth}
        \includegraphics[width=0.29\textwidth]{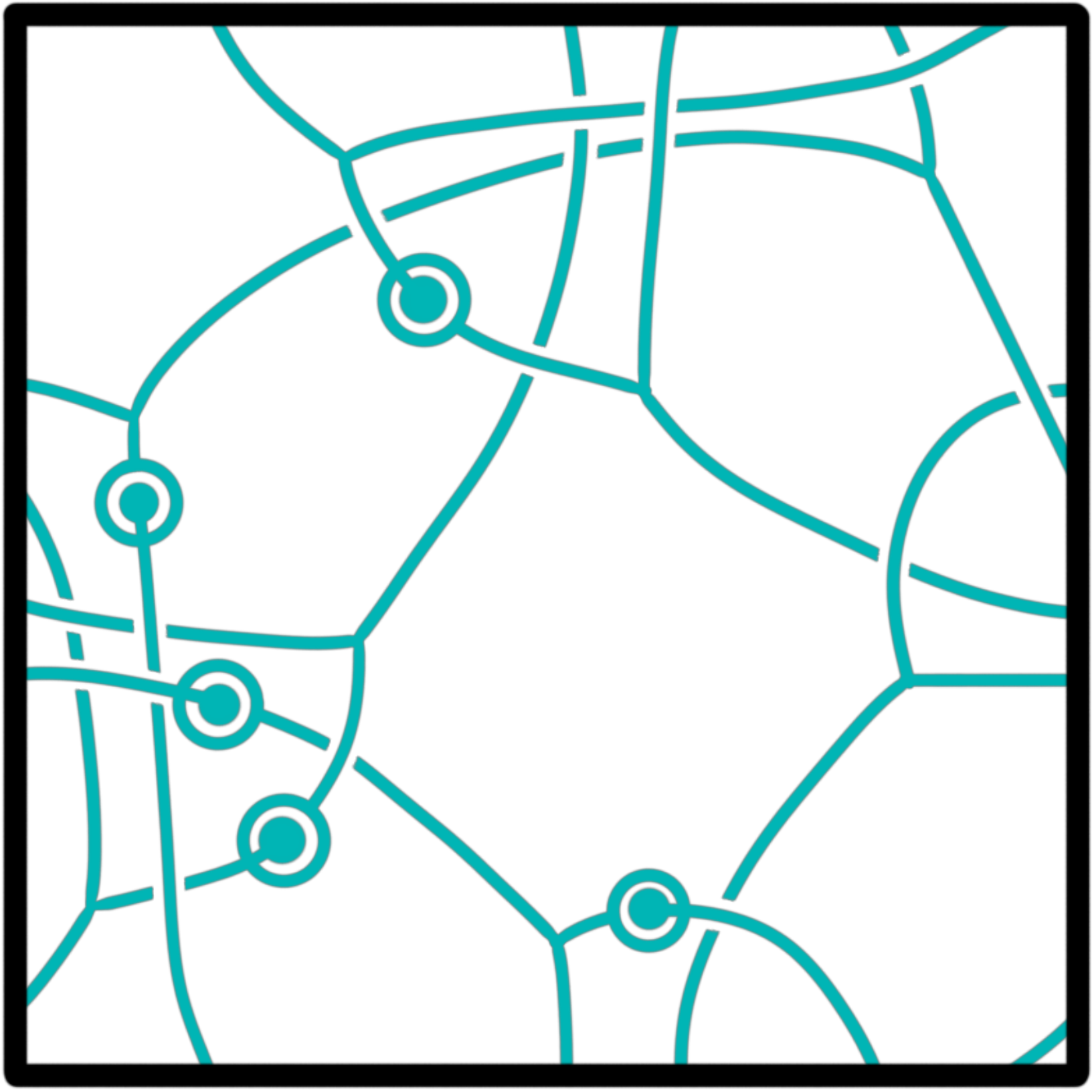}
        \hspace{0.2cm}
        \includegraphics[width=0.29\textwidth]{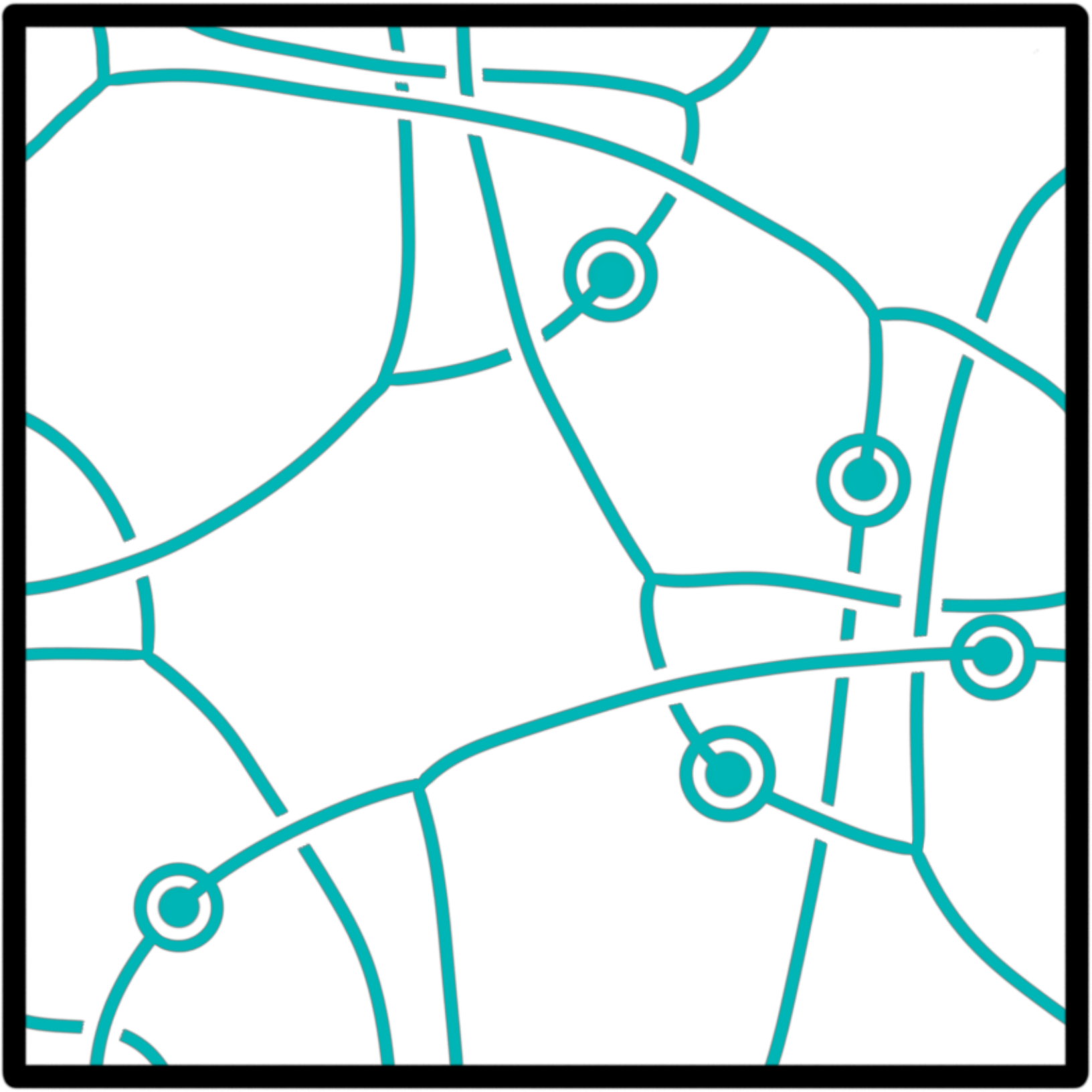}
        \hspace{0.2cm}
        \includegraphics[width=0.29\textwidth]{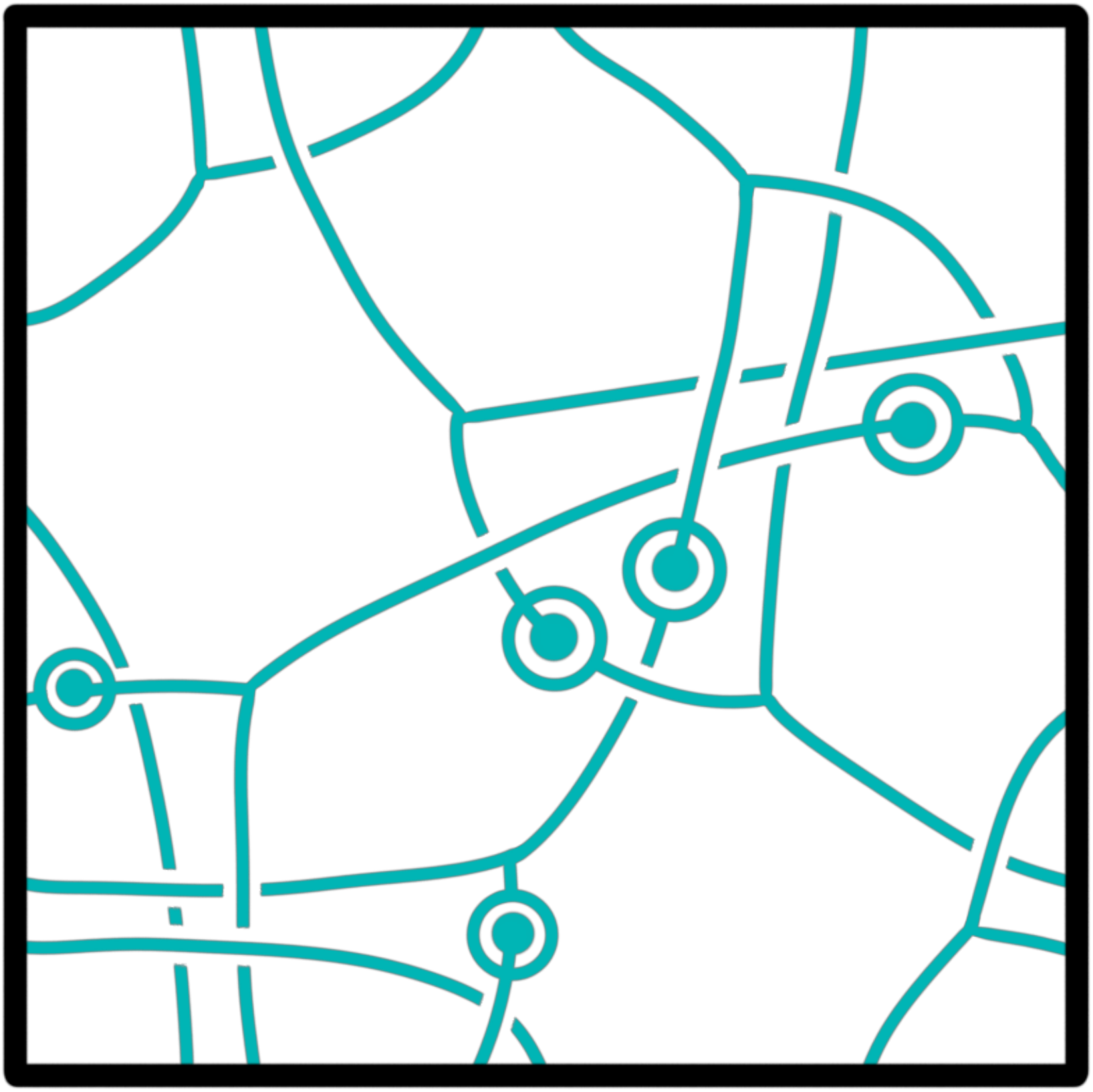}
        \caption{A tridiagram corresponding to the unit cell shown in (a).}
        \label{fig:srs-z_tridia}
    \end{subfigure}
    \caption{With respect to the unit cells shown in (a) and (c), the minimum crossing number triplets of the embeddings of \textbf{srs} displayed in Fig. \ref{fig:srs-me_and_srs-z_extended} are $(12,12,12)$ and $(16,16,16)$, respectively. Both likely have untangling number $8$, as suggested by the computations given in the Supplementary Information, indicating that one embedding is not more difficult to untangle than the other. This does not align with the hierarchy of complexity established by cycle analysis, demonstrating that different invariants capture different aspects of entanglement complexity.}
    \label{fig:srs-me_and_srs-z_uc_and_tridia}
\end{figure*}

The untangling number establishes a hierarchy of complexity between structures. For example, with respect to the unit cells shown in Fig. \ref{fig:0p6on2_to_the_6_and_n0p4on2_to_the_6_uc_and_tridia}, the untangling number of \textbf{srs-c**} and that of the embedding in Fig. \ref{fig:0p6on2_to_the_6_extended} likely are $4$ and $6$, respectively, with \textbf{srs-c*} being their nearest ground state. The computations are given in the Supplementary Information. This coincides with the intuitive idea, presented in Sect. \ref{sec:1}, that these embeddings are more tangled variants of \textbf{srs-c*}. This hierarchy also aligns with that given by the crossing number and cycle analysis: the minimum crossing number triplets are $(9,10,9)$ and $(12,12,12)$, respectively, while the strong rings of the embeddings respectively form Hopf links and $4^2_1$ links, and Hopf links, $4^2_1$ links, and $6^2_1$ links. However, it is important to note that different measures of complexity need not agree with one another. In Fig. \ref{fig:srs-me_and_srs-z_extended}, we present two embeddings of \textbf{srs} for which various measures do not agree. An analysis of the cycles of the embeddings shows that the strong rings of the embedding in Fig. \ref{fig:srs-me_extended} form Hopf links and $5^2_1$ links, whereas those of the embedding in Fig. \ref{fig:srs-z_extended}, which is \textbf{srs-z}, form only Hopf links. In addition, both possess ravelled cages, namely the $5_1$ embedding of the $\theta$ graph, which is undetected by cycle analysis. However, with respect to the unit cells shown in Fig. \ref{fig:srs-me_and_srs-z_uc_and_tridia}, the minimum crossing number triplets are $(12,12,12)$ and $(16,16,16)$, respectively, and both likely have untangling number $8$, as suggested by the computations given in the Supplementary Information. Although this may initially appear contradictory, it simply reflects the fact that different invariants quantify different aspects of entanglement complexity. Cycle analysis determines the complexity of the knots and links formed by the strong rings of a given embedding. The crossing number, via tridiagrams, indicates the `ease with which a unit cell can be drawn on three planar surfaces along three non-coplanar axes'. The untangling number captures how easily an embedding can be transformed into its least tangled state. These measures therefore need not agree. This behaviour is analogous to that of invariants of classical knots. For example, the $5_1$ knot, whose crossing number is 5, has unknotting number $2$, whereas the $6_1$ knot, whose crossing number is 6, has unknotting number $1$. This means that although the $5_1$ knot is `easier to draw on a plane', the $6_1$ knot is easier to unknot. In fact, one can construct an infinite family of knots with arbitrarily large crossing numbers ($6_1$, $7_2$, $8_1$, $9_2$, etc.), all of which have unknotting number $1$ and are therefore easier to unknot than $5_1$ \cite{Adams.book}.

Note that the cases of \textbf{srs-c*} and \textbf{srs-c**}, as well as the embedding in Fig. \ref{fig:0p6on2_to_the_6_extended}, illustrate the relevance of developing the concept of a $\mathcal{U}$-family. Indeed, it is more natural to regard \textbf{srs-c**} and the embedding in Fig. \ref{fig:0p6on2_to_the_6_extended} as more tangled variants of \textbf{srs-c*}, but not of \textbf{srs-c}, even though they are embeddings of the same graph.

\begin{figure}
    \centering
    \includegraphics[width=0.4\textwidth]{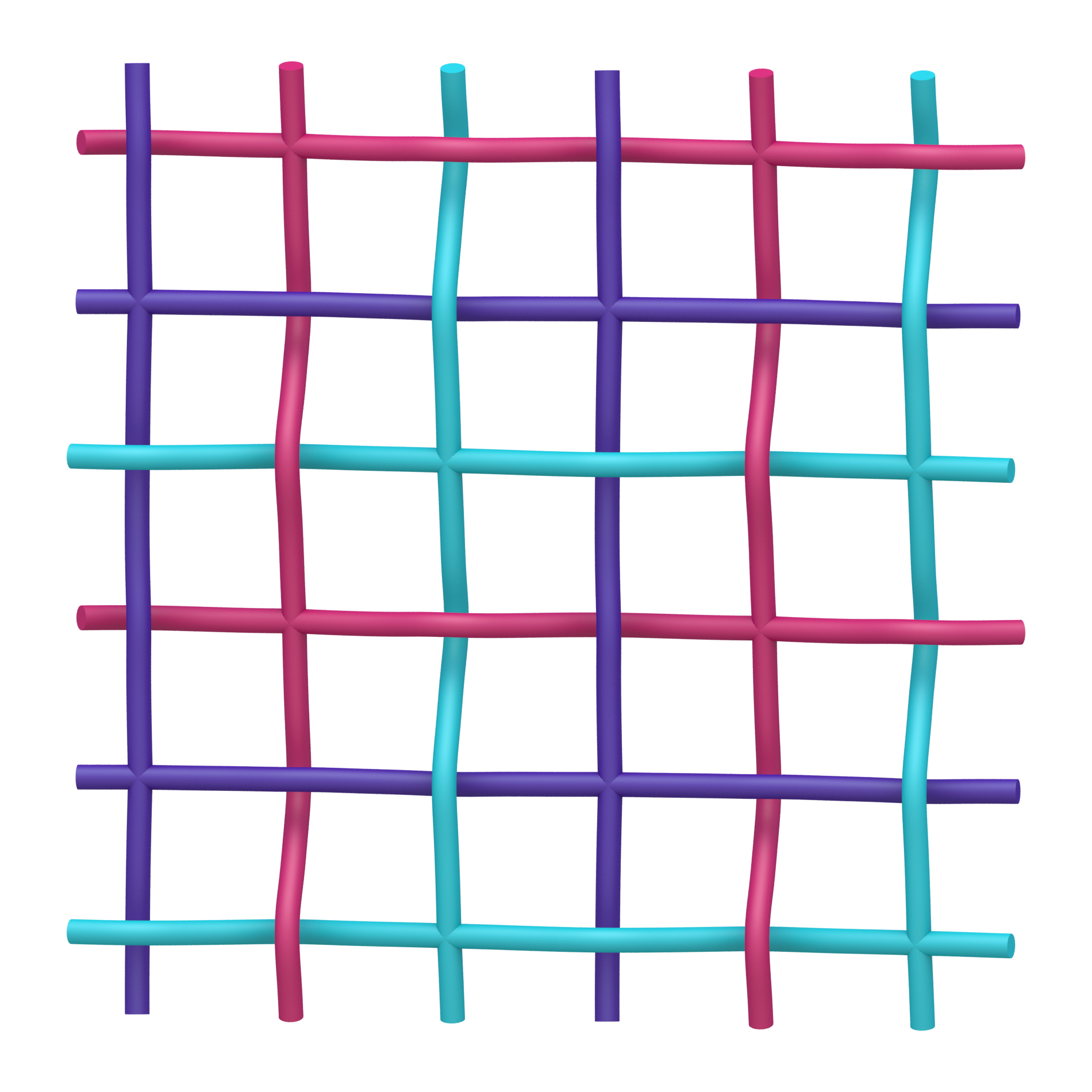}
    \caption{An example of a structure whose strong rings form Borromean rings. The network components are pairwise untangled.}
    \label{fig:borromean_sql_extended}
\end{figure}

\begin{figure*}
    \centering
    \begin{subfigure}[b]{0.18\textwidth}
        \centering
        \includegraphics[width=0.9\textwidth]{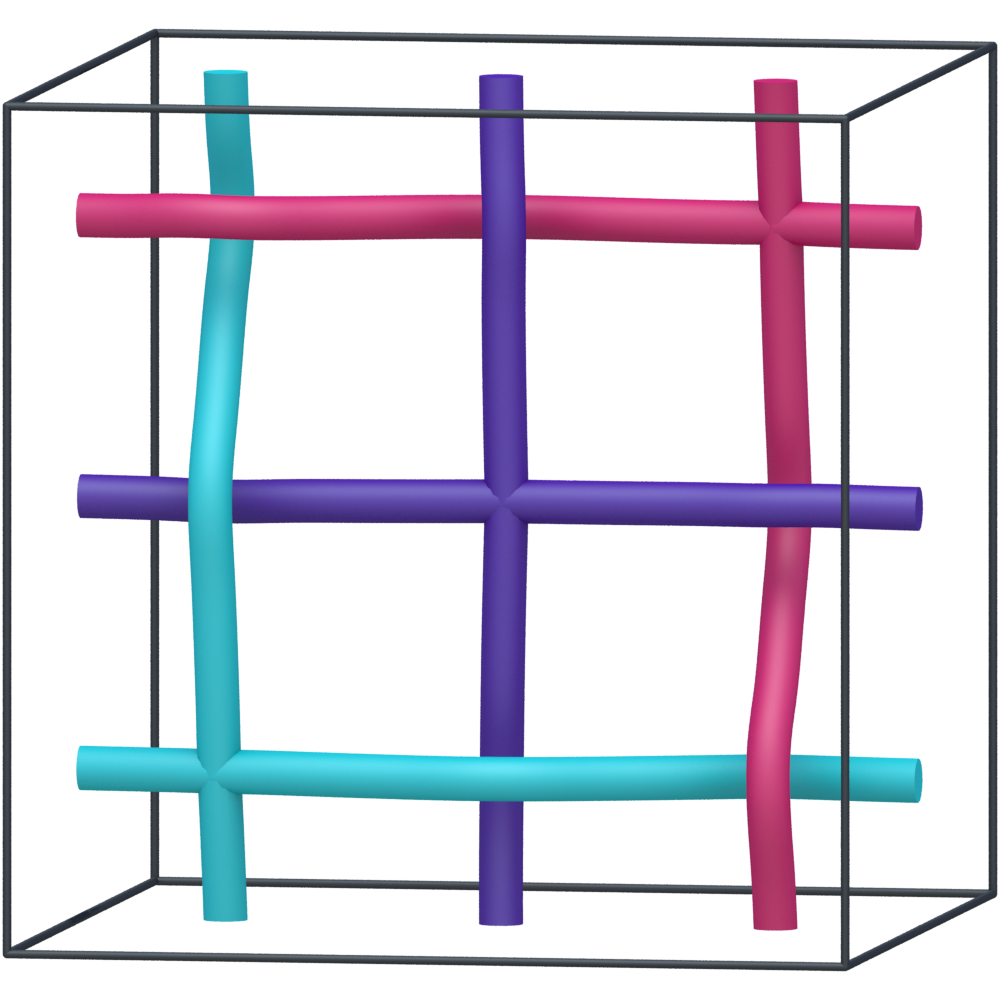}
        \caption{}
        \label{fig:borromean_sql_uc}
    \end{subfigure}
    \begin{subfigure}[b]{0.18\textwidth}
        \centering
        \includegraphics[width=0.9\textwidth]{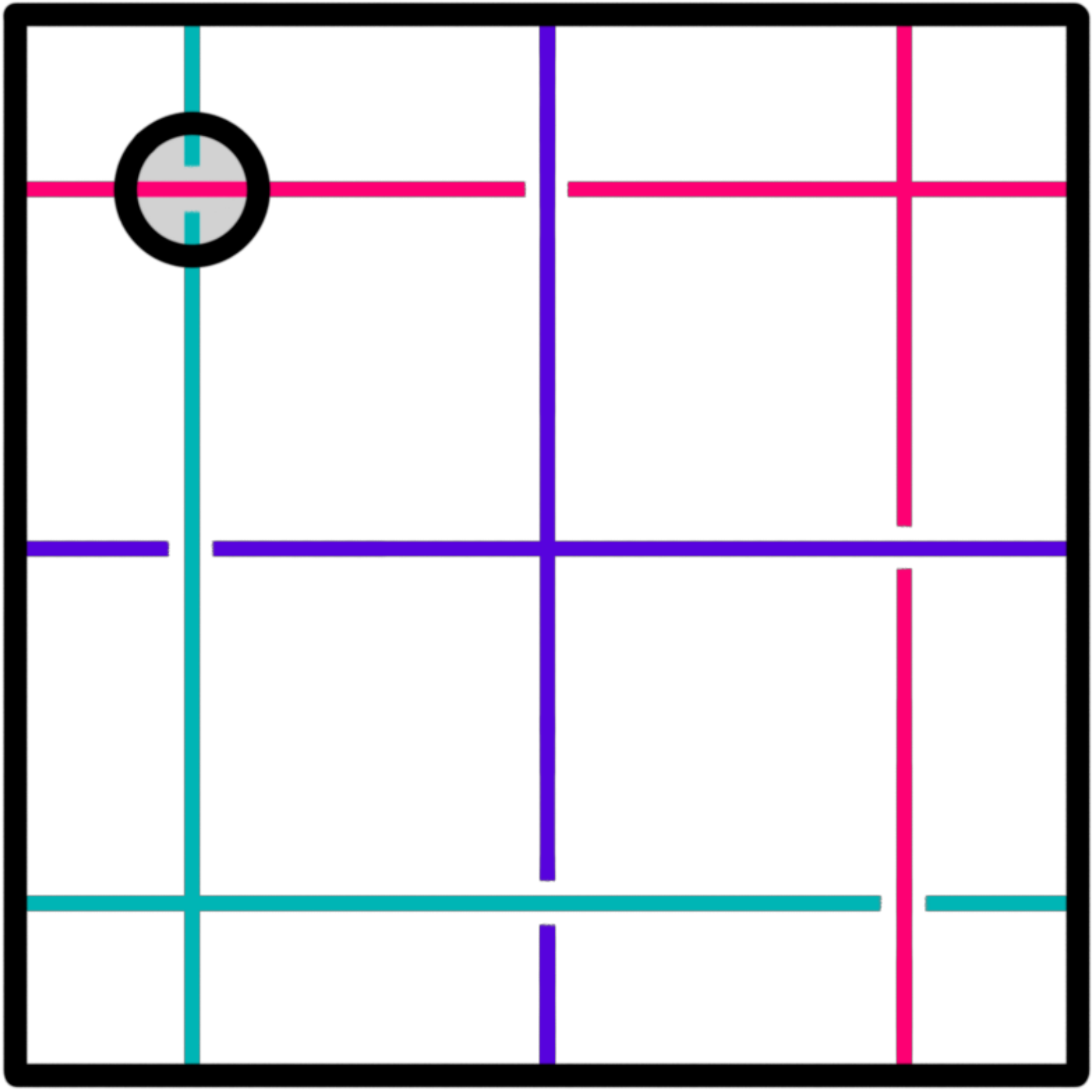}
        \caption{}
        \label{fig:untangling_borromean_sql_1}
    \end{subfigure}
    \begin{subfigure}[b]{0.18\textwidth}
        \centering
        \includegraphics[width=0.9\textwidth]{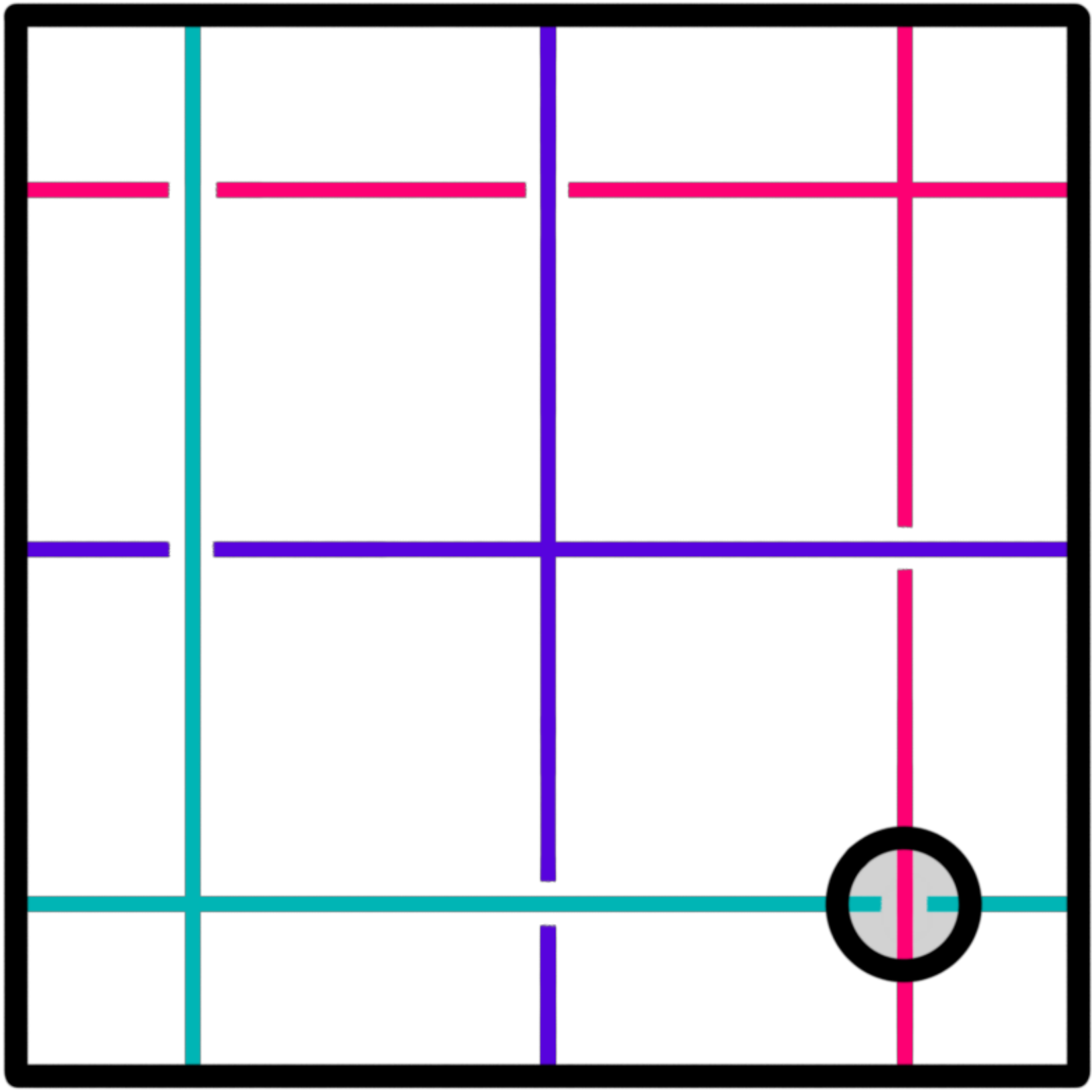}
        \caption{}
        \label{fig:untangling_borromean_sql_2}
    \end{subfigure}
    \begin{subfigure}[b]{0.18\textwidth}
        \centering
        \includegraphics[width=0.9\textwidth]{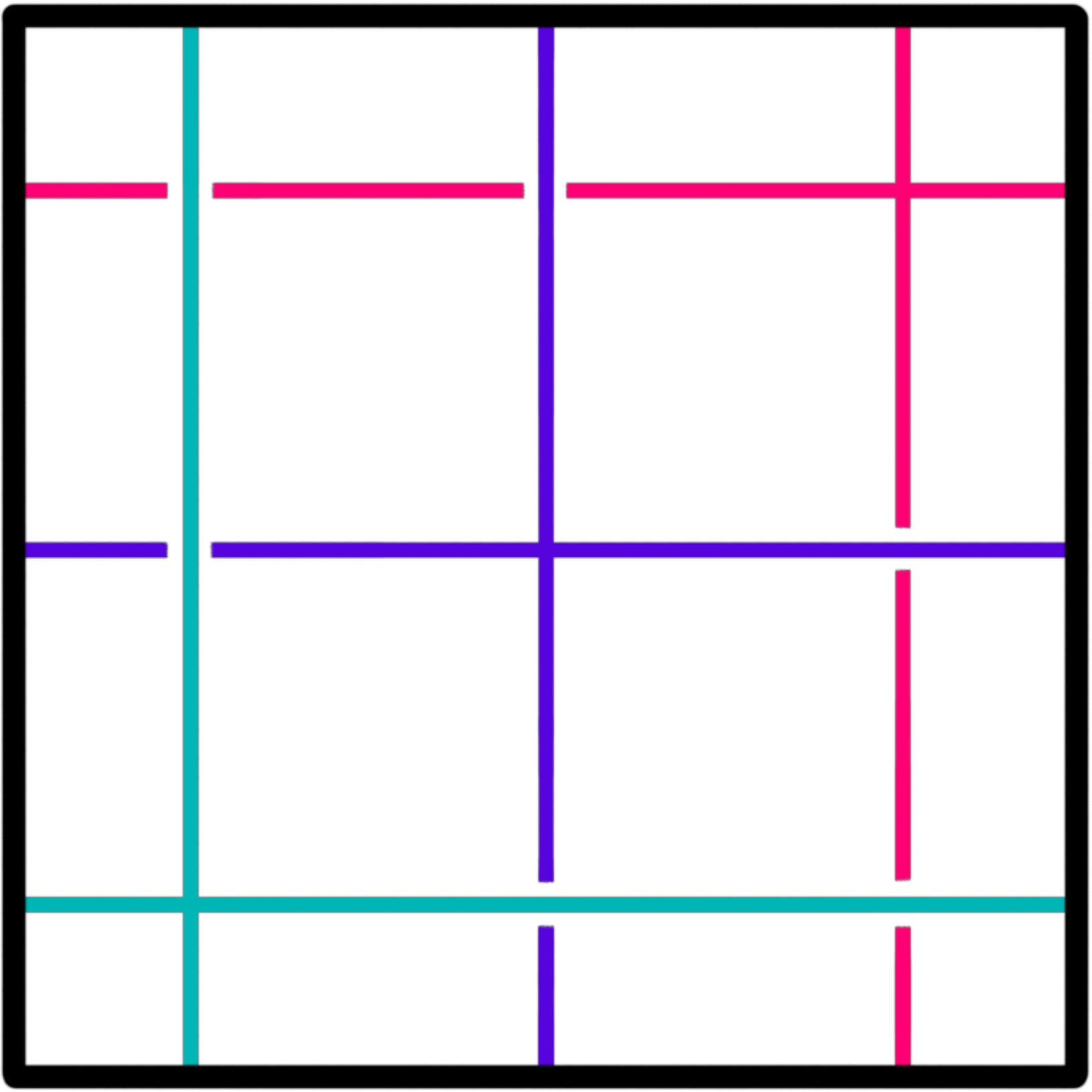}
        \caption{}
        \label{fig:untangling_borromean_sql_3}
    \end{subfigure}
    \begin{subfigure}[b]{0.18\textwidth}
        \centering
        \includegraphics[width=0.9\textwidth]{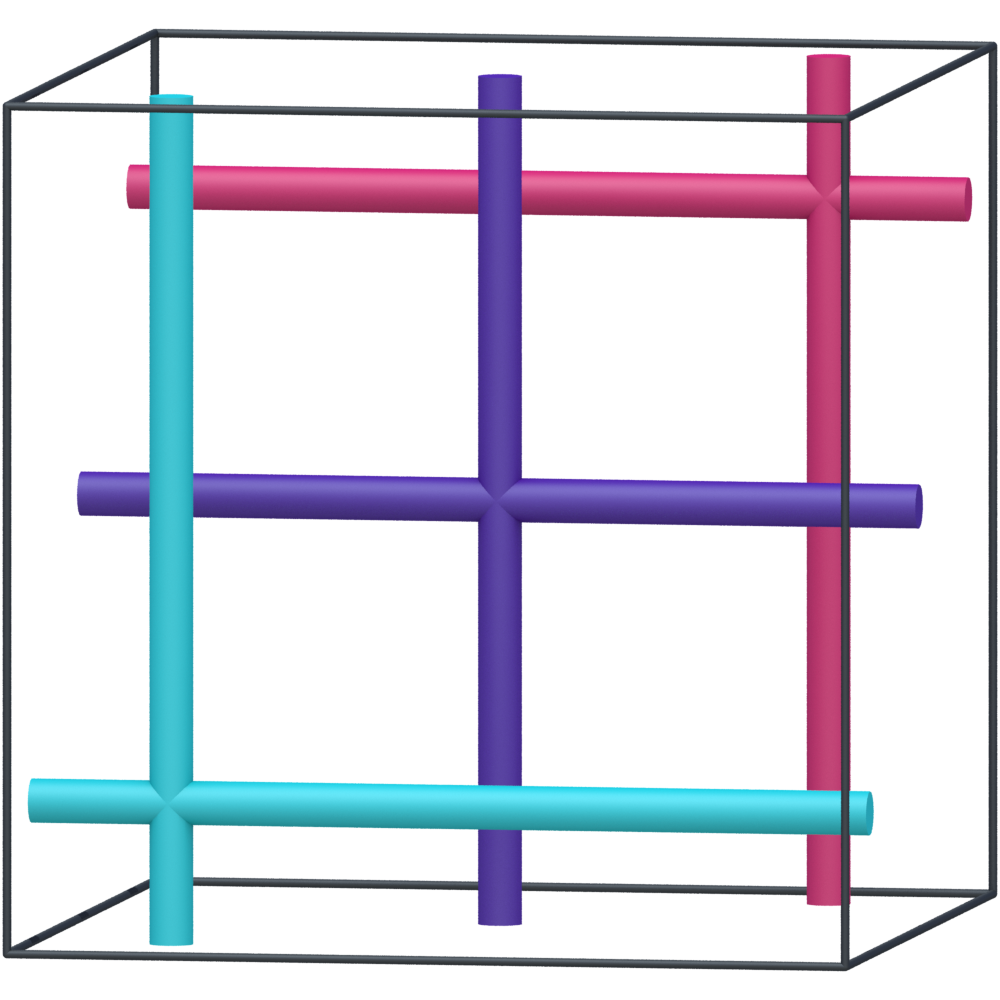}
        \caption{}
        \label{fig:triple_sql_uc}
    \end{subfigure}
    \caption{Untangling the structure in Fig. \ref{fig:borromean_sql_extended} into three layers of \textbf{sql} networks: This example shows that the concept of the untangling number defined in this paper is readily applicable on lower-periodic structures.}
    \label{fig:untangling_borromean_sql}
\end{figure*}

Not only does the untangling number seem to translate a different aspect of entanglement complexity than cycle analysis and HRN computation, but it is also based on a different perspective, one in which the knotting and linking of cycles are inherently encoded through crossing diagrams. This could in fact already be inferred from the observations in Sect. \ref{sec:least_tangled_embeddings} on \textbf{dia-c} and its primitive unit cell or on the ravelled embedding of \textbf{pcu} in Fig. \ref{fig:ravelled_pcu_extended}. Any embedding of a 3-periodic graph, regardless of the types of knots and links in its cycles, can be untangled into a ground state, provided that it is not one already as is the case for \textbf{dia-c}. By comparison, although it is possible to extend the concept of HRN to Brunnian links (links that can be transformed into unlinks by removing a single component) \cite{Alexandrov:eo5016}, the classical HRN is defined only for embeddings whose cycles are pairwise linked. The methods developed here, however, are readily defined and applicable to any embedding with any type of knotting or linking in its cycles, without additional considerations. For example, in Fig. \ref{fig:borromean_sql_extended}, we present a structure whose strong rings form Borromean rings, the simplest Brunnian link. In Fig. \ref{fig:untangling_borromean_sql}, we untangle the said structure into three layers of \textbf{sql} networks. In fact, it is this idea that the untangling number relies solely on crossing diagrams that permits its definition for 3-periodic tangles \cite{andriamanalina2025untanglingnumber3periodictangles}, a setting in which cycles are absent, and thus, cycle analysis and HRN computation are irrelevant.

\begin{figure*}[hbtp]

    \centering
    
    \begin{subfigure}[b]{0.27\textwidth}
        \includegraphics[width=0.6\textwidth]{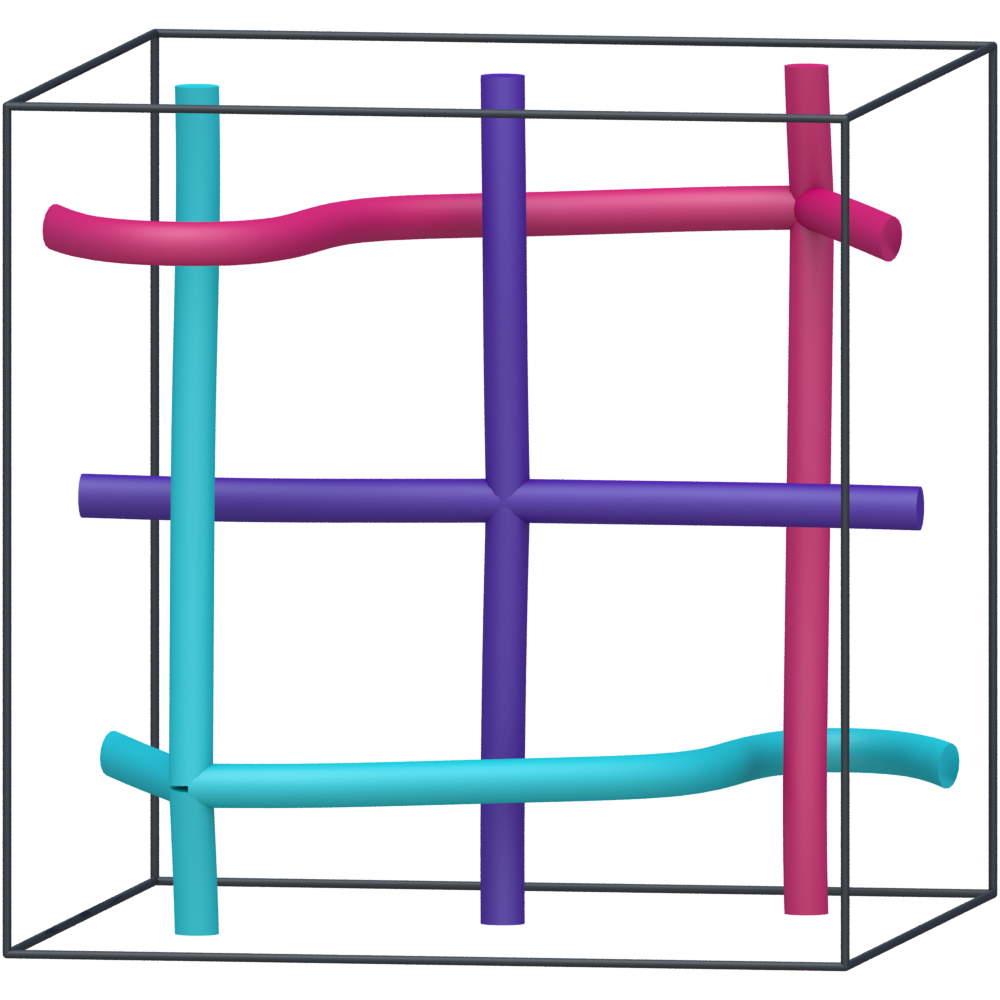}
        \caption{A unit cell of the embedding in Fig. \ref{fig:borromean_sql_extended}.}
        \label{fig:borromean_sql_minimal_uc}
    \end{subfigure}
    \begin{subfigure}[b]{0.54\textwidth}
        \includegraphics[width=0.29\textwidth]{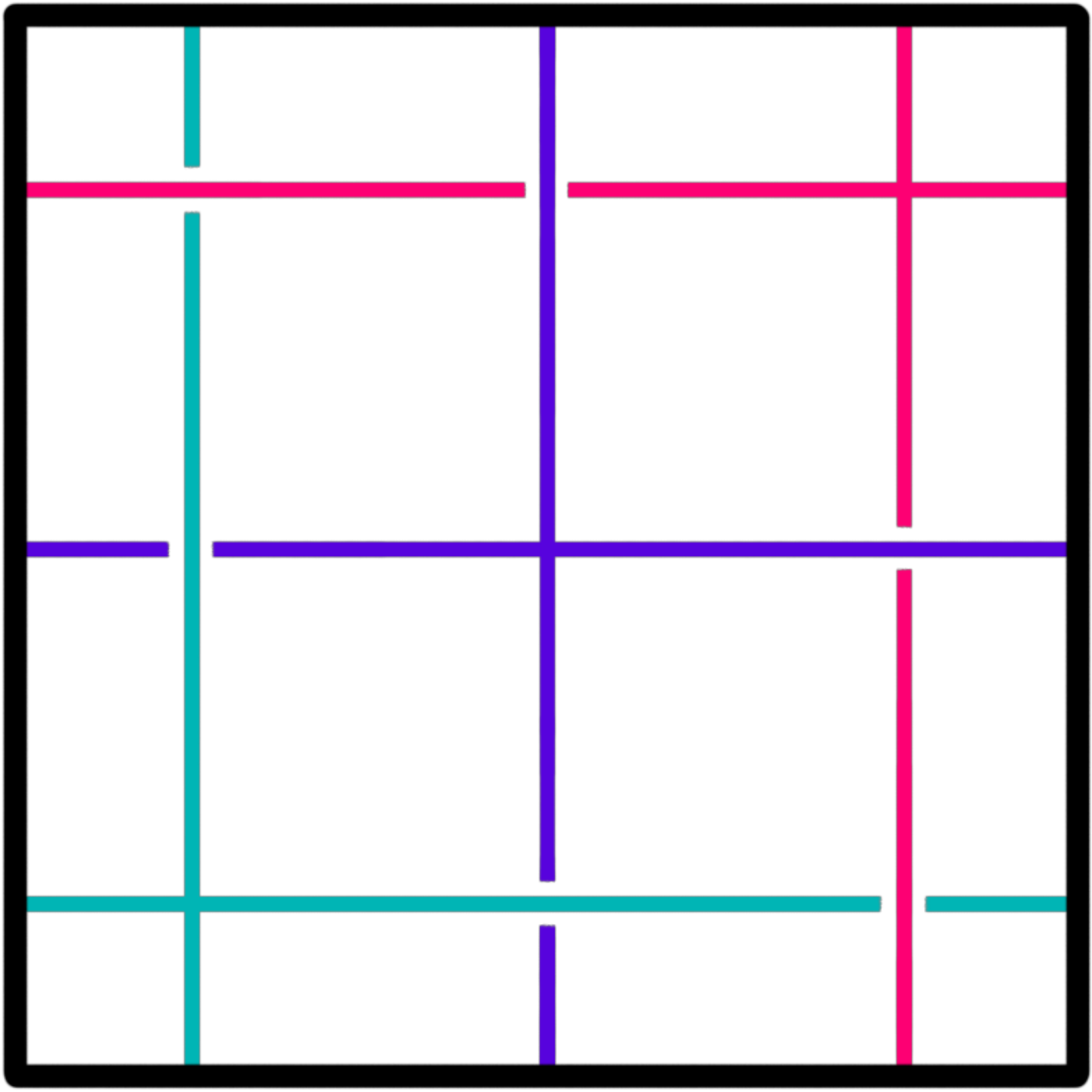}
        \hspace{0.2cm}
        \includegraphics[width=0.29\textwidth]{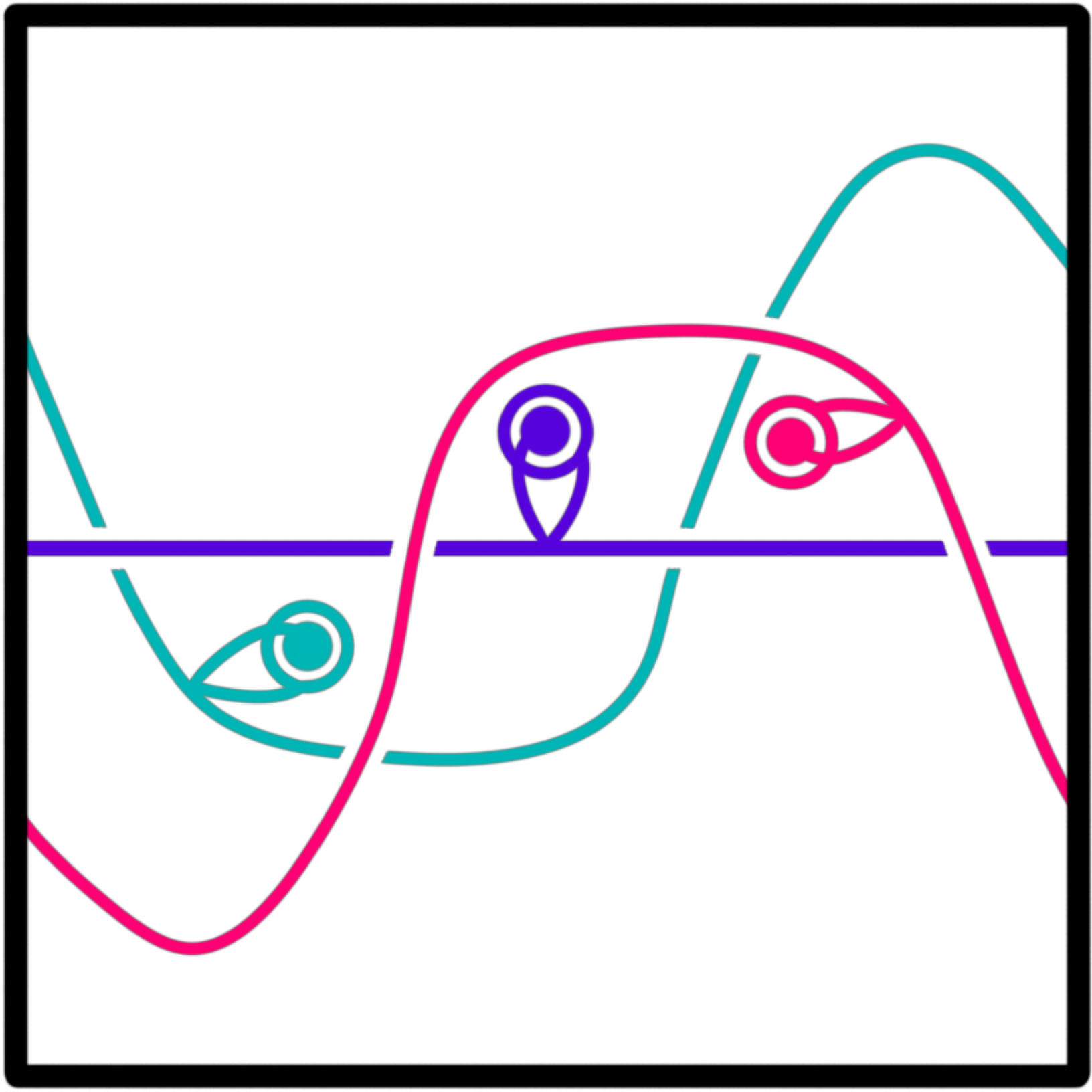}
        \hspace{0.2cm}
        \includegraphics[width=0.29\textwidth]{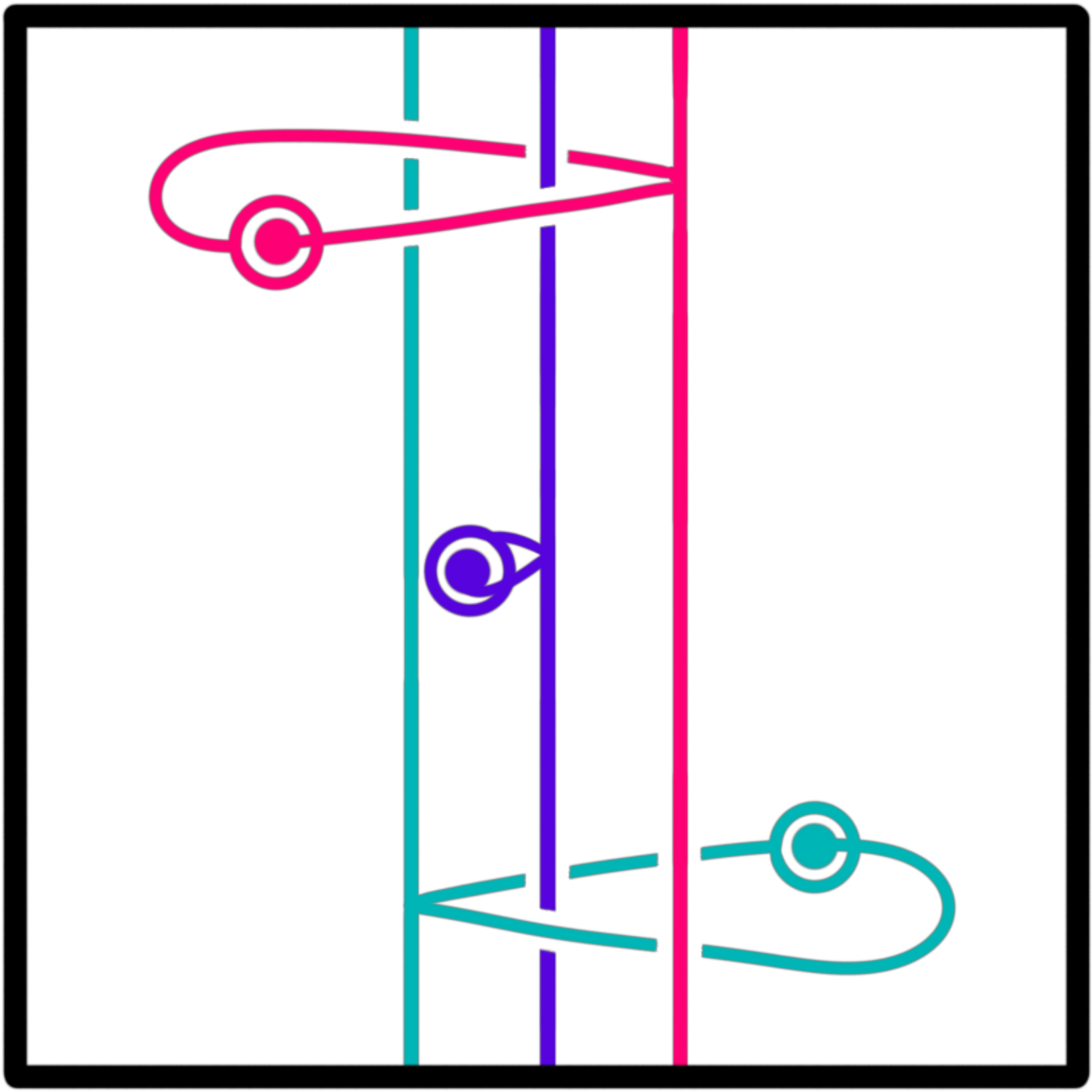}
        \caption{A tridiagram obtained from the unit cell in (a).}
        \label{fig:borromean_sql_tridia}
    \end{subfigure}
    \caption{A unit cell and a minimal tridiagram of the embedding in Fig. \ref{fig:borromean_sql_extended}: The first diagram in (b) has the same number of crossings as the diagram of the ground state shown in Fig. \ref{fig:untangling_borromean_sql_3}. This highlights the importance of working with tridiagrams, even for lower-periodic structures, through which the difference in the number of crossings becomes apparent. The minimum crossing number triplet of the tangled structure is $(6,6,8)$, whereas that of its ground state is $(6,0,0)$.}
    \label{fig:borromean_sql_uc_and_tridia}
\end{figure*}

 Note again that the structure in Fig. \ref{fig:borromean_sql_extended} and its ground state are technically 2-periodic structures. This example shows that the methods outlined in this paper can readily be applied to lower-periodic structures. It also highlights the importance of working with tridiagrams, as they provide additional information about crossings beyond that contained in single diagrams. Indeed, the two diagrams in Fig. \ref{fig:untangling_borromean_sql_1} and Fig. \ref{fig:untangling_borromean_sql_3} both have six crossings, yet the former represents a tangled structure, whereas the latter represents a ground state. The difference in the number of crossings is seen through tridiagrams and the associated minimum crossing number triplets, which are $(6,6,8)$, as shown in Fig. \ref{fig:borromean_sql_tridia}, and $(6,0,0)$, respectively.

The untangling number with respect to a given unit cell is truly dependent on the chosen unit cell. For example, in Fig. \ref{fig:untangling_dia-z}, we computed an upper bound for the untangling number of \textbf{dia-z} with respect to its primitive unit cell, and saw that two crossing changes were needed to transform it into the primitive unit cell of \textbf{dia}. Had we instead chosen a different unit cell, for example, one that untangles into the usual unit cell of \textbf{dia} shown in Fig. \ref{fig:dia_uc}, the result would have been different. This may seem inconvenient, but we argue that choosing a unit cell with which to work is not very different from considering strong rings only in a cycle analysis or in the computation of the HRN. In fact, this idea aligns with the methods outlined in \cite{Power:ib5087}. Furthermore, the untangling number can be stabilised by considering the minimum untangling number over all unit cells of a given embedding of a graph, leading to the definition of the \textit{minimum untangling number}, which is immediately an invariant of the embedding.

\section{Computability of the untangling number}\label{sec:computability}
As we have already mentioned, the untangling number is difficult to compute, and only upper bounds can be effectively given. This limitation is similar to that of the classical unknotting number \cite{Murasugi1996}, the determination of which remains to this day an open problem. Moreover, contrary to the case of 3-periodic tangles \cite{andriamanalina2025untanglingnumber3periodictangles}, there is currently no method for characterising the minimum crossing number triplets of ground states, which makes the computation even more challenging. To address this, lower bounds may be developed in future work, both for the minimum crossing number triplet of a ground state and for the untangling number itself. Inspiration can be drawn from the lower bounds for the classical unknotting number that have been developed over the past decades \cite{Murasugi1965,nakanishi_I,ma-qiu}. The majority of these lower bounds are computed via diagrams, which are now available for 3-periodic graphs. However, the extension of these techniques to the 3-periodic setting presents new challenges.

Furthermore, while the development of lower bounds is important from a theoretical perspective, we argue that, in practice, an empirical value given as an upper bound already provides insight into the structures, both for the crossing number of a ground state and for the untangling number. For the crossing number of a ground state, we assumed for example that the structure consisting of three layers of \textbf{sql} networks, whose unit cell is shown in Fig. \ref{fig:untangling_borromean_sql_3}, is a ground state, and that its minimum crossing number triplet is $(6,0,0)$. Theoretically speaking, there is no rigorous mathematical proof that confirms either of these assumptions. However, intuitively, there is no reason to believe that the structure has a tridiagram with fewer crossings, or that an embedding with fewer crossings exists. We argue that this intuition is sufficient from a practical point of view. Several ground states can in fact be empirically determined in the same manner. 
For the untangling number, developing computational methods such as those presented in \cite{Applebaum18082025} is a promising direction for determining its value. To enable this, however, a computational representation of diagrams of 3-periodic graphs must first be developed.

\section{Conclusions}\label{sec:conclusion}
In this article, we give a definition of least tangled embeddings of 3-periodic networks, that we call \textit{ground states}. The definition is given through the use of diagrammatic representations and minimum number of crossings as well as crossing changes. We gave examples of ground states, such as the barycentric embeddings of the \textbf{srs}, the \textbf{dia} and the \textbf{pcu} networks, as well as \textbf{dia-c} the standard interpenetration of two \textbf{dia} networks. We provided interesting results regarding ground states, such as the fact that an embedding, such as \textbf{dia-c}, may possess linked cycles and non-trivial HRN, yet may also possess a unit cell with no crossings.

We also defined a measure of entanglement complexity called the \textit{untangling number} of 3-periodic networks, which is analogous to the unknotting number of classical knots and links. The untangling number provides a new insight into the complexity of entanglement in 3-periodic networks, that of the ease with which an embedding can be untangled to its least tangled state. This perspective is different from that offered by other measures, such as cycle analysis and HRN computation.

Although the focus of this paper is on 3-periodic networks, the notions of ground state and untangling number extend naturally to lower-periodic networks, as well as to finite graphs. Indeed, the definitions rely solely on crossing information determined by projections along three non-coplanar axes, a construction that is well defined for any type of graph. In fact, it is for this reason that the methods here constitute a bridge between the analysis of entanglement in 3-periodic networks and that of 3-periodic filament entanglements, a setting in which cycles are absent and methods such as cycle analysis and HRN computation are irrelevant.

\section*{Conflicts of interest}
There are no conflicts to declare.




\section*{Acknowledgements}

This work is funded by the Deutsche Forschungsgemeinschaft (DFG, German Research Foundation) - Project number 468308535, and partially supported by RIKEN iTHEMS and by the JSPS KAKENHI Grant-in-Aid for Early-Career Scientists 25K17246. We also acknowledge the education licence for Houdini/SideFX, which was used for visualisation.



\balance


\onecolumn

\vspace{5cm}

\section*{Supplementary Information}

\section{Additional details on concepts introduced in the main manuscript}\label{sec:additional_details_SI}

\subsection{Equivalence of embeddings of 3-periodic graphs}
Consider a graph $\mathcal{K}$ and an embedding $K$ of $\mathcal{K}$. Essentially, deforming $K$ in space does not change its entanglement type, as long as one does not let the edges pass through each other or through vertices. Mathematically, this is the notion of ambient isotopy; in the case of graphs, one more precisely considers the notion of piecewise linear ambient isotopy (p.l. ambient isotopy), which accounts for vertices. Here, similarly to the case of 3-periodic tangles explored in \cite{ANDRIAMANALINA2025109346}, we focus on deformations that occur within unit cells to preserve the periodicity of structures. These deformations are converted into moves performed on diagrams, called \textit{$R$-moves}, which fully capture the 3-dimensional information in the unit cells.

\begin{figure}[ht]
    \centering
    \begin{subfigure}[b]{5cm}
        \includegraphics[width=\textwidth]{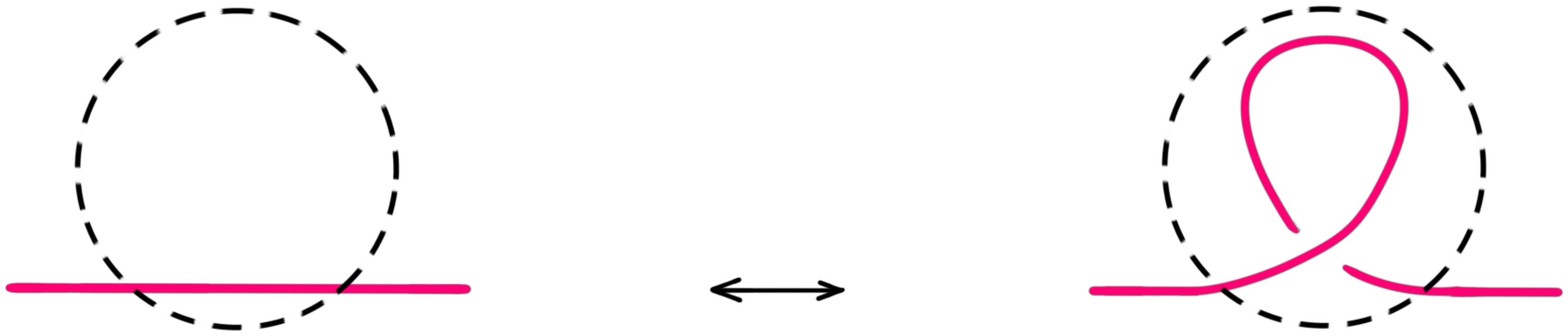}
        \caption{The $R_1$ move}
        \label{fig:Reid_move_I}
    \end{subfigure}
    \hspace{0.5cm}
    \begin{subfigure}[b]{5.5cm}
        \includegraphics[width=\textwidth]{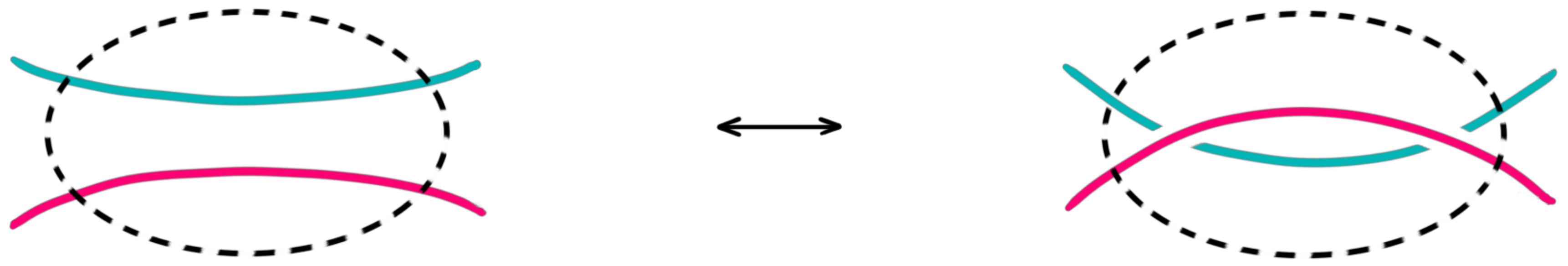}
        \caption{The $R_2$ move}
        \label{fig:Reid_move_II}
    \end{subfigure}
    \vskip\baselineskip
    \begin{subfigure}[b]{5cm}
        \includegraphics[width=\textwidth]{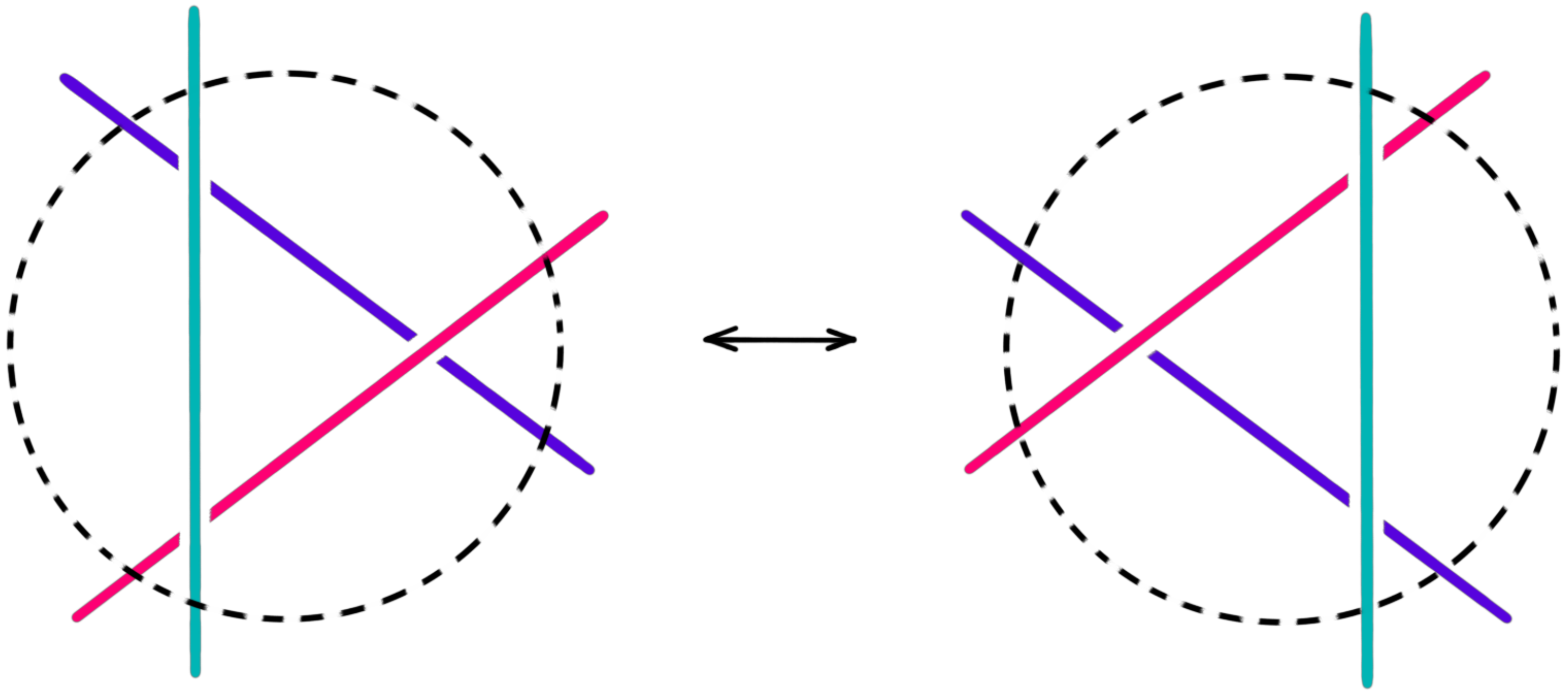}
        \caption{The $R_3$ move}
        \label{fig:Reid_move_III}
    \end{subfigure}
    \caption{The three Reidemeister moves, first defined to capture isotopies of usual knots and links in $\mathbb{R}^3$. They naturally extend to diagrams of 3-periodic graphs.}
    \label{fig:usual_Reid_moves}
\end{figure}

The \textit{$R$-moves} are a set of thirteen moves that we list below. Every isotopy deformation in a given unit cell can be converted into a finite sequence of $R$-moves. The first three moves $R_1$, $R_2$ and $R_3$ correspond to the moves originally developed to capture deformations of classical knots and links in space \cite{BurdeZieschangHeusener+2013}. They naturally extend to diagrams of 3-periodic graphs. The $R_1$ move allows one to add or remove a twist in a given edge, which also adds or removes a crossing as shown in Fig. \ref{fig:Reid_move_I}. The move $R_2$ adds or removes two crossings, as shown in Fig. \ref{fig:Reid_move_II}. The third move $R_3$ allows one to slide an edge from one side of a crossing to the other side of the crossing, as depicted in Fig. \ref{fig:Reid_move_III}.

\begin{figure}[hbtp]
    \centering
    \begin{subfigure}[b]{5.5cm}
    \centering
        \includegraphics[width=\textwidth]{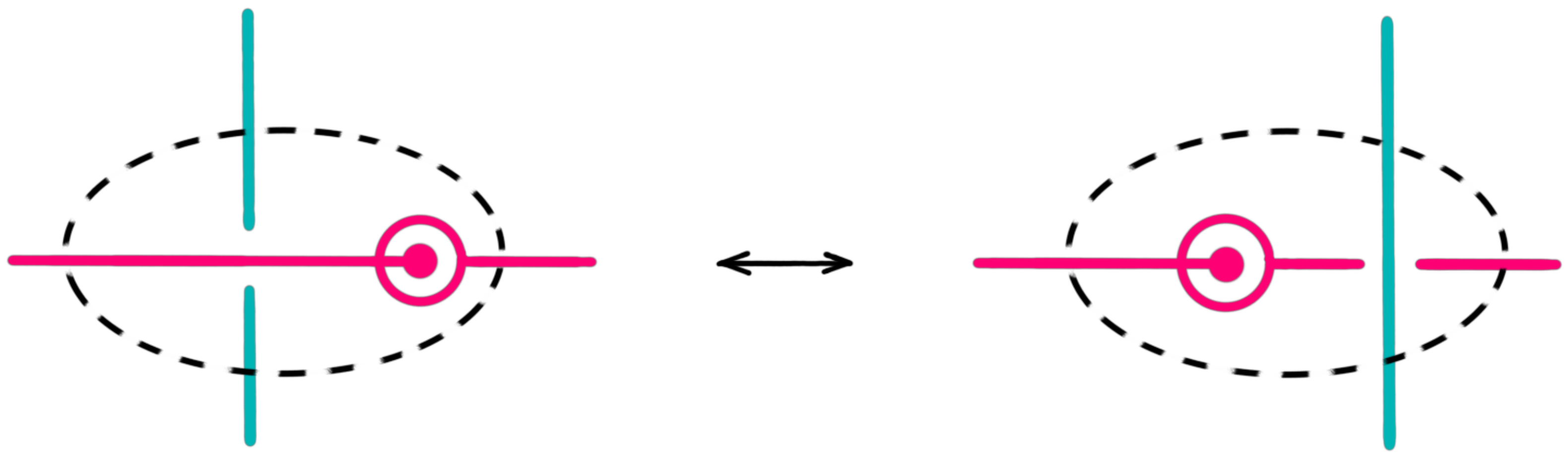}
        \caption{The $R_4$ move}
        \label{fig:Reid_move_IV}
        \vspace{0.5cm}
    \end{subfigure}
    \hspace{0.5cm}
    \begin{subfigure}[b]{5.5cm}
    \centering
        \includegraphics[width=\textwidth]{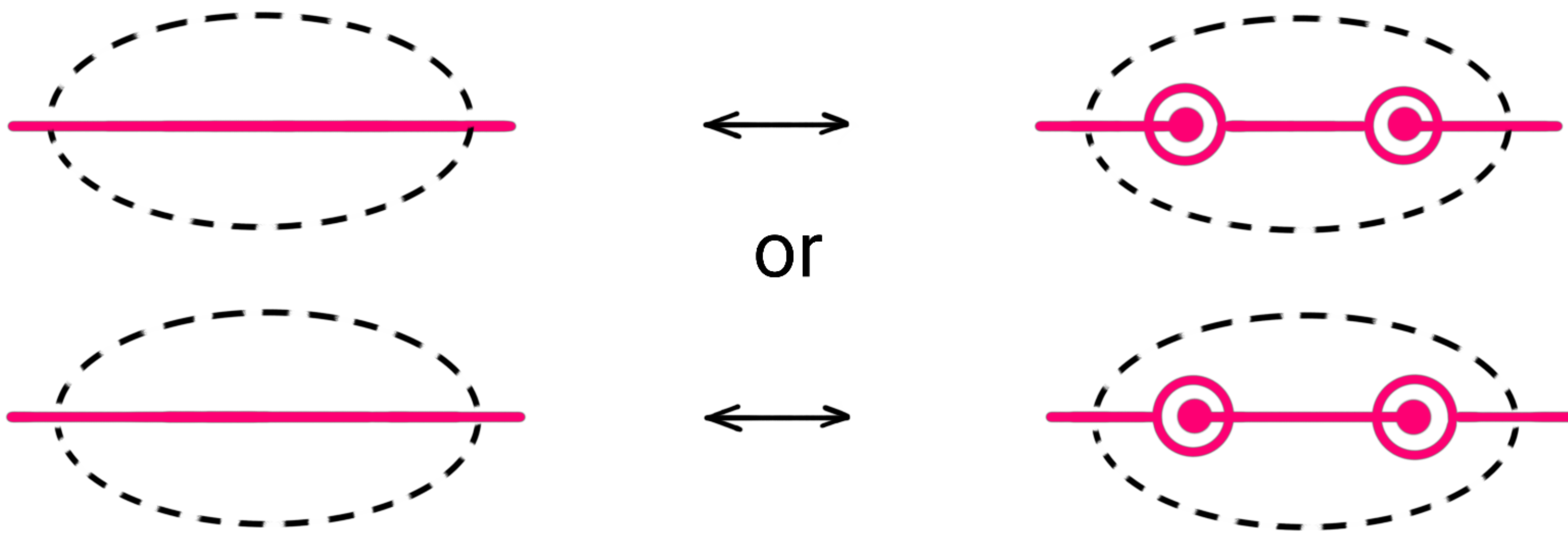}
        \caption{The $R_5$ move}
        \label{fig:Reid_move_V}
    \end{subfigure}
    
    \vskip\baselineskip
    
    \begin{subfigure}[b]{5cm}
    \centering
        \includegraphics[width=\textwidth]{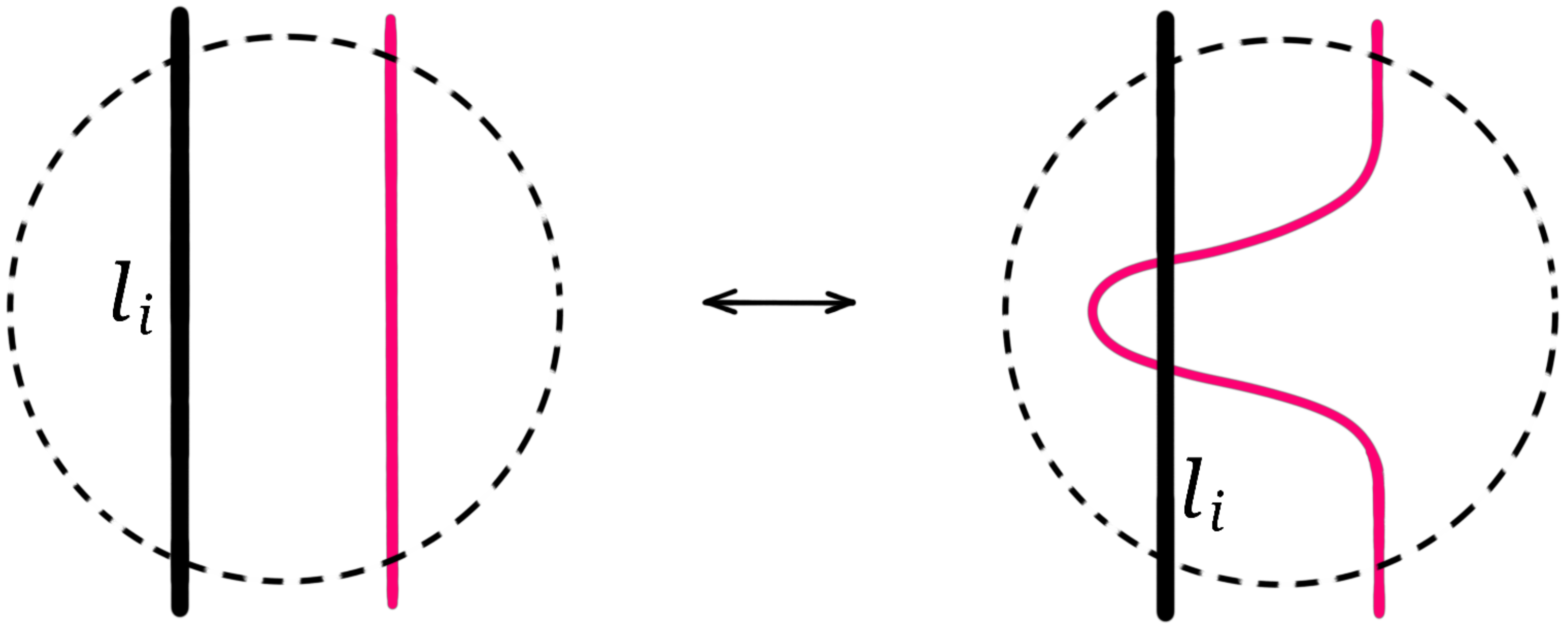}
        \caption{The $R_6$ move}
        \vspace{1.25cm}
        \label{fig:Reid_move_VI}
    \end{subfigure}
    \hspace{1cm}
    \begin{subfigure}[b]{5cm}
    \centering
        \includegraphics[width=\textwidth]{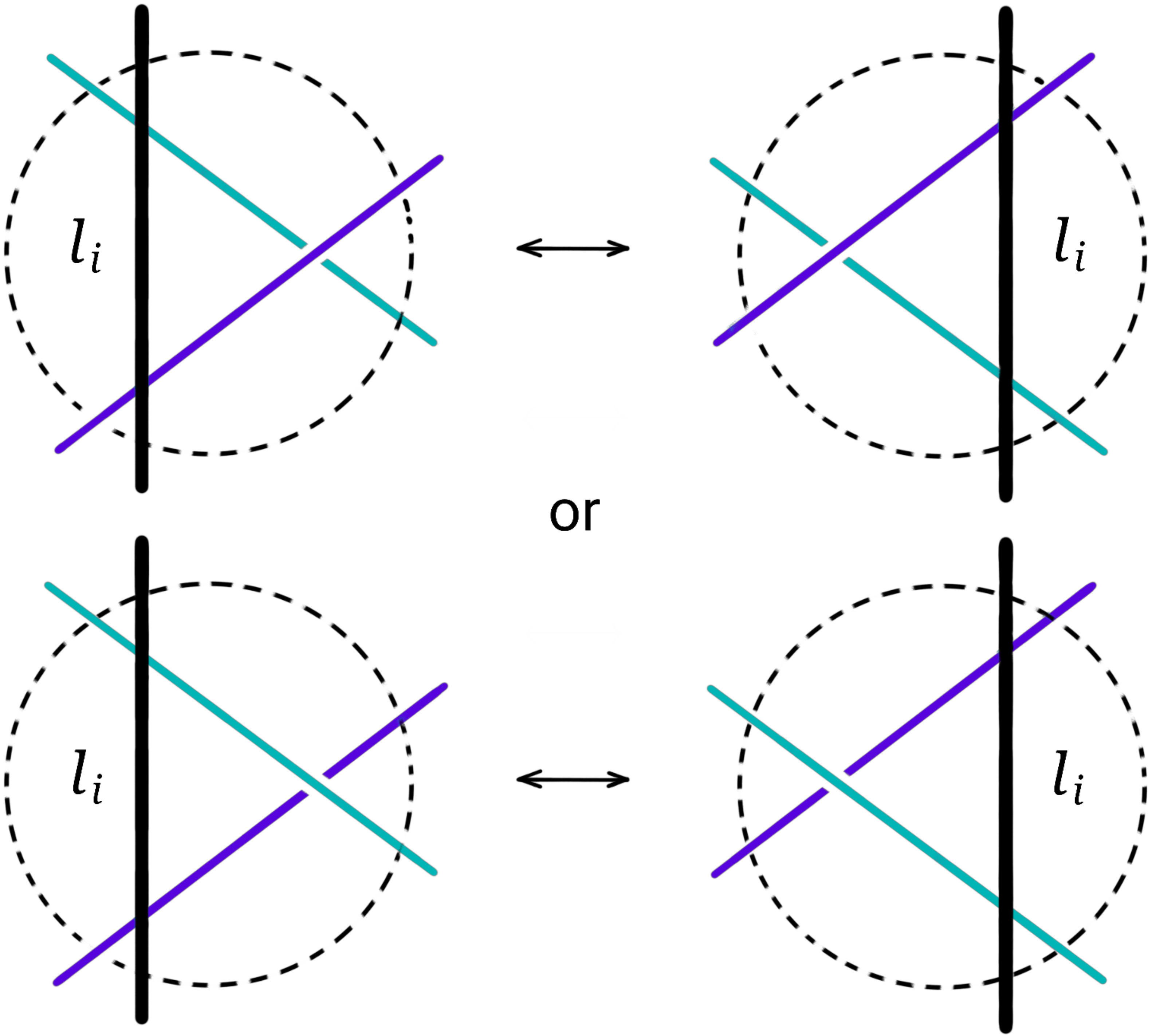}
        \caption{The $R_7$ move}
        \label{fig:Reid_move_VII}
    \end{subfigure}
    
    \vskip\baselineskip
    
    \begin{subfigure}[b]{5cm}
    \centering
        \includegraphics[width=\textwidth]{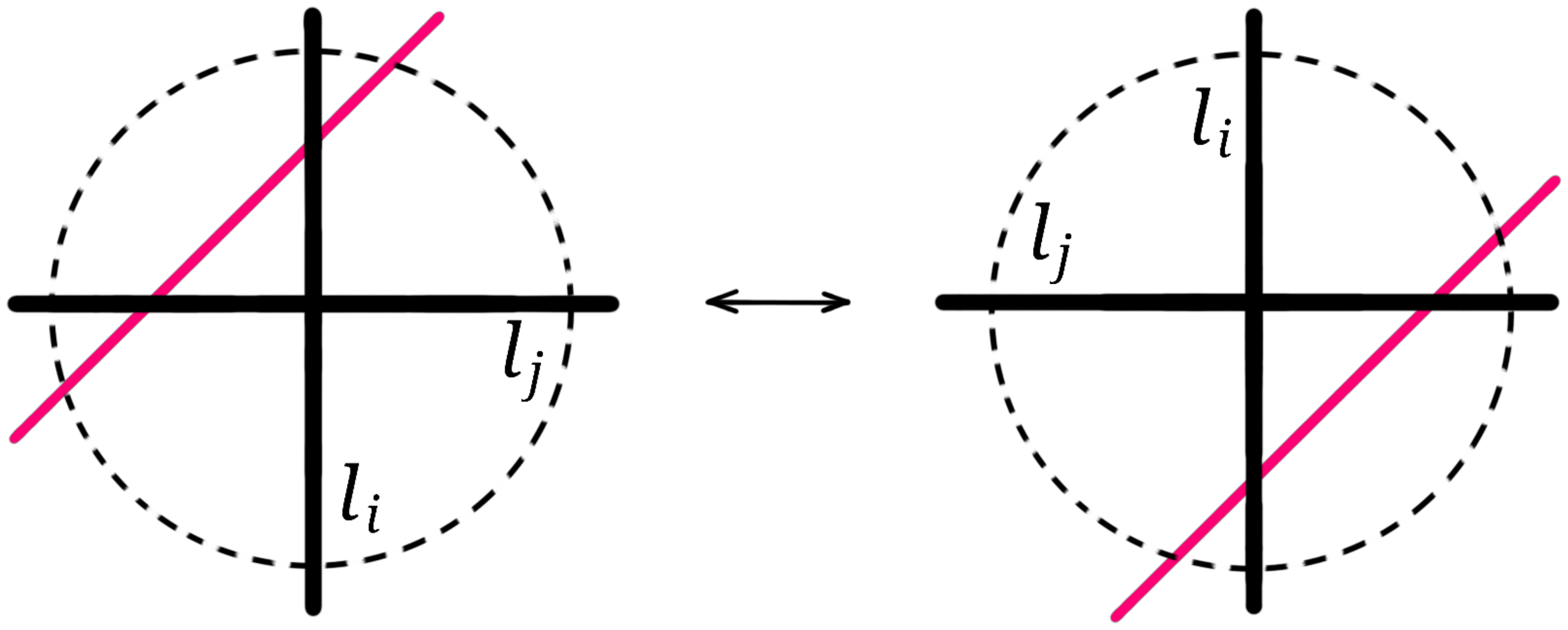}
        \caption{The $R_8$ move}
        \label{fig:Reid_move_VIII}
    \end{subfigure}
    \hspace{0.75cm}
    \begin{subfigure}[b]{5.5cm}
    \centering
        \includegraphics[width=\textwidth]{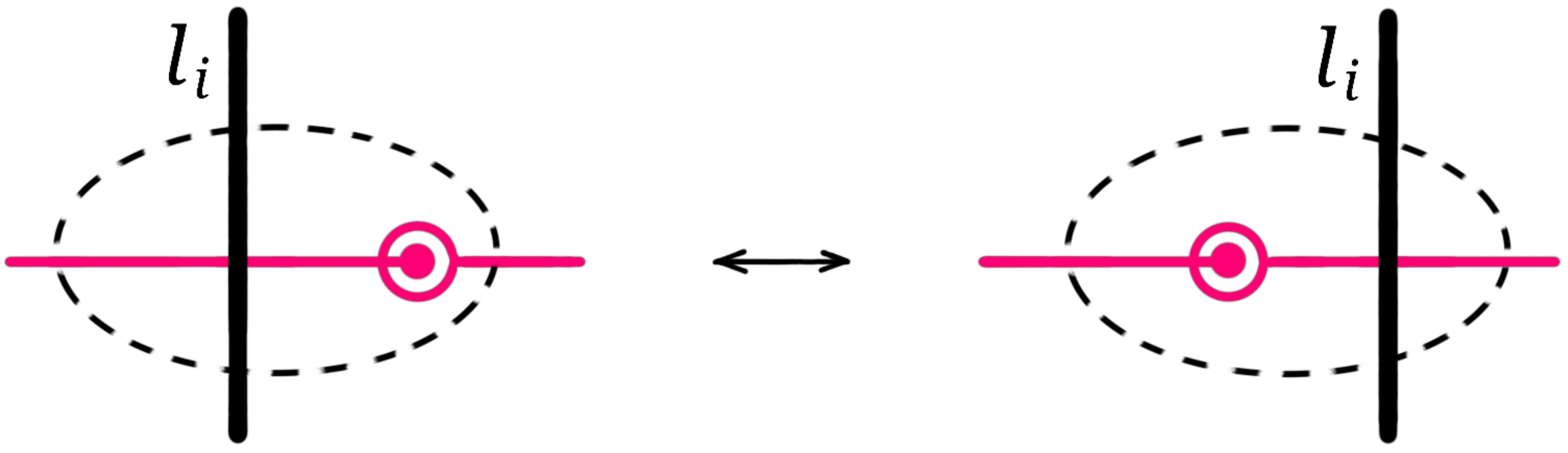}
        \caption{The $R_9$ move}
        \label{fig:Reid_move_IX}
        \vspace{0.375cm}
    \end{subfigure}
        \caption{The six new moves, originally defined to capture ambient isotopies within unit cells of 3-periodic tangles. Here $l_i$ and $l_j$, where $i,j = 1,2$, represent the edges of the squares delimiting a diagram.}
    \label{fig:new_Reid_moves}
\end{figure}

In addition to the previous three moves, there are the six moves that were first defined in \cite{ANDRIAMANALINA2025109346} to capture ambient isotopies of unit cells of 3-periodic tangles. The $R_4$ move of Fig. \ref{fig:Reid_move_IV} depicts an edge that slides between the front and back parts of another that goes through the identified front and back faces of the unit cell. The $R_5$ move of Fig. \ref{fig:Reid_move_V} allows one to pass a part of an edge through the front and back faces of a unit cell. The $R_6$ move, Fig. \ref{fig:Reid_move_VI}, depicts an edge that goes through one side of the square delimiting the diagram. The $R_7$ move shown in Fig. \ref{fig:Reid_move_VII} represents a crossing that goes through one side of the square. The $R_8$ move shown in Fig. \ref{fig:Reid_move_VIII} corresponds to an edge passing through the corners of the square. The $R_9$ move of Fig. \ref{fig:Reid_move_IX} depicts a point that is the projection of an edge intersecting the front and back faces of the unit cell, going through one side of the square.

\begin{figure}[ht]
    \centering
    \begin{subfigure}[b]{5cm}
        \includegraphics[width=\textwidth]{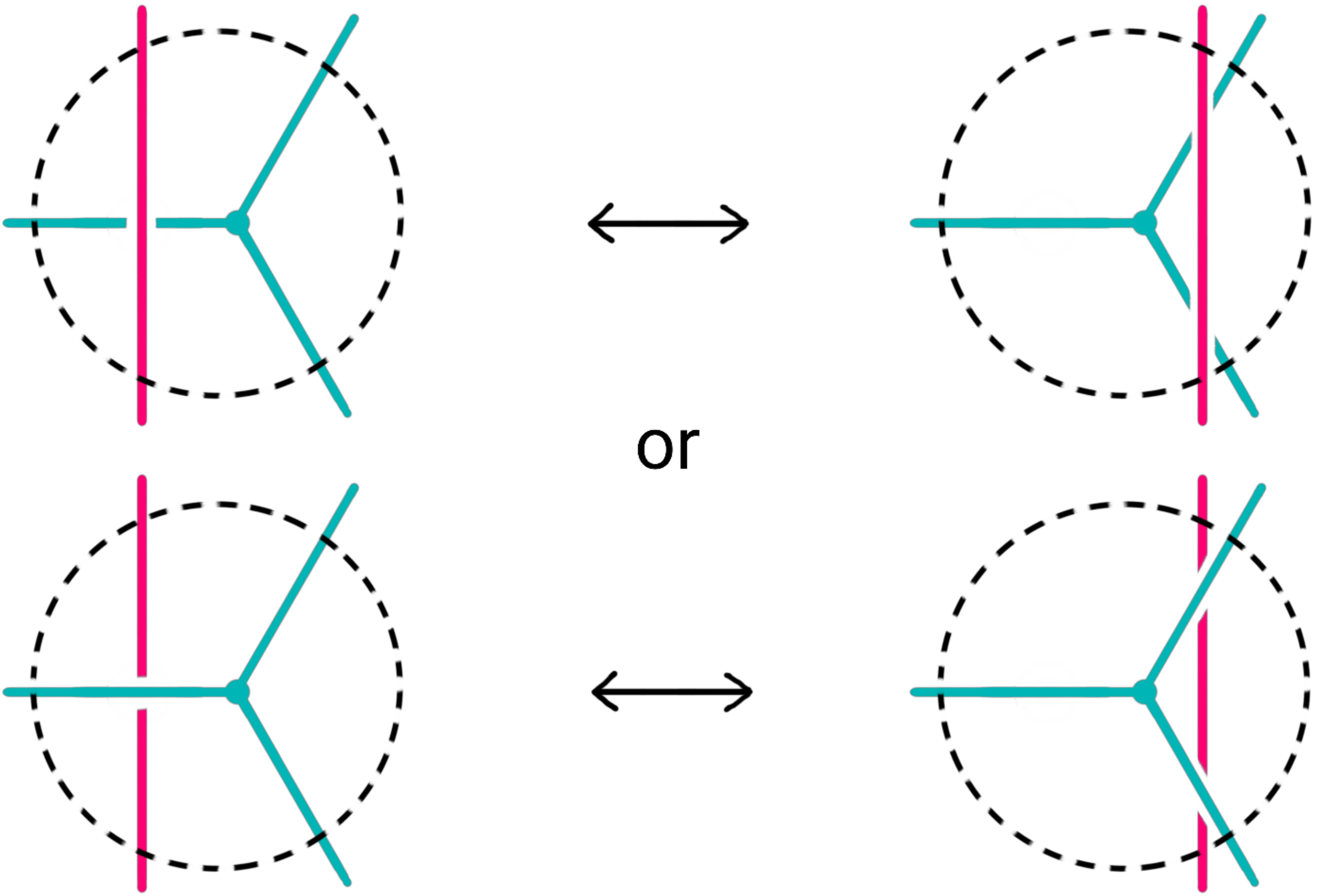}
        \caption{The $R_{10}$ move}
        \label{fig:Reid_move_X}
    \end{subfigure}
    \hspace{0.5cm}
    \begin{subfigure}[b]{5cm}
        \includegraphics[width=\textwidth]{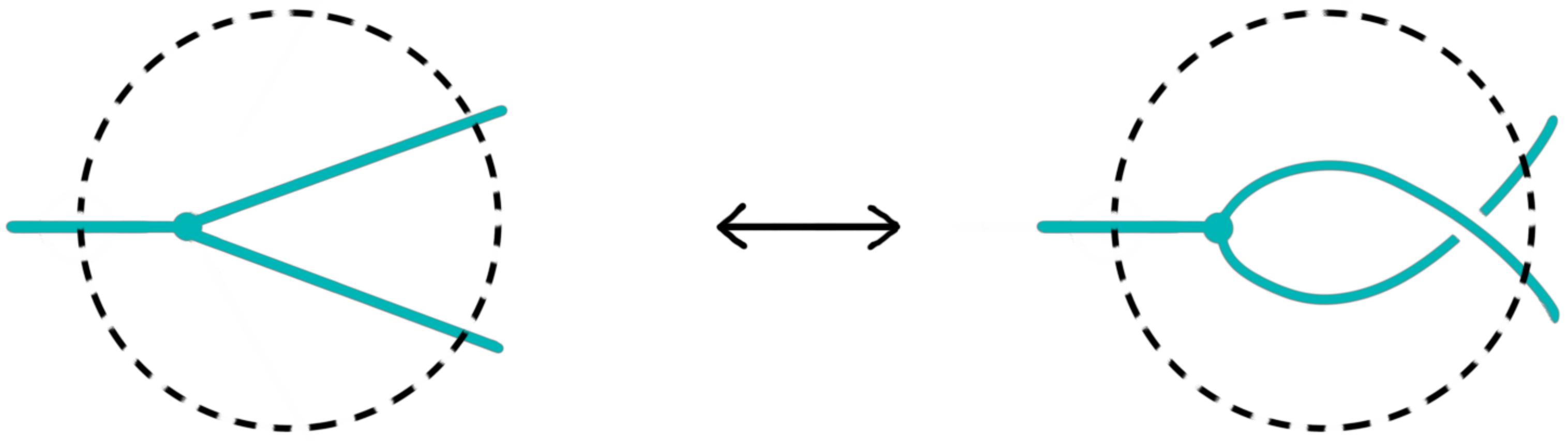}
        \caption{The $R_{11}$ move}
        \vspace{1cm}
        \label{fig:Reid_move_XI}
    \end{subfigure}
    \caption{The two moves originally defined to capture p.l. ambient isotopies of embeddings of finite graphs in space. The number of edges connected to a vertex is arbitrary.}
    \label{fig:usual_Reid_moves_graphs}
\end{figure}

In addition to these, there are also the two moves originally defined to capture p.l. ambient isotopies of finite graphs embedded in space \cite{R_moves_graph}. The $R_{10}$ move of Fig. \ref{fig:Reid_move_X} depicts an edge that slides over or under several others sharing a common vertex. The $R_{11}$ move of Fig. \ref{fig:Reid_move_XI} depicts an edge sliding over another that shares a common vertex with it.

\begin{figure}[ht]
    \centering
    \begin{subfigure}[b]{5cm}
        \includegraphics[width=\textwidth]{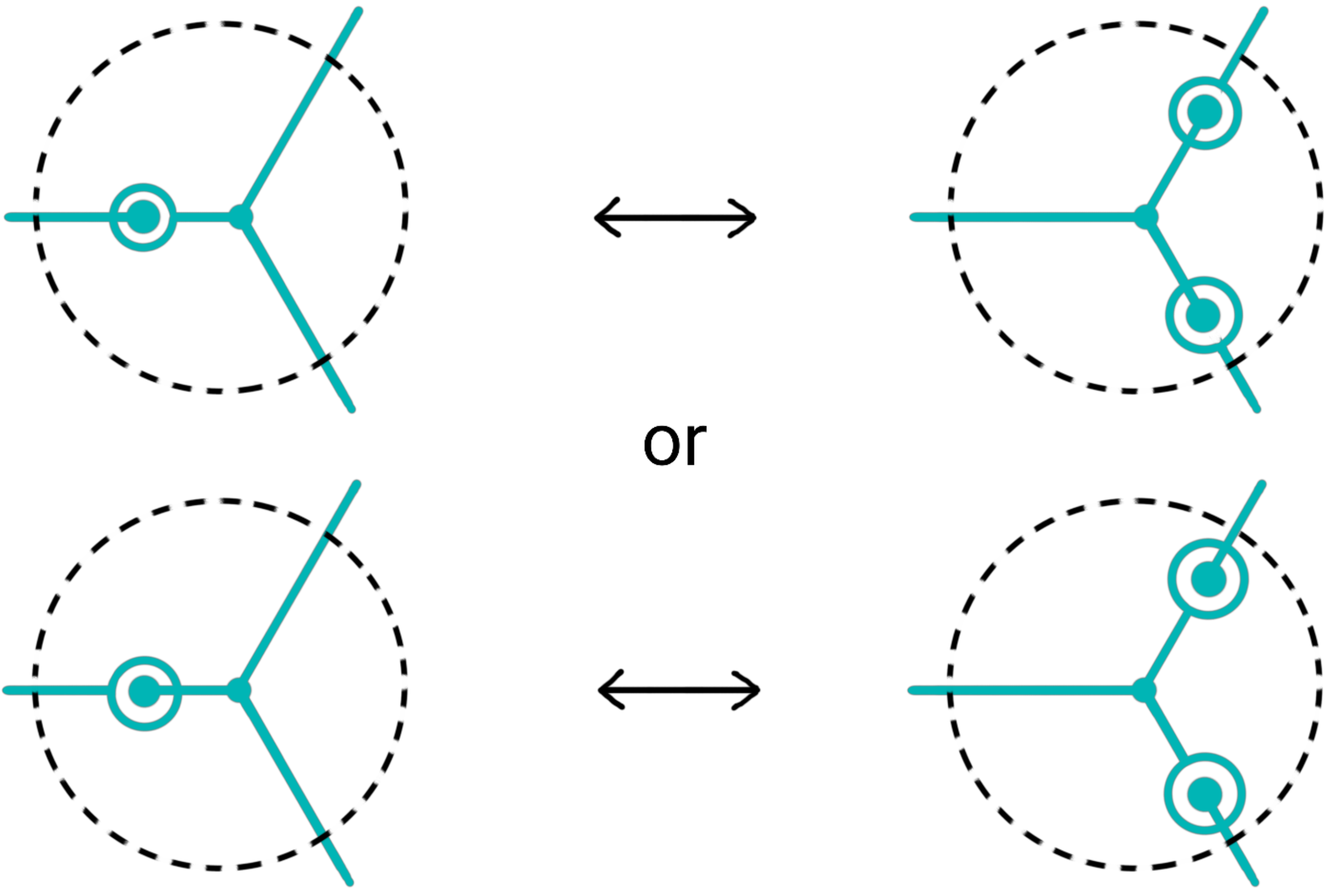}
        \caption{The $R_{12}$ move}
        \label{fig:Reid_move_XII}
    \end{subfigure}
    \hspace{0.5cm}
    \begin{subfigure}[b]{5cm}
        \includegraphics[width=\textwidth]{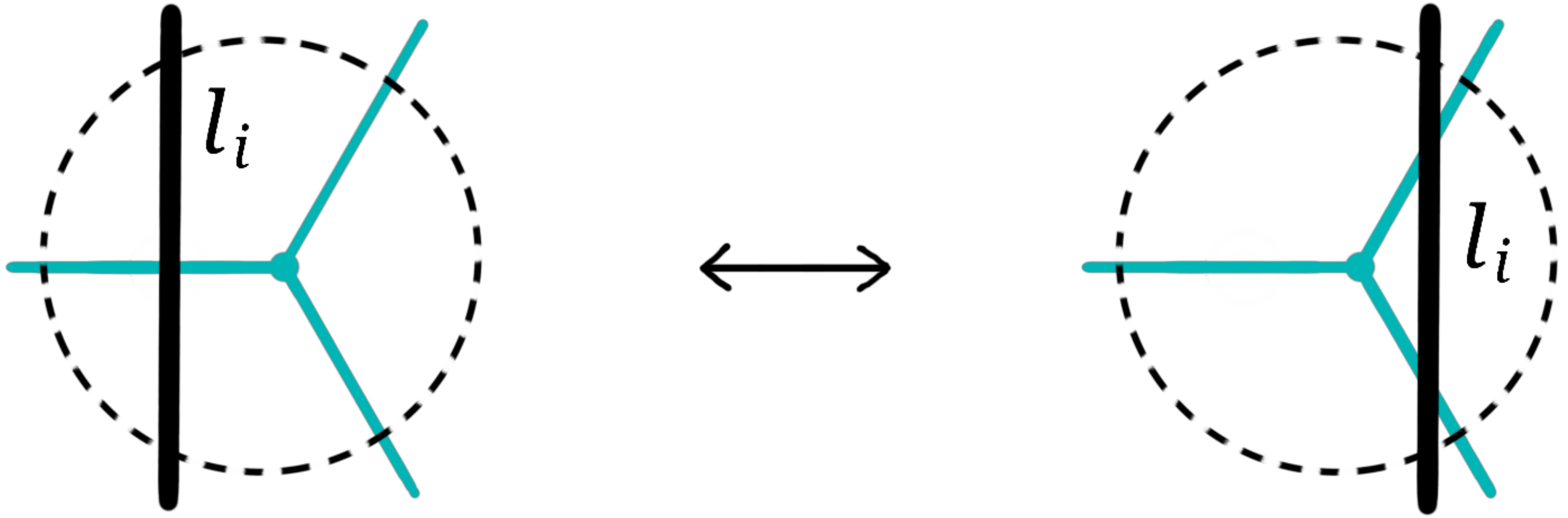}
        \caption{The $R_{13}$ move}
        \vspace{1cm}
        \label{fig:Reid_move_XIII}
    \end{subfigure}
    \caption{The two moves specific to p.l. ambient isotopies of 3-periodic embeddings of graphs. The number of edges connected to a vertex is arbitrary.}
    \label{fig:additional_Reid_move_graphs}
\end{figure}

Finally, there are additional moves specific to 3-periodic graphs. The $R_{12}$ move of Fig. \ref{fig:Reid_move_XII} corresponds to the passing of a vertex through the identified front and back faces of a unit cell. The $R_{13}$ move of Fig. \ref{fig:Reid_move_XIII} depicts a vertex going through a side of the square delimiting the diagram.

Note that in Fig. \ref{fig:usual_Reid_moves_graphs} and Fig. \ref{fig:additional_Reid_move_graphs}, the number of edges connected to a vertex is arbitrary. We choose to work with only three edges for simplicity. The moves can easily be extended to any number of edges.

\begin{figure}
    \centering
    \begin{subfigure}[b]{0.18\textwidth}
        \centering
        \includegraphics[width=0.9\textwidth]{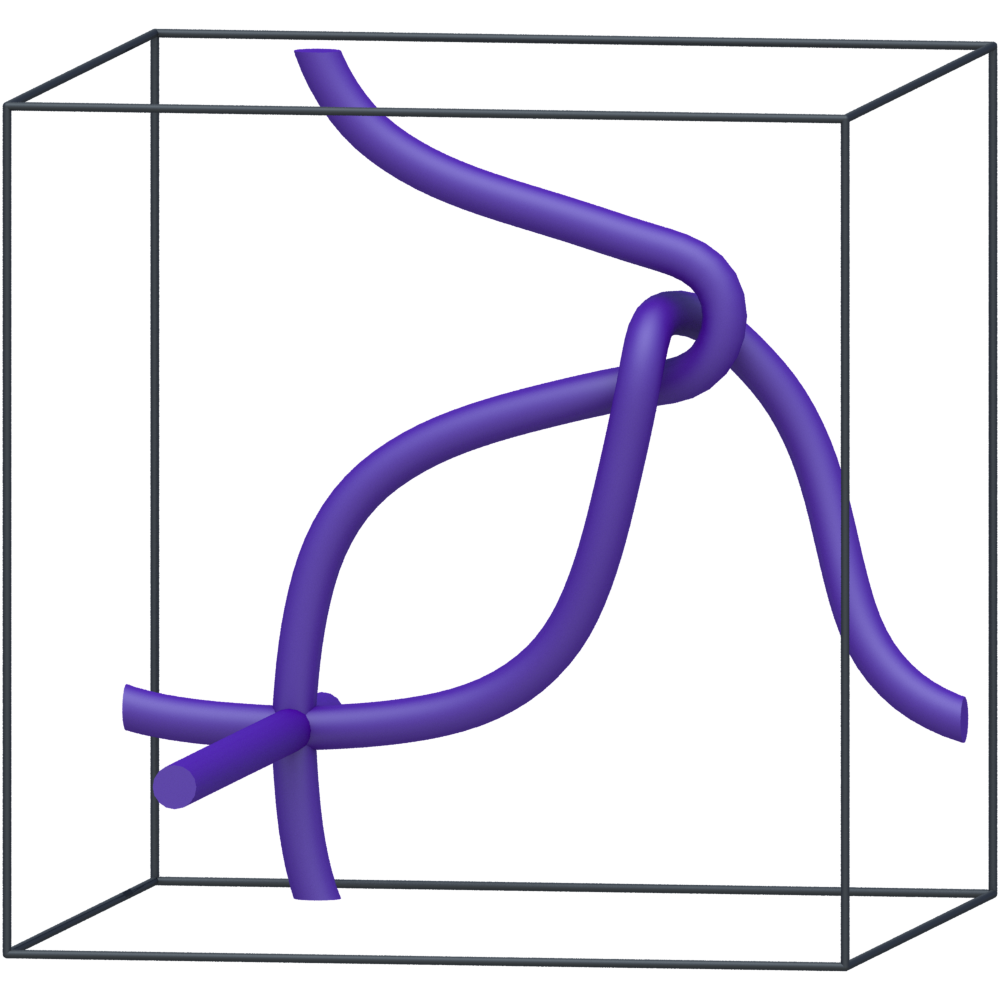}
        \caption{}
        \label{fig:pcu_1_uc}
    \end{subfigure}
    \begin{subfigure}[b]{0.18\textwidth}
        \centering
        \includegraphics[width=0.9\textwidth]{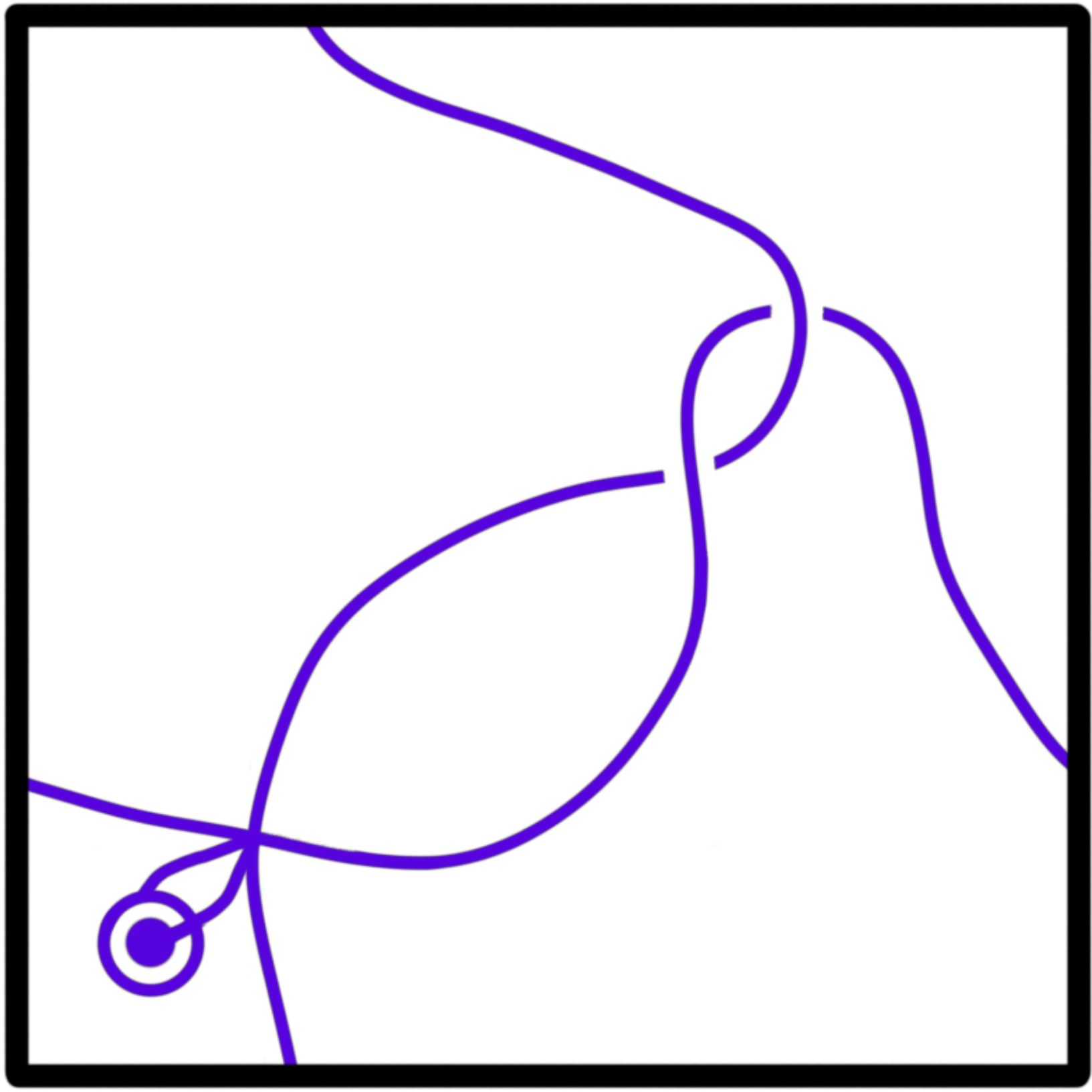}
        \caption{}
        \label{fig:pcu_1_dia}
    \end{subfigure}
    \begin{subfigure}[b]{0.18\textwidth}
        \centering
        \includegraphics[width=0.9\textwidth]{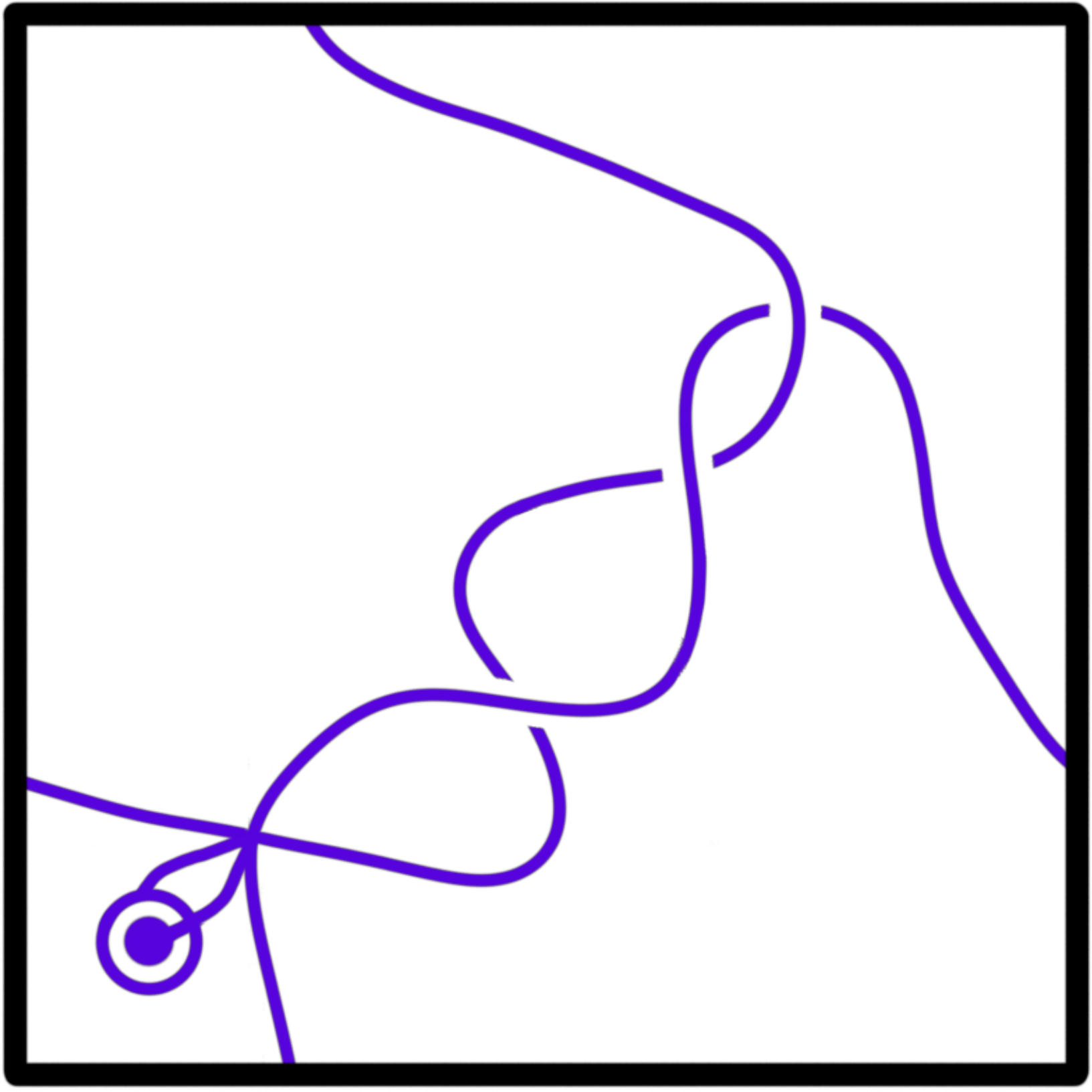}
        \caption{}
        \label{fig:pcu_2_dia}
    \end{subfigure}
    \begin{subfigure}[b]{0.18\textwidth}
        \centering
        \includegraphics[width=0.9\textwidth]{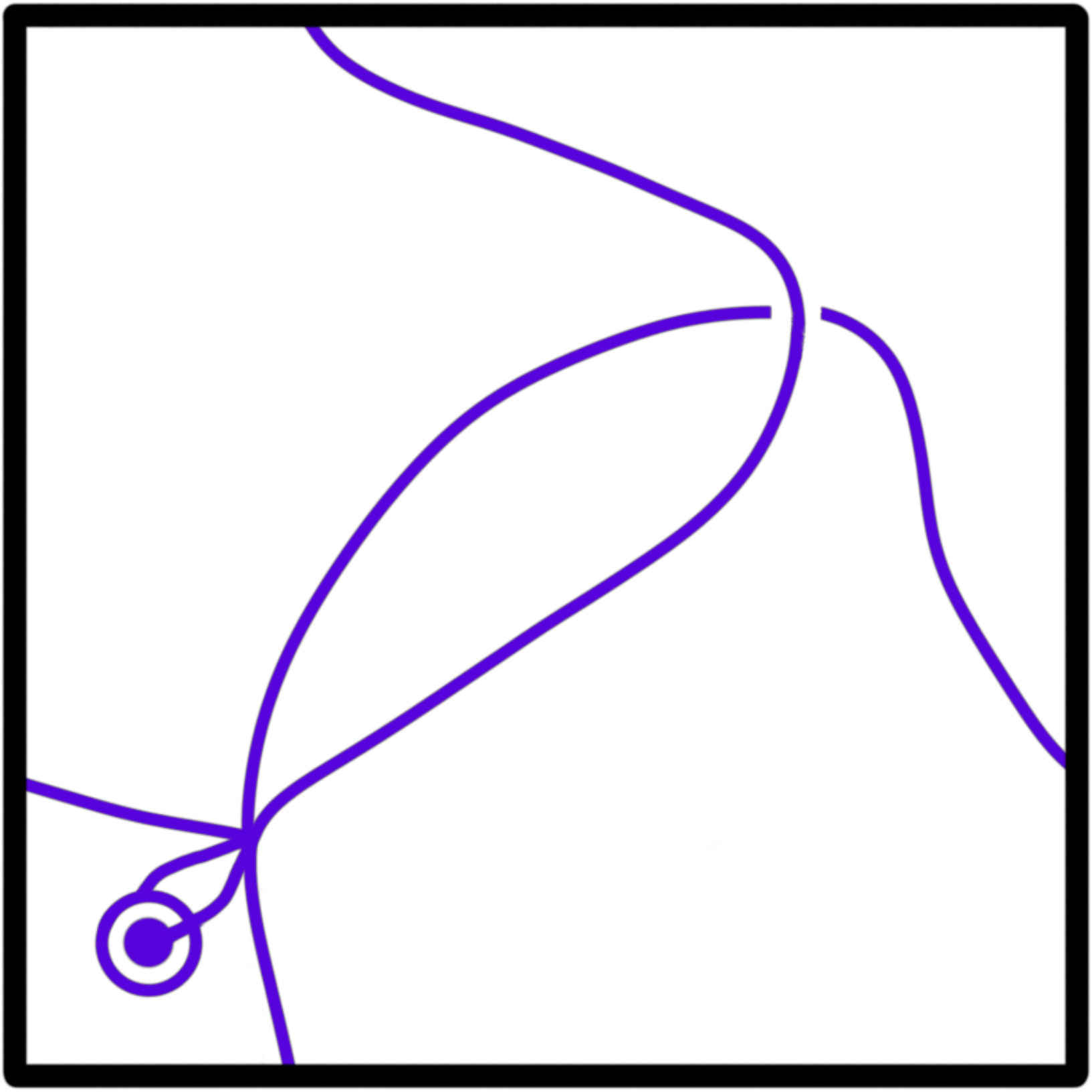}
        \caption{}
        \label{fig:pcu_3_dia}
    \end{subfigure}

    \vskip\baselineskip
    
    \begin{subfigure}[b]{0.18\textwidth}
        \centering
        \includegraphics[width=0.9\textwidth]{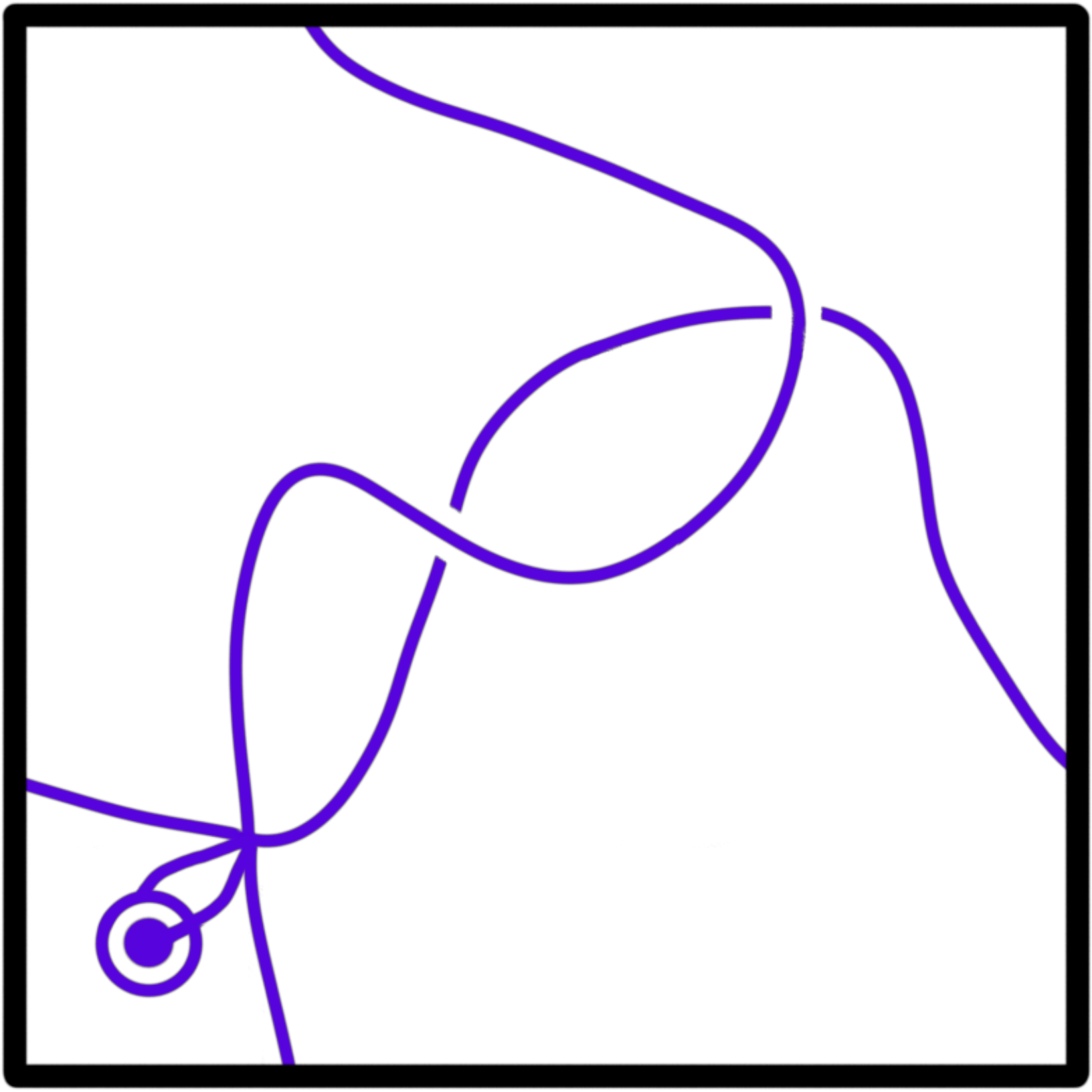}
        \caption{}
        \label{fig:pcu_4_dia}
    \end{subfigure}
    \begin{subfigure}[b]{0.18\textwidth}
        \centering
        \includegraphics[width=0.9\textwidth]{pcu_dia_F.pdf}
        \caption{}
        \label{fig:pcu_dia}
    \end{subfigure}
    \begin{subfigure}[b]{0.18\textwidth}
        \centering
        \includegraphics[width=0.9\textwidth]{pcu_uc.png}
        \caption{}
        \label{fig:pcu_uc}
    \end{subfigure}
    \caption{An example of the conversion of an ambient isotopy transformation into a finite sequence of $R$-moves: (b) A diagram obtained from the unit cell shown in (a). (b-f) Each diagram is obtained from the one that precedes it by successively applying one $R_{11}$ move, then one $R_2$ move, another $R_{11}$ move, and finally another $R_{2}$ move. (g) A unit cell represented by the diagram shown in (f).}
    \label{fig:ambient_isotopy_pcu}
\end{figure}

In Fig. \ref{fig:ambient_isotopy_pcu}, we present an example illustrating how a finite sequence of $R$-moves mimics an isotopy transformation of a unit cell. The diagram shown in Fig. \ref{fig:pcu_1_dia} is obtained from the unit cell depicted in Fig. \ref{fig:pcu_1_uc}. The sequence of $R$-moves consists of an $R_{11}$ move followed by an $R_2$ move, then another $R_{11}$ move and a final $R_2$ move. The resulting diagram, shown in Fig. \ref{fig:pcu_3_dia}, represents the unit cell shown in Fig. \ref{fig:pcu_uc}.

\subsection{Least tangled embeddings or ground states}
Let $K$ be an embedding of a graph $\mathcal{K}$ and $U$ a unit cell of $K$. We define the \textit{$\mathcal{U}$-family} of $K$ with respect to $U$, denoted by $\mathcal{U}(K,U)$, as the set containing $K$ and any embedding $K'$ that possesses a unit cell $U'$ connected to $U$ by a finitely many crossing changes.

Note that any graph embedding $K'$ belonging to $\mathcal{U}(K,U)$ is always understood to be paired with the unit cell $U'$ that is connected to $U$ by a finite series of crossing changes, and not any other unit cell representing $K'$. 

Within the $\mathcal{U}$-family $\mathcal{U}(K,U)$, we define the \textit{$\mathcal{G}$-family} of $K$ with respect to $U$, denoted by $\mathcal{G}(K,U)$, as the subset consisting of all embeddings having the crossing number that is the least among all embeddings of $\mathcal{U}(K,U)$. The elements of $\mathcal{G}(K,U)$ are called the \textit{ground states}.

Notice that if $K'$ belongs to $\mathcal{U}(K,U)$ with respect to the unit cell $U'$, then $K$ belongs to $\mathcal{U}(K',U')$ with respect to $U$, and thus $\mathcal{U}(K,U) = \mathcal{U}(K',U')$. Notice also that, by the definition of $\mathcal{G}(K,U)$, if $K'$ belongs to $\mathcal{U}(K,U)$ with respect to the unit cell $U'$, then we have $\mathcal{G}(K',U') =\mathcal{G}(K,U)$.

\subsection{Untangling number}
Let $K$ be an embedding of a graph $\mathcal{K}$ and let $U$ be a unit cell of $K$. Consider all sequences $(U_k)_{k=0,\dots,n}$  of unit cells of embeddings of $\mathcal{K}$ belonging to $\mathcal{U}(K,U)$, satisfying:
\begin{itemize}
    \item $U_0 = U$,
    \item $\forall k \in \lbrace 1, \dots,n\rbrace$,
    $U_k$ and $U_{k-1}$ are connected by one and only one crossing change,
    \item the embedding associated to $U_n$ belongs to $\mathcal{G}(K,U)$.
\end{itemize}
The \textit{untangling number of $K$ with respect to $U$}, denoted by $u(K,U)$, is the least integer $n$ over all such sequences. Any embedding belonging to $\mathcal{G}(K,U)$ that realises the untangling number of $K$ is called a \textit{nearest ground state of $K$}. \\

In the above definition, we choose to enforce the condition that every pair of consecutive unit cells in a given sequence $(U_k)_{k=0,\dots,n}$ is connected by one and only one crossing change, to ensure that the number $n$ indexing the sequence of unit cells corresponds to the total number of crossing changes.

An immediate result obtained from the above definition is that, the ground states are embeddings with untangling number 0, and that every embedding with untangling number 0 with respect to a chosen unit cell is a ground state.

It is clear that the untangling number with respect to a given unit cell varies when another unit cell is chosen. To stabilise it, one can consider the minimum of all untangling numbers among all unit cells of a given embedding of a graph, leading to the following definition of an invariant. The \textit{minimum untangling number} of an embedding $K$ of a graph $\mathcal{K}$ is defined as the minimum untangling number $u(K,U)$ over all unit cells $U$ of $K$.

\bibliography{bibliography_networks} 

@book{Adams.book,
  author = {Colin Adams},
  year = {2004},
  title = {The {K}not {B}ook: {A}n {E}lementary {I}ntroduction to the {M}athematical {T}heory of {K}nots},
  publisher = {American Mathematical Society},
  address = {Providence},
  isbn = {978-0-8218-3678-1}
}

@article{nets_MOFs,
    author = {O’Keeffe, Michael},
    title = {Nets, tiles, and metal-organic frameworks},
    journal = {APL Materials},
    volume = {2},
    number = {12},
    pages = {124106},
    year = {2014},
    issn = {2166-532X},
    doi = {10.1063/1.4901292}
}

@article{Eon:ib5036,
author = "Eon, Jean-Guillaume",
title = "{Topological features in crystal structures: a quotient graph assisted analysis of underlying nets and their embeddings}",
journal = "Acta Crystallographica Section A",
year = "2016",
volume = "72",
number = "3",
pages = "268--293",
doi = {10.1107/S2053273315022950}
}

@article{terminology_MOF_COPolymers,
title = {Terminology of metal-organic frameworks and coordination polymers ({IUPAC} {R}ecommendations 2013)},
author = {Stuart R. Batten and Neil R. Champness and Xiao-Ming Chen and Javier Garcia-Martinez and Susumu Kitagawa and Lars {Ö}hrstr{ö}m and Michael O’Keeffe and Myunghyun Paik Suh and Jan Reedijk},
pages = {1715--1724},
volume = {85},
number = {8},
journal = {Pure and Applied Chemistry},
doi = {10.1351/PAC-REC-12-11-20},
year = {2013}
}

@article{DELGADOFRIEDRICHS20052480,
title = {Crystal nets as graphs: Terminology and definitions},
journal = {Journal of Solid State Chemistry},
volume = {178},
number = {8},
pages = {2480-2485},
year = {2005},
issn = {0022-4596},
doi = {10.1016/j.jssc.2005.06.011},
author = {Olaf Delgado-Friedrichs and Michael O’Keeffe}
}

@article{alexandrov2011,
author ="Alexandrov, E. V. and Blatov, V. A. and Kochetkov, A. V. and Proserpio, D. M.",
title  ="Underlying nets in three-periodic coordination polymers: topology{,} taxonomy and prediction from a computer-aided analysis of the {C}ambridge {S}tructural {D}atabase",
journal  ="CrystEngComm",
year  ="2011",
volume  ="13",
issue  ="12",
pages  ="3947-3958",
publisher  ="The Royal Society of Chemistry",
doi  ="10.1039/C0CE00636J"
}

@article{mof_2001,
author = {Banglin Chen  and M. Eddaoudi  and S. T. Hyde  and M. O'Keeffe  and O. M. Yaghi },
title = {Interwoven {M}etal-{O}rganic {F}ramework on a {P}eriodic {M}inimal {S}urface with {E}xtra-{L}arge {P}ores},
journal = {Science},
volume = {291},
number = {5506},
pages = {1021-1023},
year = {2001},
doi = {10.1126/science.1056598}
}

@article{liu2016,
author = {Yuzhong Liu  and Yanhang Ma  and Yingbo Zhao  and Xixi Sun  and Felipe G{á}ndara  and Hiroyasu Furukawa  and Zheng Liu  and Hanyu Zhu  and Chenhui Zhu  and Kazutomo Suenaga  and Peter Oleynikov  and Ahmad S. Alshammari  and Xiang Zhang  and Osamu Terasaki  and Omar M. Yaghi },
title = {Weaving of organic threads into a crystalline covalent organic framework},
journal = {Science},
volume = {351},
number = {6271},
pages = {365-369},
year = {2016},
doi = {10.1126/science.aad4011}
}

@article{zhang2022,
author = {Zhang, Zhi-Hui and Andreassen, Bj{ö}rn J. and August, David P. and Leigh, David A. and Zhang, Liang},
title = {Molecular weaving},
journal = {Nature Materials},
volume = {21},
issue = {3},
pages = {275-283},
year = {2022},
doi = {10.1038/s41563-021-01179-w}
}

@article{rcsr,
author = {O’Keeffe, Michael and Peskov, Maxim A. and Ramsden, Stuart J. and Yaghi, Omar M.},
title = {The {R}eticular {C}hemistry {S}tructure {R}esource ({RCSR}) {D}atabase of, and {S}ymbols for, {C}rystal {N}ets},
journal = {Accounts of Chemical Research},
volume = {41},
number = {12},
pages = {1782-1789},
year = {2008},
doi = {10.1021/ar800124u}
}

@article{theta_curves,
author = {Moriuchi, Hiromasa},
title = {{A}n {E}numeration of {T}heta-curves with up to {S}even {C}rossings},
journal = {Journal of Knot Theory and Its Ramifications},
volume = {18},
number = {02},
pages = {167-197},
year = {2009},
doi = {10.1142/S0218216509006884}
}

@article{srs_net_nature,
author ="Abrahams, Brendan F. and Batten, Stuart R. and Hamit, Hasan and Hoskins, Bernard F. and Robson, Richard",
title  ="A wellsian ‘three-dimensional’ racemate: eight interpenetrating{,} enantiomorphic (10{,}3)-a nets{,} four right- and four left-handed",
journal  ="Chem. Commun.",
year  ="1996",
issue  ="11",
pages  ="1313-1314",
publisher  ="The Royal Society of Chemistry",
doi  ="10.1039/CC9960001313"
}

@article{dna_tensegrity_triangle,
author = {Lu, Brandon and Vecchioni, Simon and Ohayon, Yoel P. and Sha, Ruojie and Woloszyn, Karol and Yang, Bena and Mao, Chengde and Seeman, Nadrian C.},
title = {3{D} {H}exagonal {A}rrangement of {DNA} {T}ensegrity {T}riangles},
journal = {ACS Nano},
volume = {15},
number = {10},
pages = {16788-16793},
year = {2021},
doi = {10.1021/acsnano.1c06963}
}

@article{site_directed_placement_DNA_origami,
    author = {Martynenko, Irina V. and Erber, Elisabeth and Ruider, Veronika and Dass, Mihir and Posnjak, Gregor and Yin, Xin and Altpeter, Philipp and Liedl, Tim},
    title = {Site-directed placement of three-dimensional {DNA} origami},
    journal = {Nature Nanotechnology},
    volume = {18},
    issue = {12},
    pages = {1456-1462},
    year = {2023},
    doi = {10.1038/s41565-023-01487-z}
}

@article{C4SM01226G,
author ="Sorenson, Gregory P. and Schmitt, Adam K. and Mahanthappa, Mahesh K.",
title  ="Discovery of a tetracontinuous{,} aqueous lyotropic network phase with unusual 3{D}-hexagonal symmetry",
journal  ="Soft Matter",
year  ="2014",
volume  ="10",
issue  ="41",
pages  ="8229-8235",
publisher  ="The Royal Society of Chemistry",
doi  ="10.1039/C4SM01226G"
}

@article{DELGADOFRIEDRICHS20052533,
title = {What do we know about three-periodic nets?},
journal = {Journal of Solid State Chemistry},
volume = {178},
number = {8},
pages = {2533-2554},
year = {2005},
note = {Reticular Chemistry: Design, Synthesis, Properties and Applications of Metal-Organic Polyhedra and Frameworks},
issn = {0022-4596},
doi = {10.1016/j.jssc.2005.06.037},
author = {Olaf Delgado-Friedrichs and Martin D. Foster and Michael O’Keeffe and Davide M. Proserpio and Michael M.J. Treacy and Omar M. Yaghi}
}

@article{Hyde2008ASH,
author = {Hyde, Stephen T. and O'Keeffe, Michael and Proserpio, Davide M.},
title = {A {S}hort {H}istory of an {E}lusive {Y}et {U}biquitous {S}tructure in {C}hemistry, {M}aterials, and {M}athematics.},
journal = {Angewandte Chemie International Edition},
volume = {47},
number = {42},
pages = {7996-8000},
doi = {10.1002/anie.200801519},
year = {2008}
}

@article{Alexandrov:eo5016,
author = "Alexandrov, Eugeny V. and Blatov, Vladislav A. and Proserpio, Davide M.",
title = "{A topological method for the classification of entanglements in crystal networks}",
journal = "Acta Crystallographica Section A",
year = "2012",
volume = "68",
number = "4",
pages = "484--493",
doi = {10.1107/S0108767312019034}
}

@article{Baburin:eo5056,
author = "Baburin, Igor A.",
title = "{On the group-theoretical approach to the study of interpenetrating nets}",
journal = "Acta Crystallographica Section A",
year = "2016",
volume = "72",
number = "3",
pages = "366--375",
doi = {10.1107/S2053273316002692}
}

@article{periodic_ent_I,
author = {Evans, Myfanwy E. and Robins, Vanessa and Hyde, Stephen T.},
title = {Periodic entanglement {I}: networks from hyperbolic reticulations},
journal = {Acta Crystallographica Section A},
volume = {69},
number = {3},
pages = {241-261},
doi = {10.1107/S0108767313001670},
year = {2013}
}

@article{CARLUCCI2003247,
title = {Polycatenation, polythreading and polyknotting in coordination network chemistry},
journal = {Coordination Chemistry Reviews},
volume = {246},
number = {1},
pages = {247-289},
year = {2003},
issn = {0010-8545},
doi = {10.1016/S0010-8545(03)00126-7},
author = {Lucia Carlucci and Gianfranco Ciani and Davide M. Proserpio}
}

@article{HYDE2011676,
title = {From untangled graphs and nets to tangled materials},
journal = {Solid State Sciences},
volume = {13},
number = {4},
pages = {676-683},
year = {2011},
issn = {1293-2558},
doi = {10.1016/j.solidstatesciences.2010.10.028},
author = {Stephen T. Hyde and Olaf Delgado-Friedrichs}
}

@article{Power:ib5087,
author = "Power, Stephen C. and Baburin, Igor A. and Proserpio, Davide M.",
title = "{Isotopy classes for 3-periodic net embeddings}",
journal = "Acta Crystallographica Section A",
year = "2020",
volume = "76",
number = "3",
pages = "275--301",
doi = {10.1107/S2053273320000625}
}

@article{4NG,
author = {Tang, Yumin and Xue, Yi-nan and Yang, Shu-Gui and Zhang, Ruibin and Liu, Feng and Zeng, Xiangbing and Ungar, Goran},
title = {Gyroid {L}abyrinth of {S}upertwisted {D}ouble {H}elices in a {L}iquid {C}rystal {P}olymer},
journal = {Angewandte Chemie International Edition},
volume = {65},
number = {4},
pages = {e22314},
doi = {10.1002/anie.202522314},
year = {2026}
}

@article{ravels,
author ="Castle, Toen and Evans, Myfanwy E. and Hyde, S. T.",
title  ="Ravels: knot-free but not free. {N}ovel entanglements of graphs in 3-space",
journal  ="New J. Chem.",
year  ="2008",
volume  ="32",
issue  ="9",
pages  ="1484-1492",
publisher  ="The Royal Society of Chemistry",
doi  ="10.1039/B719665B"
}

@article{synthesised_ravel,
    author = {Li, Feng and Clegg, Jack K. and Lindoy, Leonard F. and Macquart, Ren{é} B. and Meehan, George V.},
    title = {Metallosupramolecular self-assembly of a universal 3-ravel},
    journal = {Nature Communications},
    year = {2011},
    volume = {2},
    pages = {205},
    issue = {1},
    doi = {10.1038/ncomms1208}
}

@article{sym14040822,
AUTHOR = {O’Keeffe, Michael and Treacy, Michael M. J.},
TITLE = {The {S}ymmetry and {T}opology of {F}inite and {P}eriodic {G}raphs and {T}heir {E}mbeddings in {T}hree-{D}imensional {E}uclidean {S}pace},
JOURNAL = {Symmetry},
VOLUME = {14},
YEAR = {2022},
issue = {4},
pages = {822},
DOI = {10.3390/sym14040822}
}

@article{Wells:a01232,
author = "Wells, A. F.",
title = "{The geometrical basis of crystal chemistry. {P}art 1}",
journal = "Acta Crystallographica",
year = "1954",
volume = "7",
number = "8-9",
pages = "535--544",
doi = {10.1107/S0365110X5400182X}
}

@article{Wells:a01308,
author = "Wells, A. F.",
title = "{The geometrical basis of crystal chemistry. {P}art 4}",
journal = "Acta Crystallographica",
year = "1954",
volume = "7",
number = "12",
pages = "849--853",
doi = {10.1107/S0365110X54002587}
}

@article{tesselate_decussate,
AUTHOR = {O’Keeffe, Michael and Treacy, Michael M. J.},
TITLE = {Embeddings of {G}raphs: {T}essellate and {D}ecussate {S}tructures},
JOURNAL = {International Journal of Topology},
VOLUME = {1},
YEAR = {2024},
NUMBER = {1},
PAGES = {1--10},
DOI = {10.3390/ijt1010001}
}

@article{Delgado-Friedrichs:au5000,
author = "Delgado-Friedrichs, Olaf and O'Keeffe, Michael",
title = "{Identification of and symmetry computation for crystal nets}",
journal = "Acta Crystallographica Section A",
year = "2003",
volume = "59",
number = "4",
pages = "351--360",
doi = {10.1107/S0108767303012017}
}

@article{equilibrium_placement,
  author    = {Delgado-Friedrichs, Olaf},
  title     = {Equilibrium {P}lacement of {P}eriodic {G}raphs and {C}onvexity of {P}lane {T}ilings},
  journal   = {Discrete \& Computational Geometry},
  year      = {2005},
  volume    = {33},
  number    = {1},
  pages     = {67--81},
  doi       = {10.1007/s00454-004-1147-x}
}

@article{topospro,
author = {Blatov, Vladislav A. and Shevchenko, Alexander P. and Proserpio, Davide M.},
title = {Applied {T}opological {A}nalysis of {C}rystal {S}tructures with the {P}rogram {P}ackage {T}opos{P}ro},
journal = {Crystal Growth \& Design},
volume = {14},
number = {7},
pages = {3576-3586},
year = {2014},
doi = {10.1021/cg500498k}
}

@article{myf2011_entanglement_graphs,
    author = {Castle, Toen and Evans, Myfanwy E. and Hyde, Stephen T.},
    title = {Entanglement of {E}mbedded {G}raphs},
    journal = {Progress of Theoretical Physics Supplement},
    volume = {191},
    pages = {235-244},
    year = {2011},
    issn = {0375-9687},
    doi = {10.1143/PTPS.191.235}
}

@article{myf_ideal_geo,
author = {Evans, Myfanwy E.  and Robins, Vanessa  and Hyde, Stephen T. },
title = {Ideal geometry of periodic entanglements},
journal = {Proceedings of the Royal Society A: Mathematical, Physical and Engineering Sciences},
volume = {471},
number = {2181},
pages = {20150254},
year = {2015},
doi = {10.1098/rspa.2015.0254}
}

@article{ANDRIAMANALINA2025109346,
title = {Diagrammatic representations of 3-periodic entanglements},
journal = {Topology and its Applications},
volume = {368},
pages = {109346},
year = {2025},
issn = {0166-8641},
doi = {10.1016/j.topol.2025.109346},
author = {Toky Andriamanalina and Myfanwy E. Evans and Sonia Mahmoudi}
}

@misc{andriamanalina2025untanglingnumber3periodictangles,
      title={The untangling number of 3-periodic tangles}, 
      author={Toky Andriamanalina and Sonia Mahmoudi and Myfanwy E. Evans},
      year={2026},
      eprint={2504.01747},
      archivePrefix={arXiv},
      primaryClass={math.GT},
      note = {Preprint at \url{https://arxiv.org/abs/2504.01747}}, 
}

@book{Murasugi1996,
author="Murasugi, Kunio",
title="{K}not {T}heory and {I}ts {A}pplications",
year="1996",
publisher="Birkh{\"a}user Boston",
address="Boston, MA",
isbn="978-0-8176-4719-3",
doi="10.1007/978-0-8176-4719-3"
}

@article{R_moves_graph,
 ISSN = {00029947},
 author = {Louis H. Kauffman},
 journal = {Transactions of the American Mathematical Society},
 number = {2},
 pages = {697--710},
 publisher = {American Mathematical Society},
 title = {Invariants of {G}raphs in {T}hree-{S}pace},
 volume = {311},
 year = {1989},
 doi = {10.1090/S0002-9947-1989-0946218-0}
}

@article{KAUR2017586,
title = {Gauss diagrams, unknotting numbers and trivializing numbers of spatial graphs},
author = {K. Kaur and S. Kamada and A. Kawauchi and M. Prabhakar},
journal = {Topology and its Applications},
volume = {230},
pages = {586-598},
year = {2017},
issn = {0166-8641},
doi = {10.1016/j.topol.2017.08.037}
}

@article{batten_1998,
author = {Batten, Stuart R. and Robson, Richard},
title = {Interpenetrating {N}ets: {O}rdered, {P}eriodic {E}ntanglement},
journal = {Angewandte Chemie International Edition},
volume = {37},
number = {11},
pages = {1460-1494},
doi = {10.1002/(SICI)1521-3773(19980619)37:11<1460::AID-ANIE1460>3.0.CO;2-Z},
year = {1998}
}

@book{BurdeZieschangHeusener+2013,
title = {Knots},
author = {Gerhard Burde and Heiner Zieschang and Michael Heusener},
publisher = {De Gruyter},
address = {Berlin, Boston},
doi = {10.1515/9783110270785},
isbn = {9783110270785},
year = {2013}
}

@article{k6_linked,
author = {Conway, J. H. and McA. Gordon, C.},
title = {Knots and links in spatial graphs},
journal = {Journal of Graph Theory},
volume = {7},
number = {4},
pages = {445-453},
doi = {10.1002/jgt.3190070410},
year = {1983}
}

@article{Murasugi1965,
  author    = {Kunio Murasugi},
  title     = {On a Certain Numerical Invariant of Link Types},
  journal   = {Transactions of the American Mathematical Society},
  volume    = {117},
  year      = {1965},
  pages     = {387--422},
  doi       = {10.1090/S0002-9947-1965-0171275-5}
}

@article{nakanishi_I,
author = {Nakanishi, Yasutaka},
title = {A note on unknotting number},
journal = {Mathematics Seminar Notes},
publisher = {Kobe University},
volume = {9},
number = {1},
pages = {99-108},
year = {1981},
doi = {10.24546/E0001593}
}

@article{ma-qiu,
    author = {Ma, Jiming and Qiu, Ruifeng},
    title = {A Lower Bound on Unknotting Number},
    journal = {Chinese Annals of Mathematics, Series B},
    volume = {27},
    number = {4},
    pages = {437-440},
    year = {2006},
    doi = {10.1007/s11401-004-0390-z}
}

@article{Applebaum18082025,
author = {Taylor Applebaum and Sam Blackwell and Alex Davies and Thomas Edlich and Andr{á}s Juh{á}sz and Marc Lackenby and Nenad Toma{š}ev and Daniel Zheng},
title = {The Unknotting Number, Hard Unknot Diagrams, and Reinforcement Learning},
journal = {Experimental Mathematics},
pages = {1--19},
year = {2025},
publisher = {Taylor \& Francis},
doi = {10.1080/10586458.2025.2542174},
volume = {0},
number = {0}
}
\bibliographystyle{rsc} 

\newpage


\section{Upper bounds for the untangling numbers of various structures}


\begin{figure*}[hbtp]
    \centering
    \includegraphics[width=0.95\textwidth]{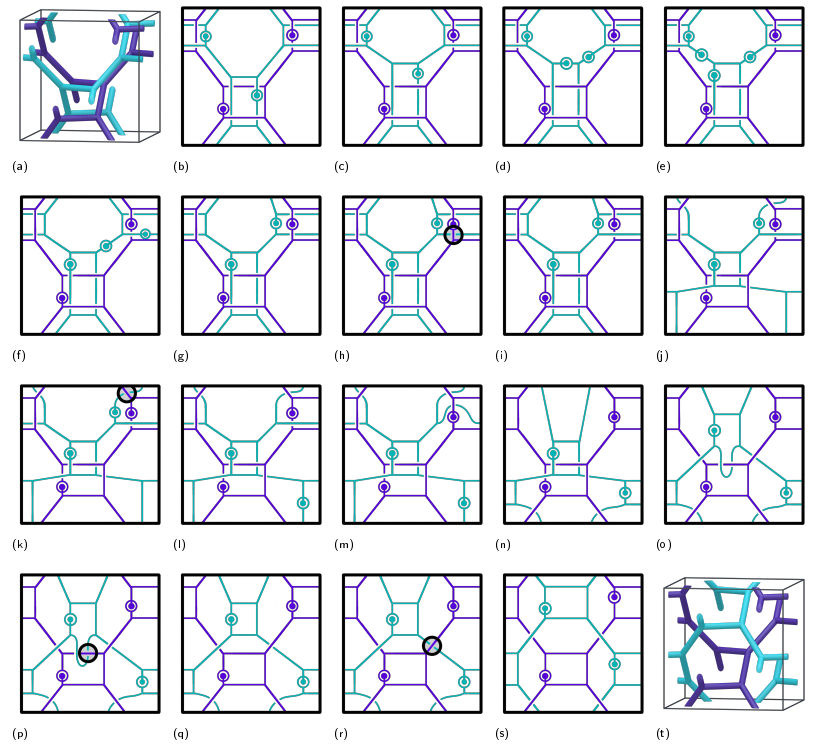}
    \caption{Computation of an upper bound for the untangling number of \textbf{utp-c***}: (a) A unit cell of \textbf{utp-c***} with respect to which the computation is done. (t) A unit cell of \textbf{utp-c**}. Four crossing changes, applied on the crossings highlighted by black circles, are needed to transform the unit cell of \textbf{utp-c***} into that of \textbf{utp-c**}. This upper bound is not as good as the one given by the computation in Fig. \ref{fig:untangling_utp-c_star_star_star_2nd_time_into_utp-c_star}.}
    \label{fig:untangling_utp-c_star_star_star_1st_time_into_utp-c_star_star}
\end{figure*}


\begin{figure*}[hbtp]
    \centering
    \includegraphics[width=0.95\textwidth]{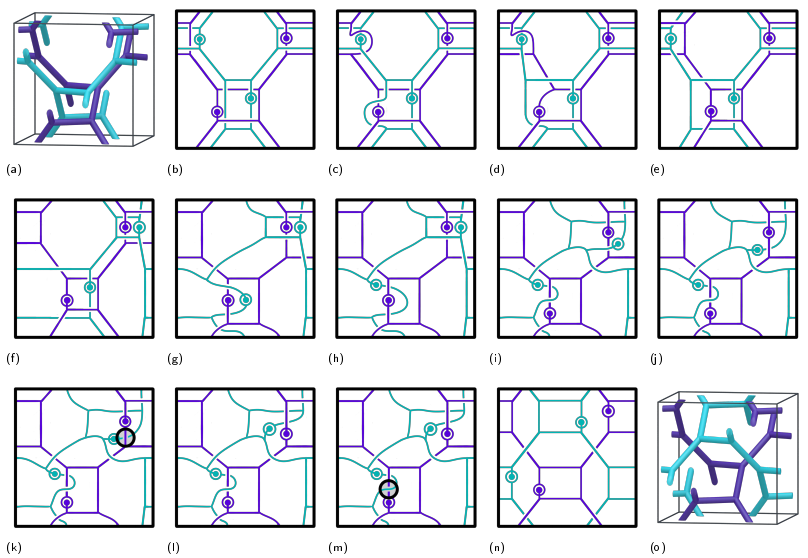}
    \caption{Computation of an upper bound for the untangling number of \textbf{utp-c***}: (a) A unit cell of \textbf{utp-c***}, the same as that shown in Fig. \ref{fig:untangling_utp-c_star_star_star_1st_time_into_utp-c_star_star}a. (o) A unit cell of \textbf{utp-c*}. Only two crossing changes, applied on the crossings highlighted by black circles, are needed to transform the unit cell of \textbf{utp-c***} into that of \textbf{utp-c*}, which is fewer than what is seen in Fig. \ref{fig:untangling_utp-c_star_star_star_1st_time_into_utp-c_star_star}, to reach the unit cell of \textbf{utp-c**}. This suggests that the untangling number of \textbf{utp-c***} with respect to the unit cell in (a) is at most 2, and its nearest ground state is \textbf{utp-c*}.}
    \label{fig:untangling_utp-c_star_star_star_2nd_time_into_utp-c_star}
\end{figure*}


\begin{figure}[hbtp]
    \centering
    \includegraphics[width=0.95\textwidth]{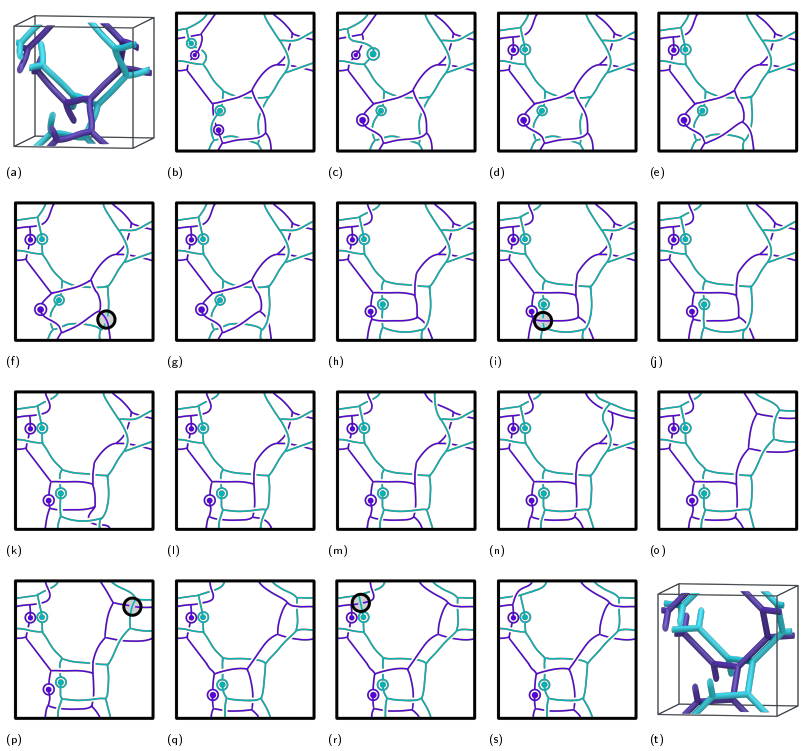}
    \caption{Computation of an upper bound for the untangling number of \textbf{srs-c**}: (a) A unit cell of \textbf{srs-c**} with respect to which the computation is done. (t) A unit cell of \textbf{srs-c*}. Four crossings changes, applied on the crossings highlighted by black circles are needed to transform the unit cell of \textbf{srs-c**} into that of \textbf{srs-c*}, which means that its untangling number is at most 4.}
    \label{fig:untangling_n0p4on2_to_the_6}
\end{figure}


\begin{figure}[hbtp]
    \centering
    \includegraphics[width=0.95\textwidth]{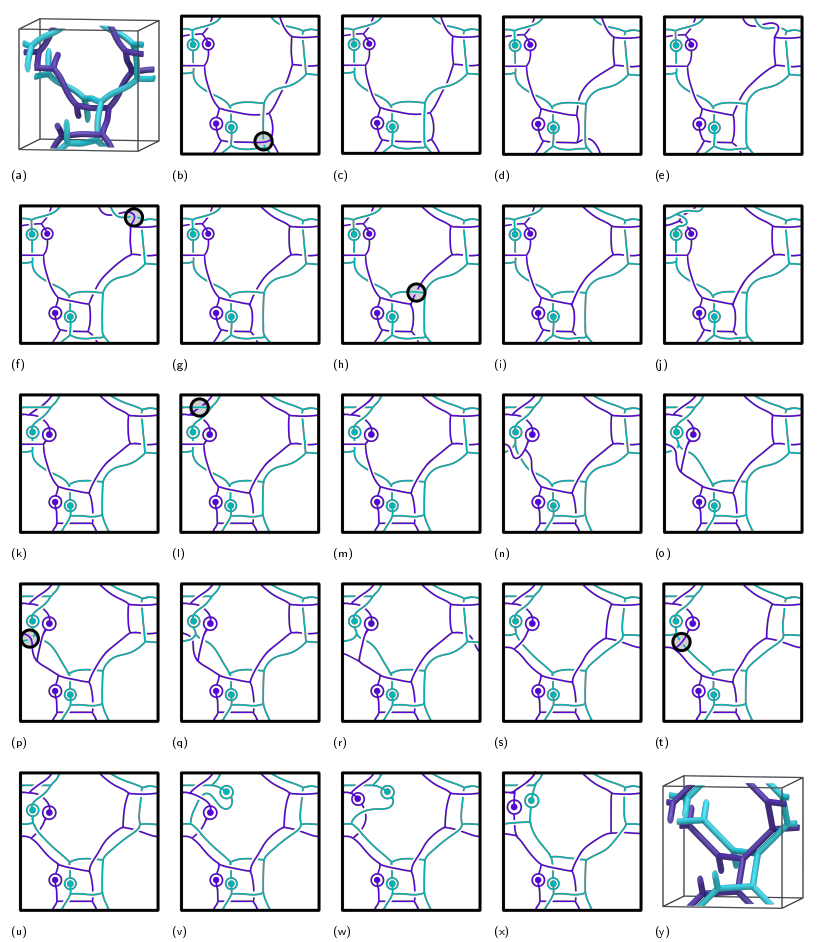}
    \caption{Computation of an upper bound for the untangling number of the embedding of a graph comprising two same-handed \textbf{srs} networks, shown in Fig. 2c of the main manuscript: (a) A unit cell of the structure with respect to which the computation is done. (y) A unit cell of \textbf{srs-c*}. Six crossing changes, applied on the crossings highlighted by black circles, are needed to transform the unit cell of the structure into that of \textbf{srs-c*}, meaning that the untangling number is at most $6$.}
    \label{fig:untangling_0p6on2_to_the_6}
\end{figure}


\begin{figure}[hbtp]
    \centering
    \includegraphics[width=0.95\textwidth]{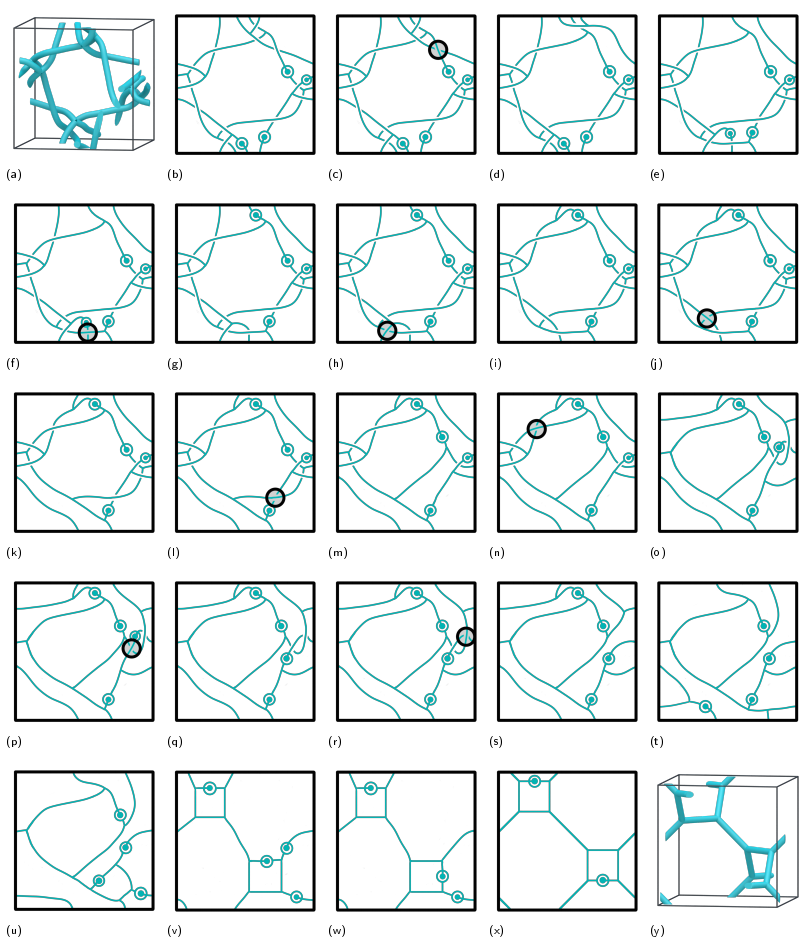}
    \caption{Computation of an upper bound for the untangling number of the tangled embedding of \textbf{srs} shown in Fig. 16a of the main manuscript: (a) A unit cell of the structure with respect to which the computation is done. (y) A unit cell of the barycentric embedding of \textbf{srs}. Eight crossing changes, applied on the crossings highlighted by black circles, are needed to transform the unit cell of the tangled embedding to that of the barycentric embedding of \textbf{srs}, which means that the untangling number is at most 8.}
    \label{fig:untangling_srs-me}
\end{figure}


\begin{figure}[hbtp]
    \centering
    \includegraphics[width=0.95\textwidth]{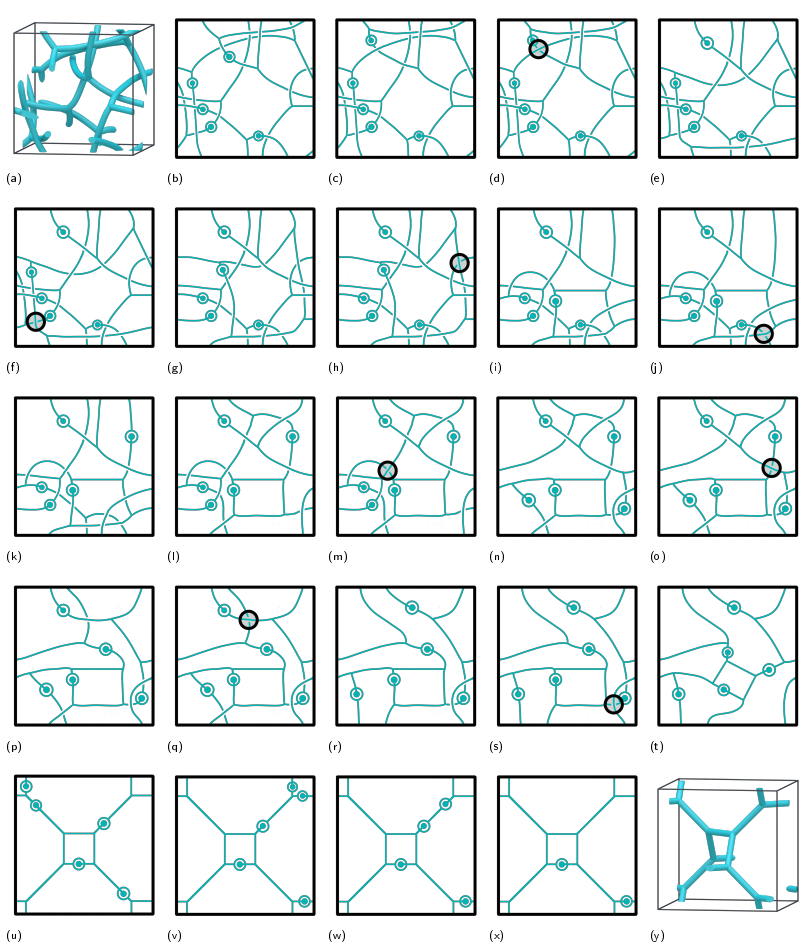}
    \caption{Computation of an upper bound for the untangling number of \textbf{srs-z}: (a) A unit cell of \textbf{srs-z} with respect to which the computation is done. (y) A unit cell of the barycentric embedding of \textbf{srs}. Eight crossing changes, applied on the crossings highlighted by black circles, are needed to transform the unit cell of \textbf{srs-z} to that of the barycentric embedding of \textbf{srs}, which means that the untangling number of \textbf{srs-z} is at most 8.}
    \label{fig:untangling_srs-z}
\end{figure}


\end{document}